\documentclass[12pt]{article}
 \usepackage[latin1]{inputenc}  
\usepackage{amssymb,amsmath,amsthm,latexsym,enumerate}
\usepackage{graphicx,subfigure}
\usepackage{url}

\usepackage{color}

\def\figpath{figures/}
\def\figpathnew{figures/}

\newtheorem{algorithm}{Algorithm}
\newtheorem{testproblem}{Test Problem}

\DeclareMathOperator{\sep}{sep}
\DeclareMathOperator{\dist}{dist}
\def\eps{\varepsilon}

\def\hu{\hat u}

\def\RR{\mathbb{R}}

\begin{document}

\title{Adaptive RBF-FD Method for Elliptic Problems with Point Singularities in
2D\thanks{This work was supported in part by Vietnam's National Foundation for
Science and Technology Development (NAFOSTED) under grant number 101.01-2014.28,
by a Natural Science Research Project of the Ministry of Education and
Training under grant number B2015-TN06-02, and by the German Academic Exchange Service (DAAD) under the Programme 
``Research Stays for University Academics and Scientists.''}}
\author{
Dang Thi Oanh\thanks{Division of Science-Technology $\&$ International Cooperation,
University of Information $\&$ Communication Technology,
Thai Nguyen University,
Quyet Thang Commune, Thai Nguyen City, Vietnam, {\tt dtoanhtn@gmail.com}},\quad 
Oleg Davydov\thanks{Department of Mathematics, University of Giessen, Arndtstrasse 2, 
35392 Giessen, Germany, {\tt oleg.davydov@math.uni-giessen.de}}\quad  
and
Hoang Xuan Phu\thanks{Institute of Mathematics, Vietnam Academy of Science and Technology, 
18 Hoang Quoc Viet Road, Hanoi, Vietnam, {\tt hxphu@math.ac.vn}} 
}
\maketitle

\begin{abstract}
We describe and test numerically an adaptive meshless generalized finite difference method based on radial basis functions that
competes well with the finite element method on standard benchmark problems with  
reentrant corners of the boundary, sharp peaks and rapid oscillations in the neighborhood of an isolated point.
This is achieved thanks to significant improvements introduced into the earlier algorithms 
of [Oleg Davydov and Dang~Thi Oanh,
Adaptive meshless centers and {RBF} stencils for {P}oisson equation,
{Journal of {C}omputational {P}hysics}, 230:287--304, 2011], 
including a new error indicator of Zienkiewicz-Zhu type.
\end{abstract}

\section{Introduction}
\label{Intro}
Let us consider the Dirichlet boundary value problem: 
Find $u: \overline{\Omega} \to \mathbb{R}$ such that
      \begin{equation}      \label{poi}
      L u=f\;\hbox{ on $\Omega$ },\quad
      u|_{\partial\Omega}=g,
      \end{equation}
where $L$ is a linear elliptic differential operator of second order, $\Omega\subset\RR^2$ is a given bounded domain, the function $f$ is defined on $\Omega$,
and the function $g$ is defined on the boundary $\partial\Omega$ of $\Omega$. 
A \emph{generalized finite difference} discretization of the Dirichlet problem (\ref{poi}) is given by the following linear system
with respect to the vector $\hu=[\hu_\xi]_{\xi\in \Xi}$:
\begin{equation}
\label{dird}
\sum_{\xi\in\Xi_\zeta}w_{\zeta,\xi}\hu_\xi=f(\zeta),
\quad \zeta\in\Xi_{\rm int};
\qquad \hu_\xi=g(\xi),\quad \xi\in\partial\Xi,
\end{equation}
where
\begin{itemize}
\item
$\Xi \subset \overline{\Omega}$ is the set of discretization centers;
\item 
$\hu$ represents the approximation of the solution $u$ of (\ref{poi}) at the
points $\xi\in \Xi$;
\item
$\partial \Xi := \Xi \cap \partial\Omega$ is the set of boundary discretization centers;

\item
$\Xi_{\rm int} := \Xi \setminus \partial \Xi$ is the set of interior discretization centers;

\item
$\Xi_\zeta$ is a set (called the \emph{stencil support} of $\zeta$) that consists of the considered center $\zeta$ and some 
selected neighbor points $\xi_i \in \Xi$;

\item
$w_{\zeta,\xi}\in\RR$ are the \emph{stencil weights} chosen such that 
$\sum_{\xi\in\Xi_\zeta}w_{\zeta,\xi}u(\xi)$ is an approximation of $Lu(\zeta)$.

\end{itemize}
To set up the system (\ref{dird}), three tasks have to be addressed: 
(a) how to generate $\Xi$, (b) how to choose the stencil supports $\Xi_\zeta$, and
(c) how to compute suitable weights $w_{\zeta,\xi}$. 

In the RBF-FD method 
the weights  $w_{\zeta,\xi}$, $\xi\in \Xi_\zeta$, are %
generated through the interpolation with radial basis functions. Referring to 
\cite{DavyOanh11,FoFlbook15} for further details and references, we briefly describe this approach. Let 
$\phi:\RR_+\to\RR$ be a positive definite radial basis function \cite{Buhmann03}, 
for example the Gaussian function
\begin{equation}
\label{G}
\phi(r)=e^{-\eps ^2 r^2},
\end{equation}
where $\eps$ is the \emph{shape parameter}. Given $\zeta\in\Xi_{\rm int}$ and
$\Xi_\zeta=\{\zeta_0,\zeta_1,\ldots,\zeta_k\}\subset\Xi$, with $\zeta_0=\zeta$, we set $\varphi_i(x)=\phi(\|x-\zeta_i\|)$,
$x\in\RR^2$, where $\|\cdot\|$ denotes the Euclidean norm in $\RR^2$. Assuming for
simplicity that the operator $L$ has the form
$$
Lu(x)=\Delta u(x)+c(x)u(x),$$
we first find the weights $w_{i}$ such that 
$$
\Delta s(\zeta)=\sum_{i=0}^k w_{i}u(\zeta_i),$$
where 
$s(x):= \sum_{i=0}^k a_i\varphi_i(x)$, $a_i\in\RR$,
satisfies the interpolation condition $s(\zeta_i)=u(\zeta_i)$, $i=0,\ldots,k$. The vector 
$w=[w_i]_{i=0}^k$ can be computed by solving the linear system
\begin{equation}
\label{ww}
\Phi_{\Xi_\zeta}w=[\Delta \varphi_i(\zeta)]_{i=0}^k,\quad\text{with}\quad
\Phi_{\Xi_\zeta}:=[\varphi_j(\zeta_i)]_{i,j=0}^k.
\end{equation}
In particular, for the Gaussian \eqref{G}, we have
\begin{equation}\label{GPhiX}
\Phi_{\Xi_\zeta}=[e^{-\eps ^2 \|\zeta_i-\zeta_j\|^2}]_{i,j=0}^k,
\quad 
\Delta \varphi_i(\zeta)=4 \eps ^2e^{-\eps^2 \|\zeta-\zeta_i\|^2}(\eps ^2\|\zeta-\zeta_i\|^2-1).%
\end{equation}
Assuming that the interpolant $s$ provides a good approximation of the function
$u$, we expect that 
$\Delta u(\zeta)\approx\Delta s(\zeta)$, and thus
$$
Lu(\zeta)\approx\sum_{i=0}^k w_{i}u(\zeta_i)+c(\zeta)u(\zeta).$$
Therefore the weights $w_{\zeta,\xi}$ in \eqref{dird} are chosen as follows:
$$
w_{\zeta,\zeta}=w_0+c(\zeta),\qquad w_{\zeta,\zeta_i}=w_i,\quad i=1,\ldots,k.$$
We refer to \cite{DS2016} for the bounds for the numerical differentiation error 
$$
\Big|Lu(\zeta)-\sum_{\xi\in\Xi_\zeta}w_{\zeta,\xi}u(\xi)\Big|.$$

The set of discretization centers $\Xi$ does not have to form a grid or mesh, therefore
RBF-FD is a meshless method \cite{NguRabBorDuf}. For more complicated problems it is
advantageous to adapt the distribution of the centers to the features of the domain $\Omega$
and/or to the singularities of the solution $u$. This can be achieved through \emph{adaptive refinement}
of $\Xi$. In \cite{DavyOanh11} we suggested a refinement algorithm and an algorithm 
for stencil support selection, leading to an effective meshless method capable of competing with the
finite element method on a number of benchmark test problems. However, further experiments have
shown certain deterioration of the approximation quality after many refinement steps, and
suboptimal performance for more difficult test problems. 

This motivated the current study, where both stencil support selection and refinement have been improved.
The new algorithms presented in Sections~2 and 3 deliver the
stencil supports $\Xi_\zeta$ with more evenly distributed 
points, and the adaptively selected sets  $\Xi$ that better reflect the singularities of
the solution. In the same time the improved method is more efficient because a costly post-processing step
 aimed at reducing the deterioration of the centers in the cause of subsequent refinements has been removed.
One of the major differences in the refinement algorithm comparing to  \cite{DavyOanh11} 
is an error indicator of Zienkiewicz-Zhu type used instead of a simple gradient estimate.
This leads to a significant improvement of the performance of the adaptive
RBF-FD method for more difficult problems. Section~4 is devoted to numerical experiments with the
test problems considered previously in \cite{DavyOanh11},  several test
problems suggested in  \cite{Mitchell2013} as benchmarks for testing adaptive grid
refinement, and a problem on a domain with a circular slit. 
In this paper we concentrate on the elliptic problems with \emph{point singularities}, such as the
reentrant corners of the boundary, sharp peaks and oscillations in the neighborhood of an isolated point.
Problems with line or curve singularities, boundary layers and wave fronts require further adjustments of the
algorithms taking into account the anisotropy of the solution $u$ and will be considered elsewhere. As in \cite{DavyOanh11}, comparison is provided
with the numerical results obtained with the adaptive finite element method of MATLAB PDE Toolbox \cite{PDEtool}.
 The results confirm the robust and competitive performance of the suggested method 
for problems with point singularities. A conclusion is given in Section~5.

\section{Meshless stencil support selection}
\label{Selec}

Given $\zeta\in\Xi_{\rm int}$, let 
$\Xi_\zeta=\{\zeta,\zeta_1,\ldots,\zeta_n\}\subset\Xi$, where the points
$\zeta_1,\ldots,\zeta_n$ are ordered counterclockwise with respect to $\zeta$. 
Following \cite[Section 5]{DavyOanh11}, we set
$$
\mu(\zeta_1,\ldots,\zeta_n):=\sum_{i=1}^n\alpha_i^2, %
\quad
\underline{\alpha}(\zeta_1,\ldots,\zeta_n):=\min_{1\le i\le n}\alpha_i,
\quad
\overline{\alpha}(\zeta_1,\ldots,\zeta_n):=\max_{1\le i\le n}\alpha_i,
$$
where $\alpha_i$ denotes the angle between the rays $\zeta\zeta_i$, $\zeta\zeta_{i+1}$ in the counterclockwise 
direction, with the cyclic identification 
$\zeta_{n+i}:=\zeta_i$. Since $\sum_{i=1}^n\alpha_i=2\pi$, 
the minimum of the expression $\sum_{i=1}^n\alpha_i^2$ is attained
for the uniformly spaced directions $\zeta\zeta_i$ such that the
angles $\alpha_i$ are all equal $\alpha_1=\cdots=\alpha_n=2\pi/n$. 
Although we cannot expect this minimum to be achieved for some subset 
$\{\zeta_1,\ldots,\zeta_n\}$ of the finite set $\Xi$, we use $\mu$ as a measure
of angle uniformity when comparing two subsets of this type. We prefer $\mu$ over the 
possible alternative measure 
$\overline{\alpha}/\underline{\alpha}$ that has the drawback that it does not `see' small
improvements in the uniformity that do not affect either minimum or maximum angle. 
It is however convenient to use it for a termination criterion
\begin{equation}\label{term1}
\overline{\alpha}(\zeta_1,\ldots,\zeta_n) \le v\,\underline{\alpha}(\zeta_1,\ldots,\zeta_n),
\end{equation}
with some tolerance $v>1.0$. Starting with $k$ nearest points, we replace $\zeta_j$'s one after
another by the more distant points $\xi\in\Xi\setminus\Xi_\zeta$ so that $\mu(\zeta_1,\ldots,\zeta_k)$ becomes smaller. On the other
hand, since nearby points are preferable, we introduce the second termination criterion
\begin{equation}\label{term2}
\|\zeta-\xi\| \ge  \frac{c}{2\,k} \sum_{j=1}^k \big(\|\zeta_j-\zeta\| + 
\|\zeta_j-\zeta_{j+1}\|\big),
\end{equation}
with a tolerance $c>1.0$ for the exchange candidate $\xi\in\Xi\setminus\Xi_\zeta$. 

Empirically chosen values
$v=2.5$ and  $c=3.0$ %
work well in all our numerical experiments, where we always choose $k=6$ to
guarantee that the density of the system matrix of \eqref{dird} is comparable to the density of
the system matrix of the conforming finite element method with linear shape functions, see 
\cite{DavyOanh11}. For the sake of efficiency we search initially for $m=50$ nearest points, and
double the size of the local cloud each time it has been exhausted.

\begin{algorithm}\label{alg1}\rm
{\bf Meshless stencil support selection} \\
\emph{Input:} $\Xi$, $\zeta$. \emph{Output:} $\Xi_\zeta$. 
\emph{Parameters:} $k$ (the number of points in $\Xi_\zeta\setminus\{\zeta\}$), 
$v>1.0$ (the angle uniformity tolerance), $c>1.0$ (distance tolerance), 
and  $m>k$  (the increment size of the local cloud).
Parameter values used in our numerical experiments: $k=6$,  $v=2.5$,  $c=3.0$
and  $m = 50$. 
\begin{enumerate}[I.]
\item\label{step1} 
 \emph{Find} $m$ nearest points  $\xi_1,\ldots,\xi_m$ in $\Xi\setminus\{\zeta\}$ to $\zeta$,
sorted by increasing distance to $\zeta$, and initialize 
$\Xi_\zeta:=\{\zeta,\zeta_1,\ldots,\zeta_k\}=\{\zeta,\xi_1,\ldots,\xi_k\}$
and  $i := k+1$.  
\item
\emph{While} $i \le m$:
\begin{enumerate}[1.]
\item 
 \emph{If} $\|\zeta-\xi_i\| \ge  \frac{c}{2\,k} \sum_{j=1}^k \big(\|\zeta_j-\zeta\| + 
\|\zeta_j-\zeta_{j+1}\|\big)$,  \emph{then} STOP and return  $\Xi_\zeta$.
\item Compute the angles $\alpha'_1,\ldots,\alpha'_{k+1}$  formed by the extended set
$\{\zeta'_1,\ldots,\zeta'_{k+1}\}=\{\zeta_1,\ldots,\zeta_k,\xi_i\}$. 
\emph{If}  both angles between $\zeta\xi_i$ \emph{and} its two neighboring rays 
are greater than the minimum angle 
$\underline{\alpha}':=\underline{\alpha}(\zeta'_1,\ldots,\zeta'_{k+1})$:
\begin{enumerate}[i.]
\item \label{step2i}
\emph{Find} $j$ such that $\alpha'_j=\underline{\alpha}'$. 
\emph{Choose} $p=j$ \emph{or} $p=j+1$ depending on whether 
$\alpha'_{j-1}<\alpha'_{j+1}$ or $\alpha'_{j-1}\ge\alpha'_{j+1}$. 
\item
\emph{If} $\mu(\{\zeta'_1,\ldots,\zeta'_{k+1}\}\setminus\{\zeta'_p\})<\mu(\zeta_1,\ldots,\zeta_k)$:
\begin{enumerate}[a.]
\item Update $\Xi_\zeta:=\{\zeta,\zeta_1,\ldots,\zeta_k\}=\{\zeta,\zeta'_1,\ldots,\zeta'_{k+1}\}\setminus\{\zeta'_p\}$.
\item \emph{If}  $\overline{\alpha}(\zeta_1,\ldots,\zeta_k) \le v\,
\underline{\alpha}(\zeta_1,\ldots,\zeta_k)$, \emph{then} STOP and return $\Xi_\zeta$.
\end{enumerate}
\end{enumerate}
\item \emph{If} $ i = m$: 
\begin{enumerate}%
\item[] 
Find the next $m$ nearest points  $\xi_{m+1},\ldots,\xi_{2m}$ in $\Xi\setminus\{\zeta\}$ to $\zeta$,
sorted by increasing distance to $\zeta$, and set $m:=2m$.
 \end{enumerate}
\item
Set $i := i+1.$
 \end{enumerate}
\end{enumerate}
\end{algorithm}

\subsubsection*{Remarks}

\begin{enumerate}[1.]

\item In the earlier version of our algorithm \cite[Algorithm 1]{DavyOanh11} we used a point
cloud with a fixed number ($m=30$) of nearest points considered for the inclusion into $\Xi_\zeta$, which sometimes led to the termination with rather non-uniform angles $\alpha_i$.
This was in part compensated by the removal of the `worst' point in this situation. In the improved
algorithm of this paper we dropped this removal step, and added %
the second termination criterion 
\eqref{term2} in Step~II.1, while allowing potentially any point of $\Xi$ to be selected into $\Xi_\zeta$. 
Theoretically, the algorithm may continue until all points in $\Xi$ are exhausted 
such that neither \eqref{term1} nor \eqref{term2} is satisfied. This however, never happened in our
experiments. Note that Algorithm~\ref{alg1} is supposed to work in the interaction with the adaptive
refinement according to Algorithm~\ref{alg2Refi}, and so its behavior may be quite different if
it is applied to a set $\Xi$ of a different nature.

\item The choice of $p$ in Step~II.2.i %
ensures that 
$$\mu\big(\{\zeta'_1,\ldots,\zeta'_{k+1}\}\setminus\{\zeta'_p\}\big)=
\min\Big\{\mu\big(\{\zeta'_1,\ldots,\zeta'_{k+1}\}\setminus\{\zeta'_j\}\big),\;
\mu\big(\{\zeta'_1,\ldots,\zeta'_{k+1}\}\setminus\{\zeta'_{j+1}\}\big)\Big\}.$$

\item For non-convex domains the points $\xi_i$ such that the segment between $\zeta$
and $\xi_i$ intersects the domain boundary are ignored in Step~I and Step~II.3.

\item \label{bsel}
Due to certain oversampling of the boundary by Algorithm~\ref{alg2Refi}, see Remark 
no.~\ref{oversamplingb} after it, many of the nearest centers to $\zeta$ may lie on a straight
segment of the boundary. Since these points are not desirable in $\Xi_\zeta$, we switch off the
termination in Step~II.1 in case when the number of such points in the current set 
$\Xi_\zeta$ is more than 3.  However, this rule is not applied if $i\ge 50$, in order to avoid 
inclusion of very far points. %

\end{enumerate}
\medskip

\begin{figure}[htbp!]
 \begin{center}
  \subfigure[]
{\includegraphics[width=7.5cm]{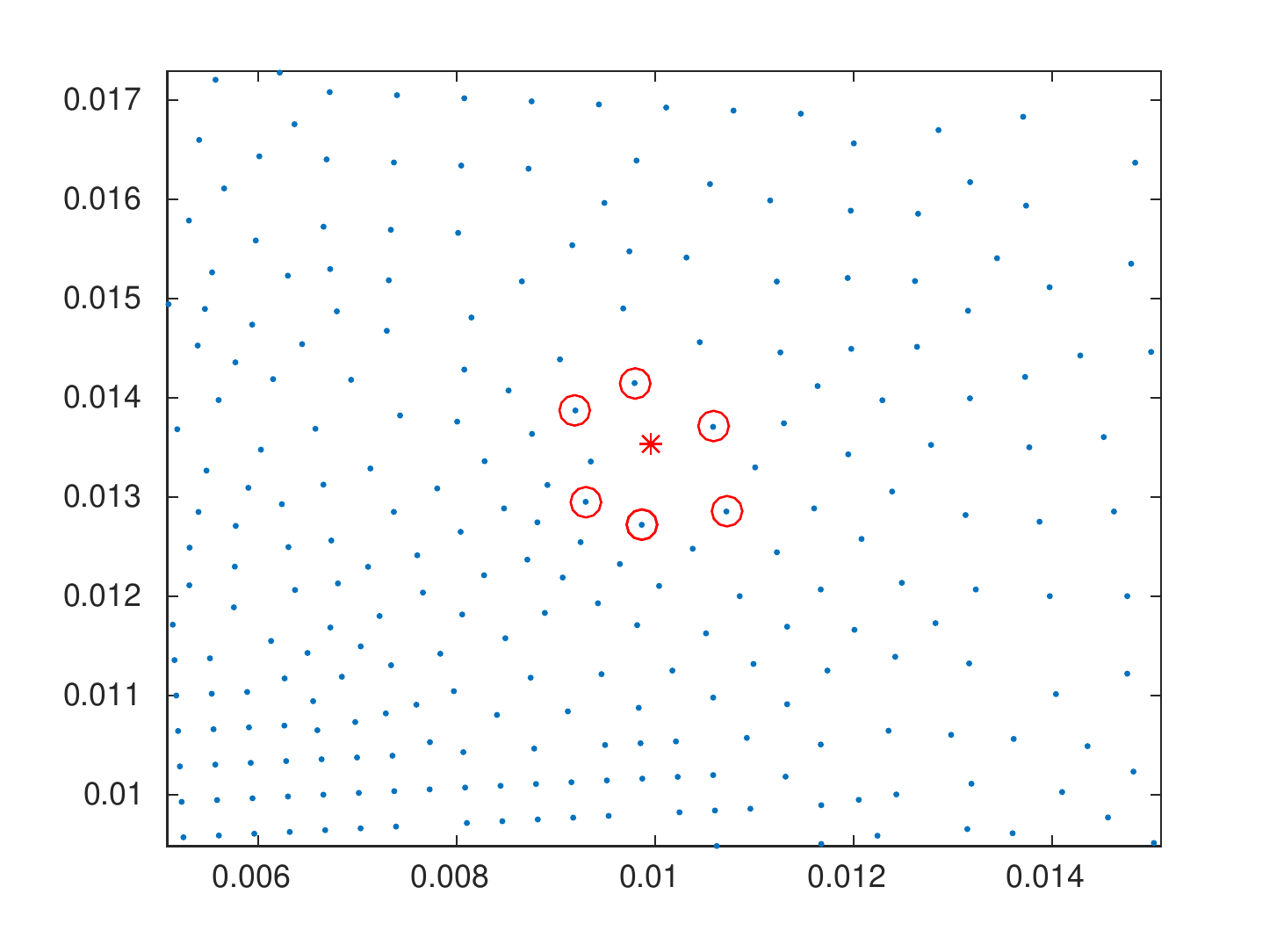}} %
 \subfigure[]
{\includegraphics[width=7.5cm]{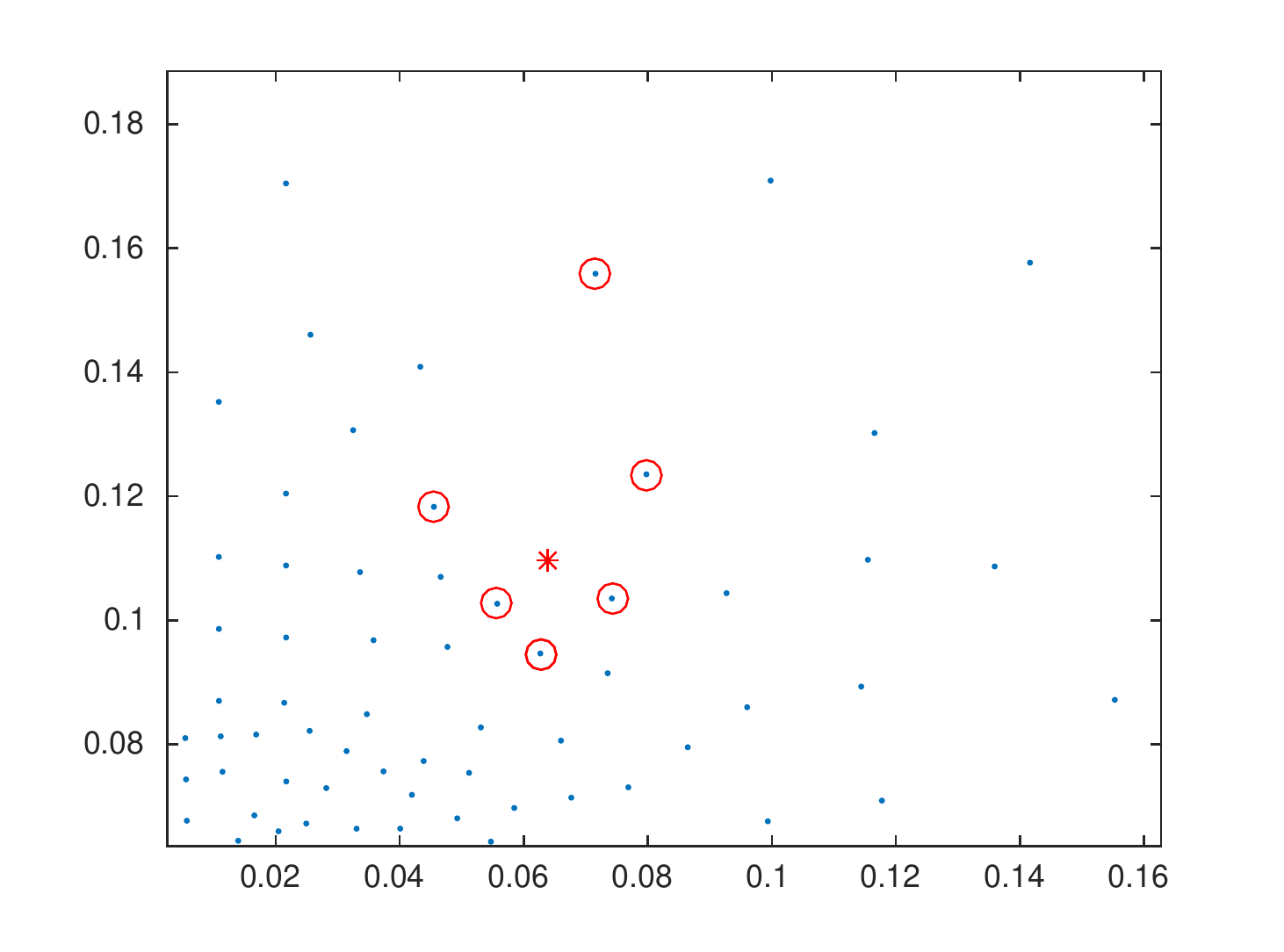}} 
 \subfigure[]
{\includegraphics[width=7.5cm]{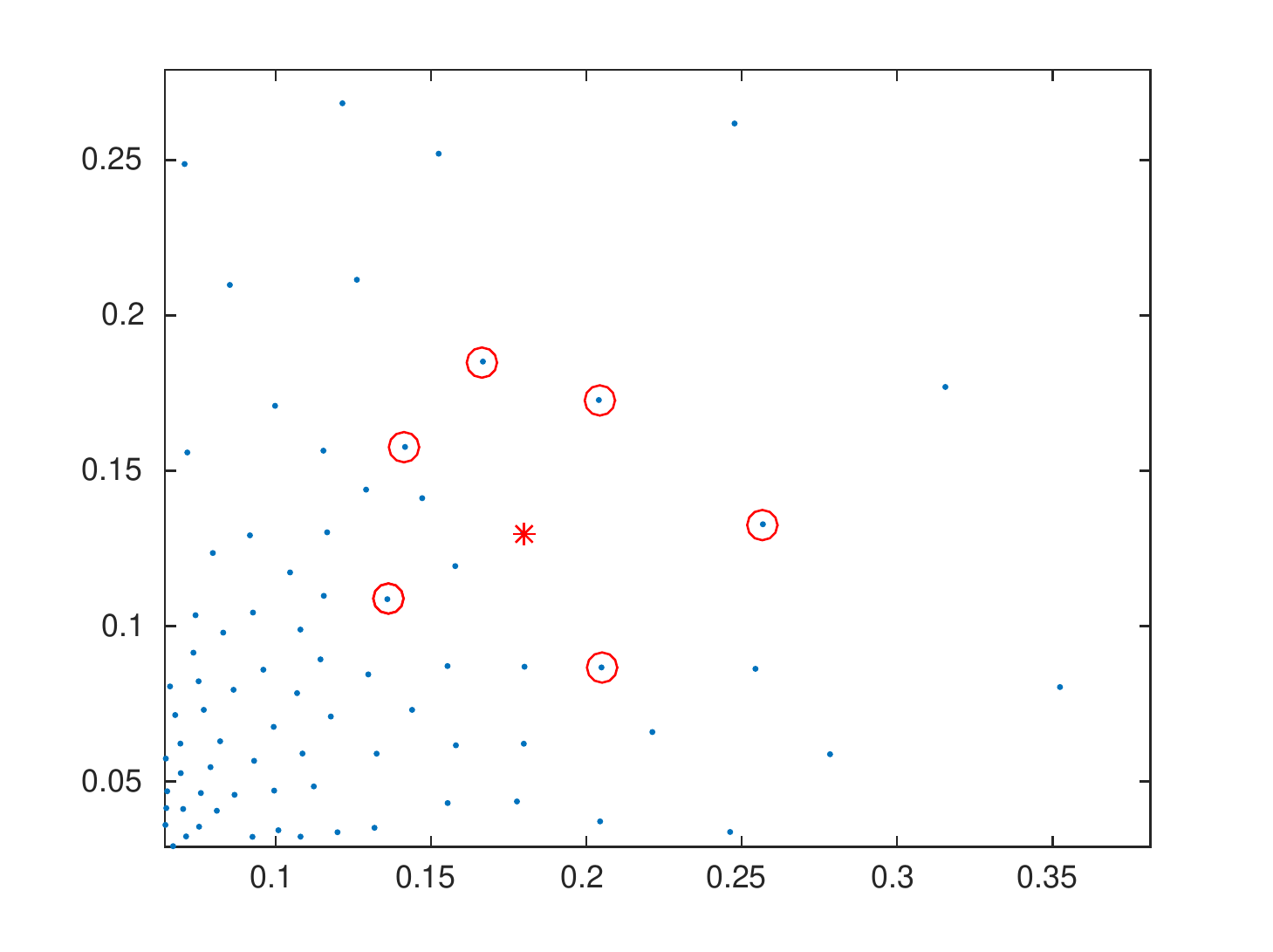}} 
 \subfigure[]
{\includegraphics[width=7.5cm]{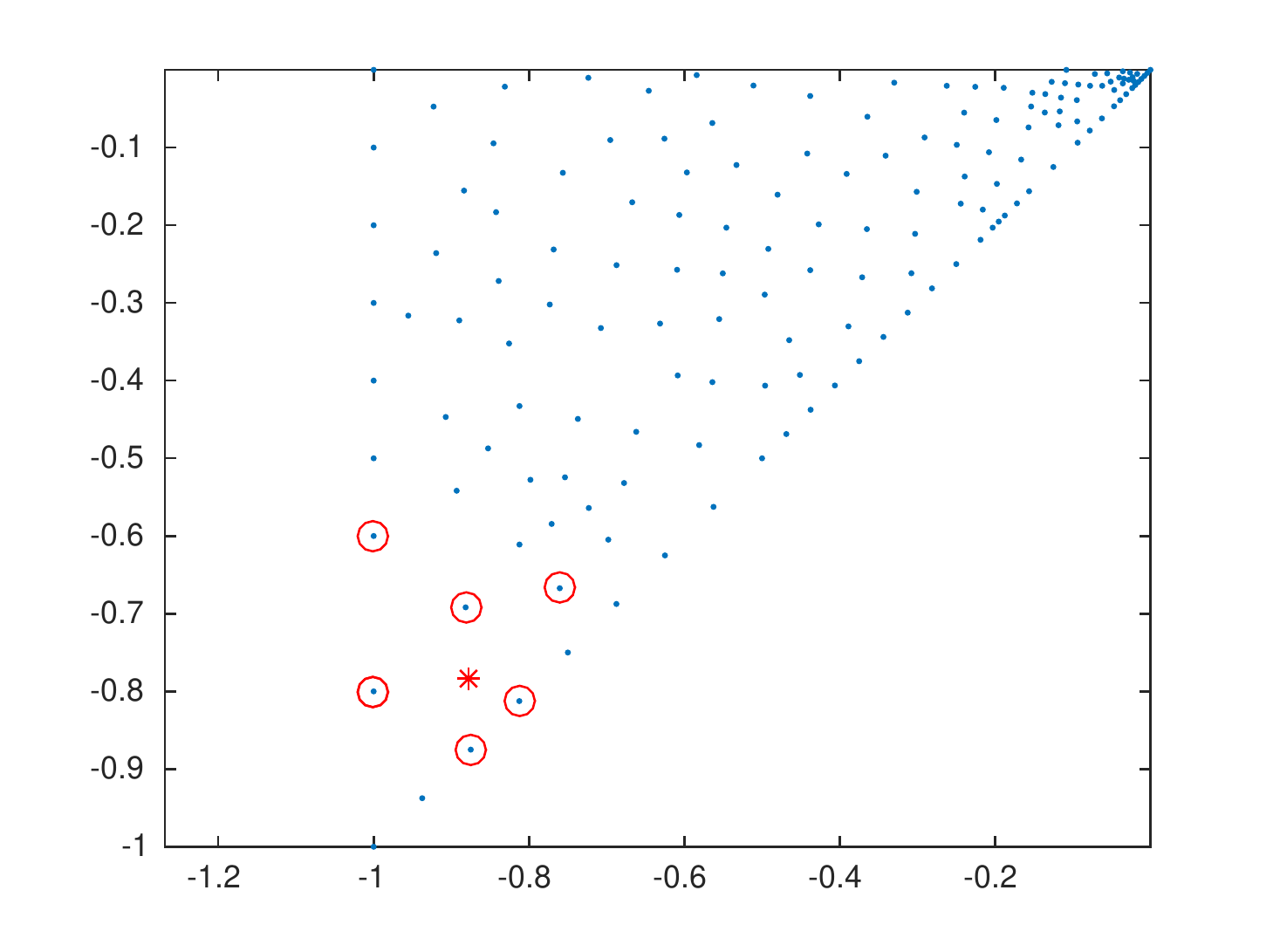}} 
 \subfigure[]
{\includegraphics[width=7.5cm]{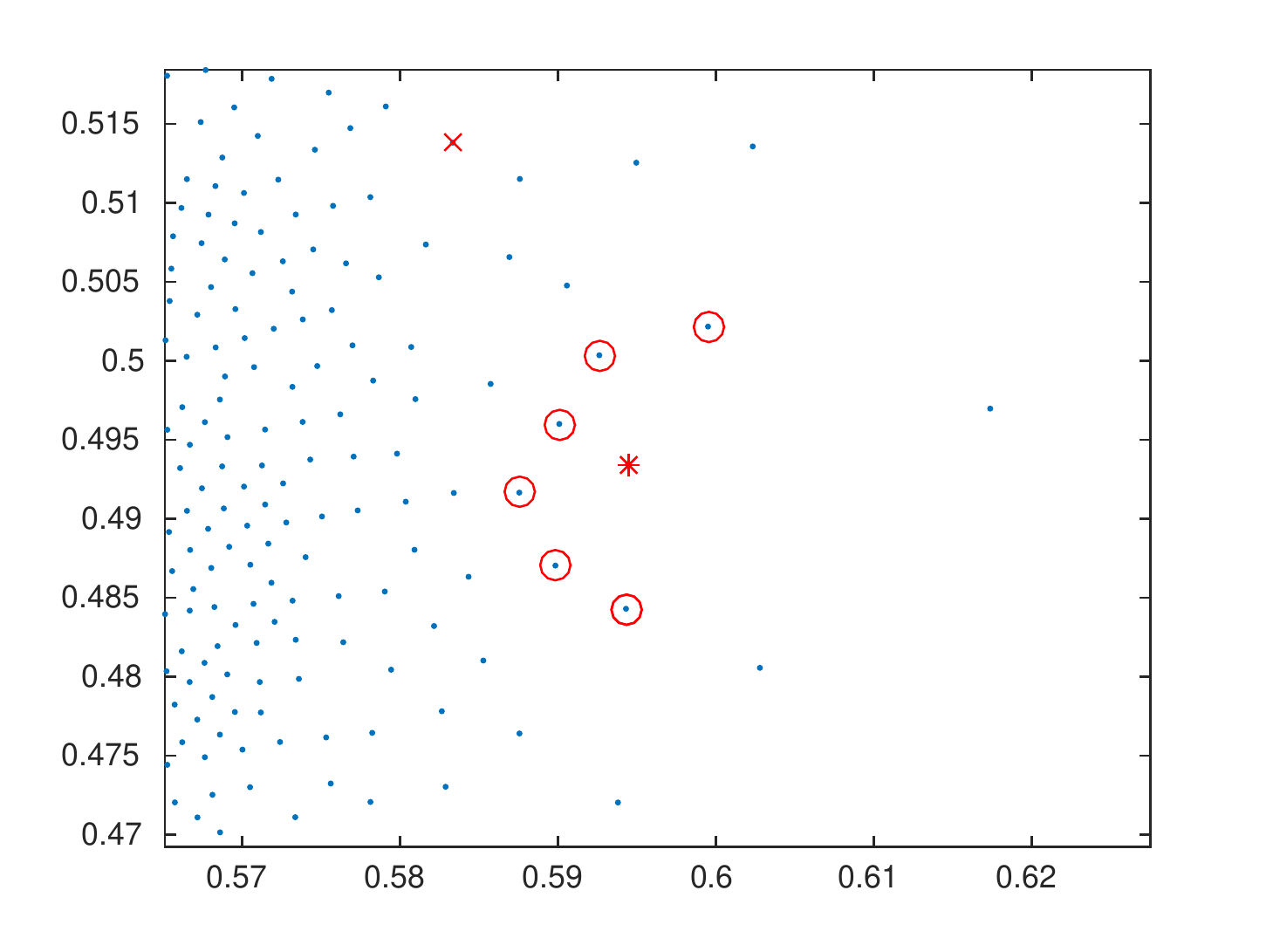}} 
 \subfigure[]
{\includegraphics[width=7.5cm]{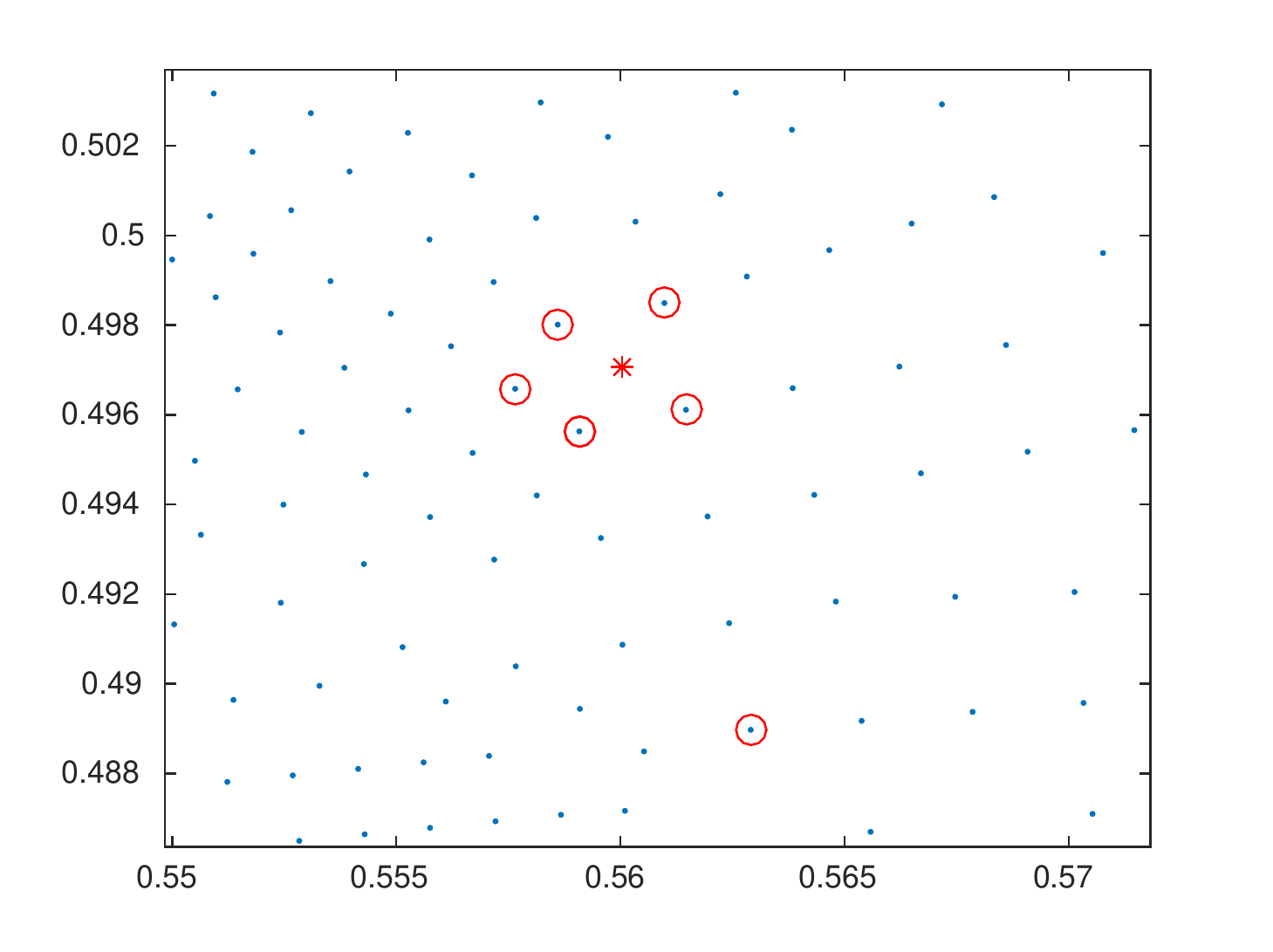}} 
        \end{center}
\caption{Stencil supports $\Xi_\zeta$ obtained by 
Algorithm~\ref{alg1}: The star shows the position of
$\zeta$ and the circles the positions of $\xi_1\ldots,\xi_6$. 
}
\label{typ_stencils}
\end{figure}

Figures~\ref{typ_stencils}(a--d) show typical stencil supports  $\Xi_\zeta$ obtained by 
Algorithm~\ref{alg1} in our experiments. The choice of parameters $c$ and $v$ determines the trade-off
between improving the uniformity of the angles $\alpha_i$ and avoiding far away points.
Figure~\ref{typ_stencils}(e) shows an example of $\Xi_\zeta$ for which Algorithm~\ref{alg1} terminates at Step~II.1. The red cross indicates
the center $\xi_i$ for which the termination condition \eqref{term2} is satisfied ($i=52$), and hence the improvement of the angles is
terminated. In fact $\overline{\alpha} /\underline{\alpha}=4.3343>v$ in this case.
Figure~\ref{typ_stencils}(f) gives an example of the other extreme, where the termination is at Step~II.2.ii.b so that
$\overline{\alpha} /\underline{\alpha}=2.4410<v$, but this is achieved at the expense of including in $\Xi_\zeta$ a rather distant center 
(the circled dot near the bottom).

In order to quantify the uniformity of the stencil supports produced by Algorithm~\ref{alg1}, we
define several measures as follows. Given a set of stencil supports $\Xi_\zeta=\{\zeta,\zeta_1,\ldots,\zeta_n\}$, we denote by
$v_\mathrm{max}$ the maximum and by $v_\mathrm{aver}$ the average values of the quotient
$\overline{\alpha}/\underline{\alpha}$ over all $\Xi_\zeta$ in this set. Similarly, 
$c_\mathrm{max}$ and $c_\mathrm{aver}$ denote the  maximum and the average values of the quotient
\begin{equation}\label{duq}
\max_{j=1,\ldots,n}\|\zeta-\zeta_j\| \Big/  \frac{1}{2\,n} \sum_{j=1}^n \big(\|\zeta_j-\zeta\| + 
\|\zeta_j-\zeta_{j+1}\|\big).
\end{equation}
The values of these measures in our experiments are provided in Table~\ref{angl_quo} at the end 
of Section~4, see there a discussion of the results as well.

\section{Refinement method}
\label{Refi}
Before giving a formal description of our algorithm we discuss its main features and changes in
comparison to \cite[Algorithm 2]{DavyOanh11}.

\bigskip

\noindent
{\bf Error indicator.}
Assuming that an approximate discrete solution $\hu$ of  the Dirichlet problem (\ref{poi})
is determined via \eqref{dird} with some $\Xi$ and stencil supports $\Xi_\zeta$ for all 
$\zeta\in\Xi_{\rm int}$, we choose an  \emph{error indicator} $\eps(\zeta,\xi)$ 
associated with each `edge' $\zeta\xi$ for all $\zeta\in\Xi_{\rm int}$ and
$\xi\in\Xi_\zeta\setminus\{\zeta\}$. In 
our earlier refinement algorithm \cite[Section 6]{DavyOanh11} we used 
$\eps(\zeta,\xi)=\eps_0(\zeta,\xi):=|\hu_\zeta-\hu_\xi|$.
However, this indicator identifies the areas where the gradient of the solution is big, which results
in some cases, see Figure~\ref{IndicF441} and related discussion in Section~4, in oversampling relatively 
flat regions and undersampling the regions of high curvature, which leads to sub-optimal solution. 
Therefore we replace it by an indicator of Zienkiewicz-Zhu type that estimates the error of the
approximation of the directional derivative along the edge $\zeta\xi$.
The error indicator used in this paper is defined as follows. For each $\zeta\in\Xi_{\rm int}$,
 let $\ell_\zeta(x)=a+b^T(x-\zeta)$ be the linear polynomial that fits the data 
$\{(\xi,\hu_\xi):\xi\in \Xi_\zeta\}$ in the least squares sense, that is its coefficients
$a\in\mathbb{R}$, $b\in\mathbb{R}^2$ are chosen such that the sum
$$
\sum_{\xi\in\Xi_\zeta}|\hu_\xi-\ell_\zeta(\xi)|^2$$
is minimized. (This minimization problem has a unique solution as soon as $\Xi_\zeta$ is not a subset
of a straight line.) We set
\begin{equation}\label{errind}
\eps(\zeta,\xi)=\eps_1(\zeta,\xi):=|(\hu_\zeta-\hu_\xi)-(\ell_\zeta(\zeta)-\ell_\zeta(\xi))|,\quad 
\zeta\in\Xi_{\rm int},\;\xi\in \Xi_\zeta\setminus\{\zeta\}.
\end{equation}
We can interprete $\eps_1(\zeta,\xi)$ as an edge-based indicator of averaging type.
Indeed, $(\hu_\zeta-\hu_\xi)/\|\zeta-\xi\|$ is an estimate of the directional derivative
of the solution based only on two points $\zeta,\xi$, whereas 
$(\ell_\zeta(\zeta)-\ell_\zeta(\xi))/\|\zeta-\xi\|$ is an estimate of the same
directional derivative based on averaging the information from all points in $\Xi_\zeta$. 
In the finite element method such error indicators are usually based on averaging the 
gradient over elements, see \cite{ZZ87,PV00}. 
Instead, in our meshless setting we naturally resort to directional derivatives
and edges. 

\bigskip

\noindent
{\bf Marking strategy.} The refinement of $\Xi$ is achieved by inserting one or more 
new centers in the vicinity of each edge $\zeta\xi$ marked for refinement. An overview of marking methods employed
in the finite element adaptive refinement algorithms can be found in \cite{PV00}. As in 
\cite{DavyOanh11} we follow the so called `maximum strategy' applied to edges rather than elements. 
Therefore the main condition for the edge $\zeta\xi$ to be refined is 
\begin{equation}\label{epsbar}
\eps(\zeta,\xi)\geq \gamma \bar\eps(\Xi), \quad\text{where} \;
\bar\eps(\Xi):=\max\{\eps(\zeta,\xi):\,\zeta\in\Xi_{\rm int},\,\xi\in\Xi_\zeta\},
\end{equation}
and $\gamma\in(0,1)$ is a prescribed tolerance. We also follow the common practice of preventing that
too few new centers are generated in one refinement step by reducing the threshold for $\eps(\zeta,\xi)$ and marking 
further edges in the event that the number of the new interior centers is
less that certain percentage of the total number of centers in $\Xi_{\rm int}$, 
see Step~III of Algorithm~\ref{alg2Refi}.

\bigskip

\begin{figure}[htbp!]
 \begin{center}
\includegraphics[width=4.5cm]{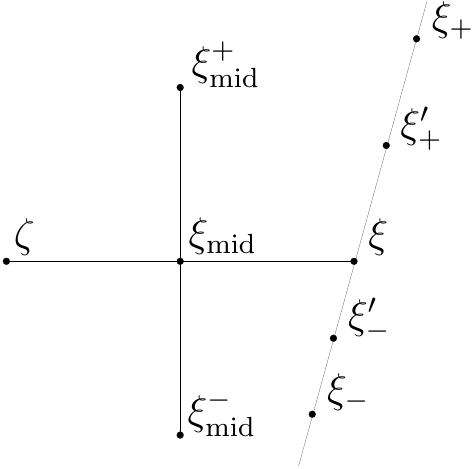} 
        \end{center}
\caption{Centers to be added in the vicinity of a marked edge. 
}
\label{refine}
\end{figure}

\noindent
{\bf Refinement in the vicinity of a marked edge.}  For each marked edge $\zeta\xi$, the candidate
new centers in \cite[Algorithm 2]{DavyOanh11} are the middle point $\xi_{\rm mid}=(\zeta+\xi)/2$ of the edge, as well as two
neighboring points on the boundary should $\xi$ be a boundary center. To make the local distribution of
centers more uniform and isotropic, we  now also consider the points 
$\xi_{\rm mid}^\pm:=\xi_{\rm mid}\pm d\bar\nu$, where $d=\|\zeta-\xi\|/2$ and $\bar\nu$ is the unit vector perpendicular to the edge
$\zeta\xi$, as candidate new centers, see Figure~\ref{refine}.

As in  \cite{DavyOanh11}, we only  proceed with the refinement if 
\begin{equation}\label{mu}
\dist(\xi',\Xi')\ge\mu\sep_{\xi'}(\Xi'),\quad\text{with a prescribed tolerance $\mu\in(0,1)$,}
\end{equation}
that is, if the insertion of a new center $\xi'$ into the current set $\Xi'$ does not significantly reduce the 
\emph{local separation} defined as
\begin{equation}\label{sep}
\sep_{\xi'}(\Xi'):=\frac{1}{4}\sum_{i=1}^4\dist(\xi_i,\Xi'\setminus\{\xi_i\}),
\end{equation}
where $\xi_1,\ldots,\xi_4$ are the four closest points in $\Xi'$ to $\xi'$, and 
$\dist(x,Y):=\inf\{\|x-y\|:\,y\in Y\}$
is the distance from a point $x$ to a set $Y$. Moreover, we also check that
$\xi'\notin\partial\Omega$ is not placed too close to the boundary, 
see Steps II.2 and II.3.i of Algorithm~\ref{alg2Refi}.

\begin{algorithm}\label{alg2Refi}\rm {\bf Adaptive meshless refinement} \\
\emph{Input:} The set of centers $\Xi$ and stencil supports $\{\Xi_\zeta:\,\zeta\in\Xi_{\rm int}\}$. \\
\emph{Output:} The refined set of  centers $\Xi'$. \\
\emph{Parameters:}   $\gamma = 0.5$ (error indicator tolerance),
$\mu = 0.8$  (separation tolerance) and $n=15$ (percentage of added centers). 
\begin{enumerate}[I.] 
\item\label{2step1} Compute the error indicator threshold $\bar\eps=\gamma\bar\eps(\Xi)$ 
and initialize  $\Xi':=\Xi$.
\item \label{2step2}
\emph{For each} edge $\zeta\xi$, $\zeta\in\Xi_{\rm int}$, $\xi\in\Xi_\zeta\setminus\{\zeta\}$,
such that $\eps(\zeta,\xi)\ge\bar\eps$: %
\begin{enumerate}[1.]

\item Compute $\xi_{\rm mid}:=(\zeta+\xi)/2$, $\xi_{\rm mid}^+:=\xi_{\rm mid}+ d\bar\nu$ and
$\xi_{\rm mid}^-:=\xi_{\rm mid}- d\bar\nu$, where $d:=\|\zeta-\xi\|/2$ and $\bar\nu$ 
is the unit vector perpendicular to the edge $\zeta\xi$.\\ Initialize  $\Xi_C:=\emptyset$.

\item \emph{If} $\xi\in\Xi_{\rm int}$, 
\emph{then for each} $\xi'\in\{\xi_{\rm mid},\xi_{\rm mid}^+,\xi_{\rm mid}^-\}$:
\begin{enumerate}[a.] 
\item[] \emph{If} $\dist(\xi',\partial\Omega) \ge d/2$ \emph{and}  
$\dist(\xi',\Xi')\ge\mu\sep_{\xi'}(\Xi')$, \emph{then} set  $\Xi_C:=\Xi_C\cup\{\xi'\}$.
\end{enumerate}

\item \emph{ElseIf} $\xi\in\partial\Xi$:
\begin{enumerate}[i.] 
\item \emph{For each} $\xi'\in\{\xi_{\rm mid},\xi_{\rm mid}^+,\xi_{\rm mid}^-\}$:\\
\emph{If} $\dist(\xi',\partial\Omega) \ge d/2$ \emph{and}  
$\dist(\xi',\Xi')\ge d/2$, \emph{then} set  $\Xi_C:=\Xi_C\cup\{\xi'\}$.
\item\label{2stepII3b} \emph{If} $\Xi_C \neq \emptyset$ \emph{or} $\dist(\xi_{\rm mid},\partial\Omega) < d/2$:\\
Find two neighbors $\xi_-,\xi_+$ of $\xi$ in $\partial\Xi$, 
 one in each direction from $\xi$ along the boundary, and 
compute two middle points
$\xi'_-,\xi'_+\in\partial\Omega$  defined by the pairs $\xi,\xi_-$ 
and $\xi,\xi_+$, respectively. Set  $\Xi_C:=\Xi_C\cup\{\xi'_+, \xi'_-\}$.
\end{enumerate}

\item Set $\Xi':=\Xi'\cup\Xi_C$.
\end{enumerate}
\item \emph{If} the number of centers in $\Xi'_{\rm int}\setminus\Xi_{\rm int}$ is less than $n$\% 
of the number of centers in $\Xi_{\rm int}$,
\emph{then} set $\bar\eps:= \gamma\bar\eps$ and \emph{goto} Step II.\\
\emph{Else} STOP and return $\Xi'$.

\end{enumerate}

\end{algorithm}

\subsubsection*{Remarks}

\begin{enumerate}[1.]
\item Algorithm~\ref{alg2Refi} is applied recursively starting from an initial
non-adaptive set of centers $\Xi$.   In our experiments in Section 4 it is the set of
centers of the initial finite element mesh created by MATLAB PDE Toolbox with default parameters.

\item The middle points in Step~II.3.ii are found using the parameterizations of the respective boundary
components connecting $\xi$ with either $\xi_-$ or $\xi_+$. If the middle point $\xi'_-$ or $\xi'_+$
is already in $\Xi'$, it is not added again, which is easy to avoid by keeping track of the pairs of
centers on the boundary that have been refined.

\item \label{oversamplingb}
Generally, our algorithm leads to a slight oversampling of the boundary, which can be seen in the
numerical results in Section 4. This is generally harmless because the boundary centers do not bear any
degrees of freedom and hence the size of the system matrix does not increase. However, this leads
to certain difficulties for Algorithm~\ref{alg1}, see Remark no.~\ref{bsel} after it. To reduce the
oversampling, in the experiments below we avoid adding one or both points $\xi'_-,\xi'_+$ in Step~II.3.ii in
the following situations:
\begin{itemize}
\item $\xi'_-$ is not added to $\Xi_C$ if $\|\xi'_--\xi_{\rm mid}\| \ge\|\xi'_+-\xi_{\rm mid}\| $, 
$\|\xi'_--\xi\| \le\min\{d,2\|\xi'_+-\xi\|\} $ and 
$\|\xi-\xi'_-\|  + \|\xi-\xi'_+\| \le 2 \|\xi'_+ -\xi'_-\| $.
\item $\xi'_+$ is not added to $\Xi_C$ if $\|\xi'_+-\xi_{\rm mid}\| \ge\|\xi'_--\xi_{\rm mid}\| $, 
$\|\xi'_+-\xi\| \le\min\{d,2\|\xi'_--\xi\|\} $ and 
$\|\xi-\xi'_-\|  + \|\xi-\xi'_+\| \le 2 \|\xi'_+ -\xi'_-\| $.
\end{itemize}

\item 
To get faster convergence of the numerical solution for easier problems, we enforce a reduction of the 
error indicator threshold $\bar\eps$ between subsequent applications of Algorithm~\ref{alg2Refi}. Let 
$\bar\eps_{\rm prev}$ be the threshold used for the previous refinement, and let $\bar\eps$ be the value
computed at Step~I of the current refinement. If $\bar\eps_{\rm prev}< \bar\eps$,
then we replace $\bar\eps$  by $\bar\eps:=\bar\eps_{\rm prev}/2$, and proceed to Step~II as usual. 
This procedure is used in the numerical experiments for Test Problems \ref{f4}, \ref{f8} and \ref{f42} with 
$\omega=\pi+0.01, 5\pi/4, 7\pi/4$. Moreover, in this case we use $n=5$ instead of the default $n=15$
which helps to get smoother error curves. The value $n=5$ for the percentage of the added centers 
is also used for Test Problem~\ref{f42} with $\omega=2\pi$ and Test Problem~\ref{curved_slit}. %

\item Note that the current algorithm no longer requires an (expensive) post-processing as in Step~III of  
\cite[Algorithm 2]{DavyOanh11}.

\end{enumerate}

\section{Numerical results}
To illustrate the performance of the improved algorithms, we consider a number of benchmark test
problems where adaptively refined centers are known to be of great advantage, and compare our 
results with those obtained with the help of the PDE Toolbox~\cite{PDEtool} 
using the finite element method with piecewise linear shape functions and default parameters of the adaptive
refinement. This comparison is fair because the density of the system matrices resulting from our method is close to the
density of the system matrices of such finite element method thanks to the fact that we choose $k=6$ points in the stencil supports of
Algorithm~\ref{alg1}. We refer to \cite{DavyOanh11} for a detailed discussion and numerical results
on the density of the system matrices.

We use the RBF-FD weights $w_{\zeta,\xi}$, $\xi\in\Xi_\zeta$, with
Gaussian function \eqref{G} described in the introduction.  Our investigation in \cite{DavyOanhSP} has
shown that small values of the shape parameter $\eps$ in \eqref{G} lead to reliable, albeit
not always optimal results. For simplicity we use a fixed small value $\eps =10^{-5}$ and compute the
weights $w_{\zeta,\xi}$ by the Gauss-QR method as described in \cite{DavyOanhSP}, 
see also \cite{FLP13}. Note that in \cite{DavyOanh11,DavyOanhSP} we added a constant term to the
Gaussian form to ensure that the scheme is exact for constants. This however turns out to be
insignificant for the performance of the method as can be seen from the numerical results in this
paper. The constant term causes certain technical complications for the computation of the weights by  
the Gauss-QR method, see \cite{DavyOanhSP}.

As a measure of the quality of the results  we use the root mean square (rms) error $E_c$ of the 
solution $\hu$ of (\ref{dird}) against the exact solution $u$ of (\ref{poi}) on the
centers,
\begin{equation}\label{rmser}
E_c=\Big(\frac{1}{N}\sum_{\zeta\in\Xi_{\rm int}}(\hu_\zeta-u(\zeta))^2\Big)^{1/2}.
\end{equation} 
where $N=\#\Xi_{\rm int}$ is the number of interior discretization centers. For the finite element
method the set $\Xi_{\rm int}$ consists of all interior vertices of the triangulation. Since these
sets are different for both methods, we also compute the rms error $E_g$  on a fixed uniform
grid with step size $0.001$, where the values of the 
approximate solution $\hu$ at the grid points are obtained by evaluating the piecewise linear interpolant
with respect to the Delaunay triangulation of the centers. 
  
\medskip
	 
The first two test problems have	been considered in \cite{DavyOanh11}, see Test Problems 2 and 3 of that
paper.

\begin{testproblem}\label{f4}\rm 
Laplace equation $\Delta u=0$ in the circle sector $\Omega$ given by the inequalities $r<1$, $-3\pi/4<\varphi< 3\pi/4$ 
in polar coordinates, with Dirichlet
boundary conditions defined by $u(r,\varphi)=\cos (2\varphi/3)$ along the arc, and $u(r,\varphi)=0$ along the
straight lines.
The exact solution is $u(r,\varphi)=r^{2/3}\cos (2\varphi/3)$.
\end{testproblem}
 
\begin{testproblem}\label{f8}\rm 
Laplace equation $\Delta u=0$ in the domain $\Omega=(0.01,1.01)^2$ with Dirichlet boundary conditions chosen
such that the exact solution is $u(x,y)=\log(x^2+y^2)$.
\end{testproblem}

\begin{figure}[htbp!]
 \begin{center}
  \subfigure%
{\includegraphics[width=6cm,height=4cm]{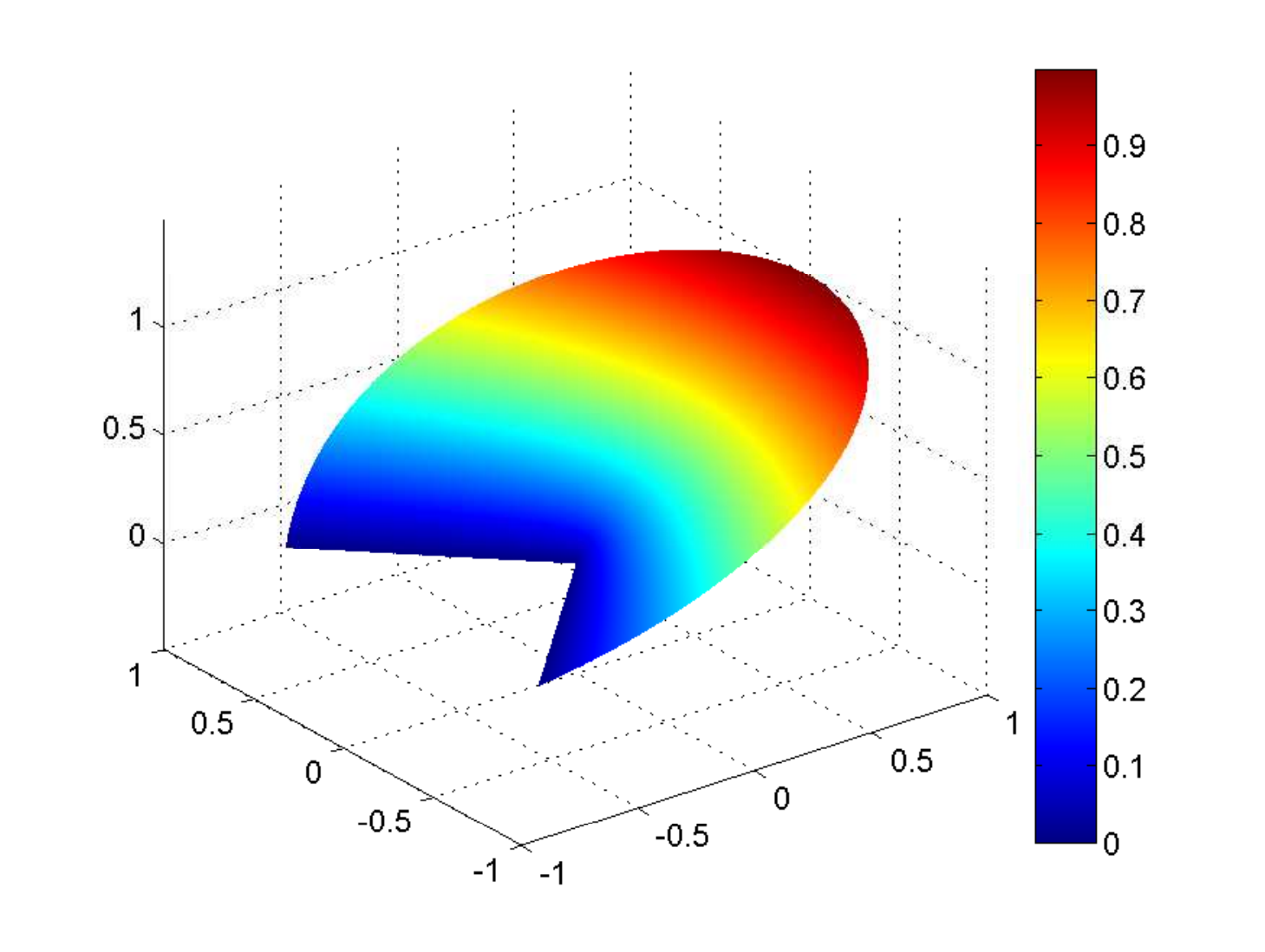}} \qquad
 \subfigure%
{\includegraphics[width=6cm,height=4cm]{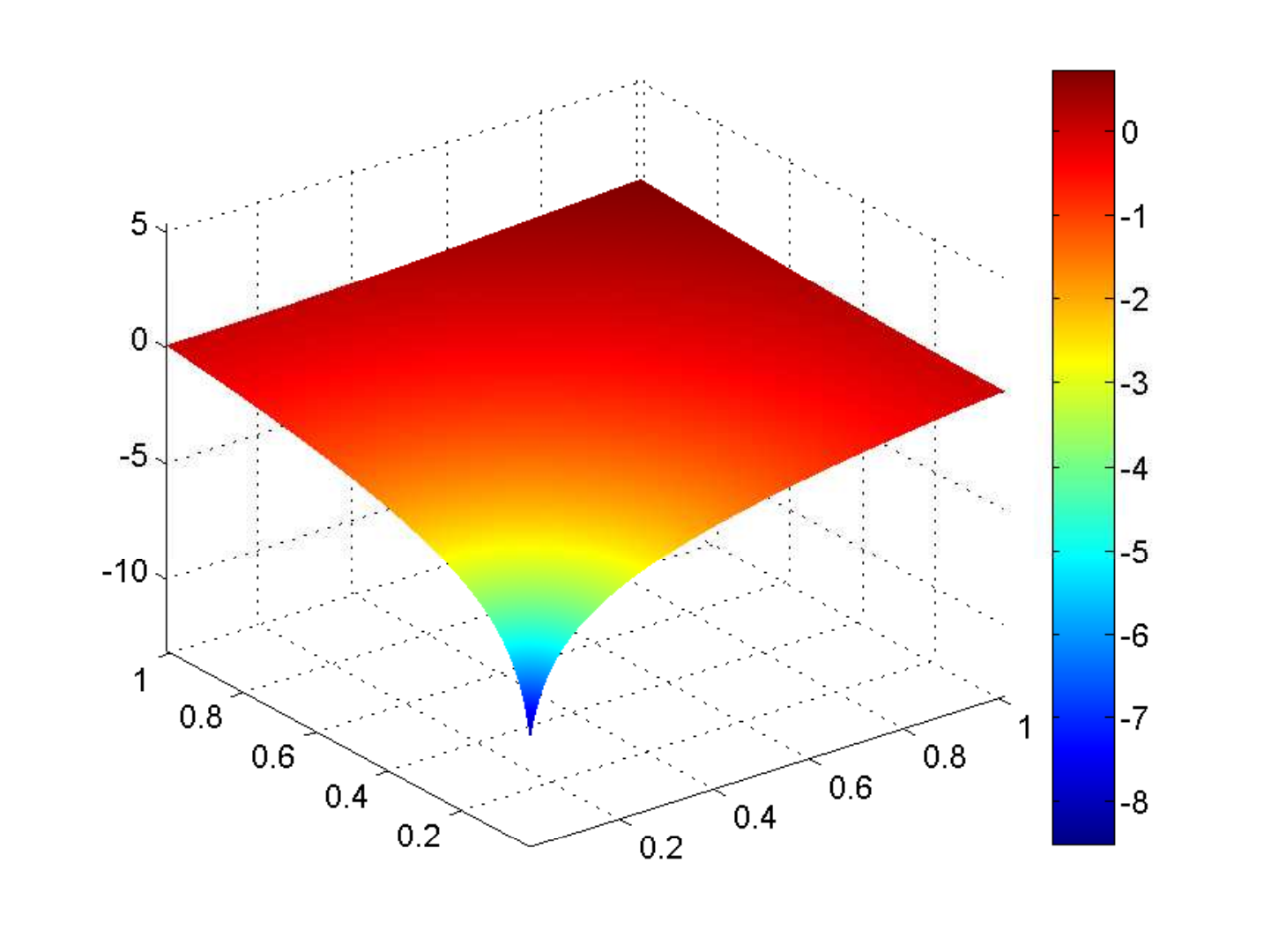}} 
        \end{center}
\caption{Exact solutions of Test Problems~\ref{f4} (left) and \ref{f8} (right). 
}
\label{F4F8exact}
\end{figure}

\begin{figure}[htbp!]
 \begin{center}
 \subfigure[Errors on centers]%
{\includegraphics[width=6cm,height=4cm]{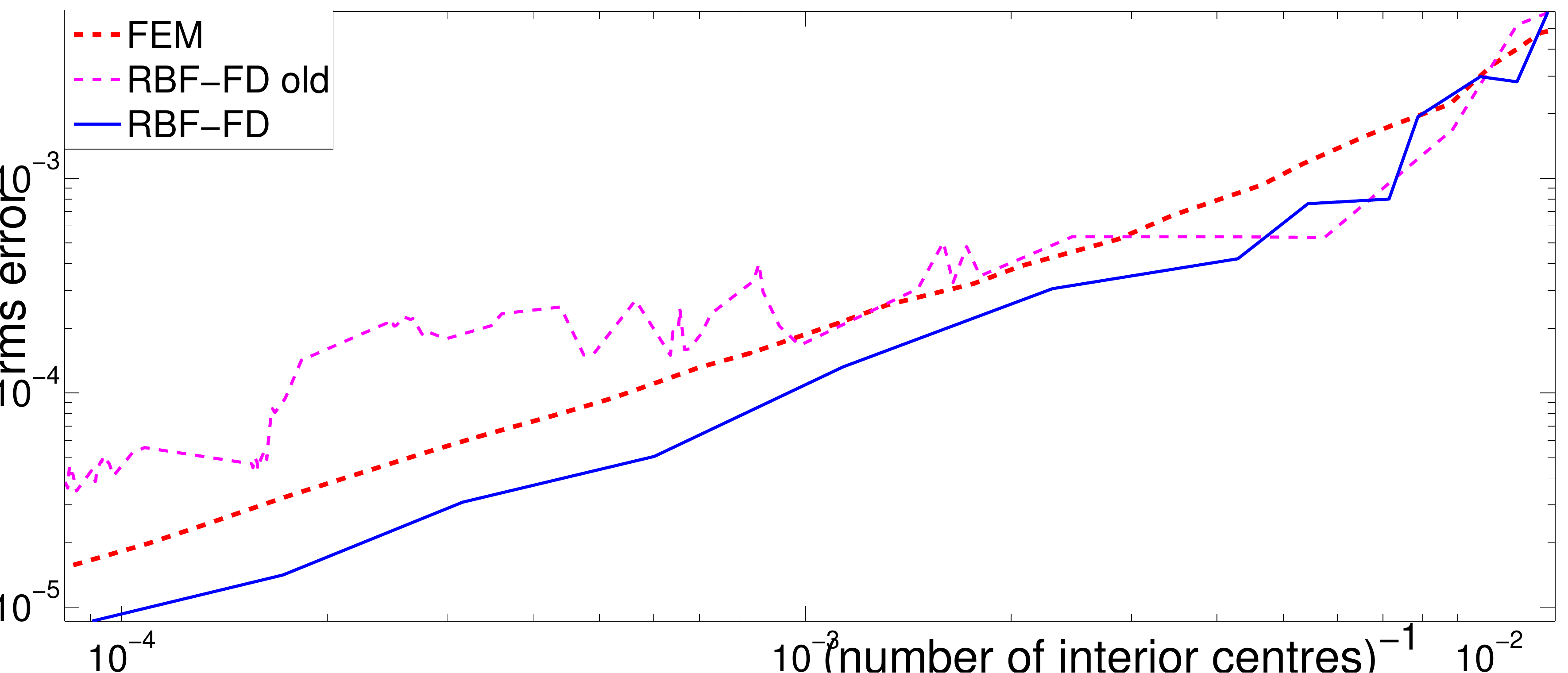}} \qquad 
 \subfigure[Errors on grid]
{\includegraphics[width=6cm,height=4cm]{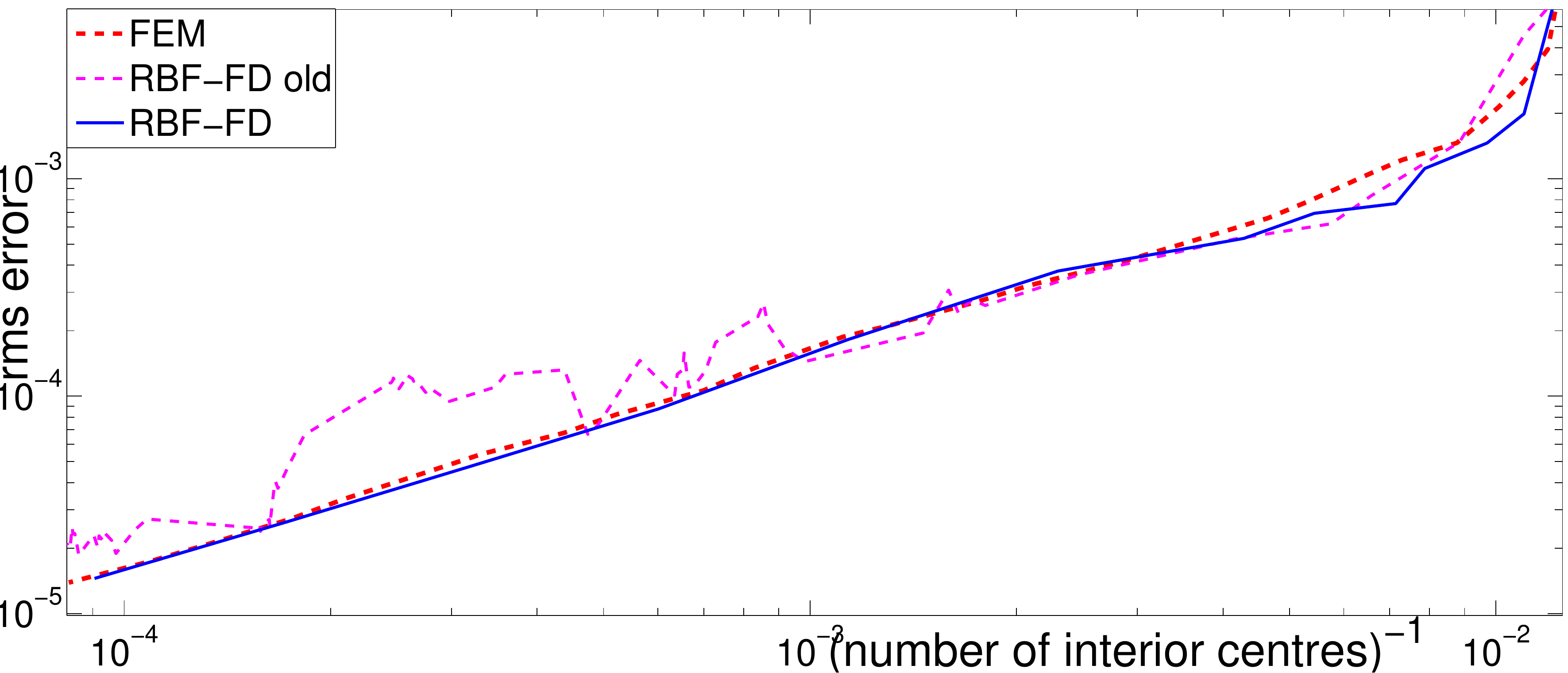}} 
 \subfigure[RBF-FD error]
{\includegraphics[width=6cm,height=4cm]{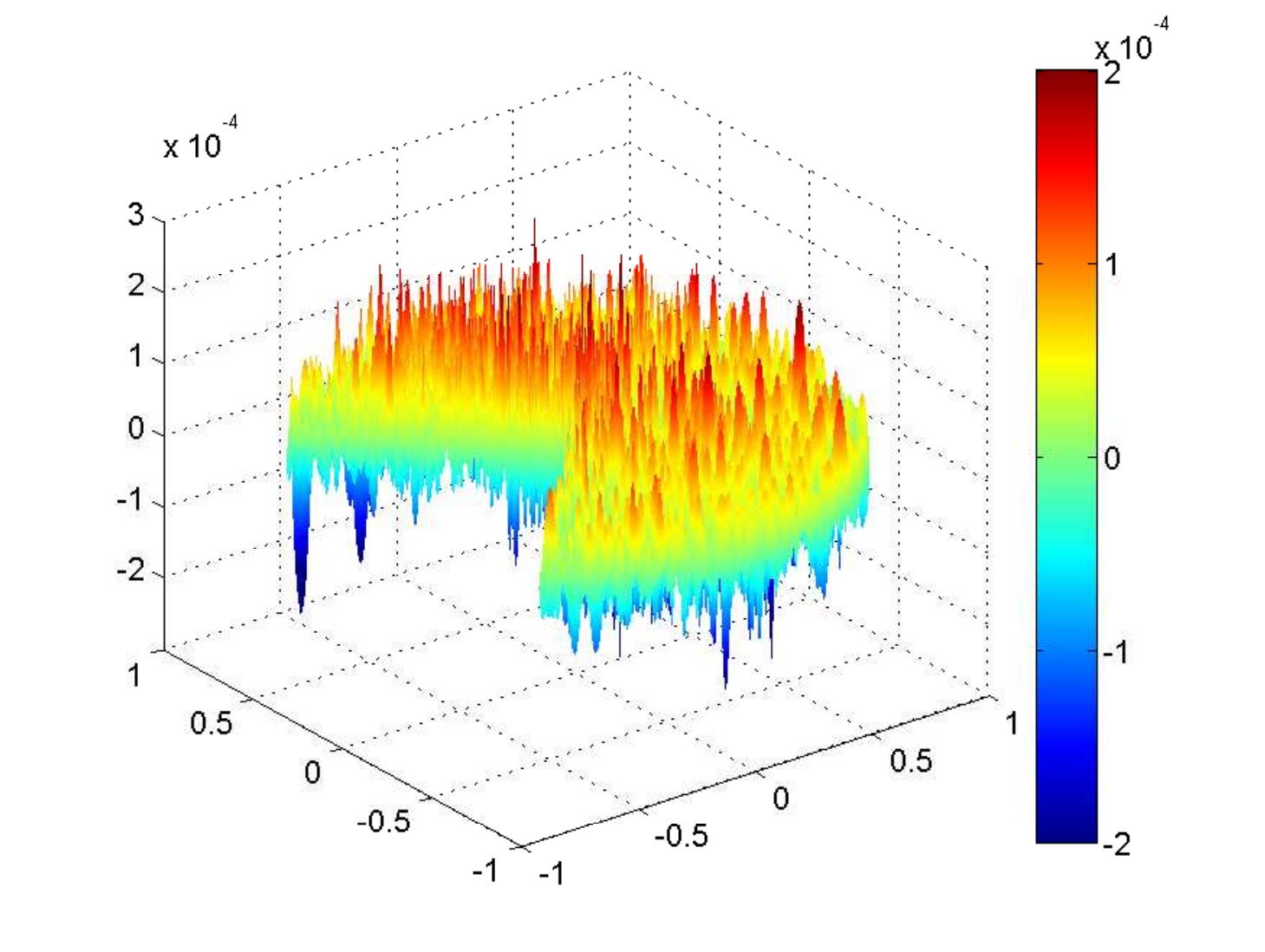}}\qquad 
 \subfigure[FEM error]
{\includegraphics[width=6cm,height=4cm]{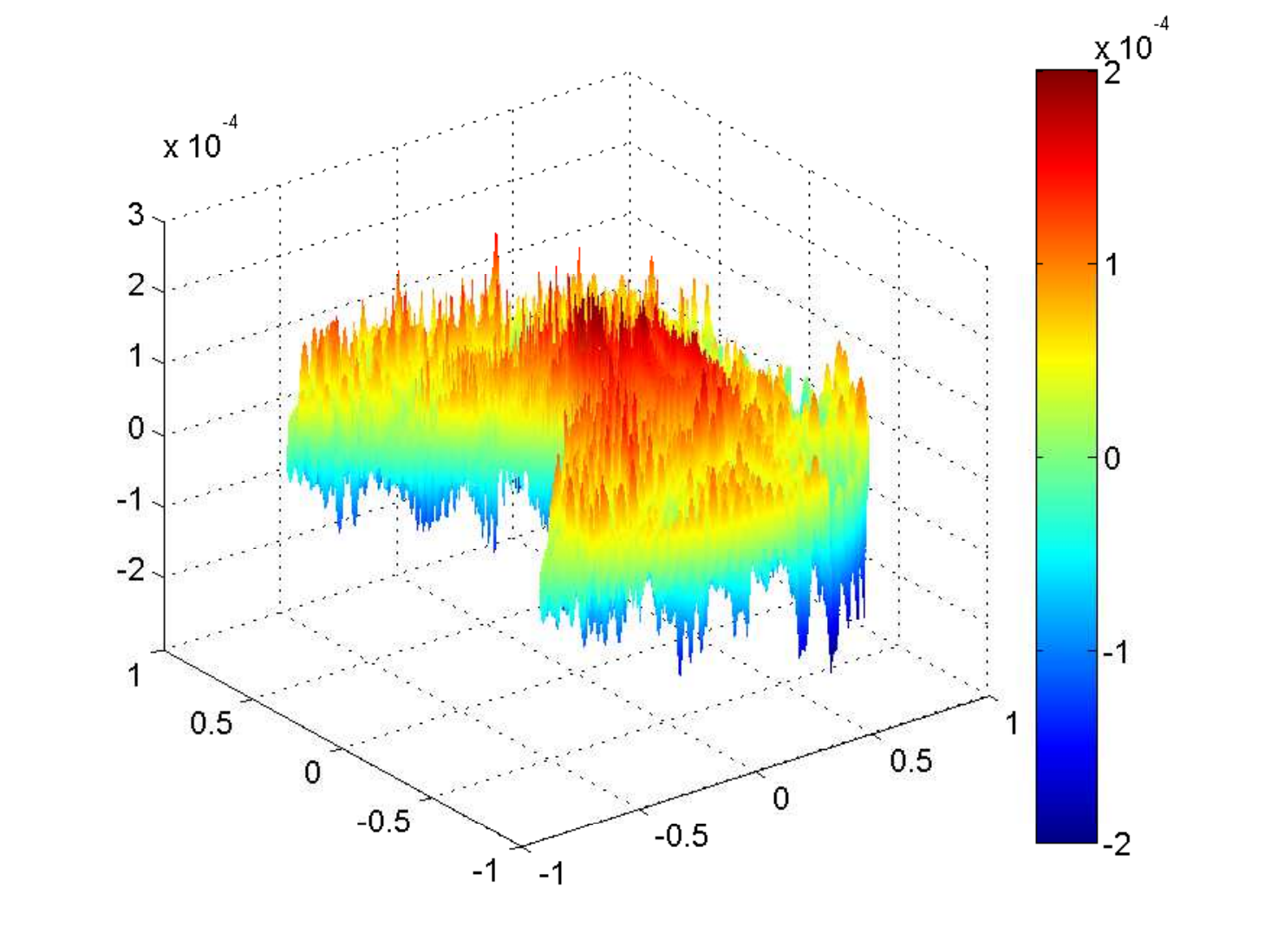}} 
 \subfigure[RBF-FD centers (3169)]
{\includegraphics[width=6cm,height=4cm]{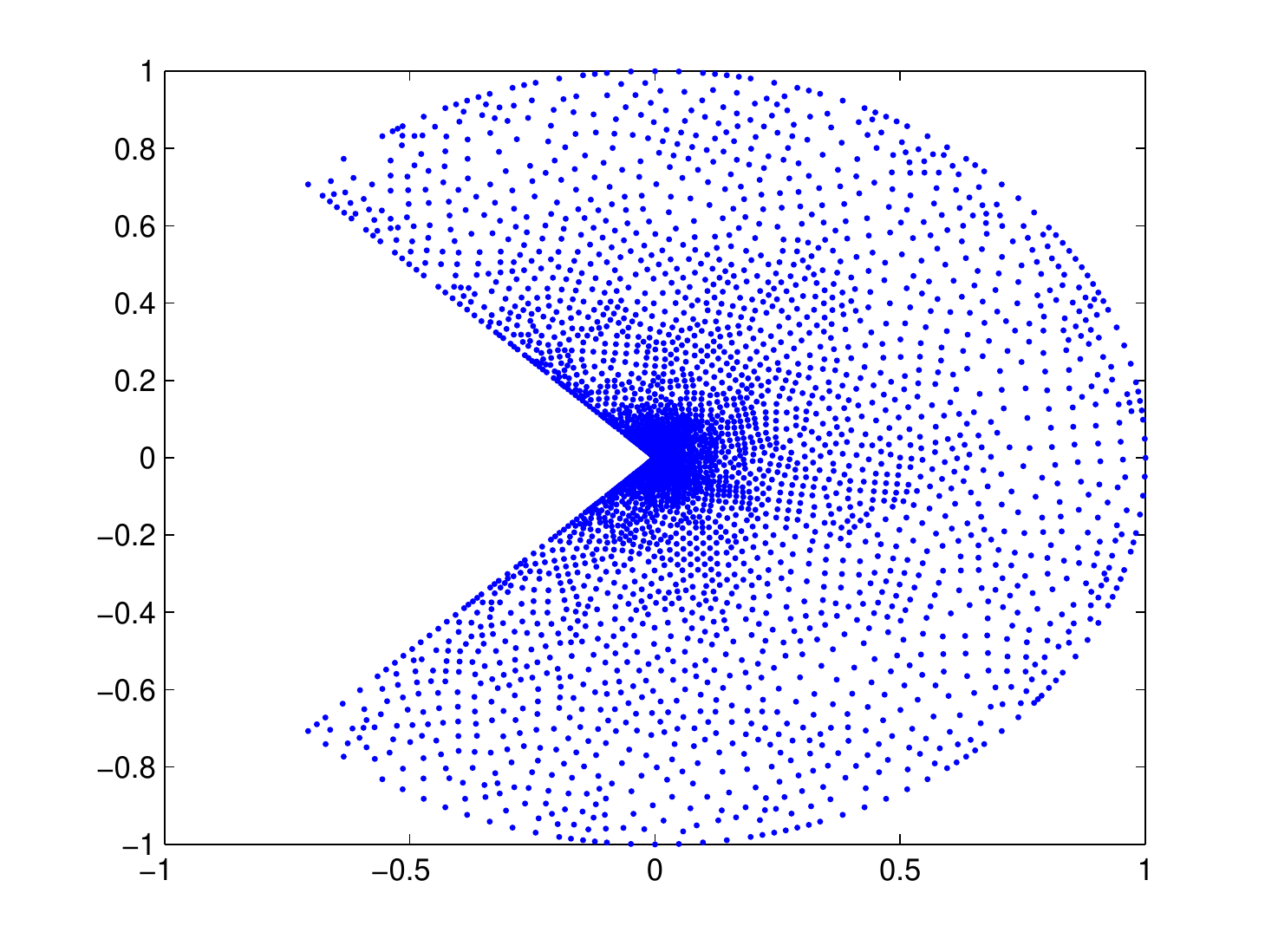}} \qquad 
 \subfigure[FEM centers (3009)]
{\includegraphics[width=6cm,height=4cm]{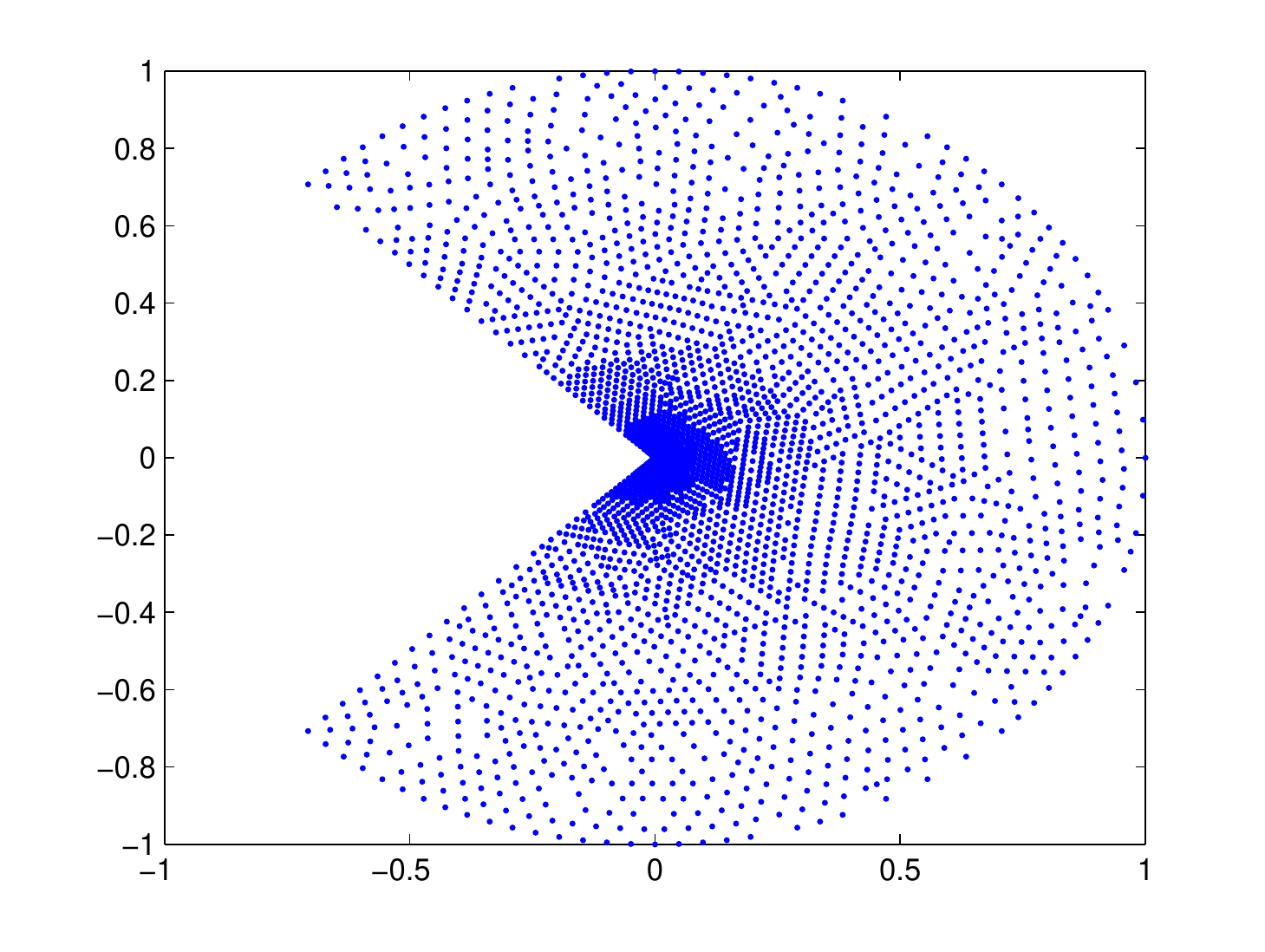}} 
 \subfigure[RBF-FD centers: zoom]%
{\includegraphics[width=6cm,height=4cm]{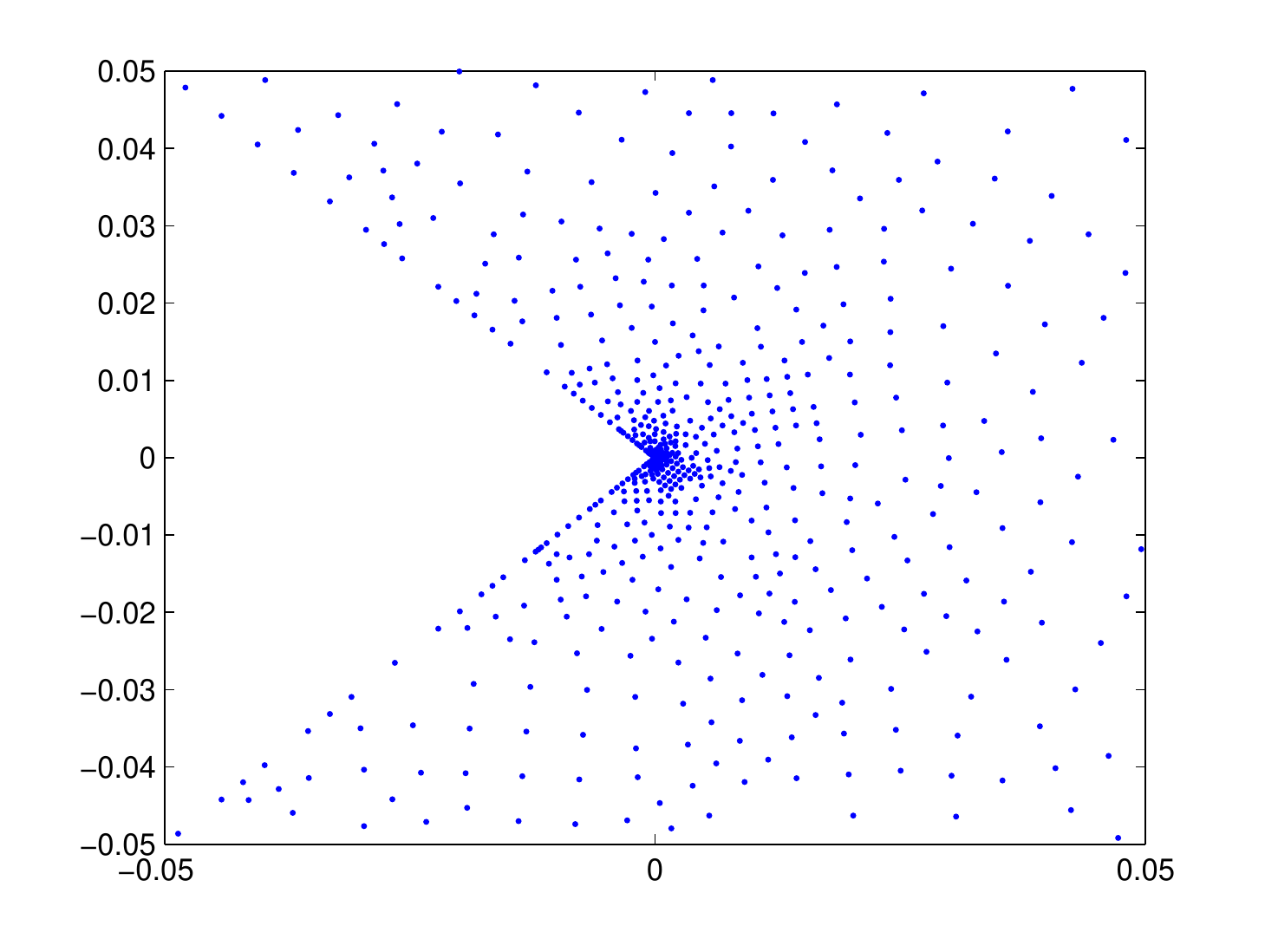}}  \qquad 
 \subfigure[FEM centers: zoom]%
{\includegraphics[width=6cm,height=4cm]{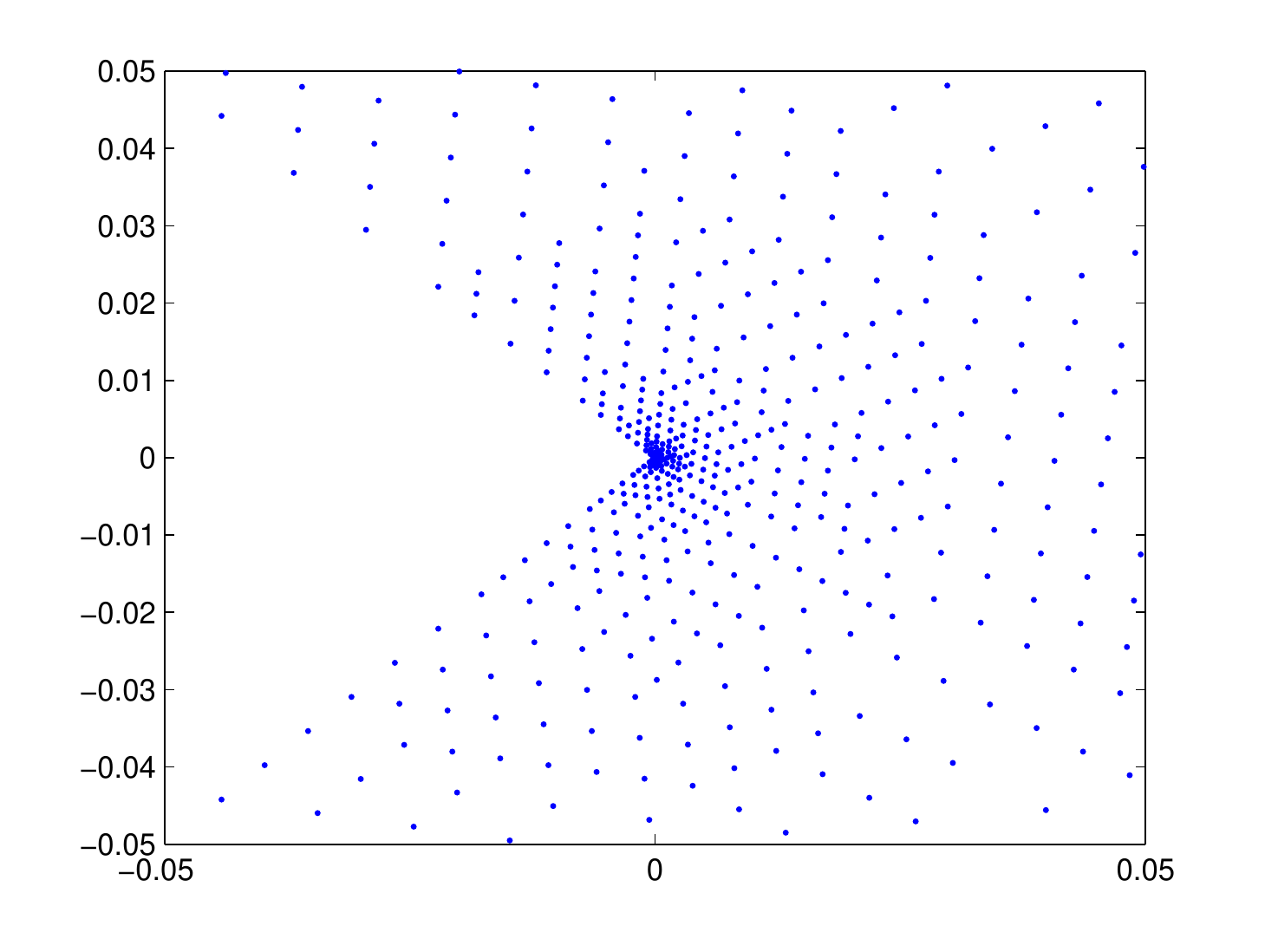}} 
        \end{center}
\caption{Test Problem~\ref{f4}: (a) Error $E_c$ of the discrete solution on the centers 
generated by consecutive refinements with the respective method, using FEM and two versions of the RBF-FD method 
({\tt RBF-FD old}: the method of \cite{DavyOanh11}, {\tt RBF-FD}: the method of this paper) 
as function of the reciprocal of the number of interior centers.
(b) Error $E_g$ of the interpolated solutions on a uniform  grid.
(cd) Error function $u-\hu$ for the RBF-FD solution of this paper on 3169 interior centers %
and   the FEM solution on 3009  interior vertices. %
(ef) The centers used for the respective solutions. (gh) Zooms into both sets of centers.
} 
\label{figF4}
\end{figure}

\begin{figure}[htbp!]
 \begin{center}
 \subfigure[Errors on centers]%
{\includegraphics[width=6cm,height=4cm]{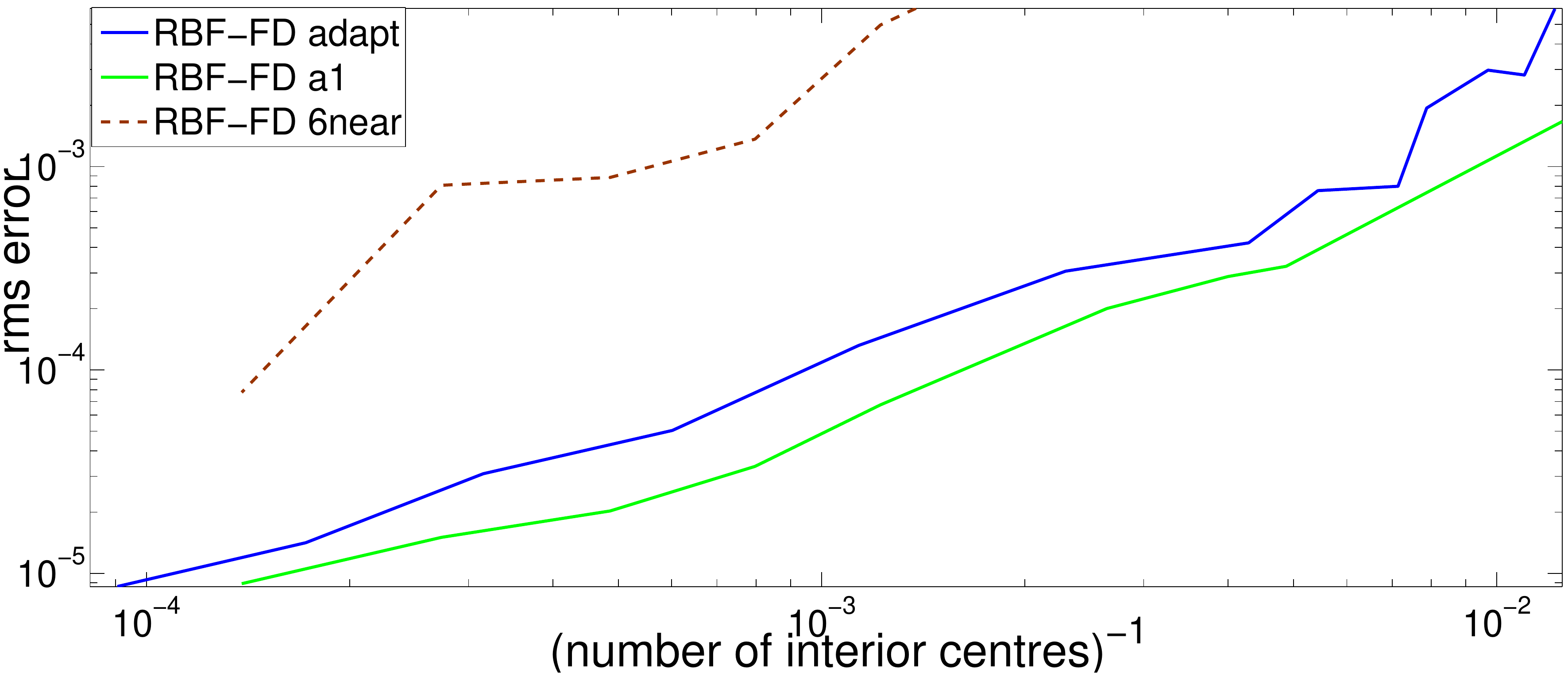}} \qquad 
 \subfigure[Errors on grid]
{\includegraphics[width=6cm,height=4cm]{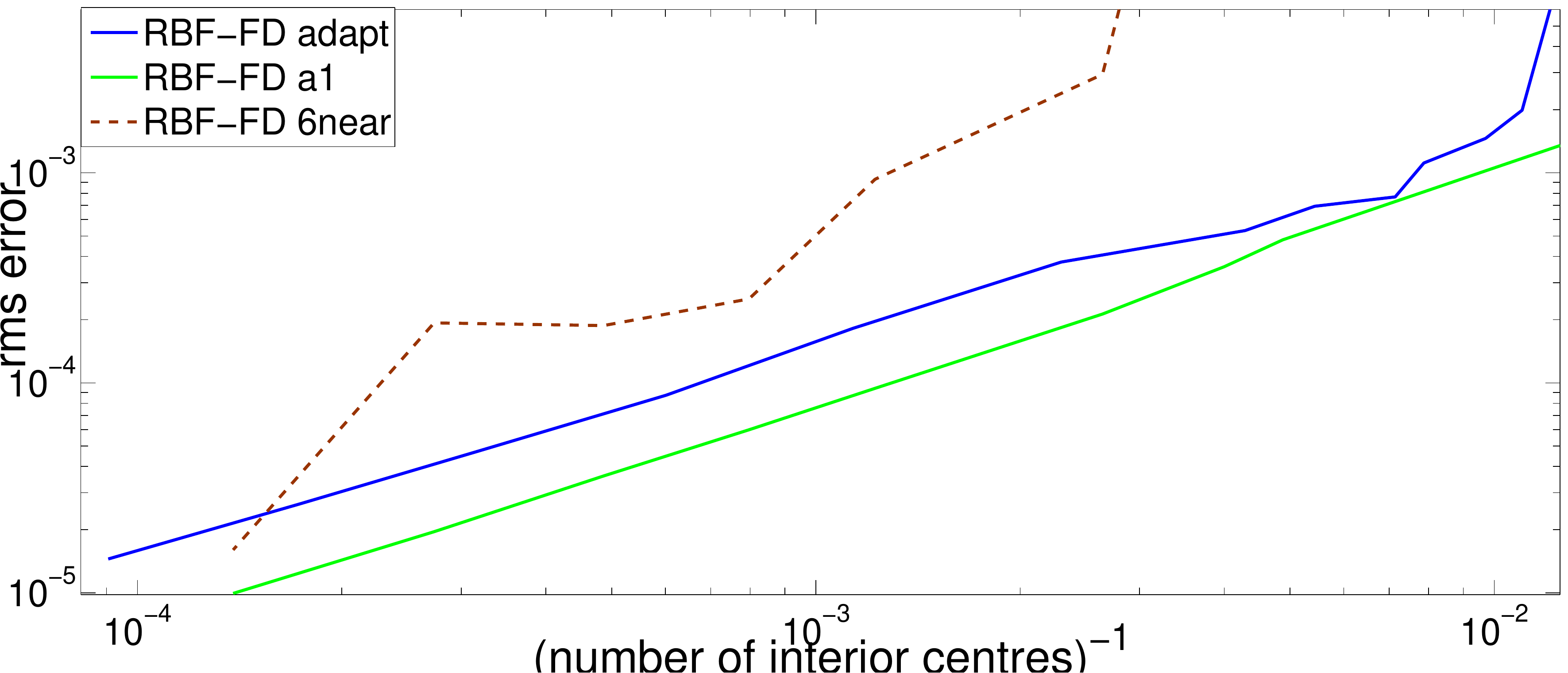}} 
 \subfigure[Centers (3654)]
{\includegraphics[width=6cm,height=4cm]{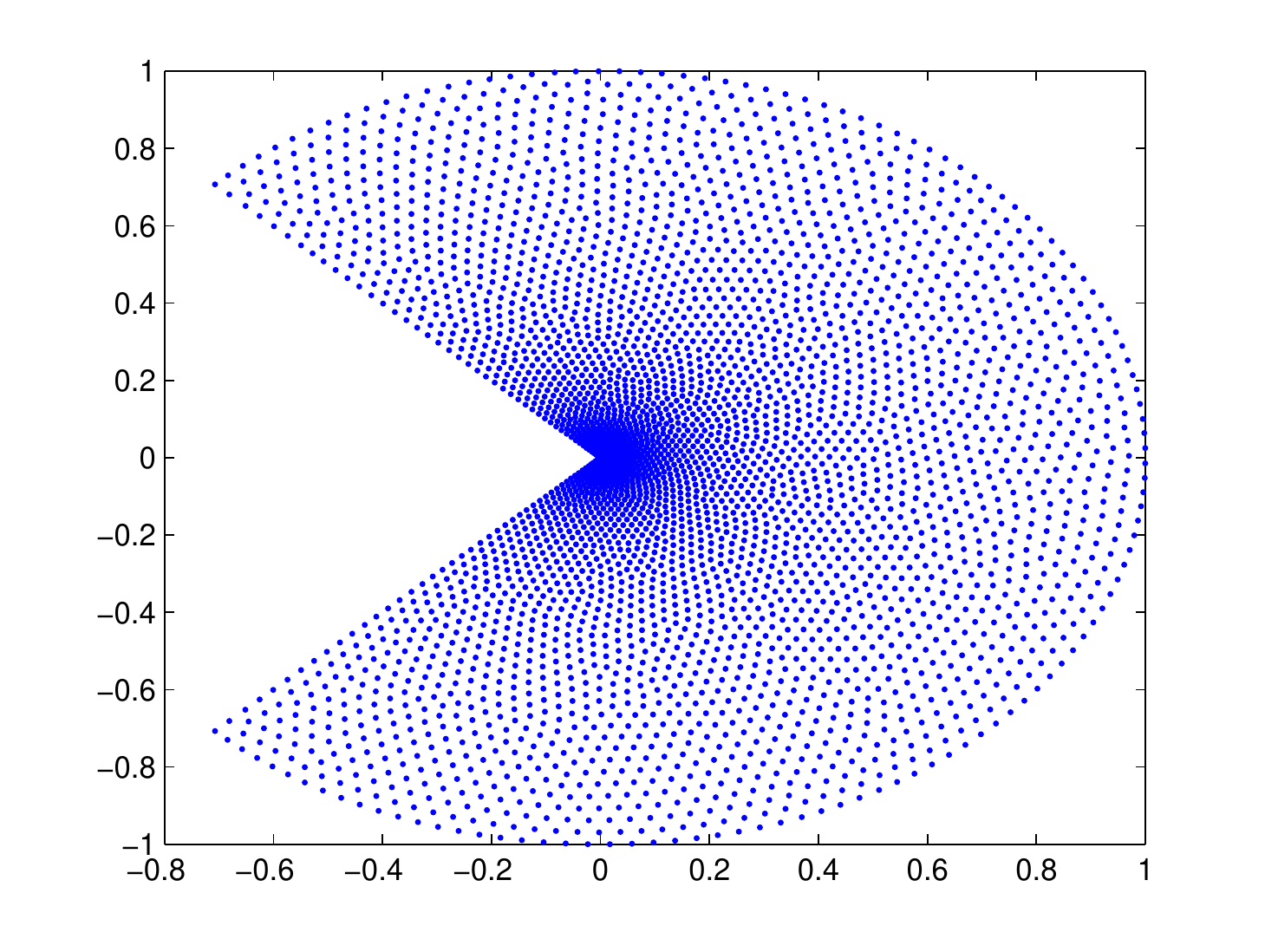}} \qquad 
 \subfigure[Centers: zoom]%
{\includegraphics[width=6cm,height=4cm]{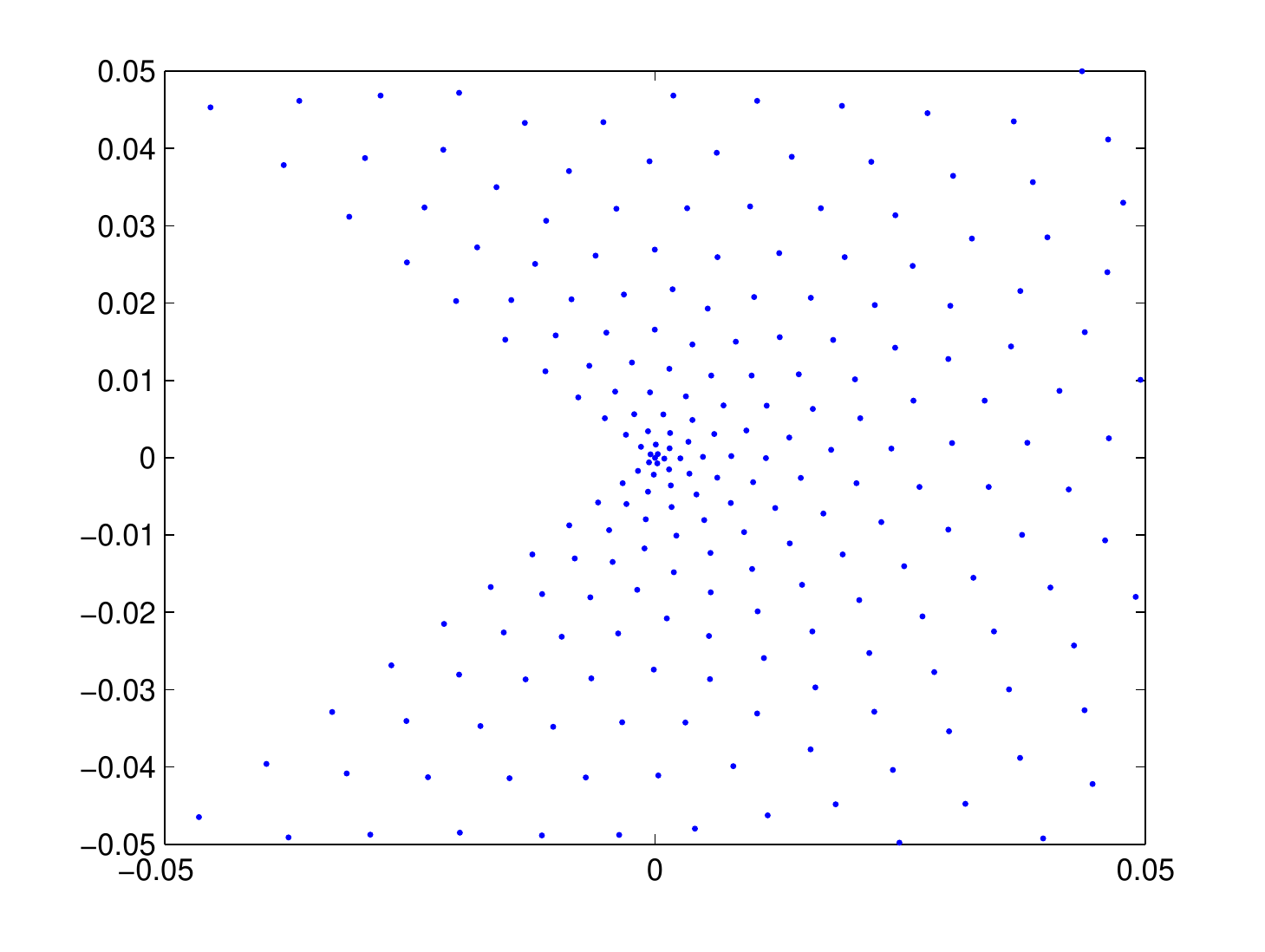}}  \qquad 
\subfigure[$\Xi_\zeta$ for {\tt RBF-FD 6near}]
{\includegraphics[width=6cm,height=4cm]{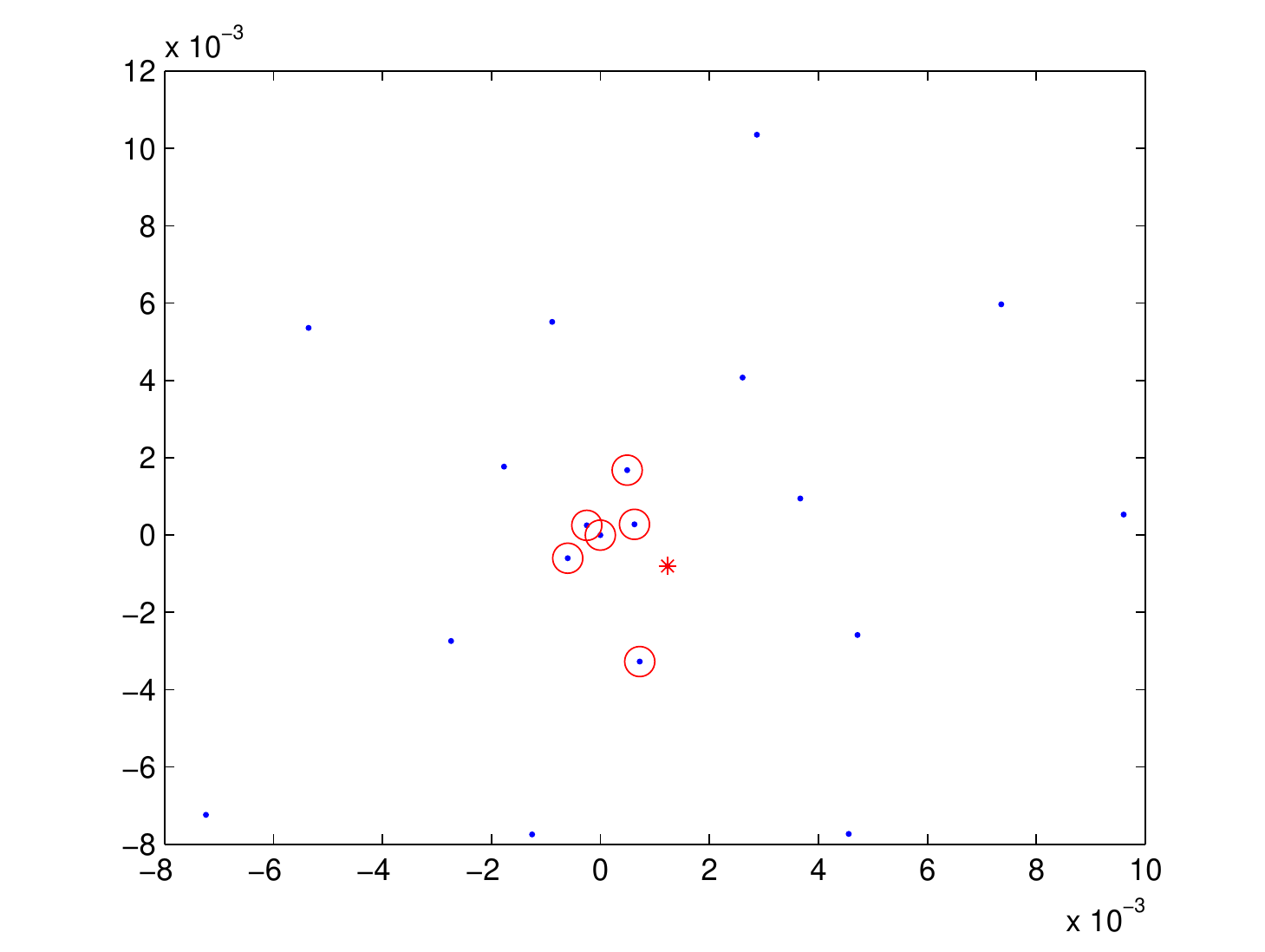}} \qquad 
 \subfigure[$\Xi_\zeta$ by Algorithm~\ref{alg1}]%
{\includegraphics[width=6cm,height=4cm]{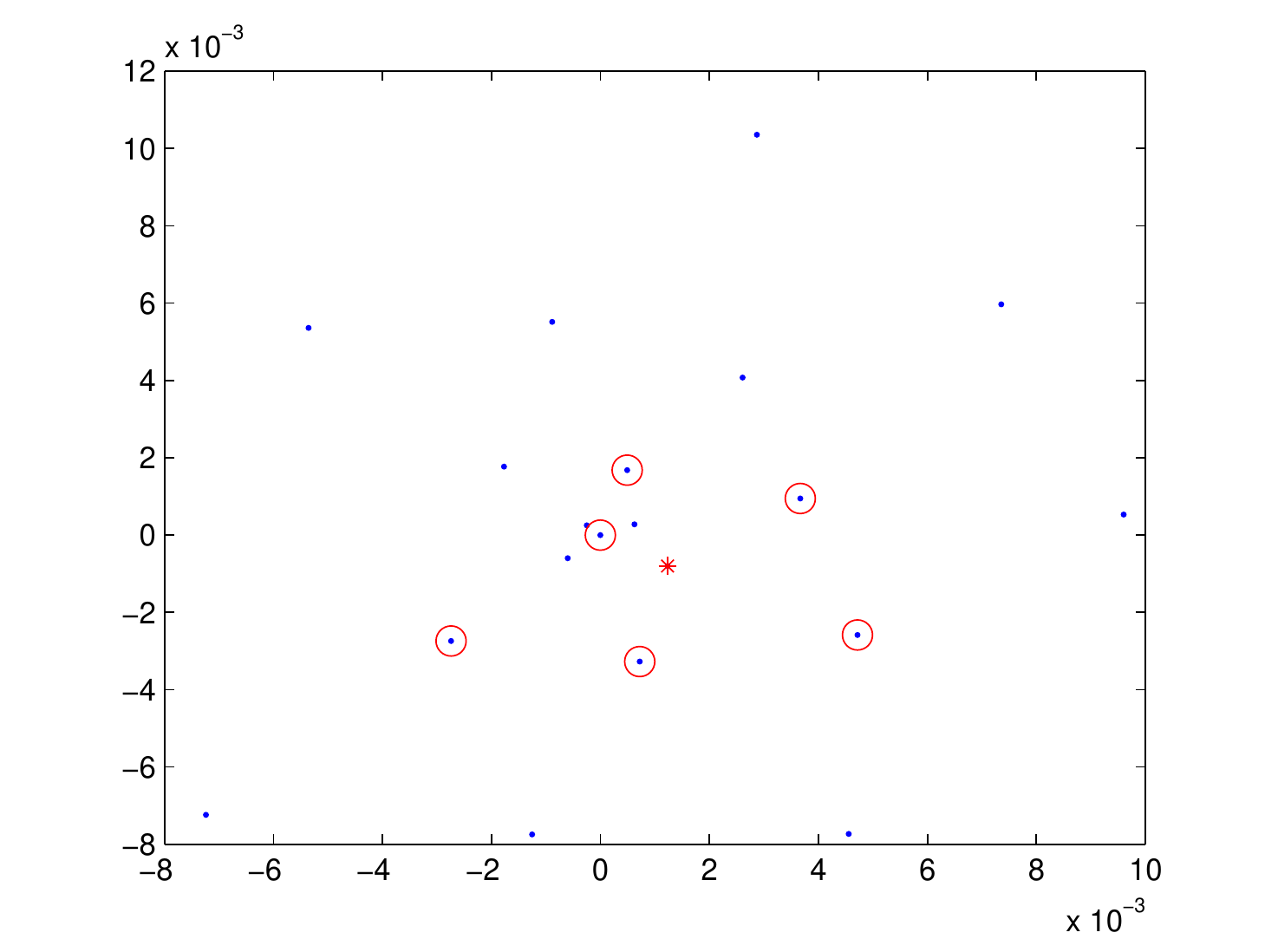}}  \qquad 
        \end{center}
\caption{Test Problem~\ref{f4}: Results on smoothly distributed centers 
obtained by \emph{a priori} refinement using {\tt DistMesh}. (ab) Errors on the centers and on a uniform grid as in 
Figure~\ref{figF4}. Two methods on smooth centers, RBF-FD method with stencil
support selection according to Algorithm~\ref{alg1} ({\tt RBF-FD a1}) and RBF-FD method with stencil supports $\Xi_\zeta$ 
obtained by choosing $\zeta$ and its 6 nearest points in $\Xi$ ({\tt RBF-FD 6near}).
The plots also include for comparison the error curves for the adaptive RBF-FD method of Figure~\ref{figF4}
({\tt RBF-FD adapt}). (cd) The 3654 centers obtained by {\tt DistMesh} with $h_0=0.0007$ and $\beta=0.5$, and a zoom.
(ef) Example where $\Xi_\zeta$ obtained by Algorithm~\ref{alg1} is significantly different from 6 nearest points.
}
\label{figF4Repul}
\end{figure}

\begin{figure}[htbp!]
 \begin{center}
 \subfigure[Errors on centers]%
{\includegraphics[width=6cm,height=4cm]{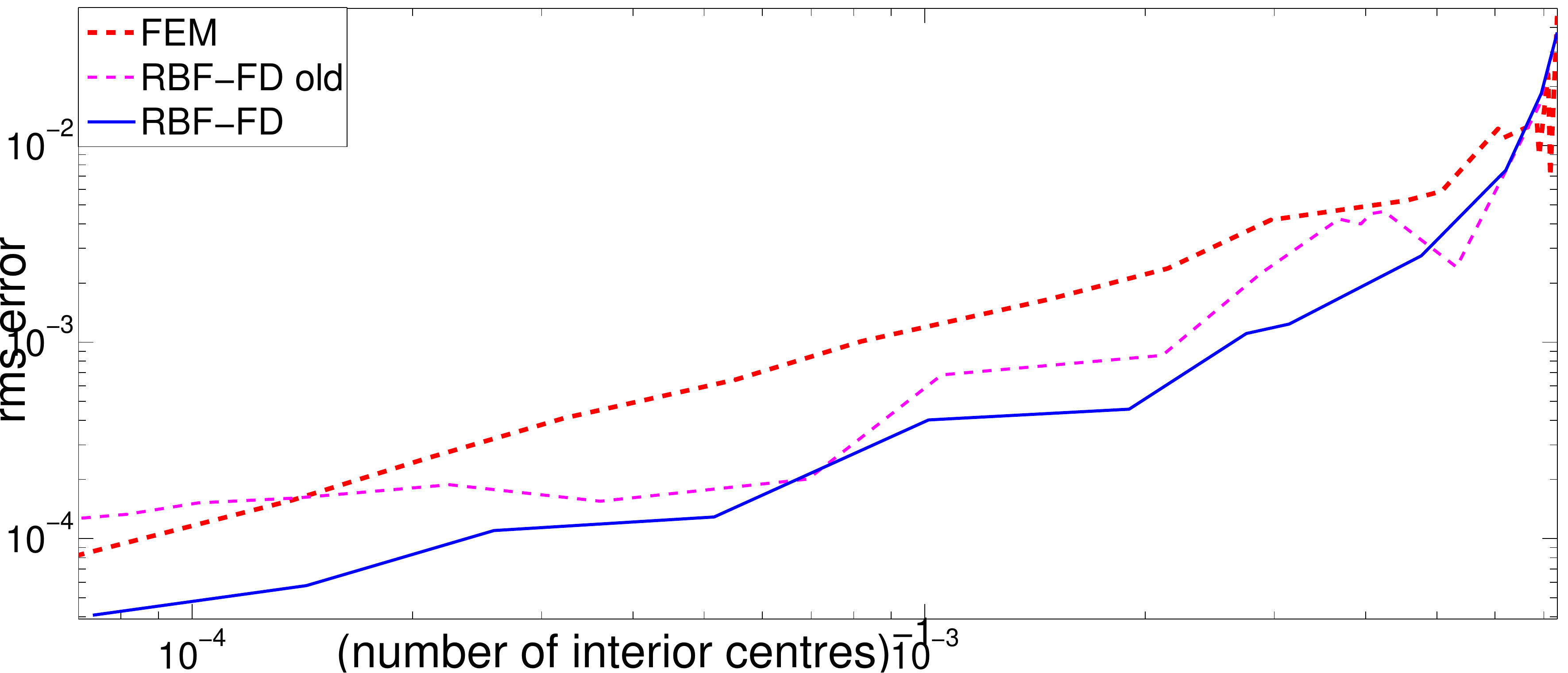}} \qquad
 \subfigure[Errors on grid]%
{\includegraphics[width=6cm,height=4cm]{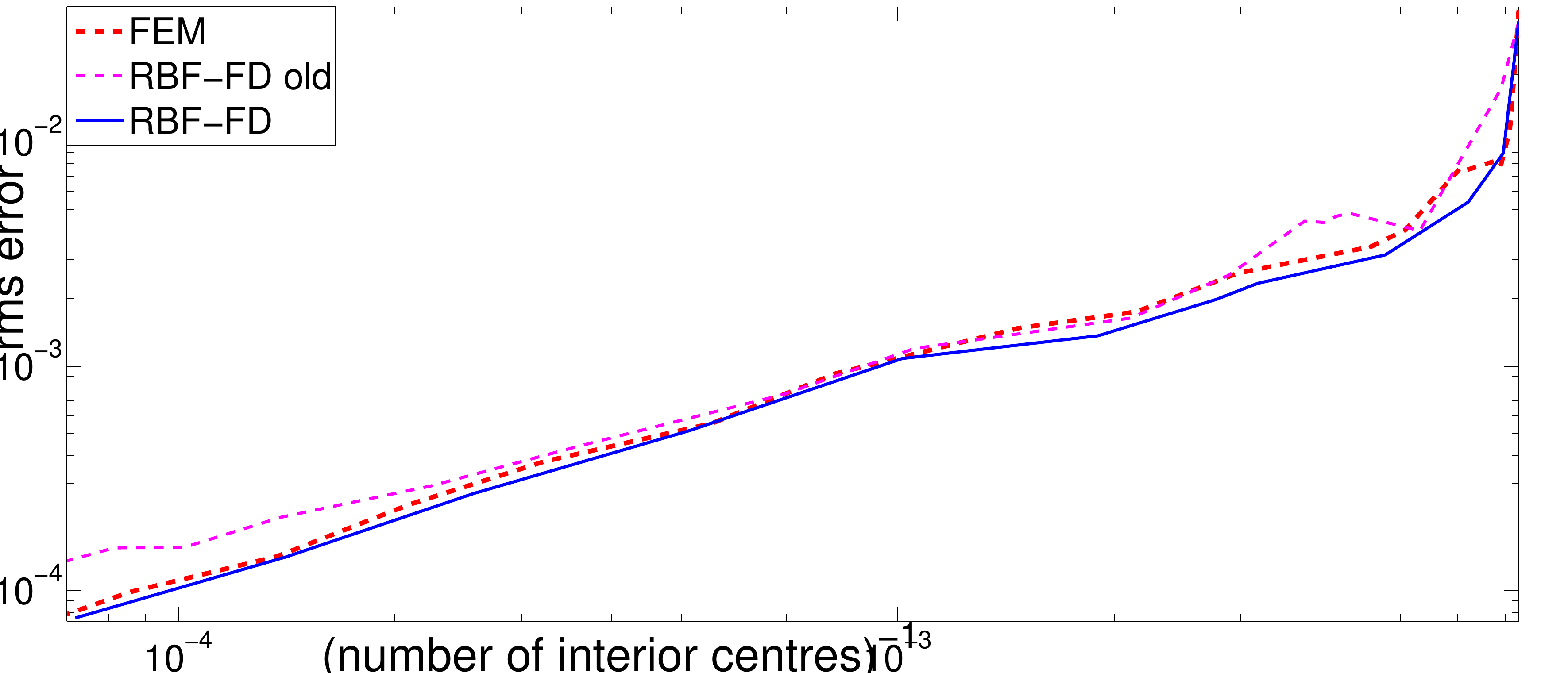}} 
 \subfigure[RBF-FD error]%
{\includegraphics[width=6cm,height=4cm]{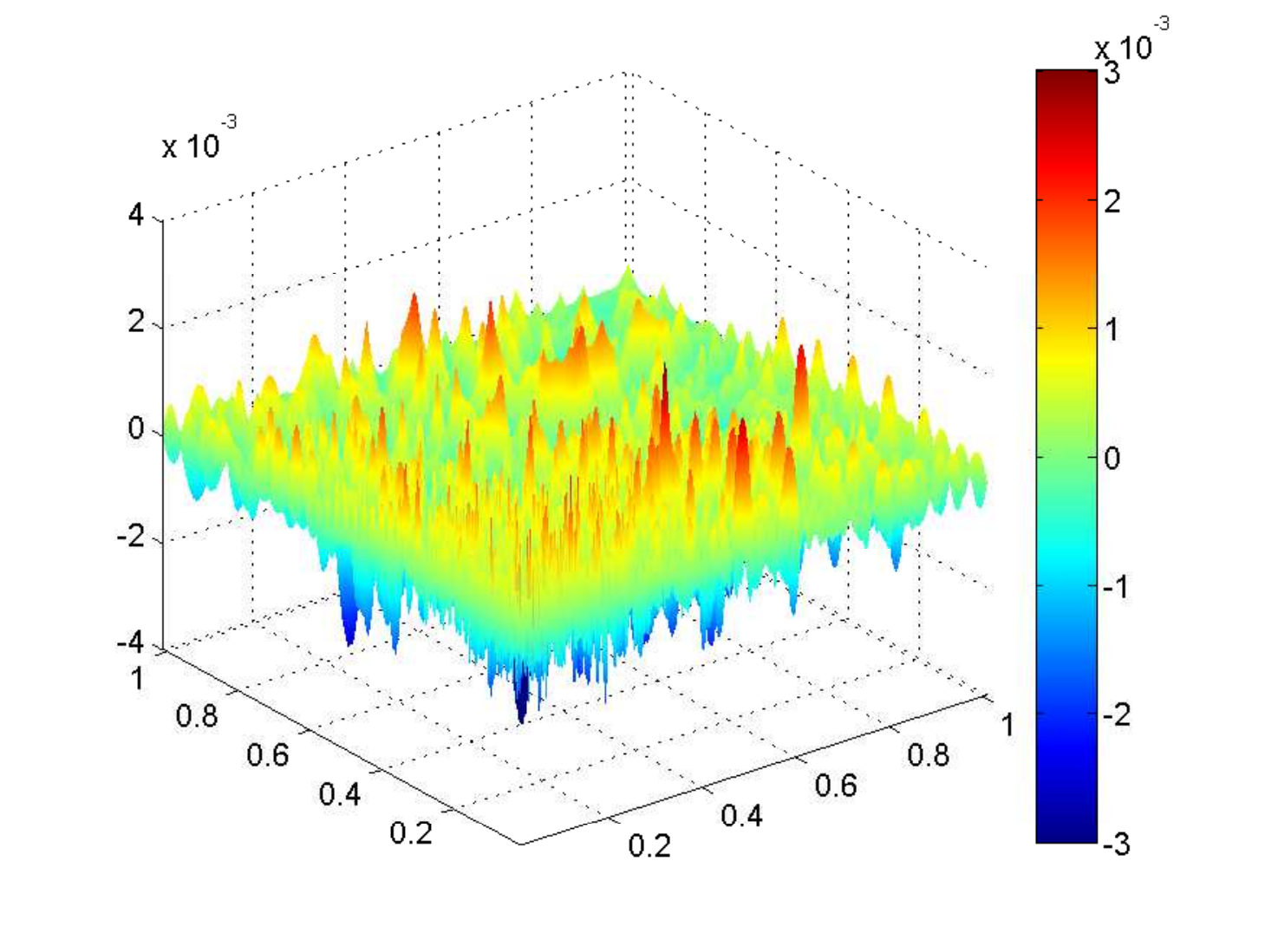}} \qquad
 \subfigure[FEM error]%
{\includegraphics[width=6cm,height=4cm]{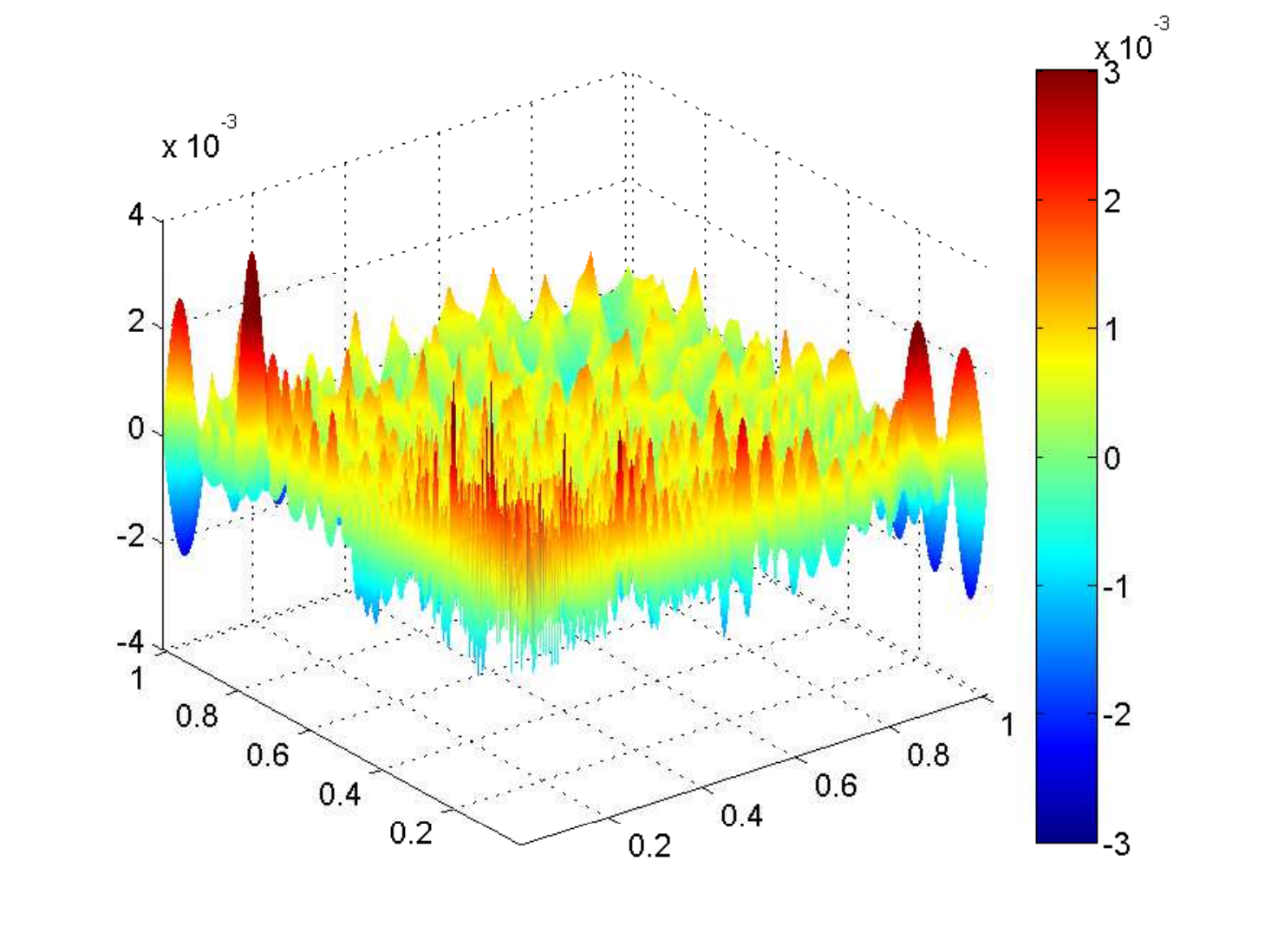}} 
 \subfigure[RBF-FD centers (1938)]
{\includegraphics[width=6cm,height=4cm]{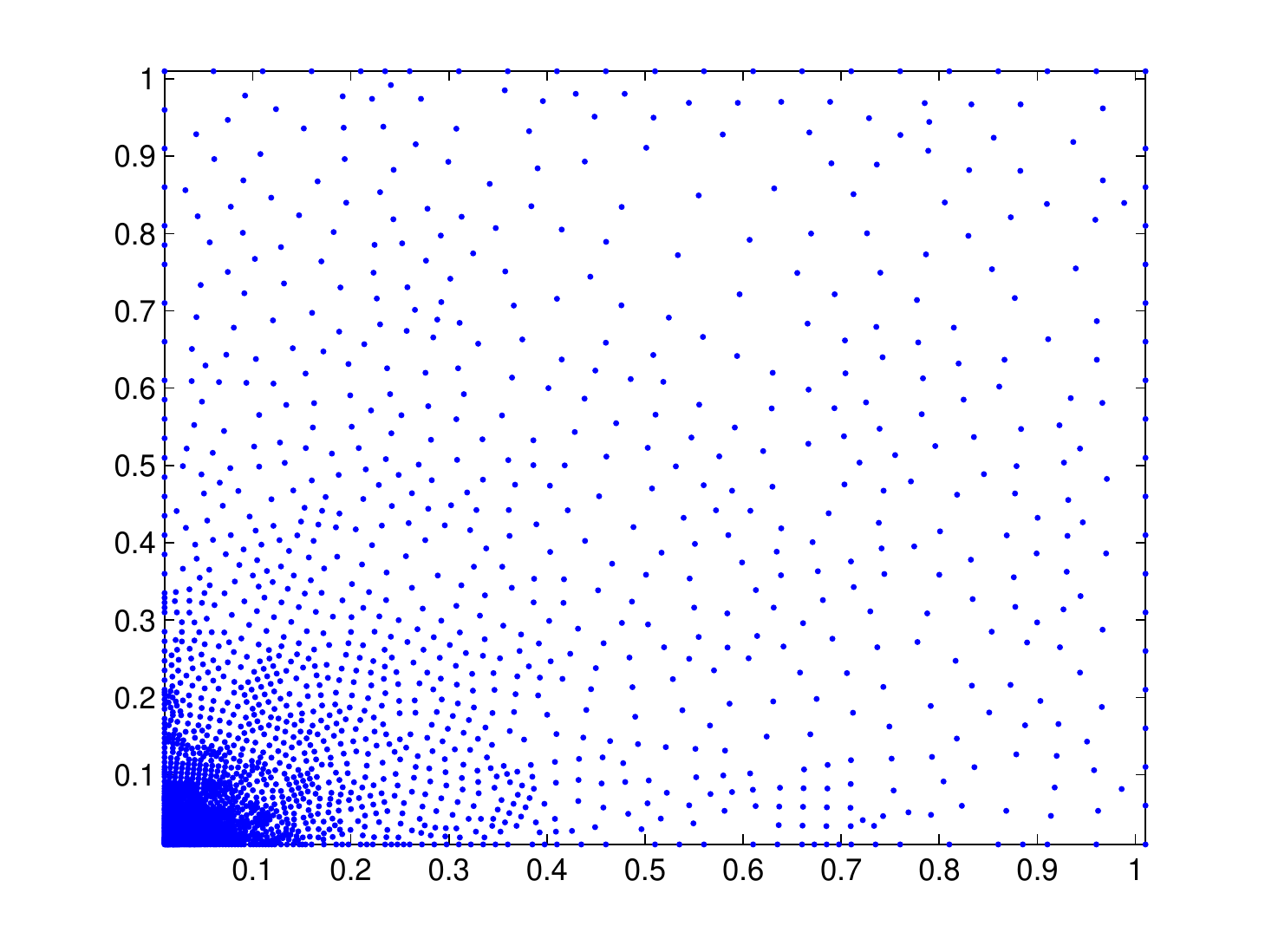}} \qquad
 \subfigure[FEM centers (1811)]
{\includegraphics[width=6cm,height=4cm]{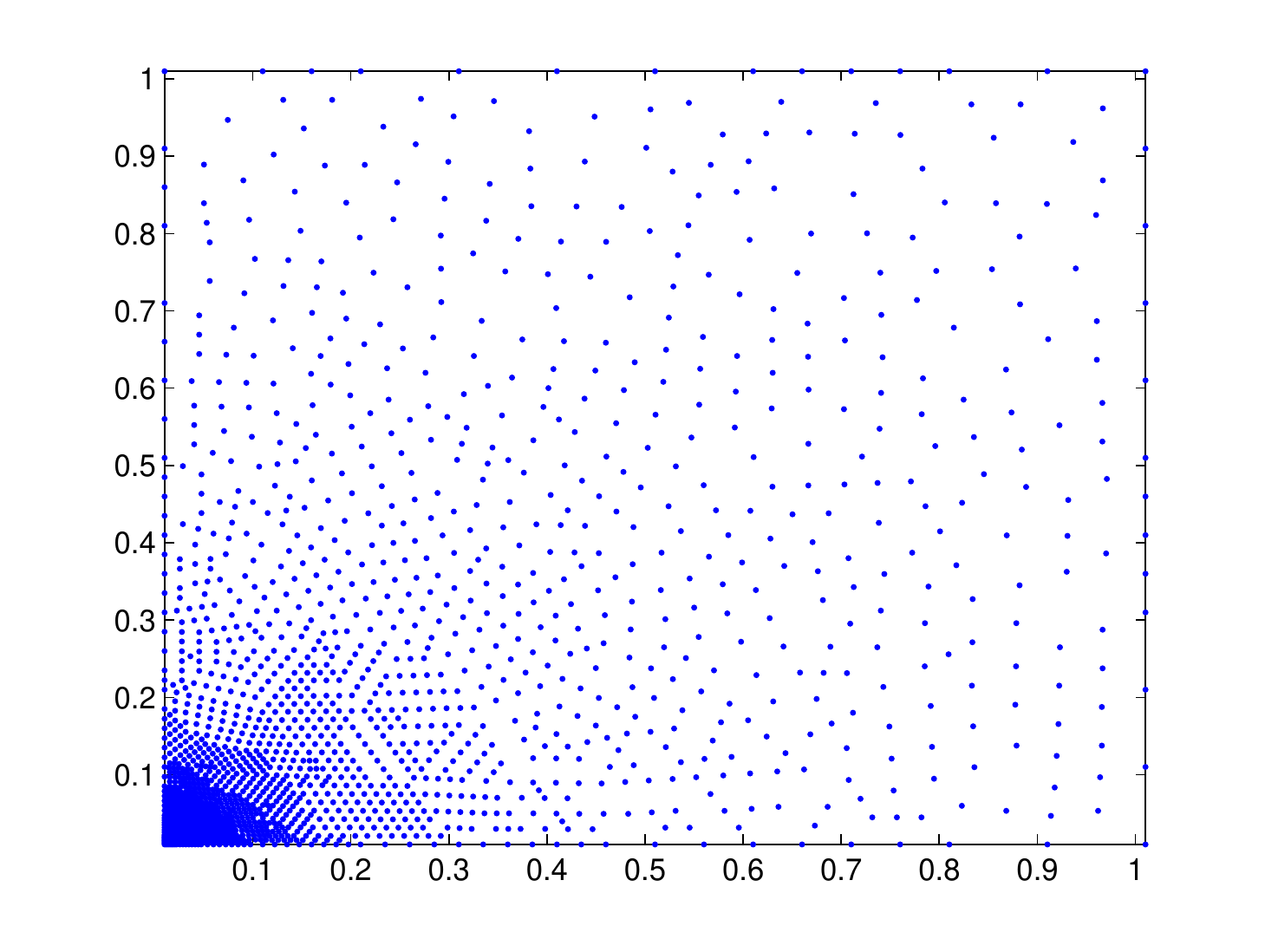}} 
 \subfigure[RBF-FD centers: zoom]%
{\includegraphics[width=6cm,height=4cm]{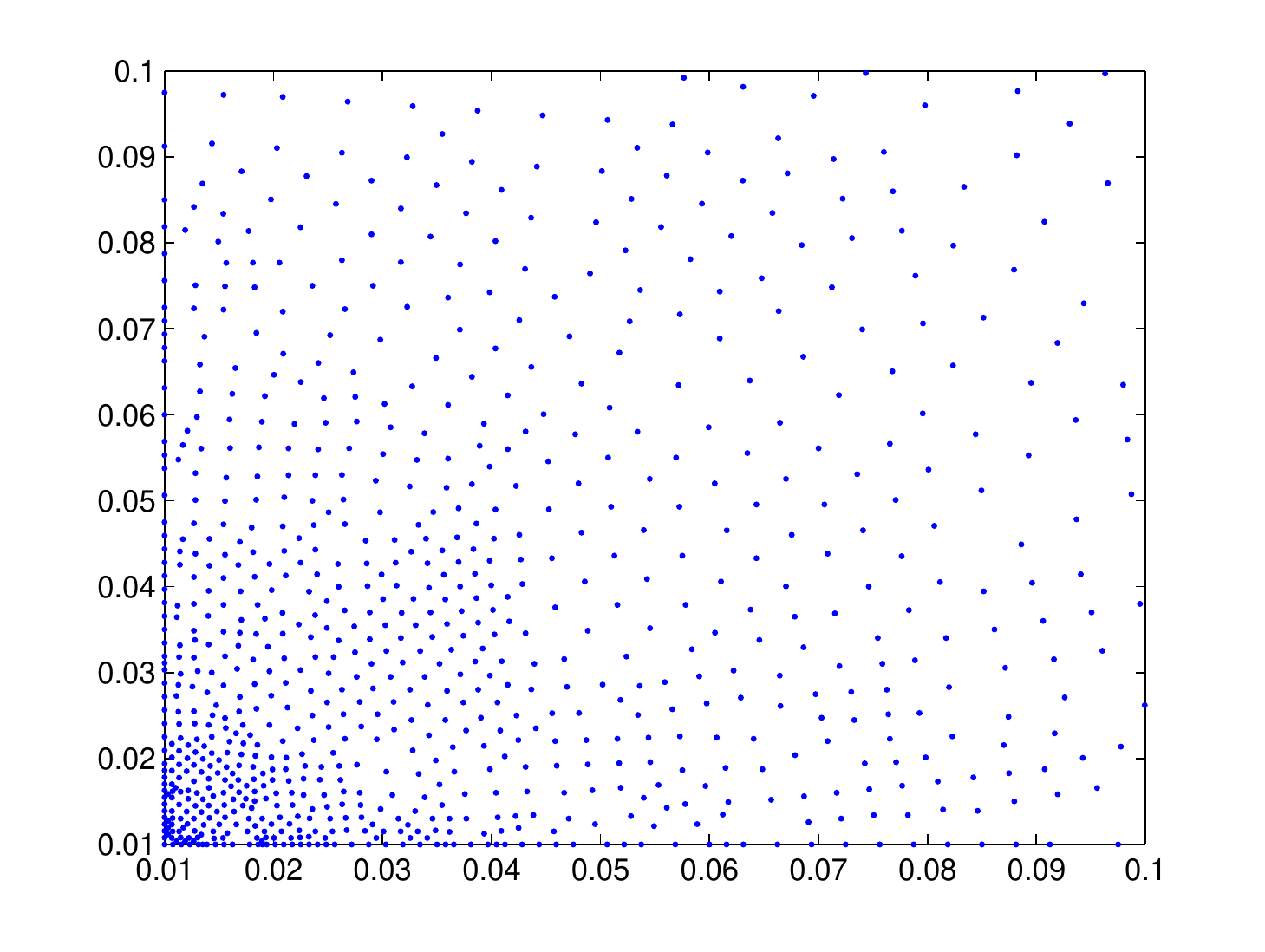}} \qquad
 \subfigure[FEM centers: zoom]%
{\includegraphics[width=6cm,height=4cm]{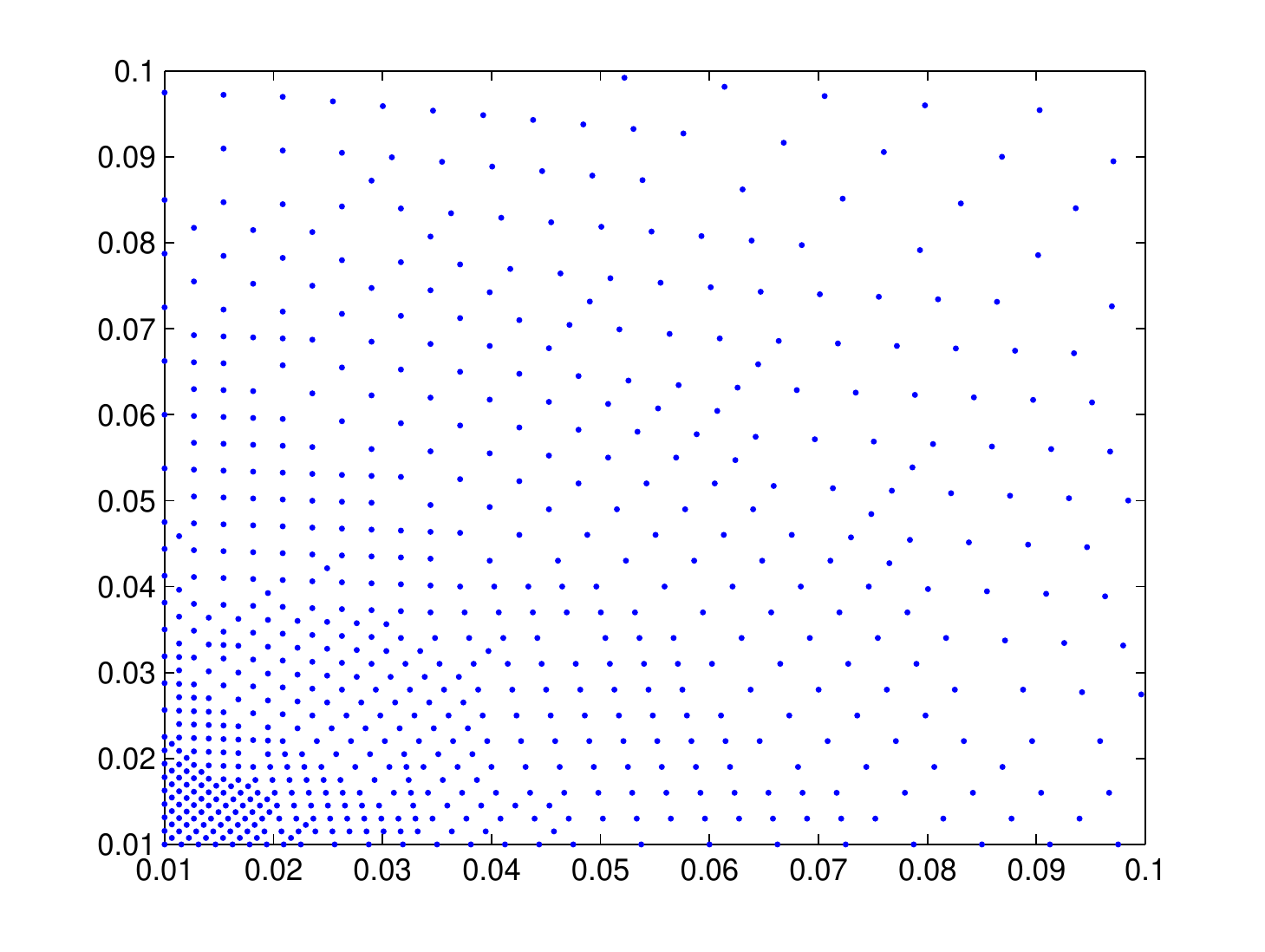}} 
        \end{center}
\caption{Test Problem~\ref{f8}: Errors and centers as in Figure~\ref{figF4}. 
The plots in (cd) are based on the RBF-FD 
solution on 1938 interior centers shown in (e) and the FEM 
solution on 1811 interior vertices shown in (f).
}
\label{figF8}
\end{figure}

The solutions of both problems have a singularity at the origin, as illustrated in 
Figure~\ref{F4F8exact}.
Our numerical results for these problems are presented in Figures~\ref{figF4}--\ref{figF8}.

\begin{table}[htbp!]
\begin{center}\renewcommand{\arraystretch}{1.2}\small
\begin{tabular}{|c|c|c|c|c|c|c|c|c|}
\hline
\hline
$\#\Xi_{\rm int}$  &205  &    250&      378   &    818&      1255   &   2057 &     3654    &  7226 \\ 
\hline
$h_0$&  0.0007 &0.0007 &0.001&  0.0012& 0.00065& 0.0012 & 0.0007&  0.0013   \\
\hline
$\beta$ & 0.700 &   0.675  &  0.650   & 0.600&    0.575&    0.525&    0.500 &  0.425\\
 \hline
\end{tabular}
\caption{Parameters of {\tt distmesh2d.m} used to produce the \emph{a priori} sets of centers 
for  Test Problem~\ref{f4}: $\#\Xi_{\rm int}$ is the number of interior centers,  $h_0$
the initial edge length, and  $\beta$ the power of the scaled edge 
length function in the form $\sigma(r)=r^{\beta}$. The set with $\#\Xi_{\rm int}=3654$
is illustrated in Figure~\ref{figF4Repul}(cd).
}
\label{h0Pow}
\end{center}
\end{table}

Figures~\ref{figF4}(ab) and \ref{figF8}(ab) show the graphs of the rms errors $E_c$ and $E_g$ of the RBF-FD method in
comparison to the finite element method, as function of the reciprocal of the number of interior
centers. The curves labeled {\tt FEM} stand for the rms error of the finite element solution computed using PDE Toolbox
with default parameters as in the example presented in \cite[function {\tt adaptmesh}]{PDEtool}.  
The curve {\tt RBF-FD old} gives the error of RBF-FD method as described in \cite{DavyOanh11},
where stencil selection and adaptive refinement are performed according to  
\cite[Algorithms 1 and 2]{DavyOanh11}, whereas {\tt RBF-FD} follows 
Algorithms~\ref{alg1} and \ref{alg2Refi} of this paper. We see that the error of the current RBF-FD
method is generally smaller on its centers than the error of the finite element solution on the vertices of its
triangulations. The errors of both FEM and RBF-FD on the grid are very close. 
The figures also show that  {\tt RBF-FD} is significantly more robust and accurate after repeated
refinements than {\tt RBF-FD old}.  Note that neither adding a constant
term to the Gaussian sum, nor using the `safe' shape
parameter as in \cite{DavyOanh11} changes the RBF-FD results significantly, so that the improvement
in the performance should be mostly attributed to the improved stencil selection and refinement
algorithms.

Figures~\ref{figF4}(cd) and \ref{figF8}(cd) compare the error functions for RBF-FD and FEM 
for the sets of centers of comparable size, whereas \ref{figF4}(e-h) and \ref{figF8}(e-h) illustrate the distribution of these
centers. The results are strikingly similar, which shows that the RBF-FD method generates
reasonably placed adaptive centers and rather uniformly distributed error, without 
significant outliers near the singularity.

In addition, for Test Problem~\ref{f4} we have checked how the results compare if 
Algorithm~\ref{alg2Refi} is replaced by an \emph{a priori} refinement method that 
produces smoothly distributed centers with the help of  a prescribed local separation 
distance function. Such centers can be obtained  by simulating the movement of small particles 
under electrostatic repulsion forces until they reach an equilibrium. 
We generated them by using MATLAB software {\tt DistMesh} \cite{PerssonStrang04}
available from \url{http://persson.berkeley.edu/distmesh/}. 
The function {\tt distmesh2d} from this package produces the centers as vertices of a triangulation 
and requires as input the initial edge length and a scaled edge length function which controls the local
separation of the centers. As suggested in \cite{Wahlbin91}, nearly optimal convergence of the 
adaptive finite element method for a problem with the behavior of the exact solution as 
$r^\alpha$, $0<\alpha<1$, where $r$ is the distance to a singular point, is expected if the edge
length is proportional to $r^{1-\alpha/k}$, where $k$ is the order of the method, that is $k=2$ for
the piecewise linear finite elements. For Test Problem~\ref{f4} this corresponds to the scaled edge 
length function $\sigma(r)=r^{2/3}$. 
We obtained several sets of centers by running 
{\tt distmesh2d} with various values of the initial edge length $h_0$ and choosing the scaled edge 
length function in the form $\sigma(r)=r^{\beta}$ for some $\beta$, see Table~\ref{h0Pow} that shows
the number of interior centers obtained this way. Note that the same $h_0$ may lead to significantly
different sizes of $\Xi_{\rm int}$ depending on $\beta$. 
The results are presented in Figure~\ref{figF4Repul}. They show the improvement of the RBF-FD errors with a 
factor of about 2 in comparison to our adaptive method using Algorithms~\ref{alg1} and \ref{alg2Refi} 
can be achieved on the smooth centers, compare the curves marked {\tt RBF-FD a1} and {\tt RBF-FD adapt} in 
Figure~\ref{figF4Repul}(ab). On the other hand, even on the smoothly distributed centers we had to use our 
Algorithms~\ref{alg1} to do the computations for {\tt RBF-FD a1} because a simple algorithm that 
obtains $\Xi_\zeta$ by choosing 6 nearest points in addition to $\zeta$ does not perform well ({\tt RBF-FD 6near}),
see also Figure~\ref{figF4Repul}(ef) that illustrates the difference in  $\Xi_\zeta$ obtained by both methods.
Figure~\ref{figF4Repul}(cd) confirms that the centers are indeed smoothly distributed. 
We conclude that although somewhat better errors can be obtained by designing smooth \emph{a priori} sets of centers if the
strength of the singularity is known, the performance of our adaptive method that does not rely on such knowledge is
nevertheless competitive.

\medskip

In what follows, Test Problems~\ref{f42}, \ref{f441} and \ref{f43}
are borrowed from \cite{Mitchell2013} and cover all problems with
isolated point singularities suggested there for testing adaptive algorithms. 
We added Test Problem~\ref{curved_slit} to consider a domain with a curved slit. 

 \begin{testproblem}\label{f42}\rm 
\cite[Section 2.2: Reentrant Corner]{Mitchell2013} 
Dirichlet problem for the Laplace equation $\Delta u=0$ 
in the domain $\Omega_\omega=(-1,1)^2\cap\{(r,\varphi) : 0<\varphi<\omega\}$, where 
$r,\varphi$ are the polar coordinates, for several values of $\omega\in(0,2\pi]$.  
The boundary conditions are chosen such that the exact solution is 
$r^\alpha \sin(\alpha \phi)$ in polar coordinates, where %
$\alpha = \pi/\omega$.
\end{testproblem}
 
\begin{figure}[htbp!]
 \begin{center}
 \subfigure[$\omega=\pi+0.01$: errors on centers]%
{\includegraphics[width=6cm,height=4cm]{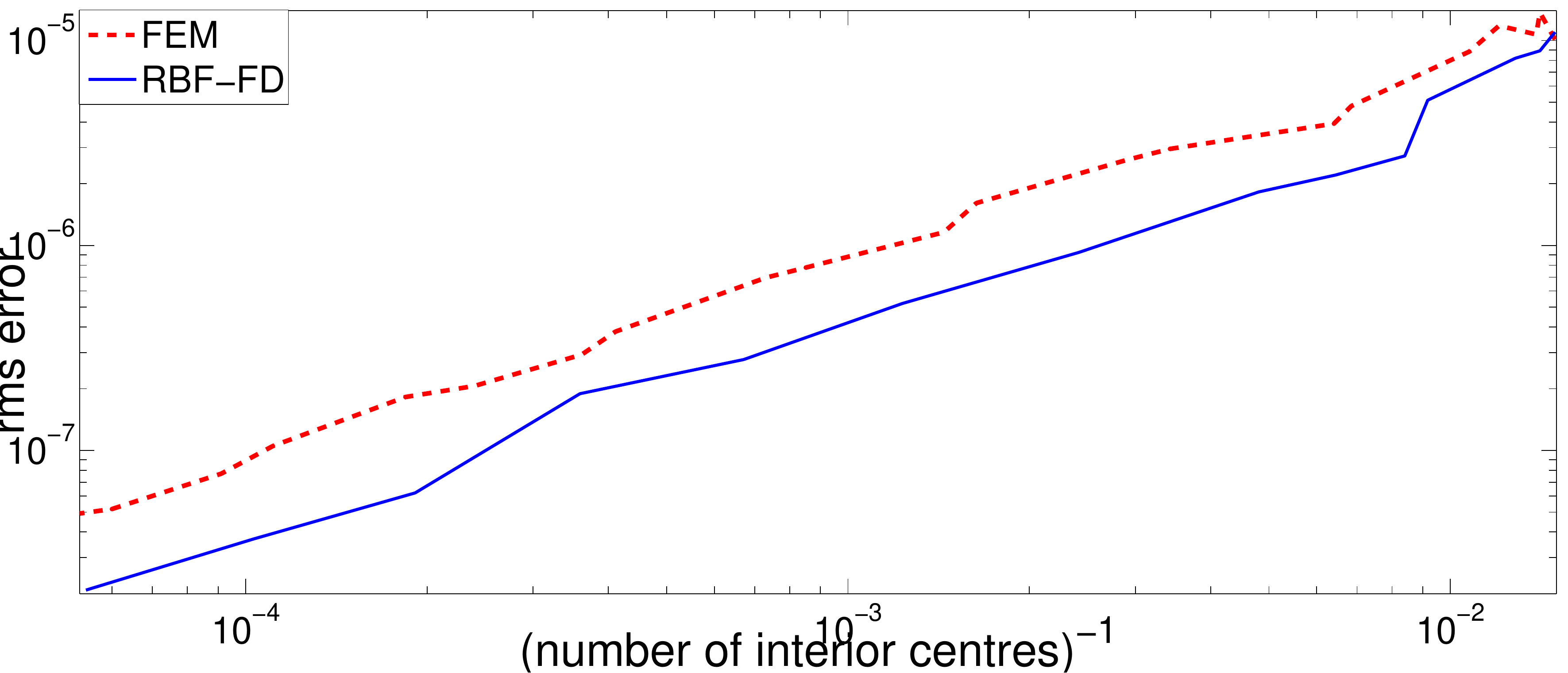}}\qquad 
 \subfigure[$\omega=\pi+0.01$: errors on grid]
{\includegraphics[width=6cm,height=4cm]{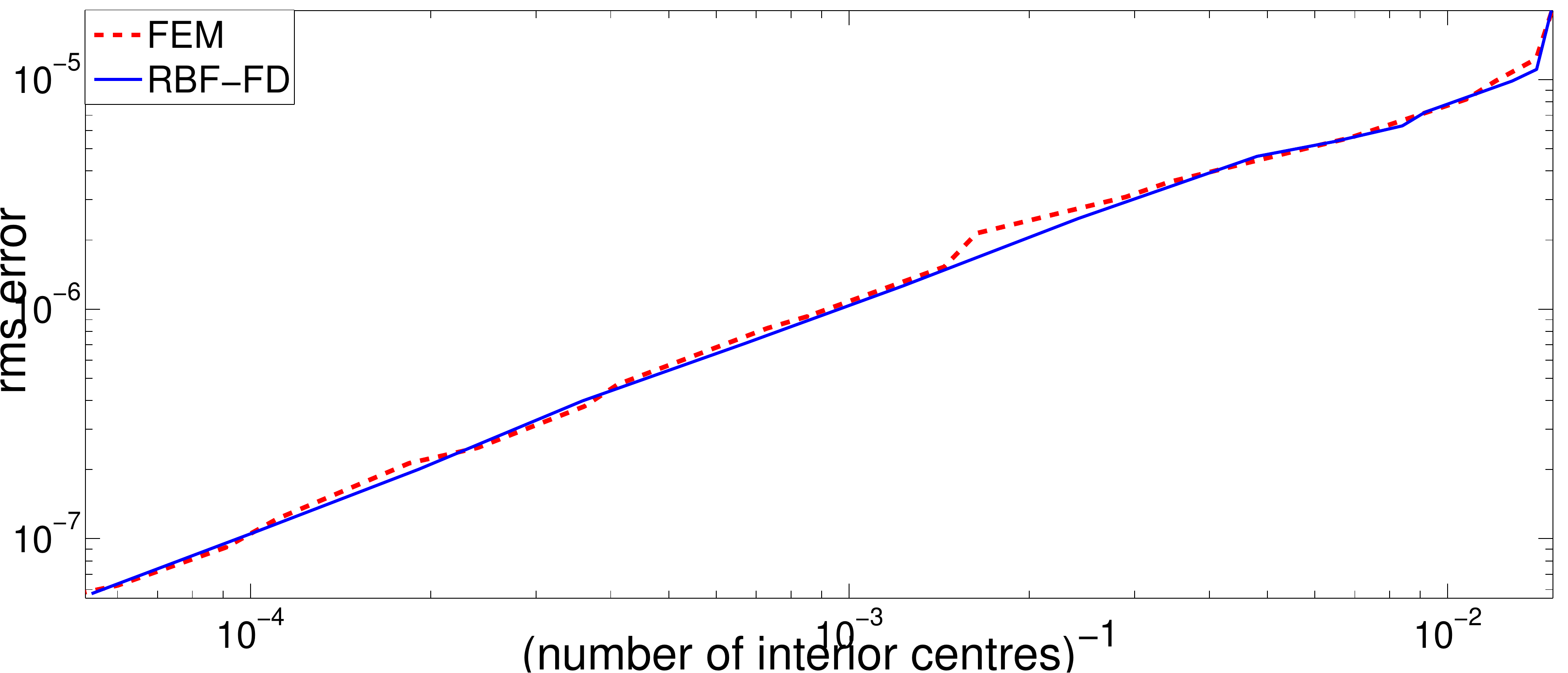}}  
 \subfigure[$\omega=5 \pi/4$: errors on centers]%
{\includegraphics[width=6cm,height=4cm]{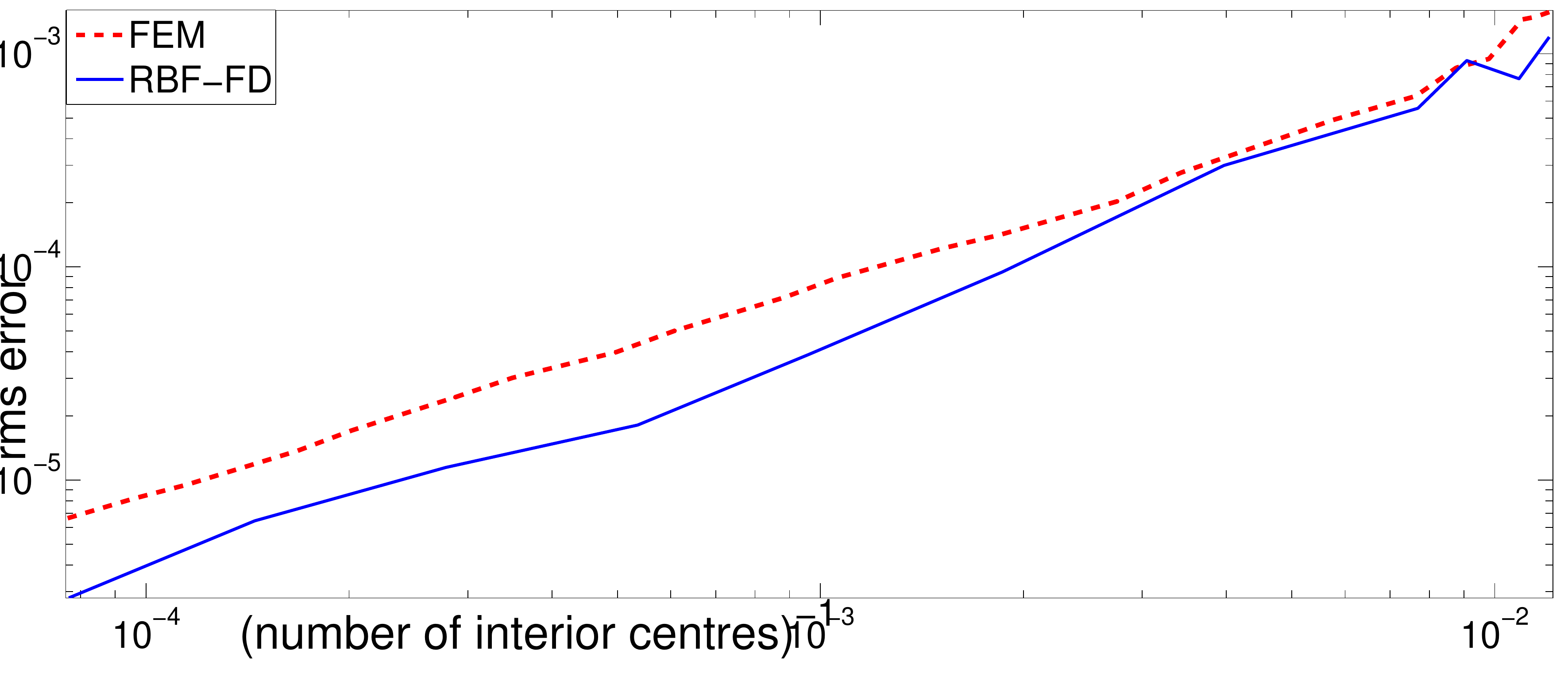}} \qquad 
 \subfigure[$\omega=5 \pi/4$: errors on grid]
{\includegraphics[width=6cm,height=4cm]{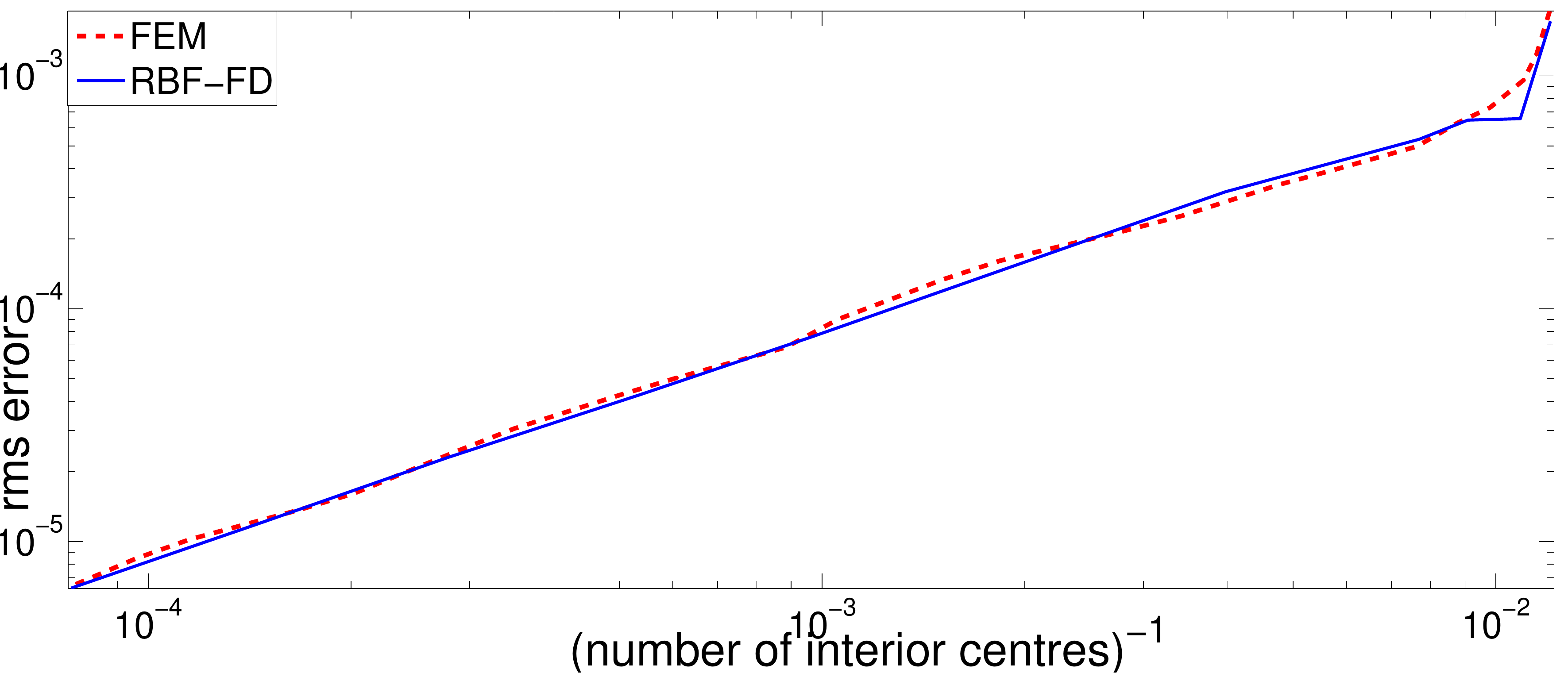}}  
  \subfigure[$\omega=7 \pi/4$: errors on centers]%
{\includegraphics[width=6cm,height=4cm]{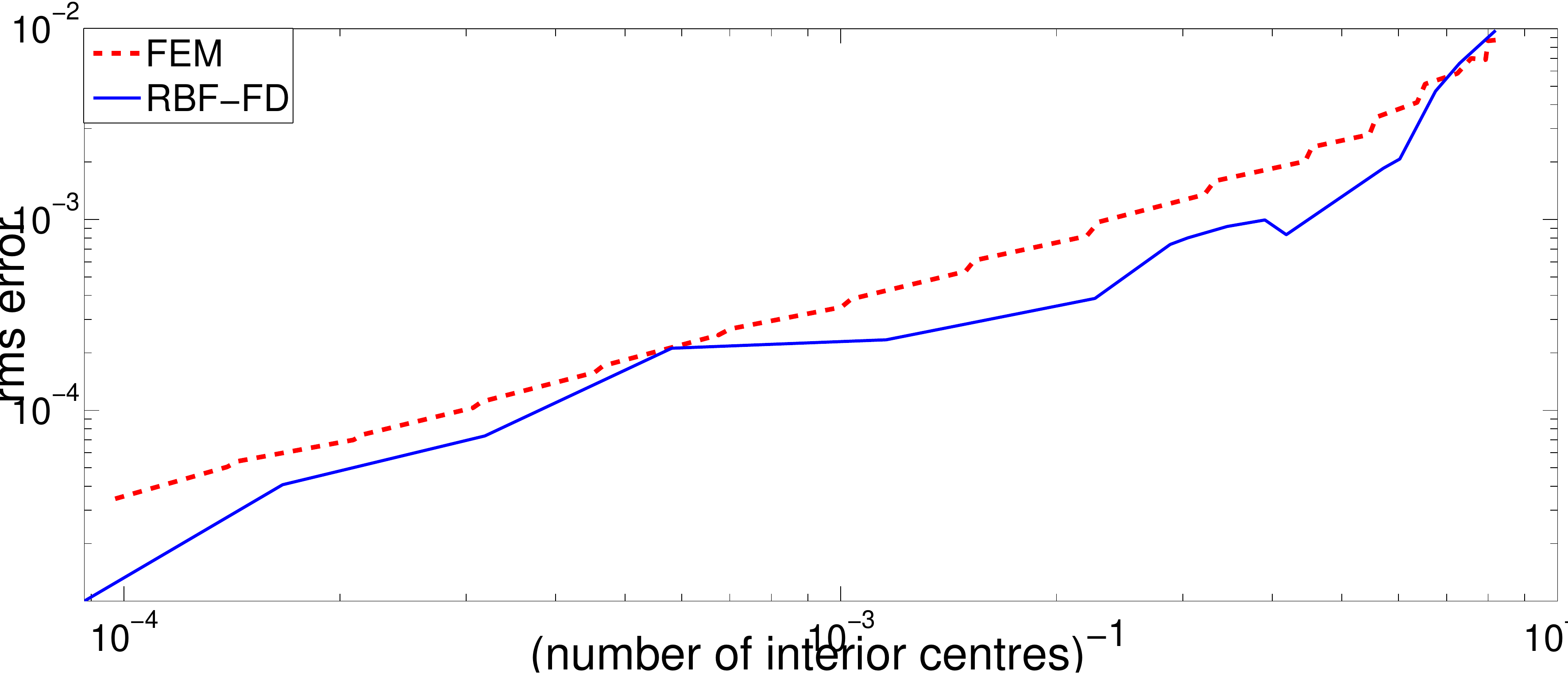}} \qquad 
 \subfigure[$\omega=7 \pi/4$: errors on grid]
{\includegraphics[width=6cm,height=4cm]{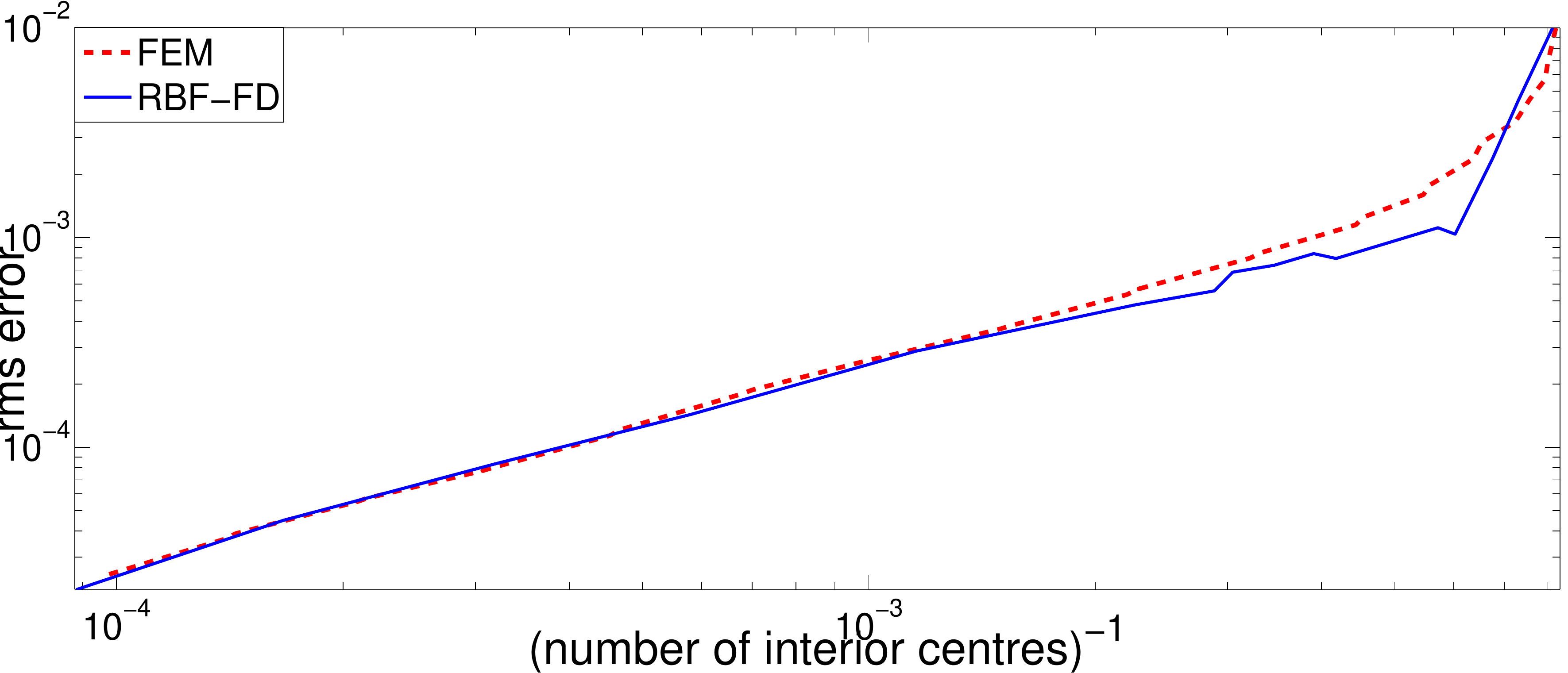}} 
 \subfigure[$\omega=2 \pi$:  errors on centers]%
{\includegraphics[width=6cm,height=4cm]{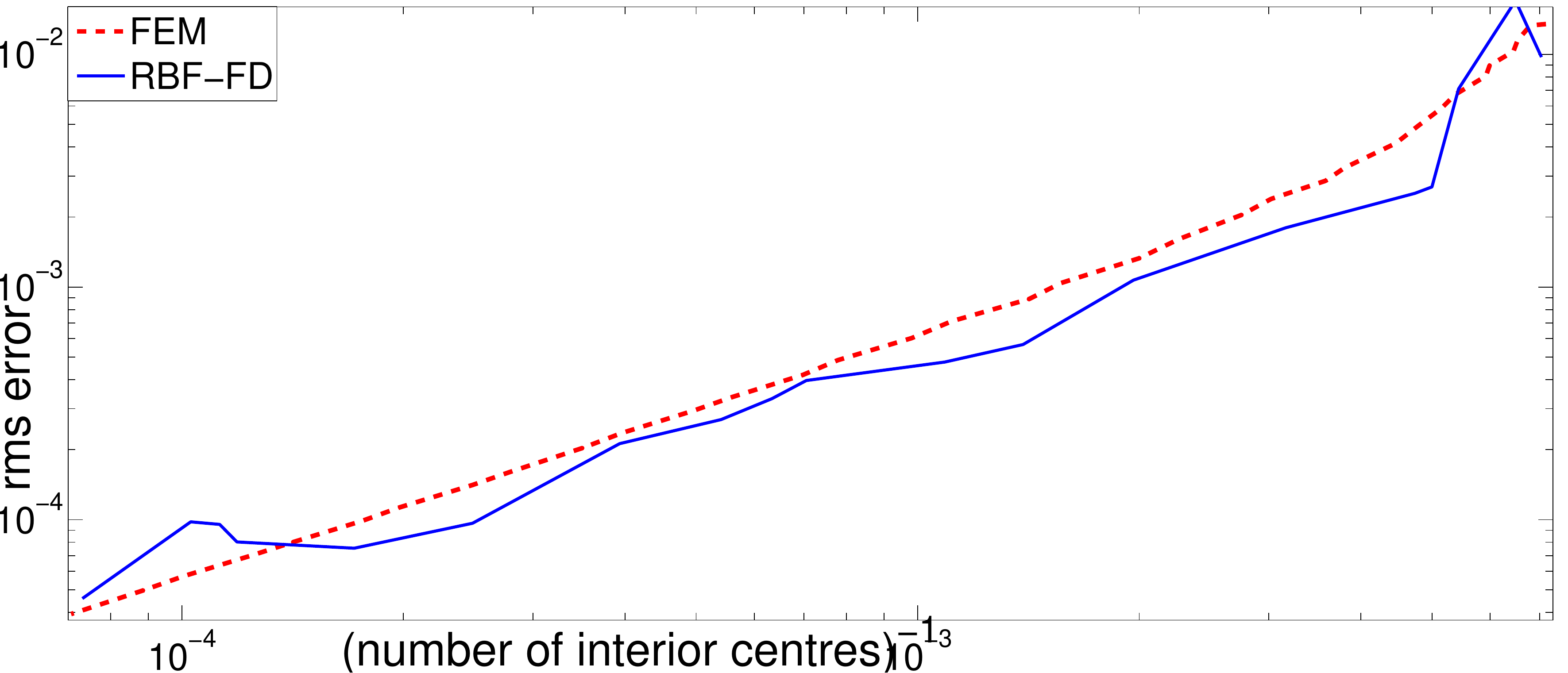}} \qquad 
 \subfigure[$\omega=2 \pi$: errors on grid]
{\includegraphics[width=6cm,height=4cm]{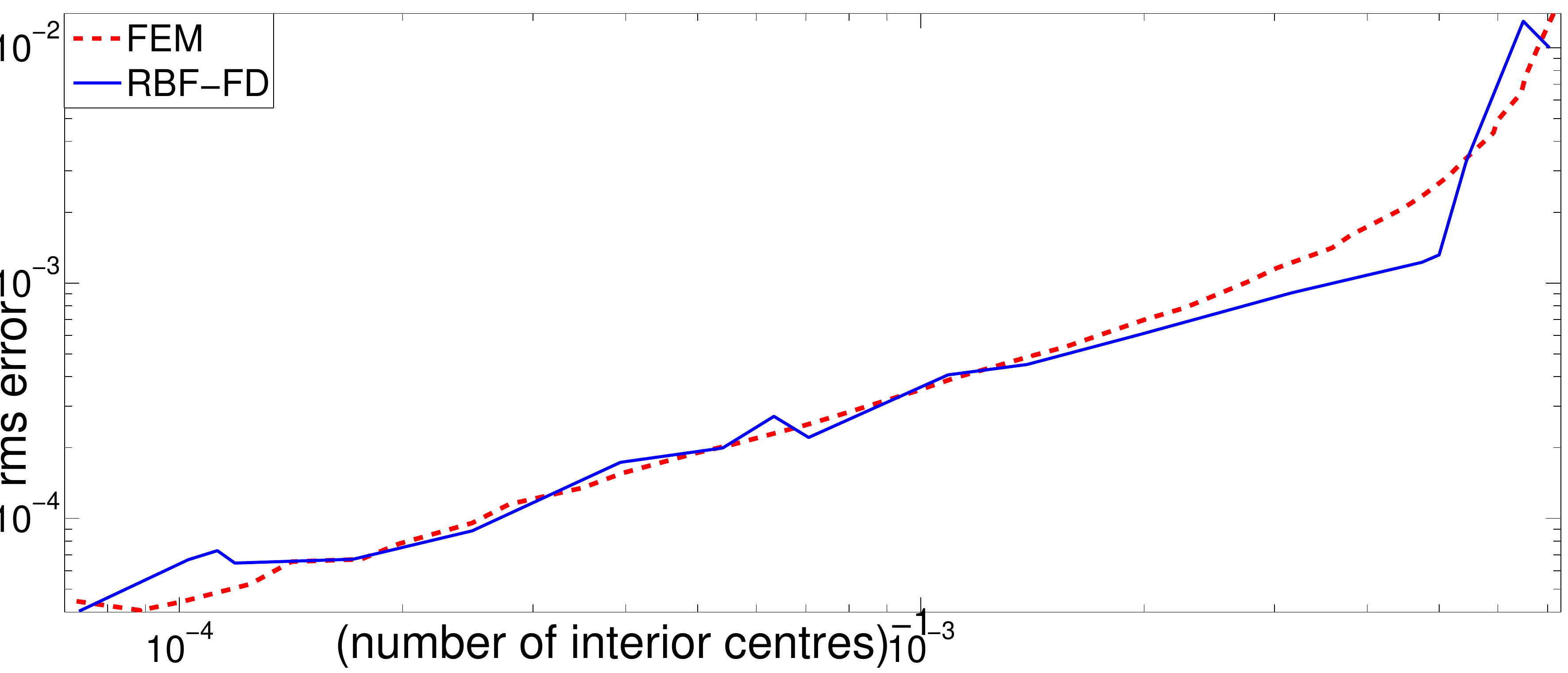}} 
        \end{center}
\caption{Test Problem~\ref{f42}: Errors for various
values of $\omega$ on the centers (left) and on a  grid (right). {\tt RBF-FD}:  the method of this paper,
{\tt FEM}: finite element method with  piecewise linear shape functions, where the solution is computed by 
using MATLAB PDE Toolbox with default parameters of the adaptive refinement.
}
\label{errF42}
\end{figure}

\begin{figure}[htbp!]
 \begin{center}
\subfigure[$\omega=\pi+0.01$: RBF-FD]
{\includegraphics[width=6cm,height=4cm]{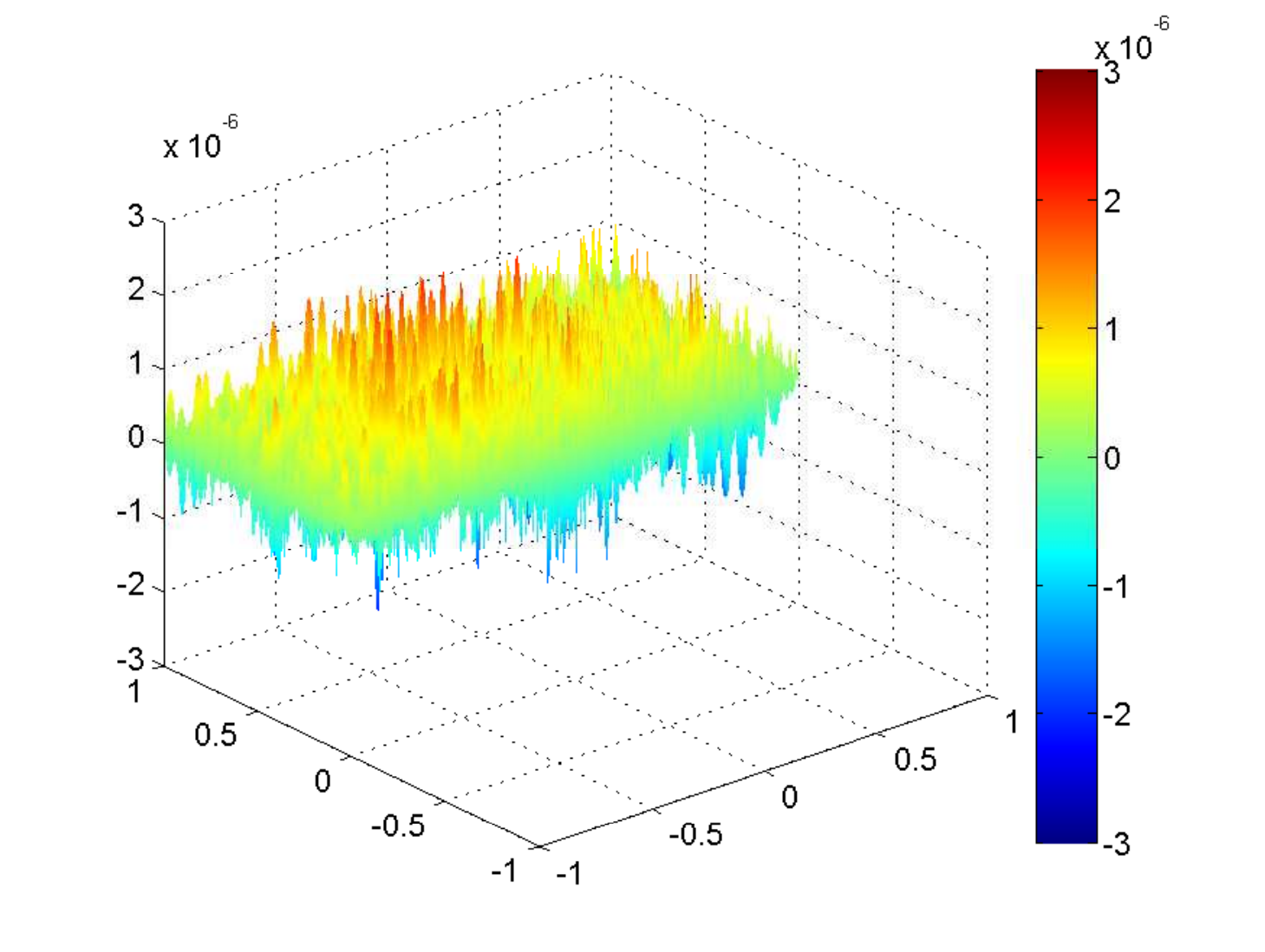}} \qquad
\subfigure[$\omega=\pi+0.01$: FEM]
{\includegraphics[width=6cm,height=4cm]{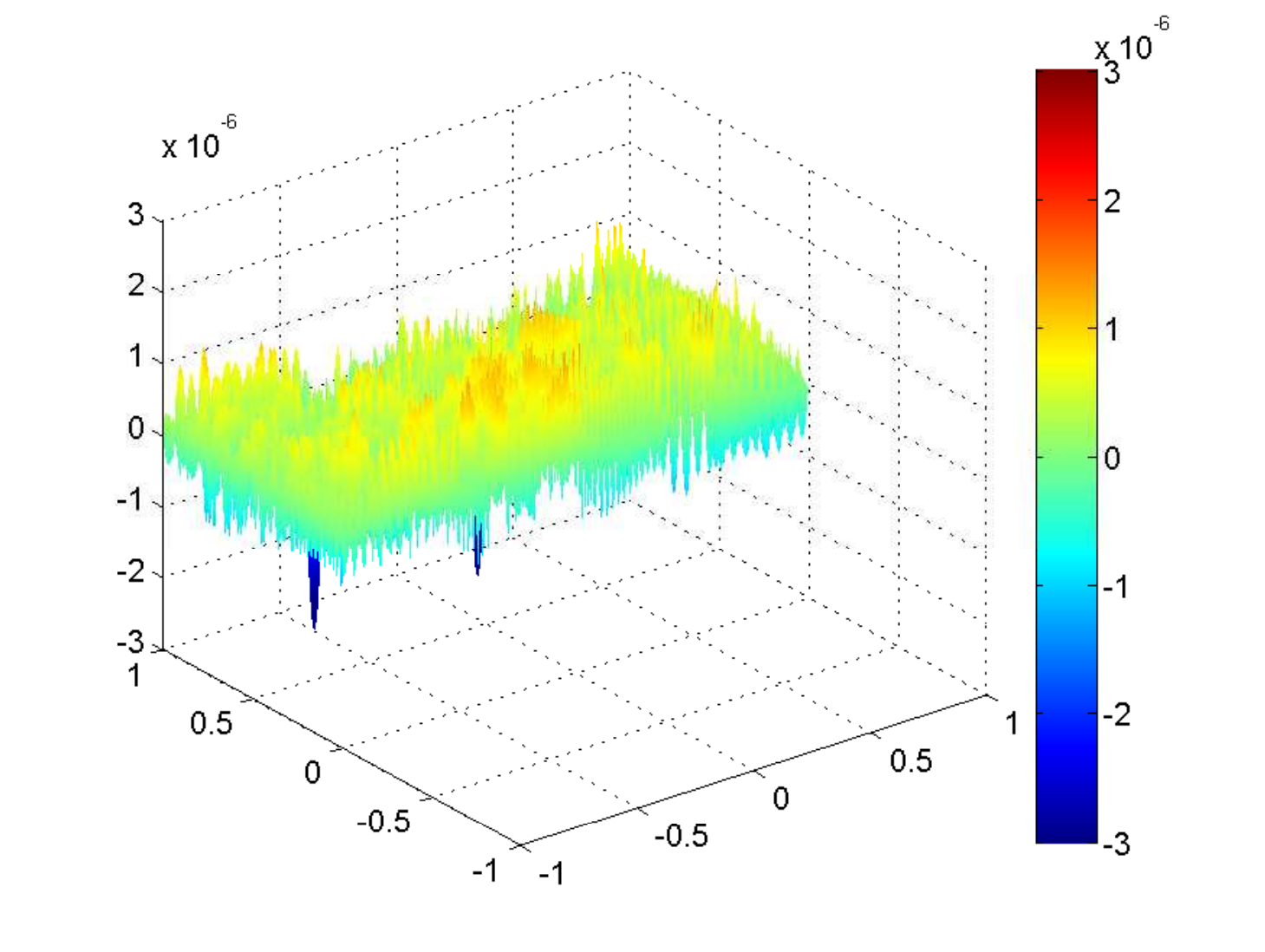}} 
 \subfigure[$\omega=5 \pi/4$: RBF-FD]
{\includegraphics[width=6cm,height=4cm]{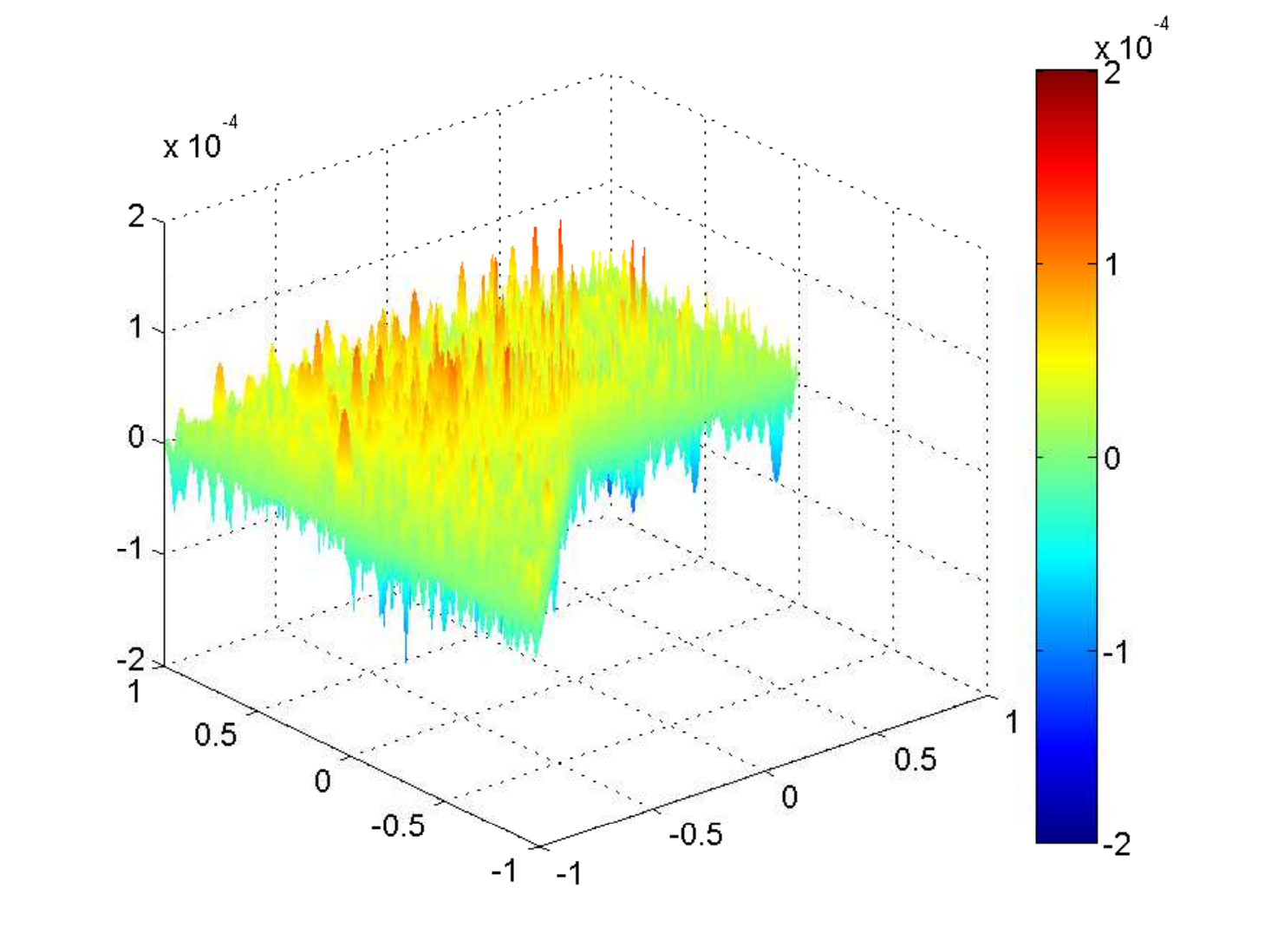}} \qquad
 \subfigure[$\omega=5 \pi/4$: FEM]
{\includegraphics[width=6cm,height=4cm]{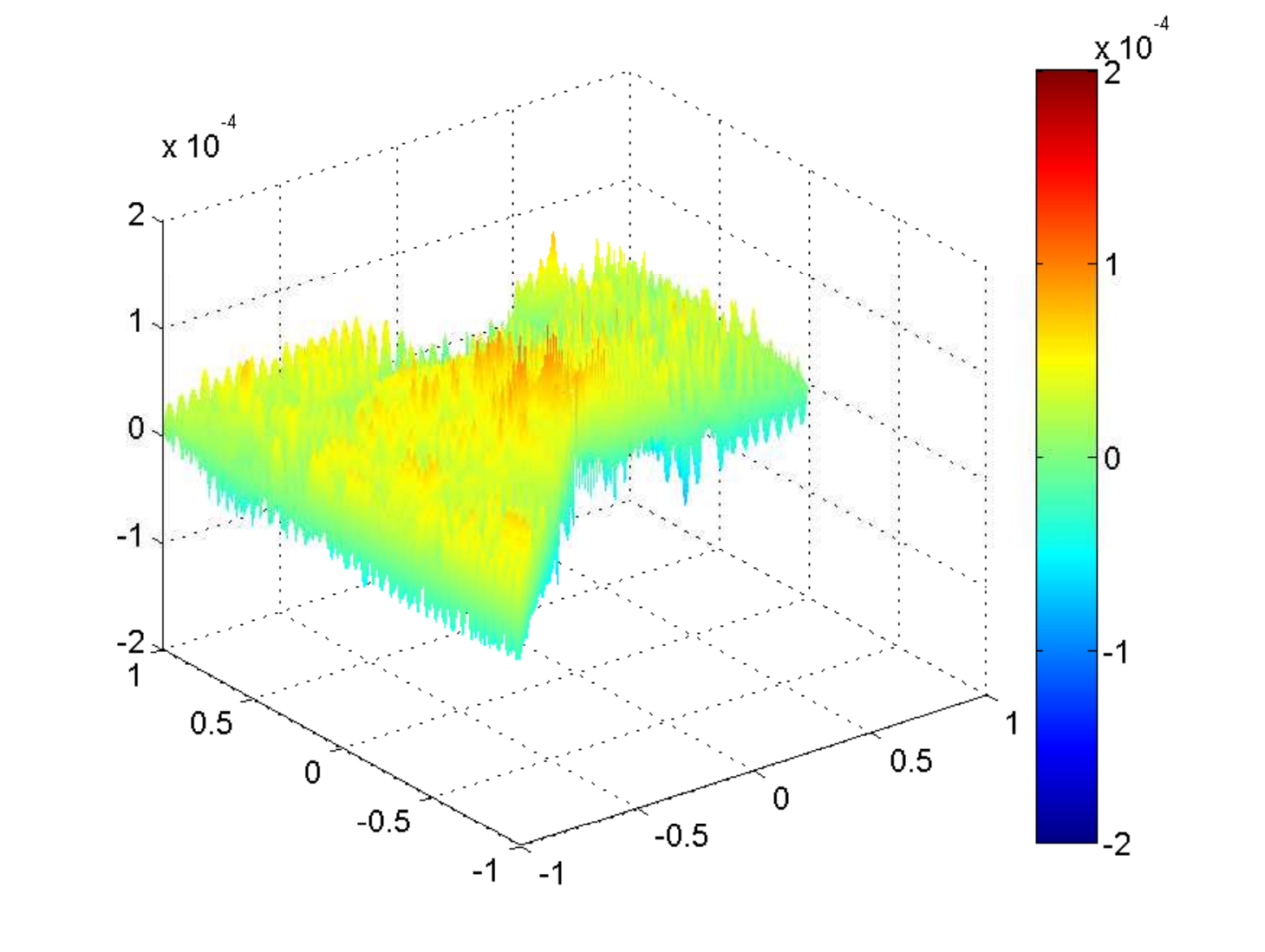}} 
  \subfigure[$\omega=7 \pi/4$: RBF-FD]
{\includegraphics[width=6cm,height=4cm]{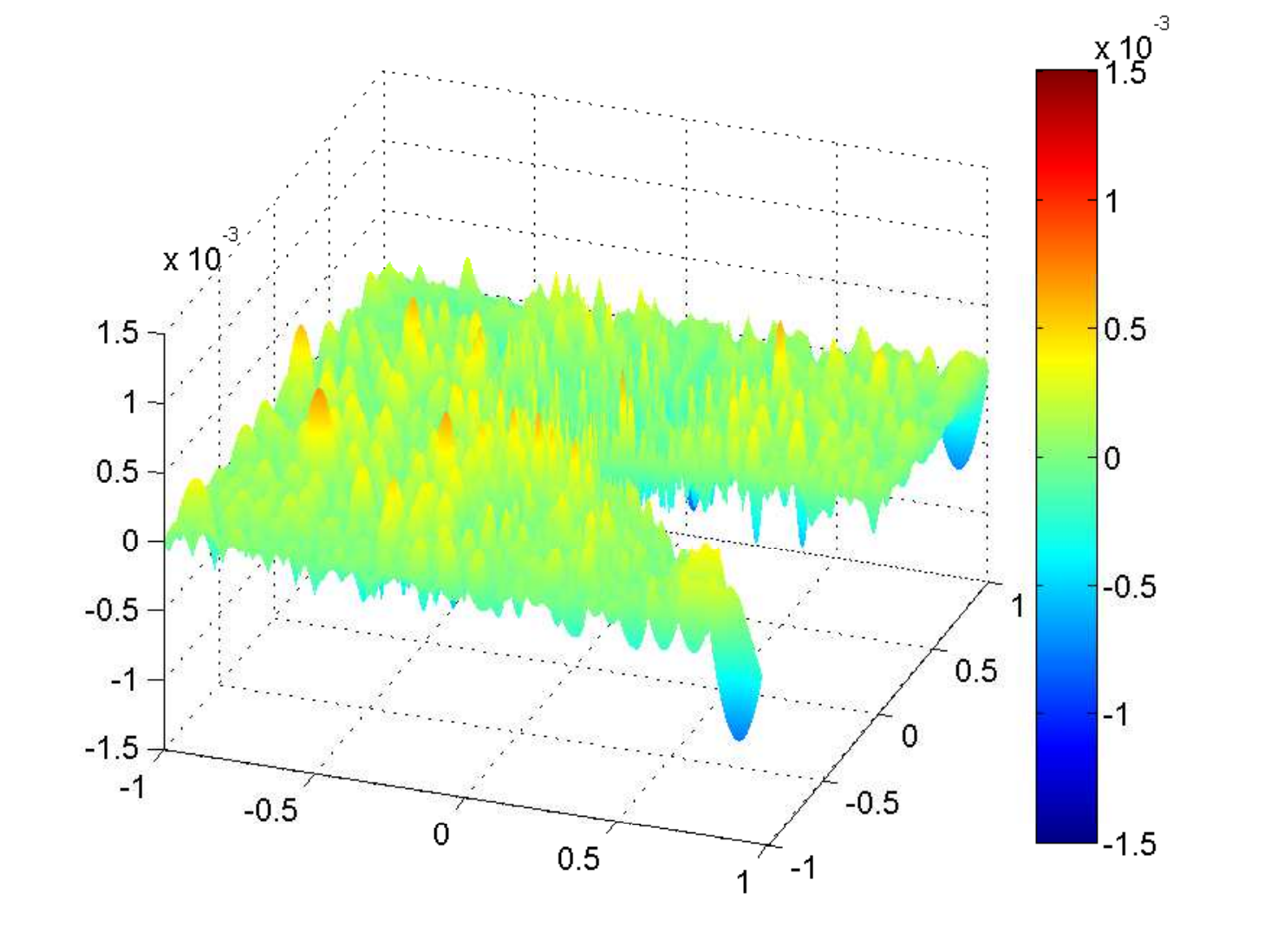}} \qquad
  \subfigure[$\omega=7 \pi/4$: FEM]
{\includegraphics[width=6cm,height=4cm]{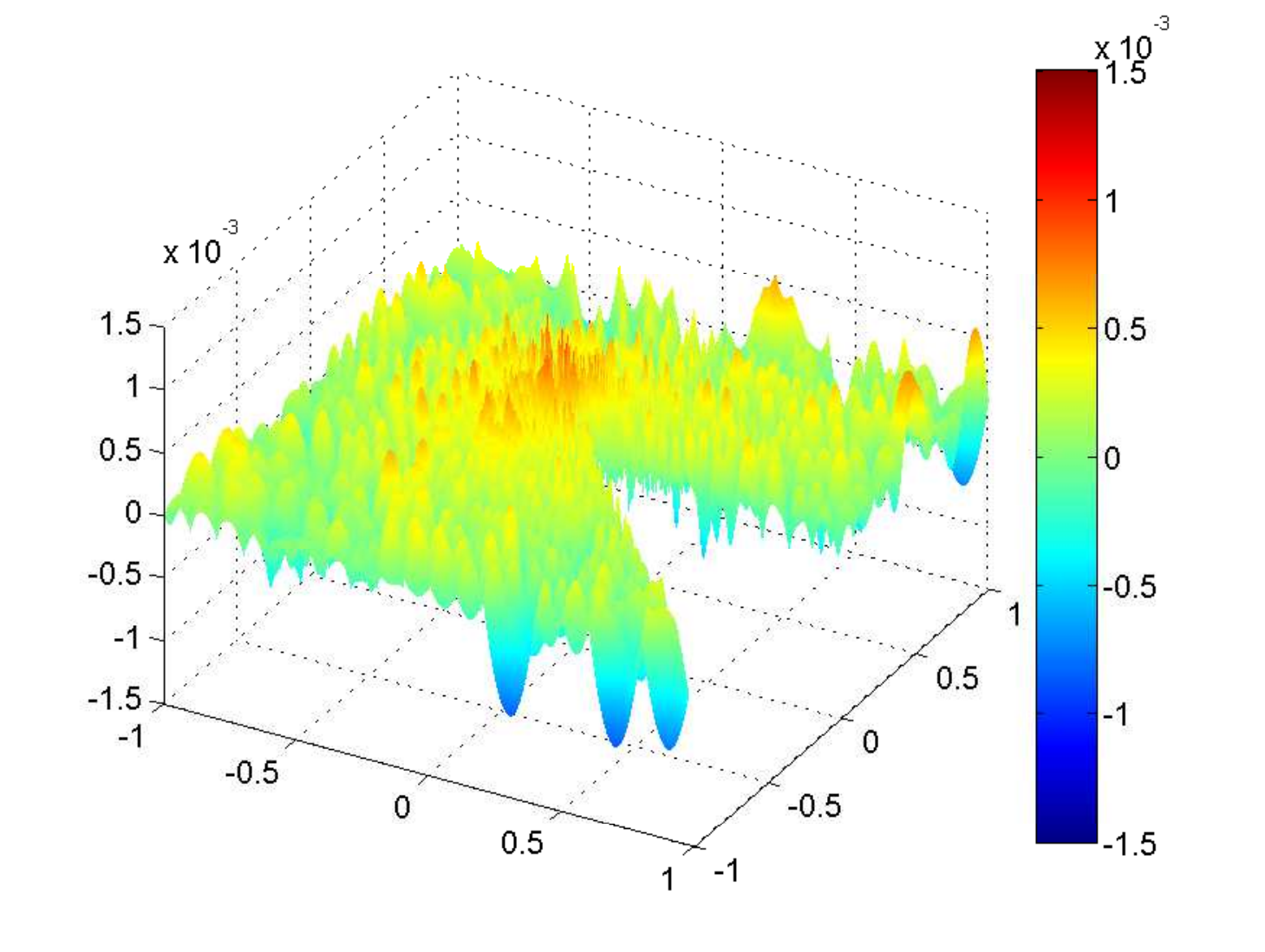}} 
 \subfigure[$\omega=2 \pi$: RBF-FD]
{\includegraphics[width=6cm,height=4cm]{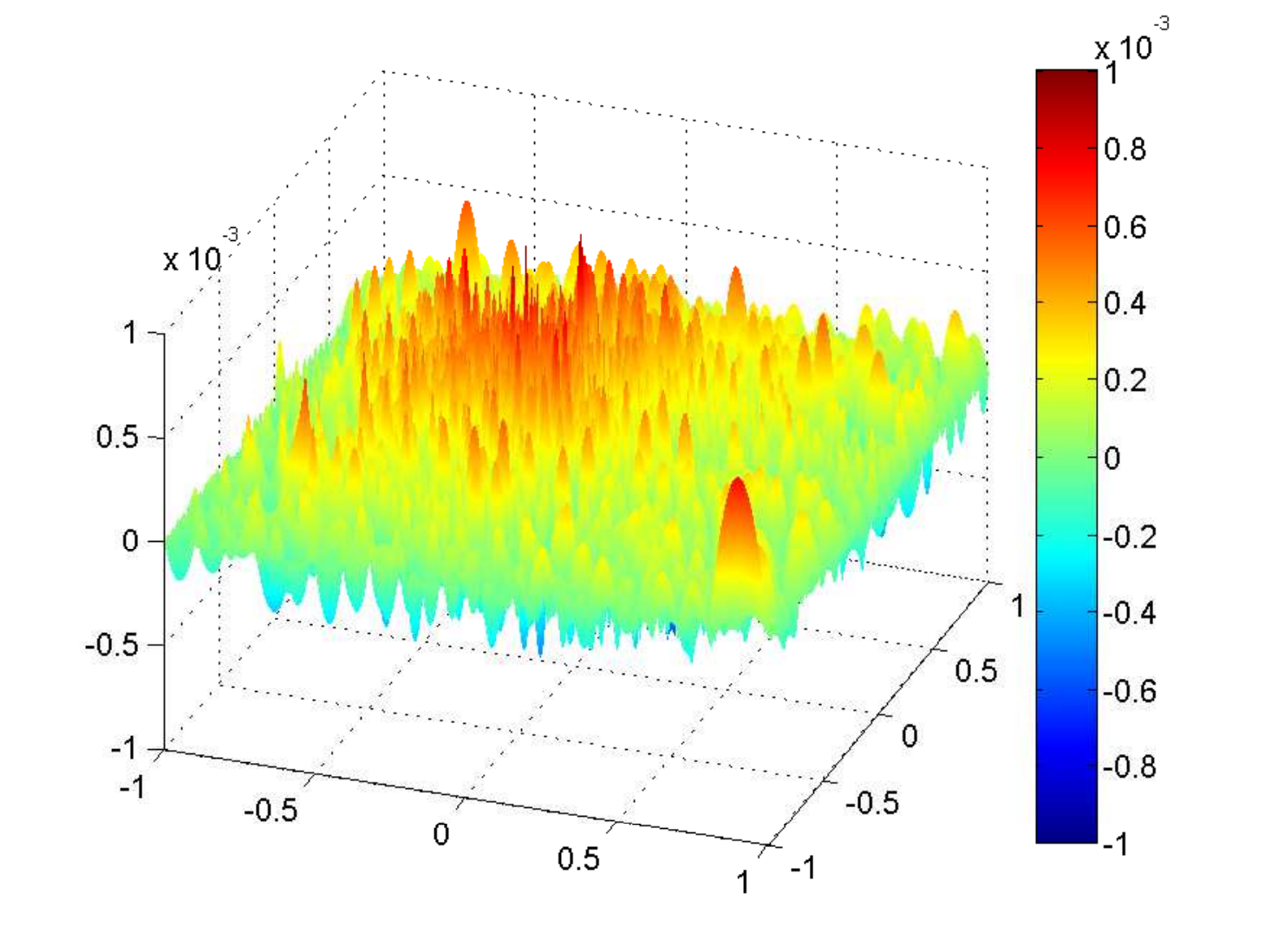}} \qquad
 \subfigure[$\omega=2 \pi$: FEM]
{\includegraphics[width=6cm,height=4cm]{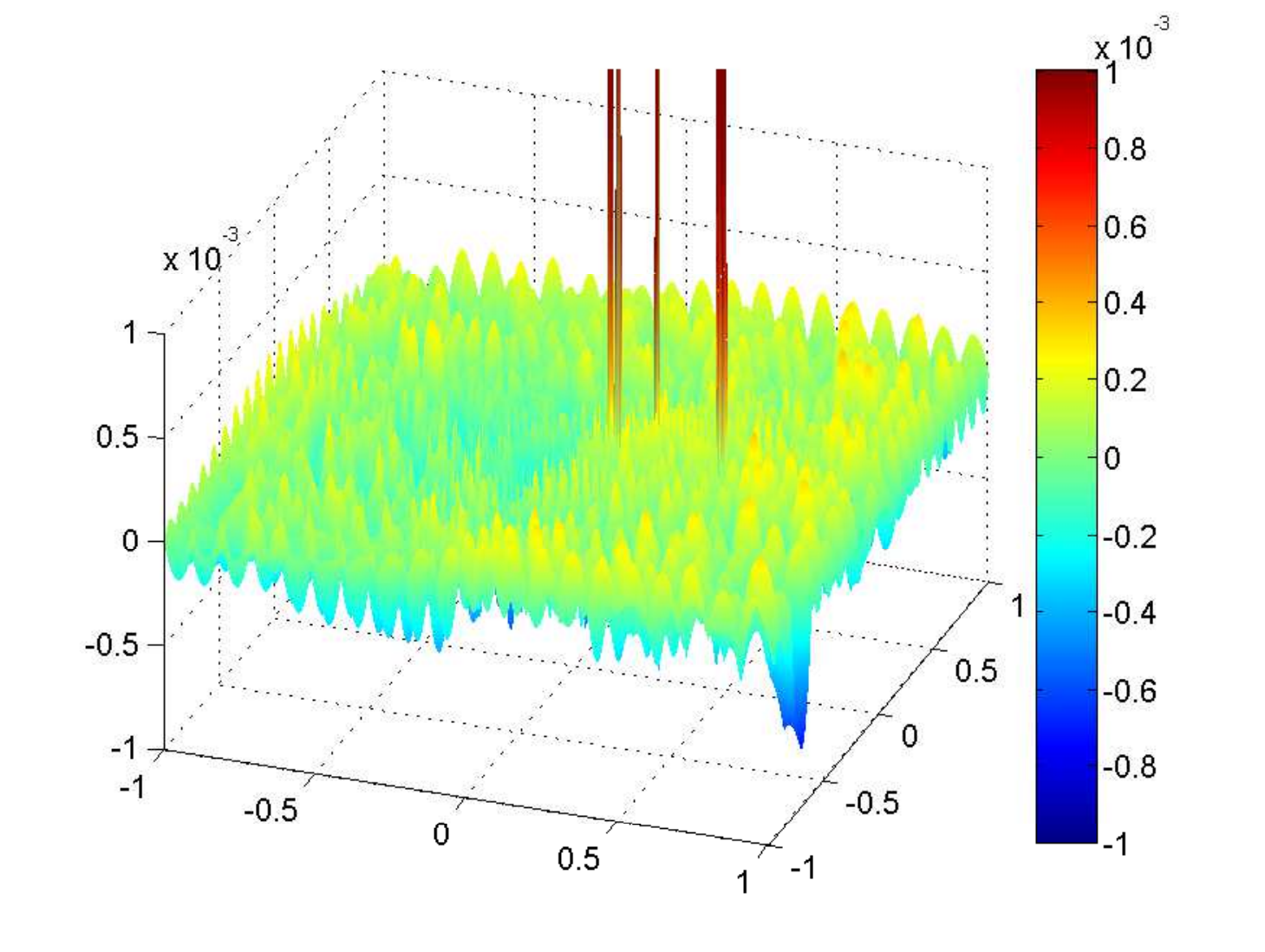}} 
        \end{center}
\caption{Test Problem~\ref{f42}: Error functions for RBF-FD (left) and FEM (right). 
The number of interior centers/vertices: 
(a) {2786}, %
(b) 2768, %
(c) {3592}, %
(d) 3478, %
(e) {1721}, %
(f) 1427, %
(g) {2553}, %
(h) 2521. %
These centers and vertices are shown in Figure~\ref{centresF42}.
}
\label{mapF42}
\end{figure}

\begin{figure}[htbp!]
  \begin{center}  
\subfigure[$\omega=\pi+0.01$: RBF-FD ]
{\includegraphics[width=6cm,height=4cm]{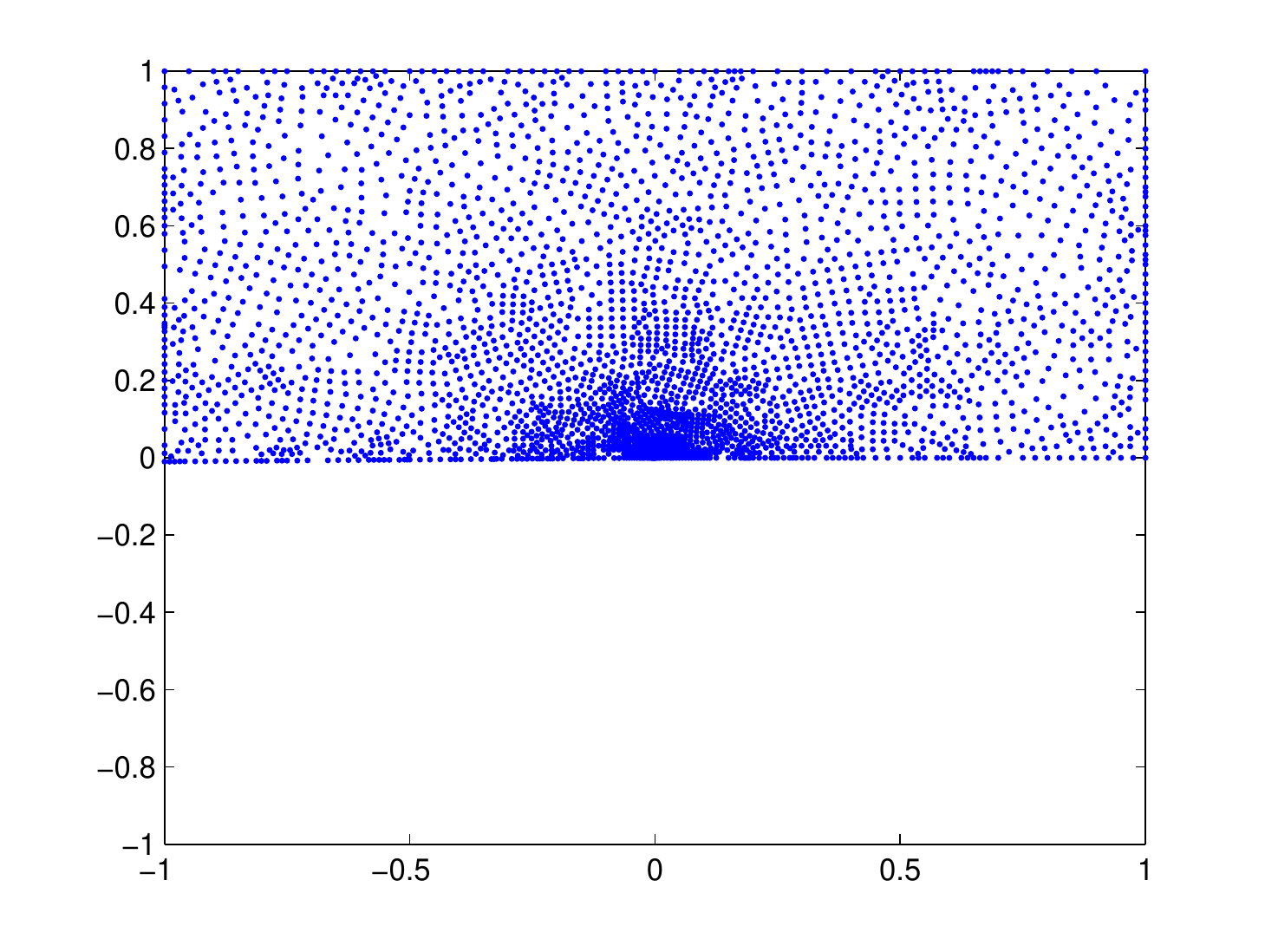}} \qquad
\subfigure[$\omega=\pi+0.01$: FEM]
{\includegraphics[width=6cm,height=4cm]{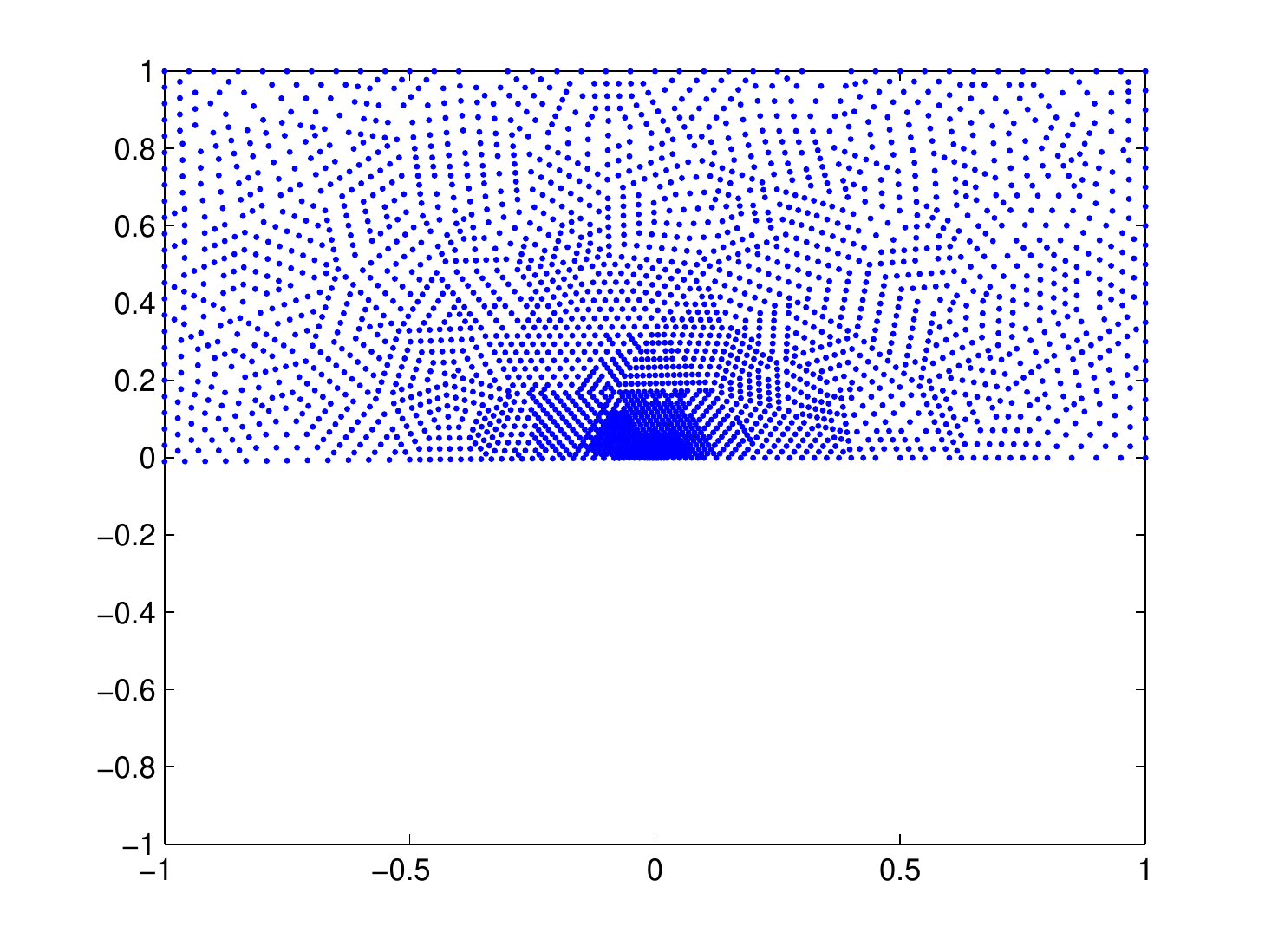}} 
 \subfigure[$\omega=5 \pi/4$: RBF-FD]
{\includegraphics[width=6cm,height=4cm]{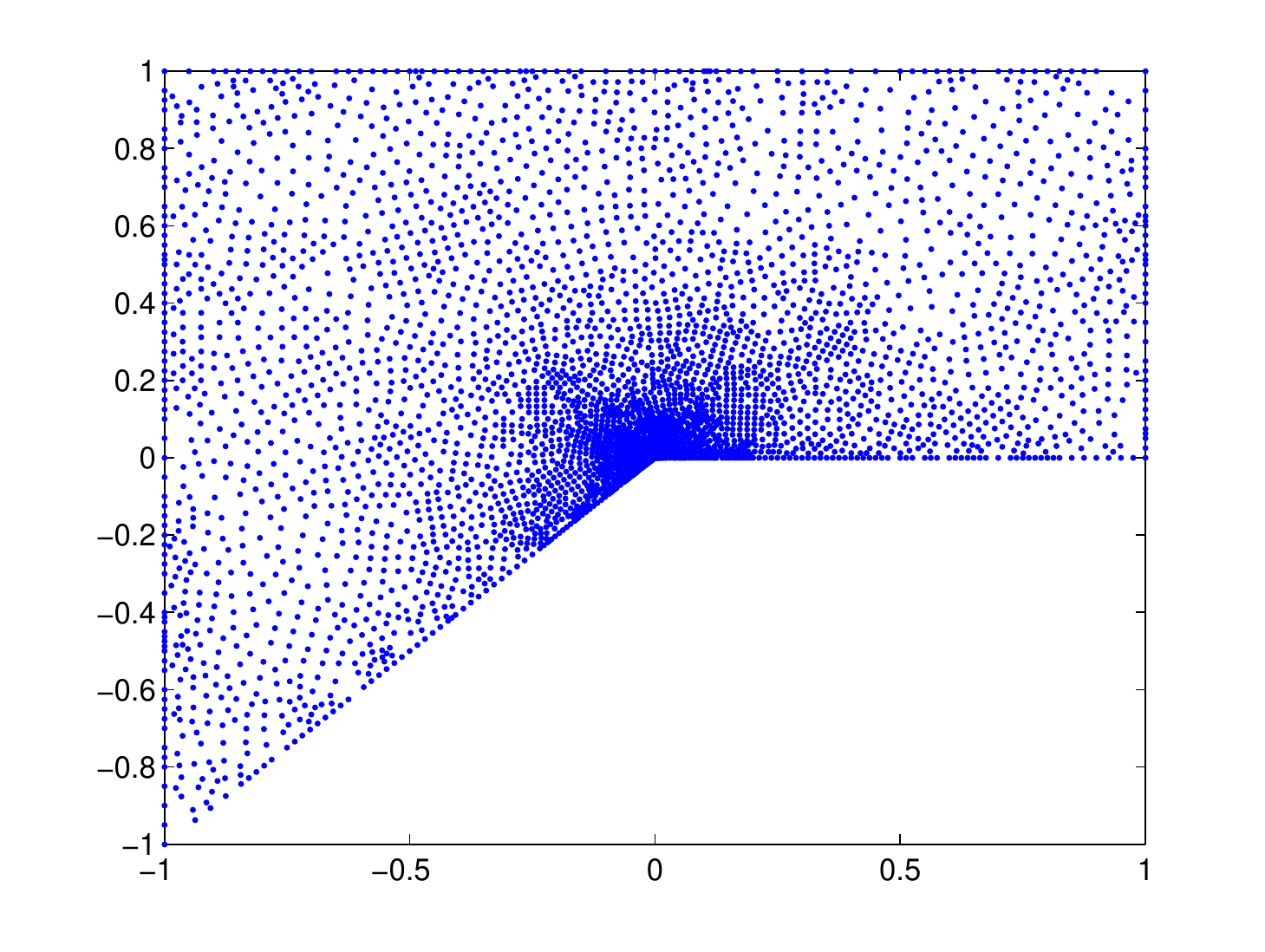}} \qquad
 \subfigure[$\omega=5 \pi/4$: FEM ]
{\includegraphics[width=6cm,height=4cm]{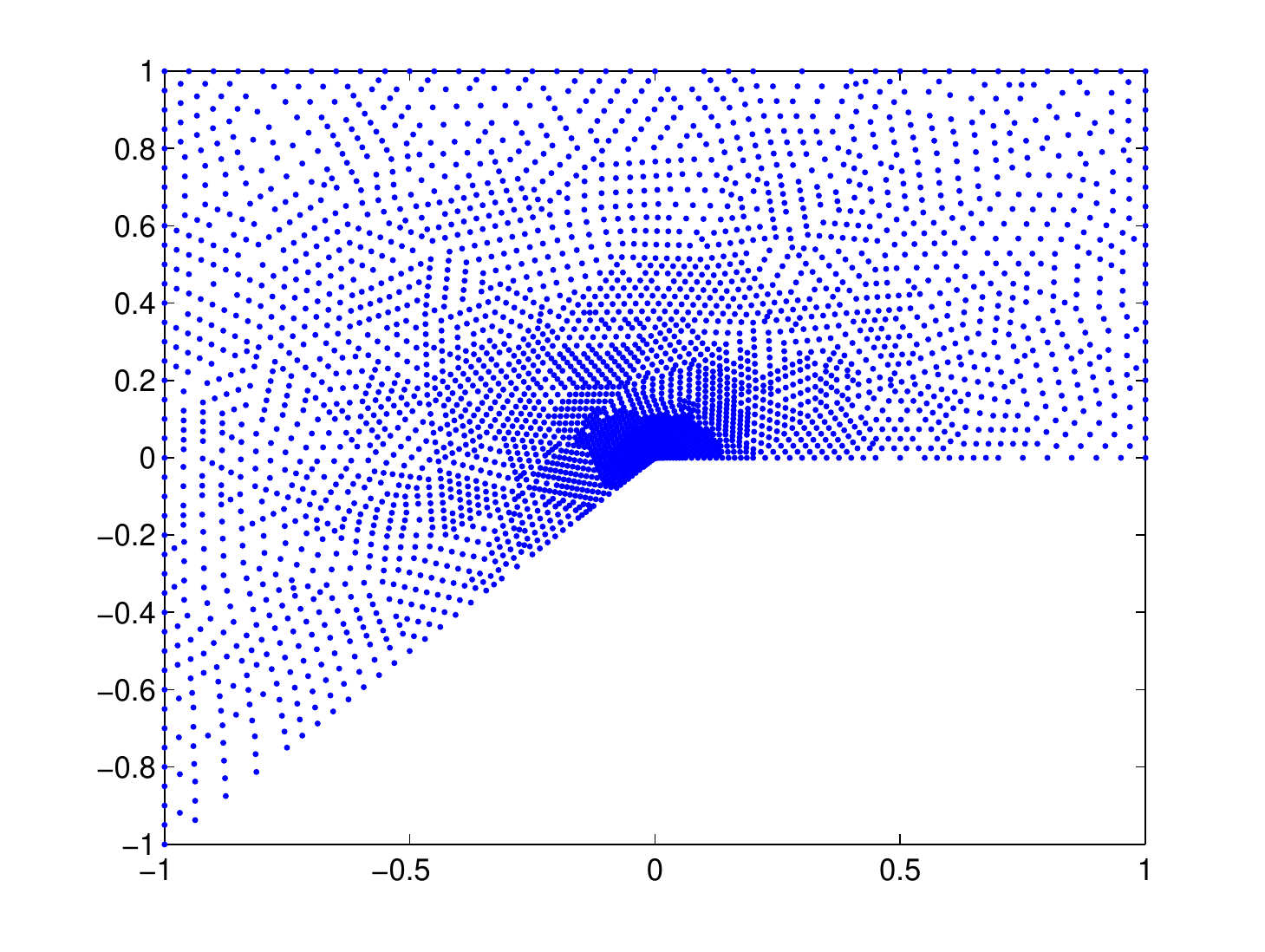}} 
  \subfigure[$\omega=7 \pi/4$: RBF-FD]
{\includegraphics[width=6cm,height=4cm]{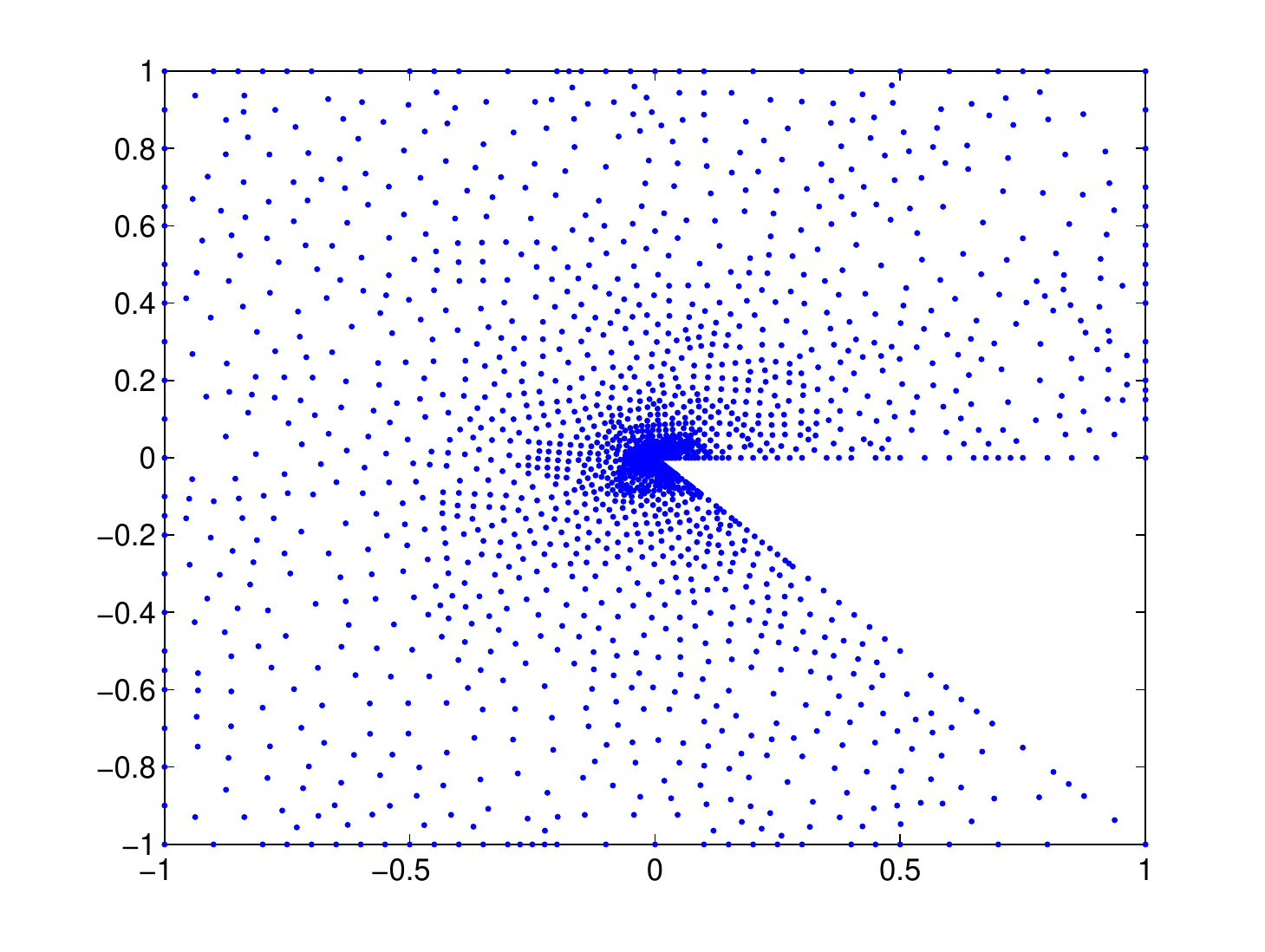}} \qquad
  \subfigure[$\omega=7 \pi/4$: FEM]
{\includegraphics[width=6cm,height=4cm]{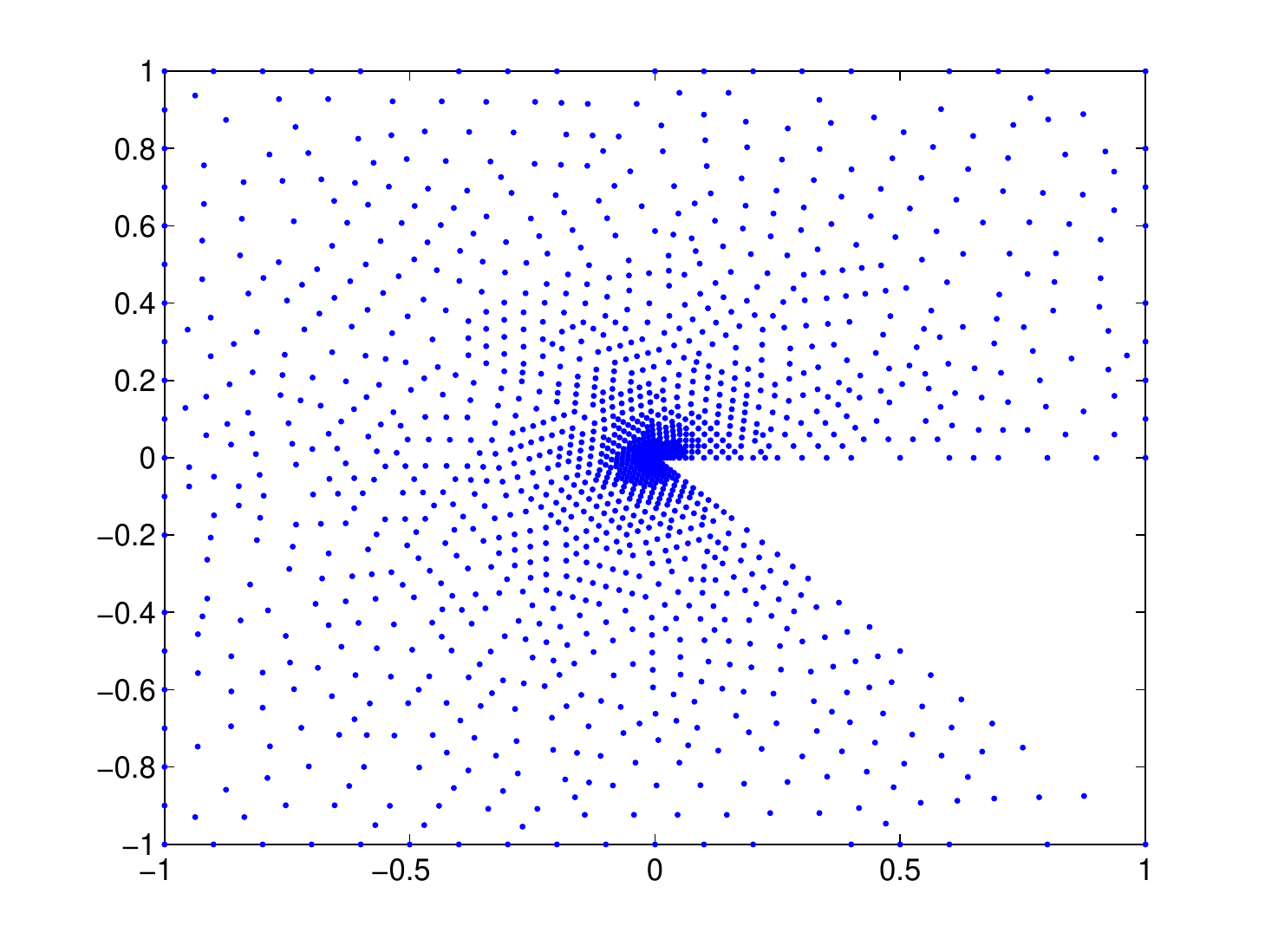}} 
 \subfigure[$\omega=2 \pi$: RBF-FD ]
{\includegraphics[width=6cm,height=4cm]{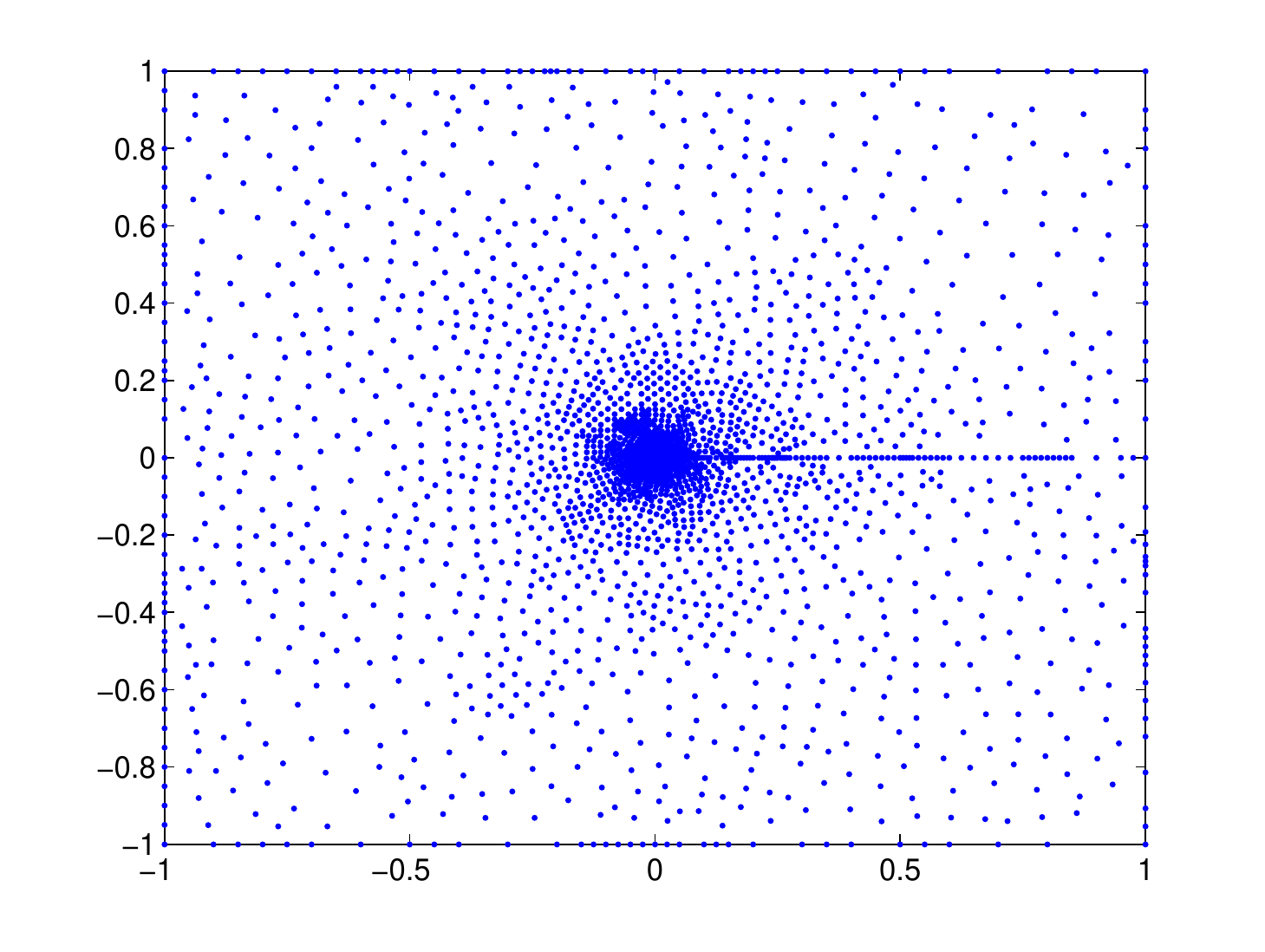}} \qquad
  \subfigure[$\omega=2 \pi$: FEM]
{\includegraphics[width=6cm,height=4cm]{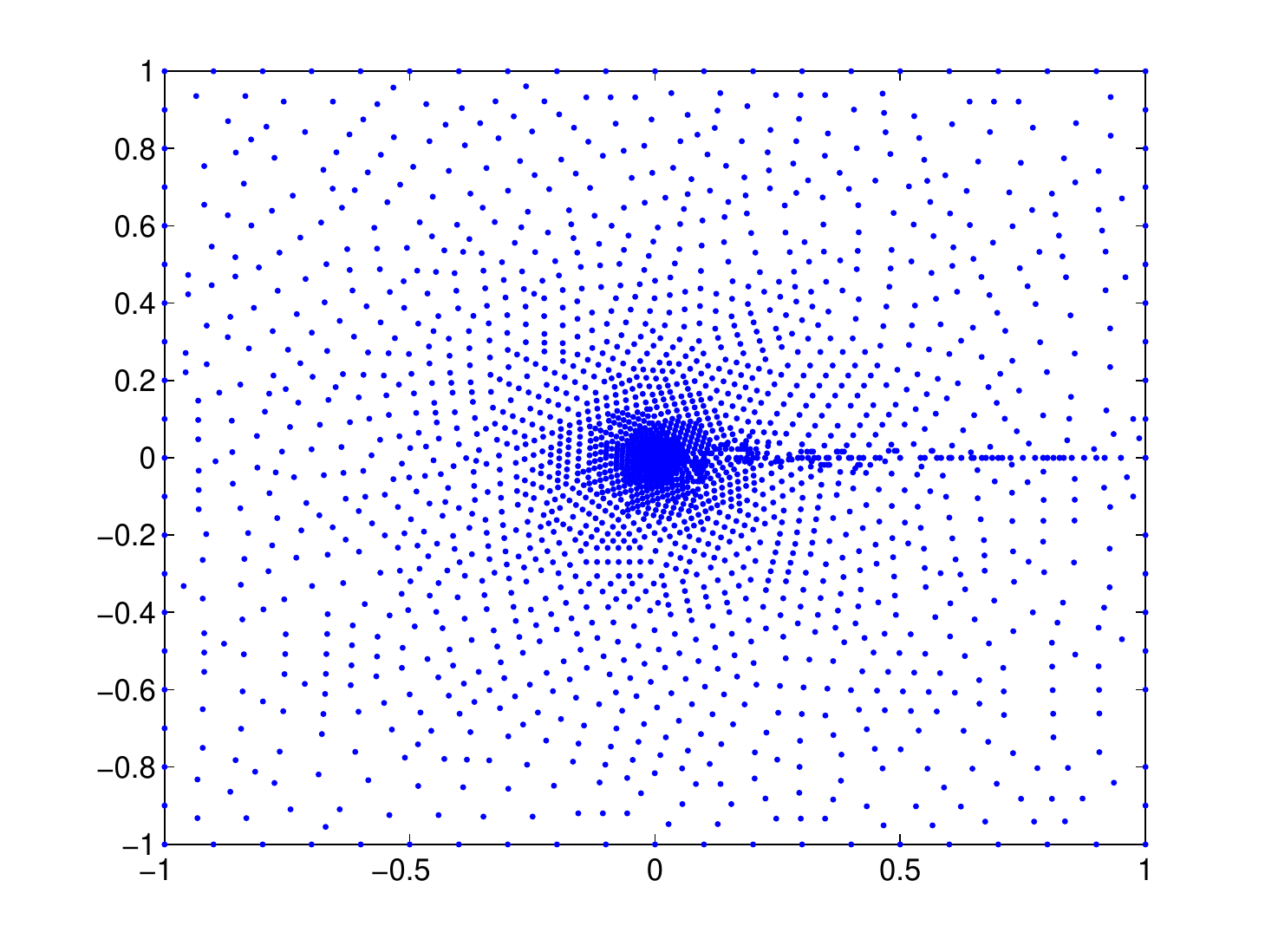}} 
     \end{center}
\caption{Test Problems~\ref{f42}: Centers generated by the adaptive RBF-FD method (left) and
the vertices of the triangulations generated by the adaptive FEM (right) for the solutions whose error plots are
shown in Figure~\ref{mapF42}.
} 
 \label{centresF42}
\end{figure}

 Numerical results for this problem are presented in Figures~\ref{errF42}--\ref{centresF42} 
which show for $\omega=\pi+0.01, 5\pi/4, 7\pi/4,2\pi$ the graphs of the rms errors of RBF-FD and FEM solutions on the centers and on a grid, 
the error functions for both methods, and the distribution of the centers. We do not consider 
$\omega=3\pi/2$ because we already had this reentrant corner in Test Problem~\ref{f4}. These results are in 
agreement with the above observations for Test Problems~\ref{f4} and \ref{f8}. 
Note that for the slit
domain with $\omega=2 \pi$ we removed some centers located too close to the slit in the initial
triangulation generated by PDE Toolbox. It would have been cumbersome to also modify the triangulation for the FEM
because the inner structures of the triangulation would need to be repaired. 
Several rather badly shaped triangles in the initial triangulation is the most plausible reason for the
 distortion of the error of the FEM solution near the slit tip, leading to the sharp peak in the error seen in
Figure~\ref{mapF42}(h). (We plot the function $\hu-u$ rather than $u-\hu$ in 
\ref{mapF42}(gh) to make the peaks clearly visible.)

\begin{testproblem}\label{curved_slit}\rm (Curved Slit)
Dirichlet problem for the Laplace equation $\Delta u=0$ 
in the domain $\Omega$ obtained from $(-1,1)^2$ by removing the arc 
of the circle with center at $(1,-0.75)$ and radius $1.25$ between the points $(0,0)$ and
$(1,0.5)$.
The boundary conditions are chosen such that the exact solution is 
$u(x,y)=\operatorname{Re}\sqrt{(3-4i)z/(z-2)}$, with $z=x+iy$.
\end{testproblem}

\begin{figure}[htbp!]
 \begin{center}
 \subfigure%
{\includegraphics[width=6cm,height=4cm]{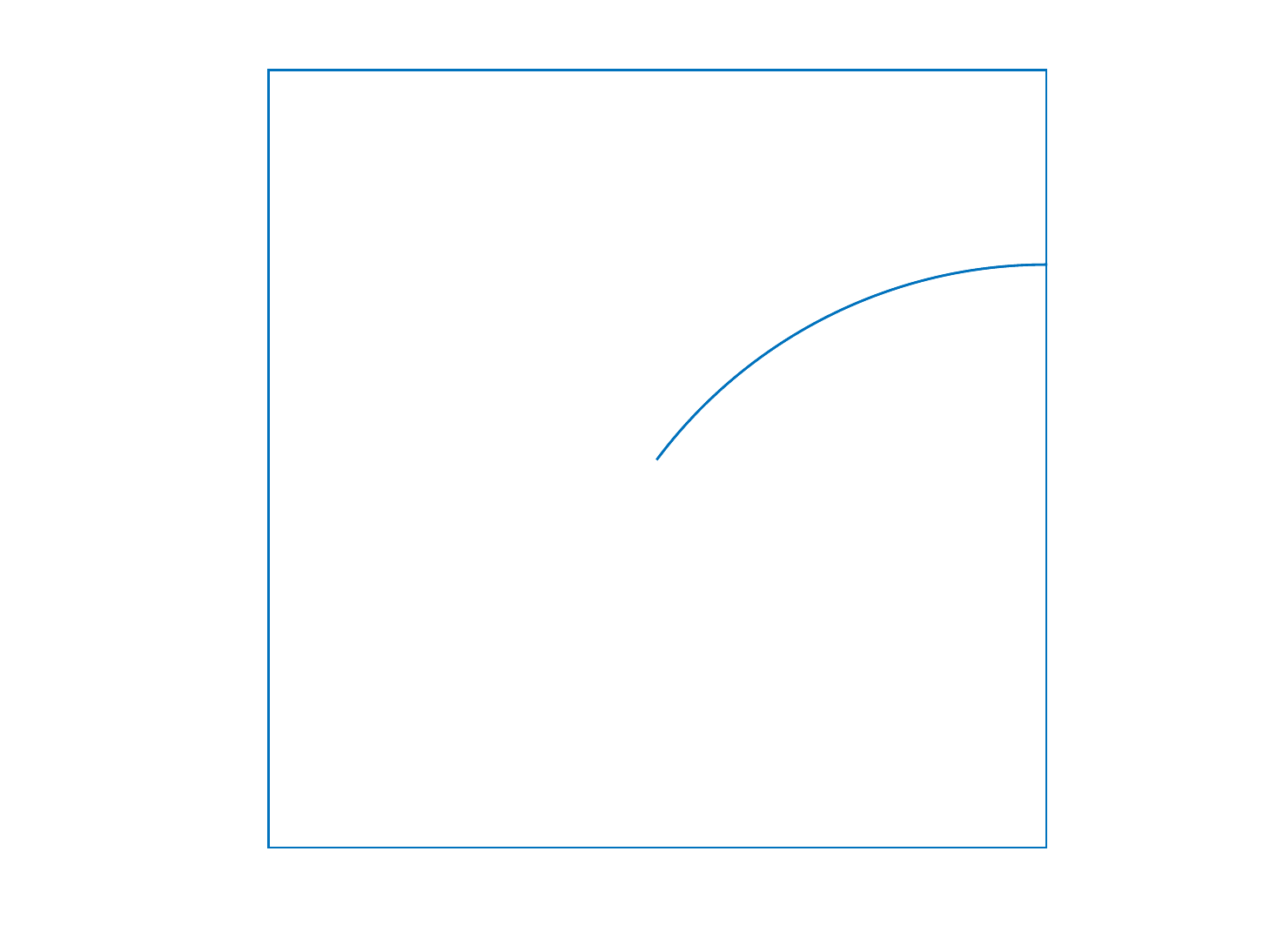}}  \qquad
  \subfigure%
{\includegraphics[width=6cm,height=4cm]{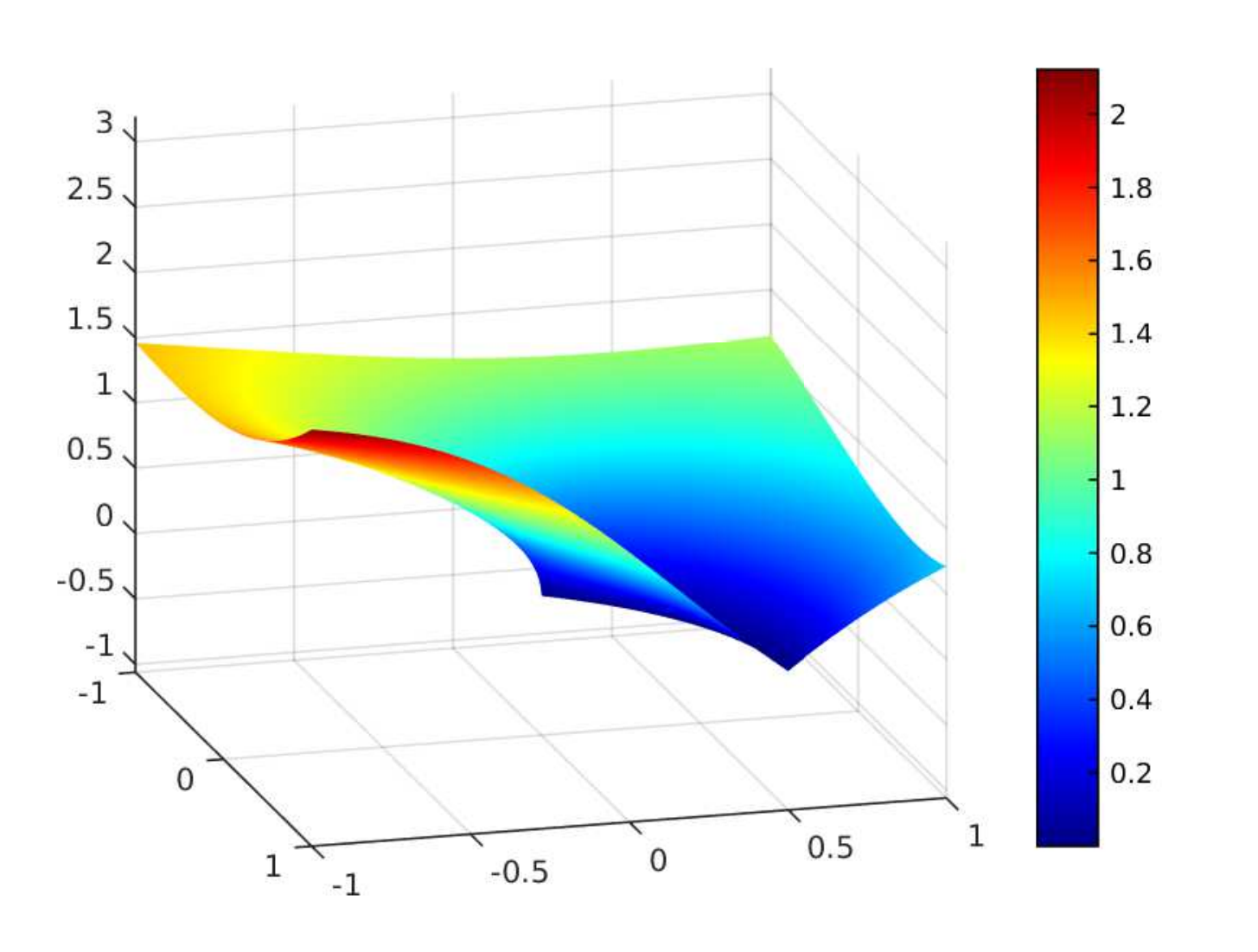}}
        \end{center}
\caption{Test Problem~\ref{curved_slit}: Domain with a curved slit (left) and exact solution (right).}
\label{CSexact}
\end{figure}

\begin{figure}[htbp!]
 \begin{center}
 \subfigure[Errors on centers]%
{\includegraphics[width=6cm,height=4cm]{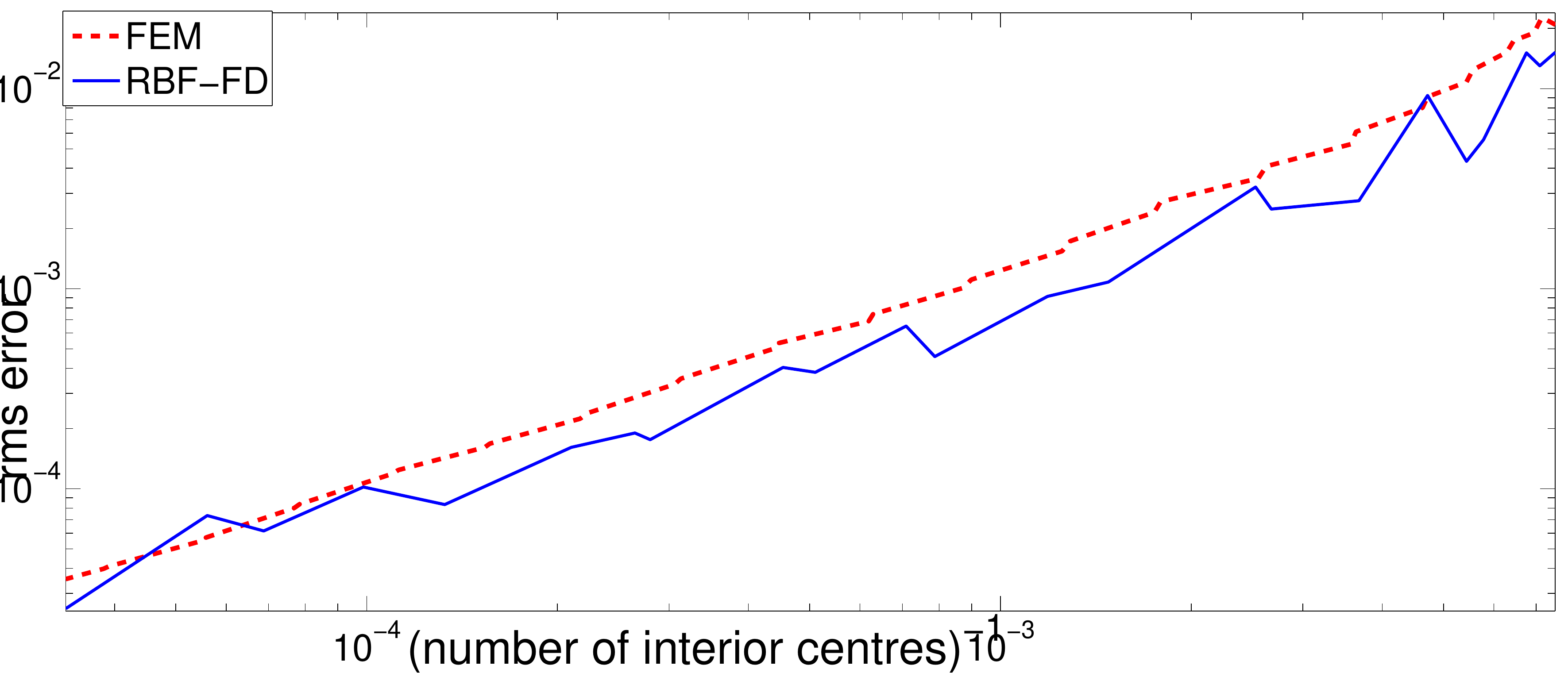}} \qquad 
 \subfigure[Errors on grid]
{\includegraphics[width=6cm,height=4cm]{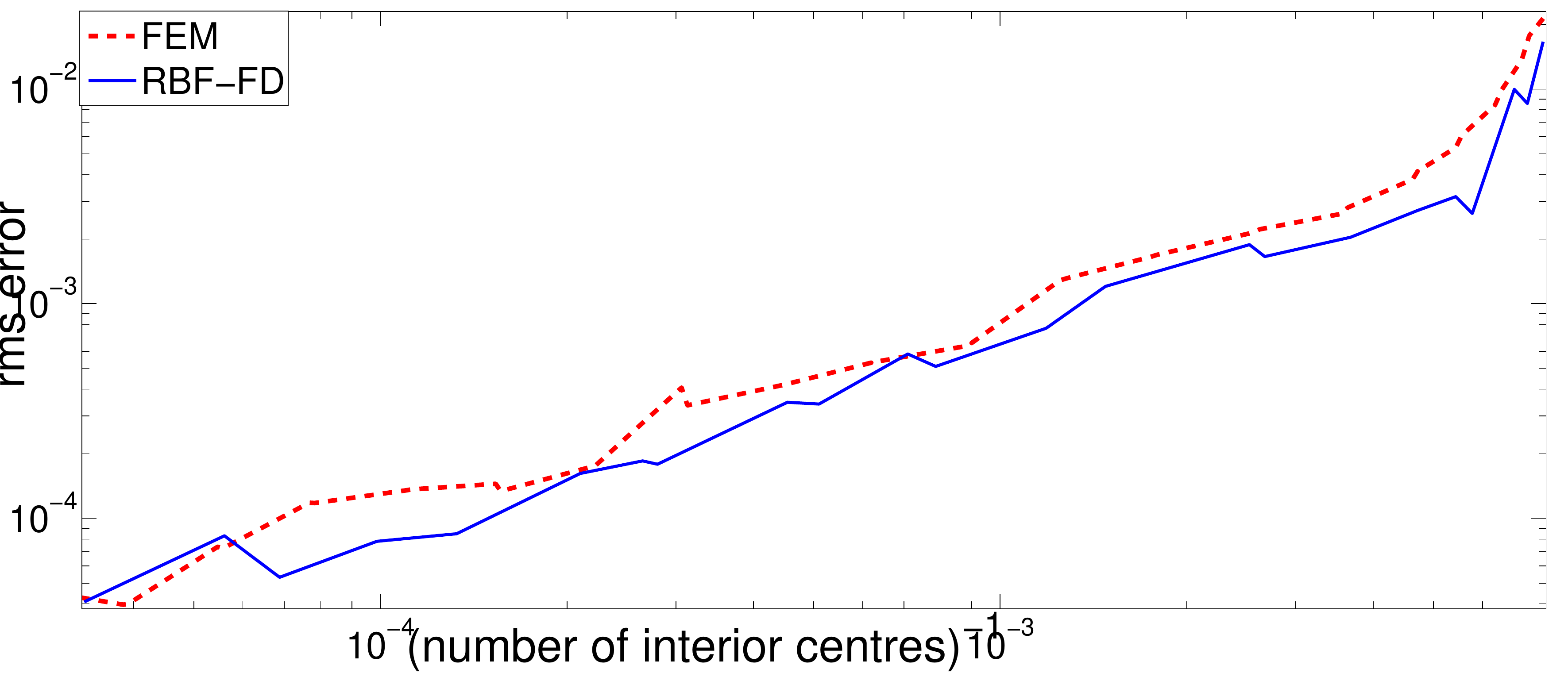}} 
 \subfigure[RBF-FD error (2204) ]
{\includegraphics[width=6cm,height=4cm]{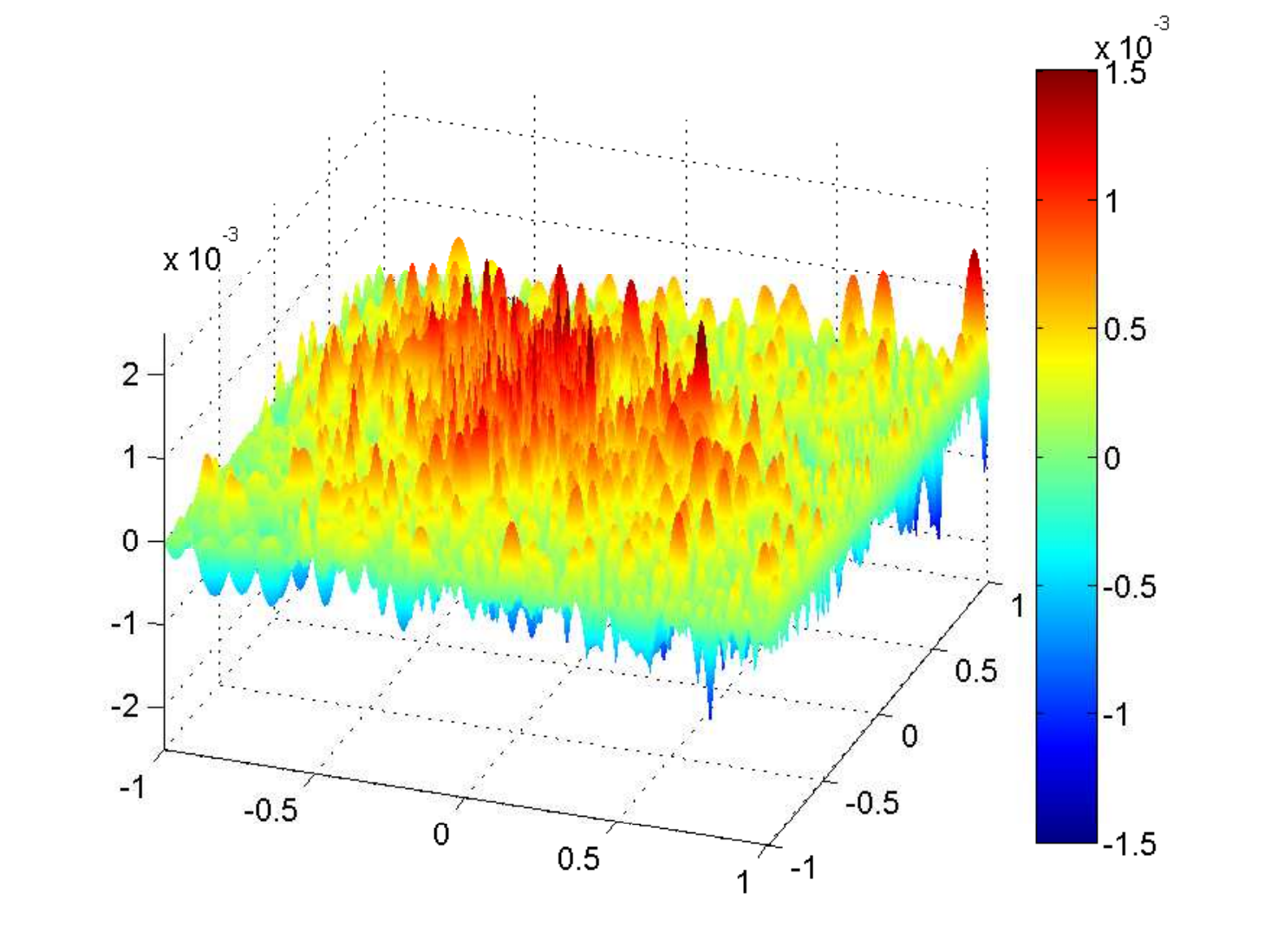}}\qquad 
 \subfigure[FEM error (2236)] 
{\includegraphics[width=6cm,height=4cm]{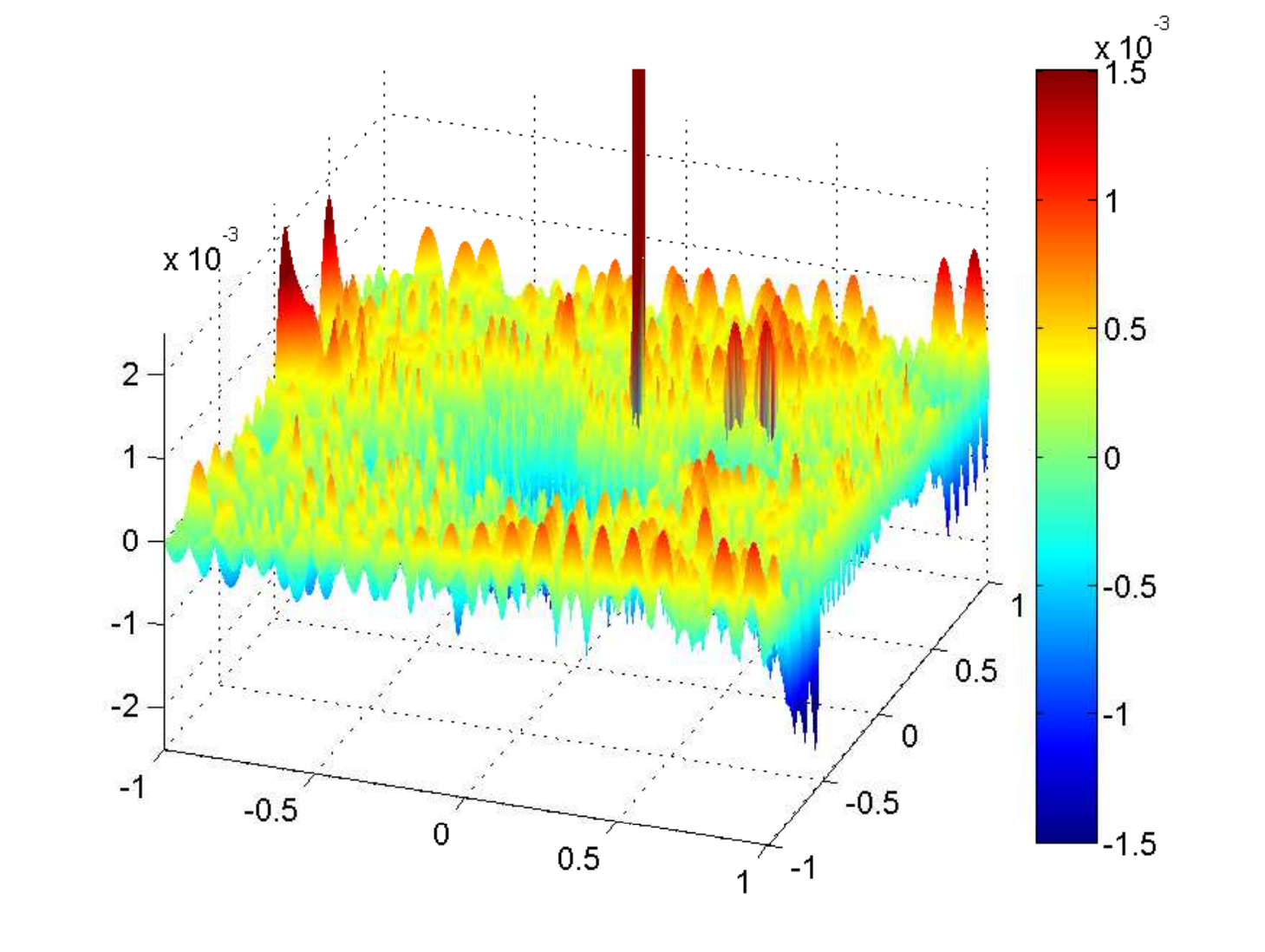}} 
 \subfigure[RBF-FD centers (2204) ]
{\includegraphics[width=6cm,height=4cm]{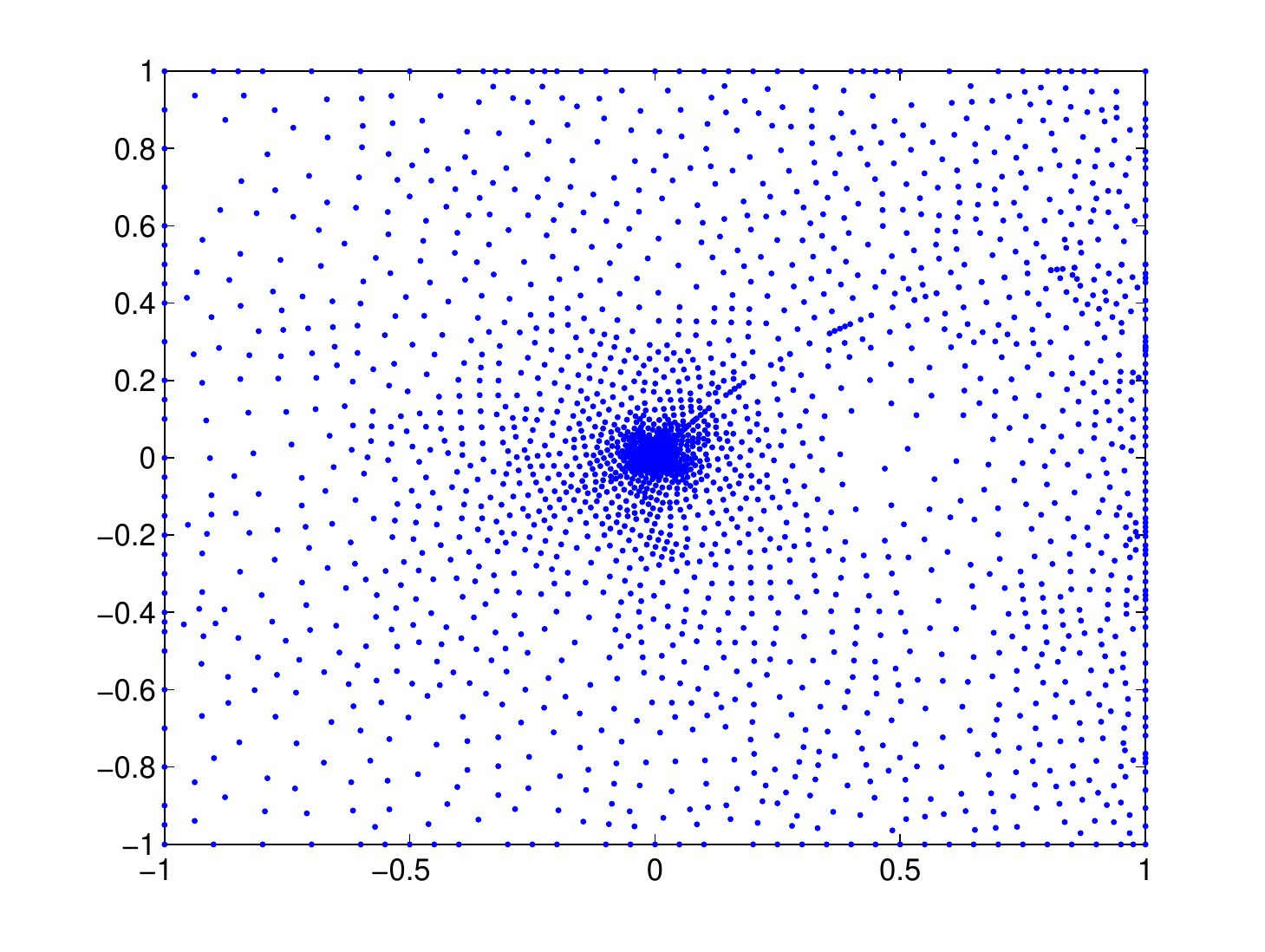}} \qquad 
 \subfigure[FEM centers (2236)]
{\includegraphics[width=6cm,height=4cm]{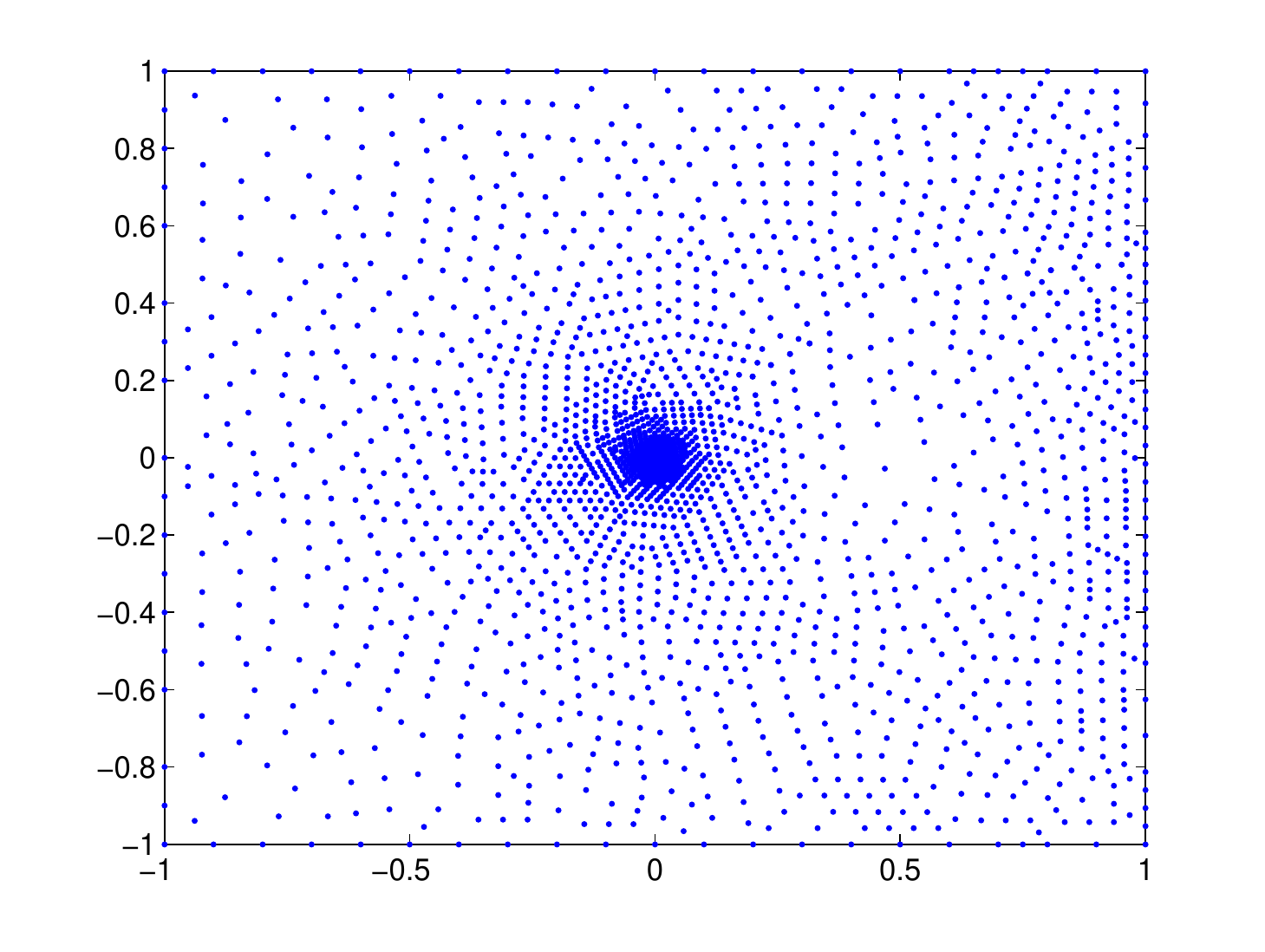}} 
 \subfigure[RBF-FD centers: zoom ]
{\includegraphics[width=6cm,height=4cm]{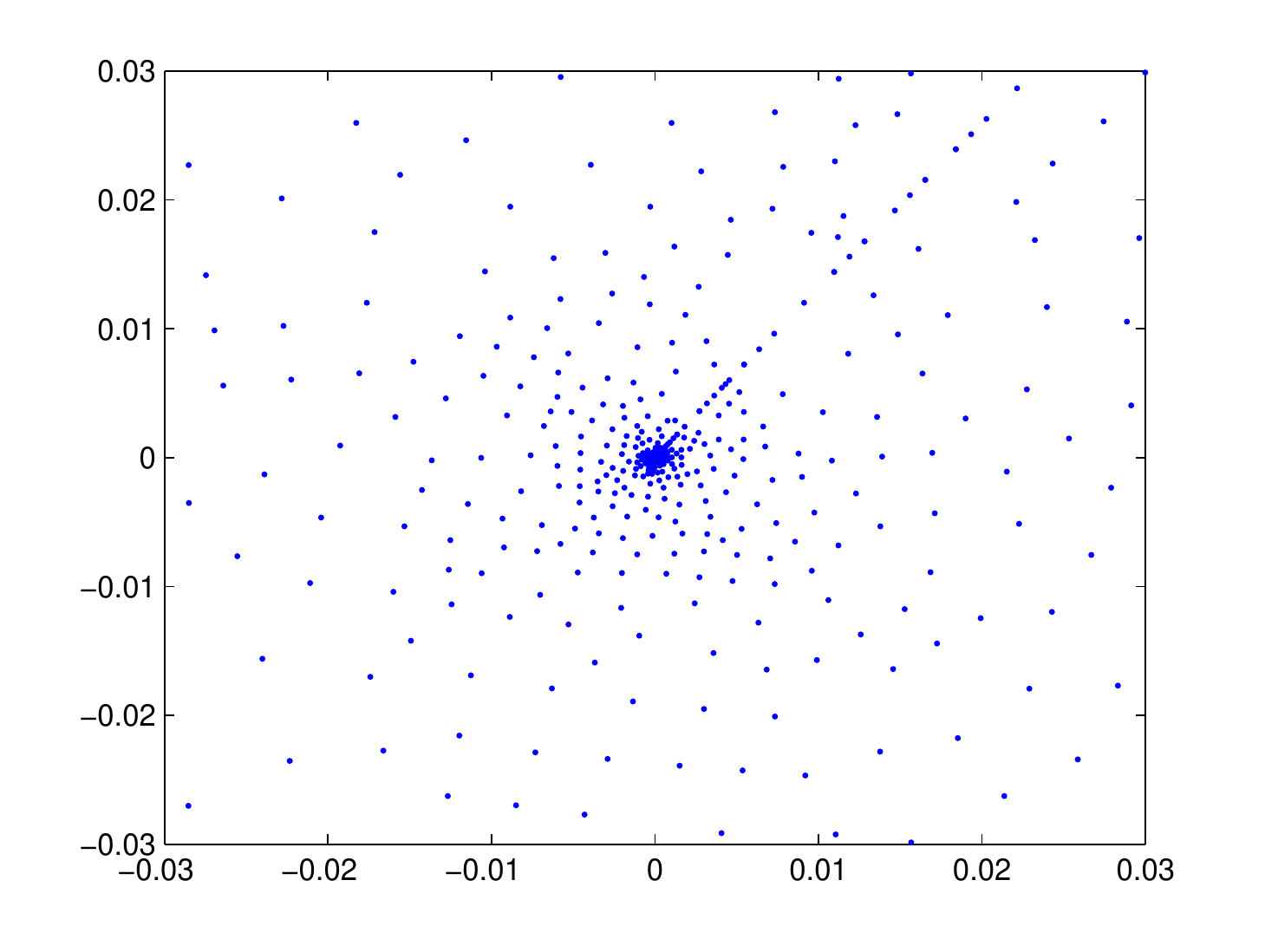}} \qquad 
 \subfigure[FEM centers: zoom ]
{\includegraphics[width=6cm,height=4cm]{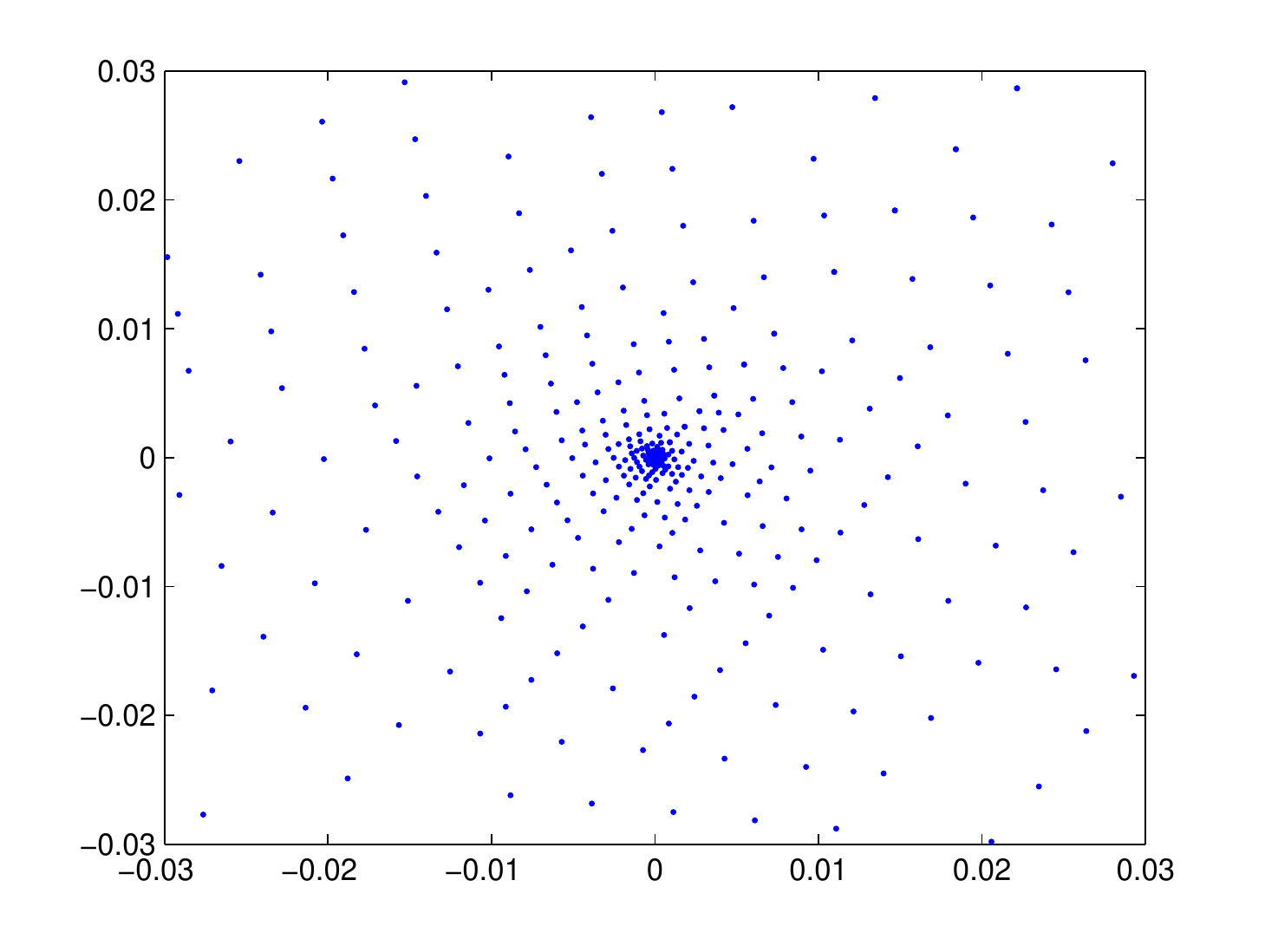}} 
        \end{center}
\caption{Test Problem~\ref{curved_slit}: Errors and centers as in Figure~\ref{figF4}. 
The plots in (cd) are based on the RBF-FD 
solution on 2204 interior centers shown in (e) and the FEM 
solution on 2236 interior vertices shown in (f).
} 
\label{figF426}
\end{figure}

The domain and the exact solution $u$ are illustrated in Figure~\ref{CSexact}. The exact
solution behaves as $r^{1/2}$ at the origin, that is the strength of the singularity is the
same as for the slit domain of Test Problem~\ref{f42} with $\omega=2\pi$. The numerical results are
presented in Figure~\ref{figF426} and are similar to the results for Test Problem~\ref{f42} 
with $\omega=2\pi$, which shows that the curvature of the slit does not present any significant additional
challenge.

\begin{testproblem}\label{f441}\rm 
 \cite[Section~2.8: Oscillatory]{Mitchell2013}
Dirichlet problem (\ref{poi})  for the Helm\-holtz equation
$-\Delta u-\frac{1}{(\alpha+r)^4}u=f$, $r = \sqrt{x^2+y^2}$,
in the domain $\Omega=(0,1)^2$, where the right hand side and the boundary conditions are chosen
such that the exact solution $u$ is $\sin(\frac{1}{\alpha +r})$, with $\alpha=\frac{1}{10 \pi}$
or $\frac{1}{50 \pi}$.
\end{testproblem}

\begin{figure}[htbp!]
 \begin{center}
\subfigure[$\alpha=\frac{1}{10 \pi}$]
{\includegraphics[width=6cm,height=4cm]{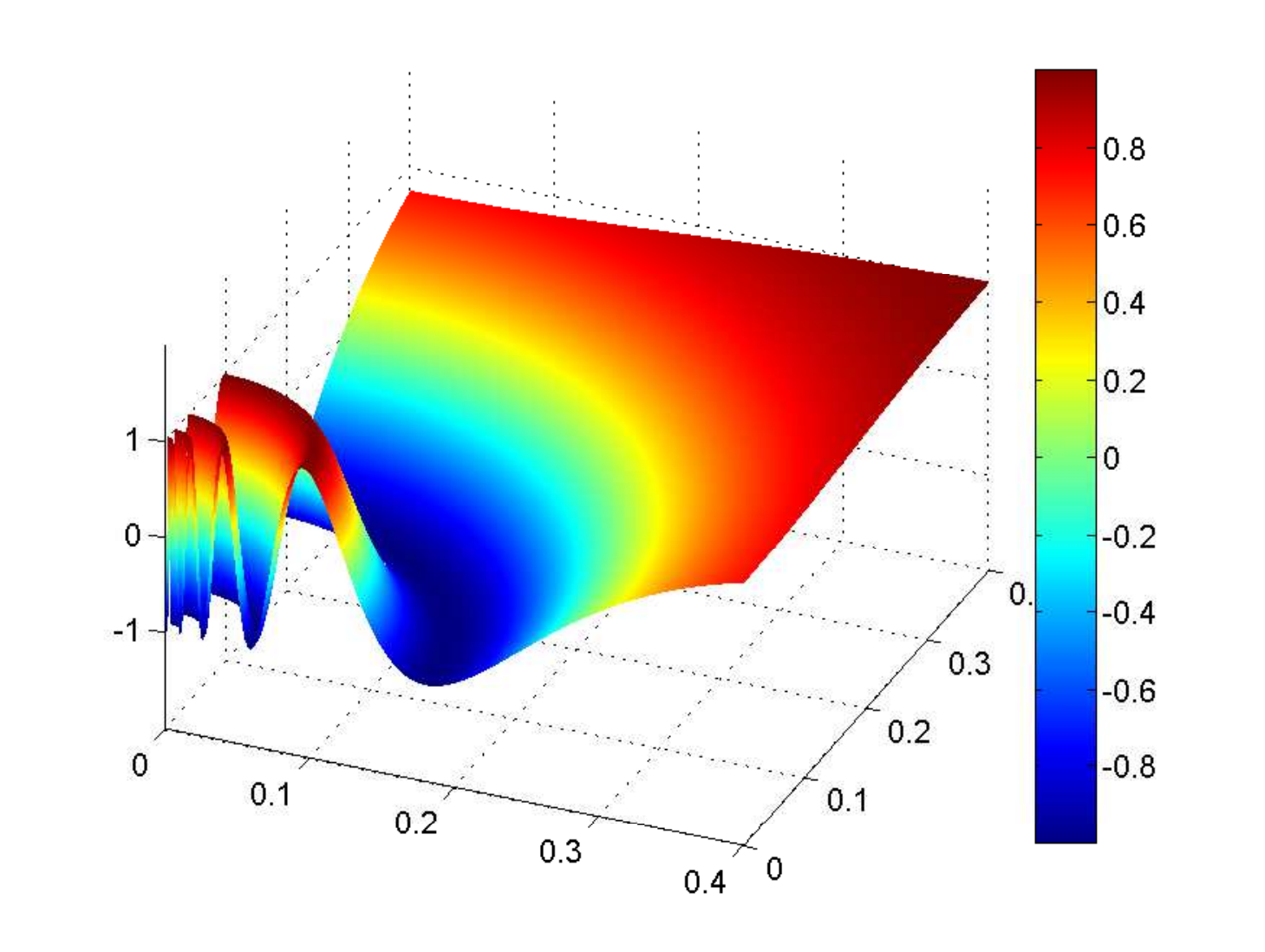}} \qquad
 \subfigure[$\alpha=\frac{1}{50 \pi}$]
{\includegraphics[width=6cm,height=4cm]{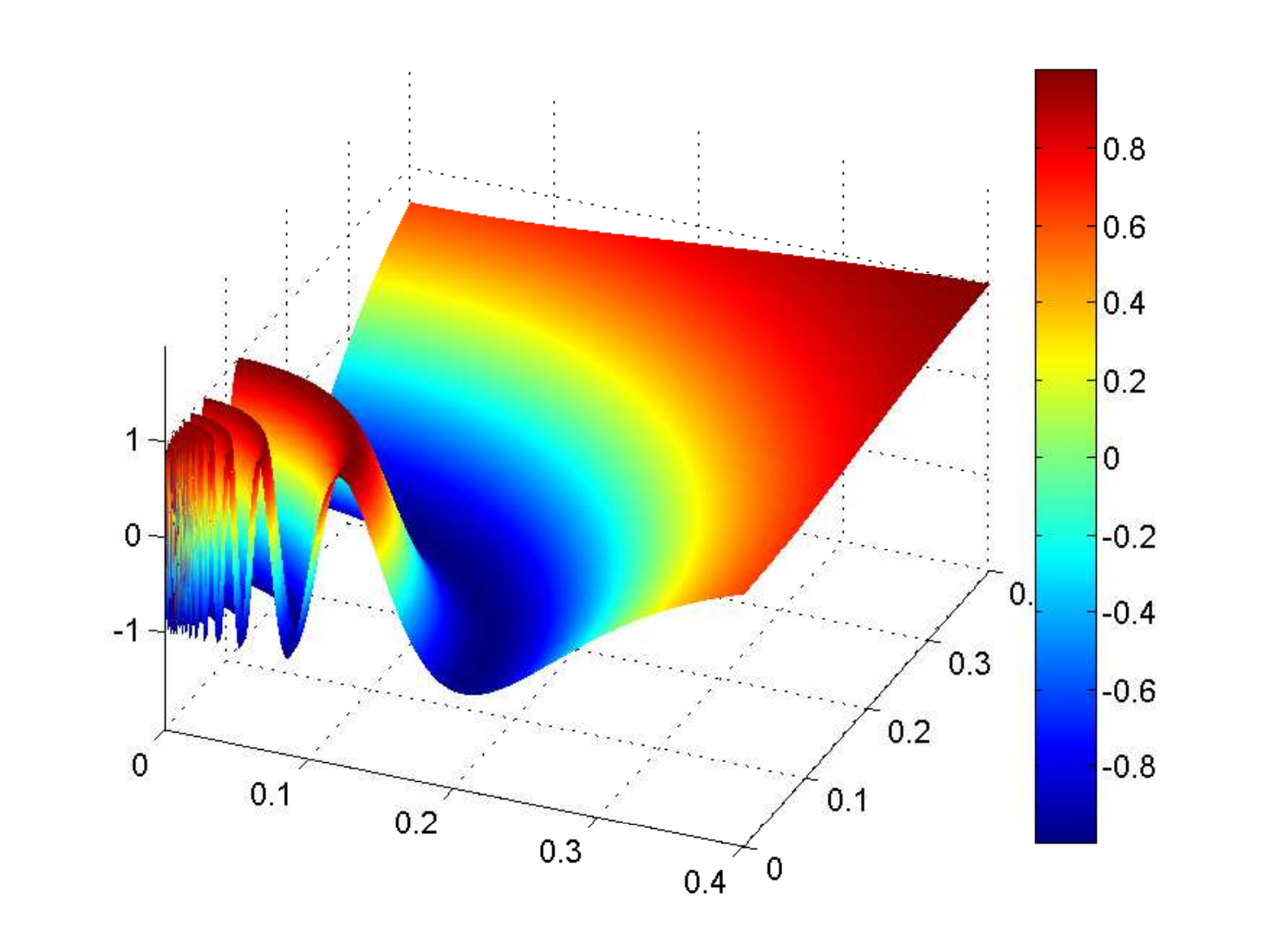}} 
        \end{center}
\caption{Test Problem~\ref{f441}: Exact solutions for $\alpha=\frac{1}{10 \pi}$ and 
$\frac{1}{50 \pi}$.
}
\label{exactF44}
\end{figure}

The solutions for both $\alpha=\frac{1}{10 \pi}$ and $\frac{1}{50 \pi}$ are highly oscillatory near
the origin, with increasing frequency closer to the origin, see Figure~\ref{exactF44}.

\begin{figure}[htbp!]
 \begin{center}
 \subfigure[Errors on centers]%
{\includegraphics[width=6cm,height=4cm]{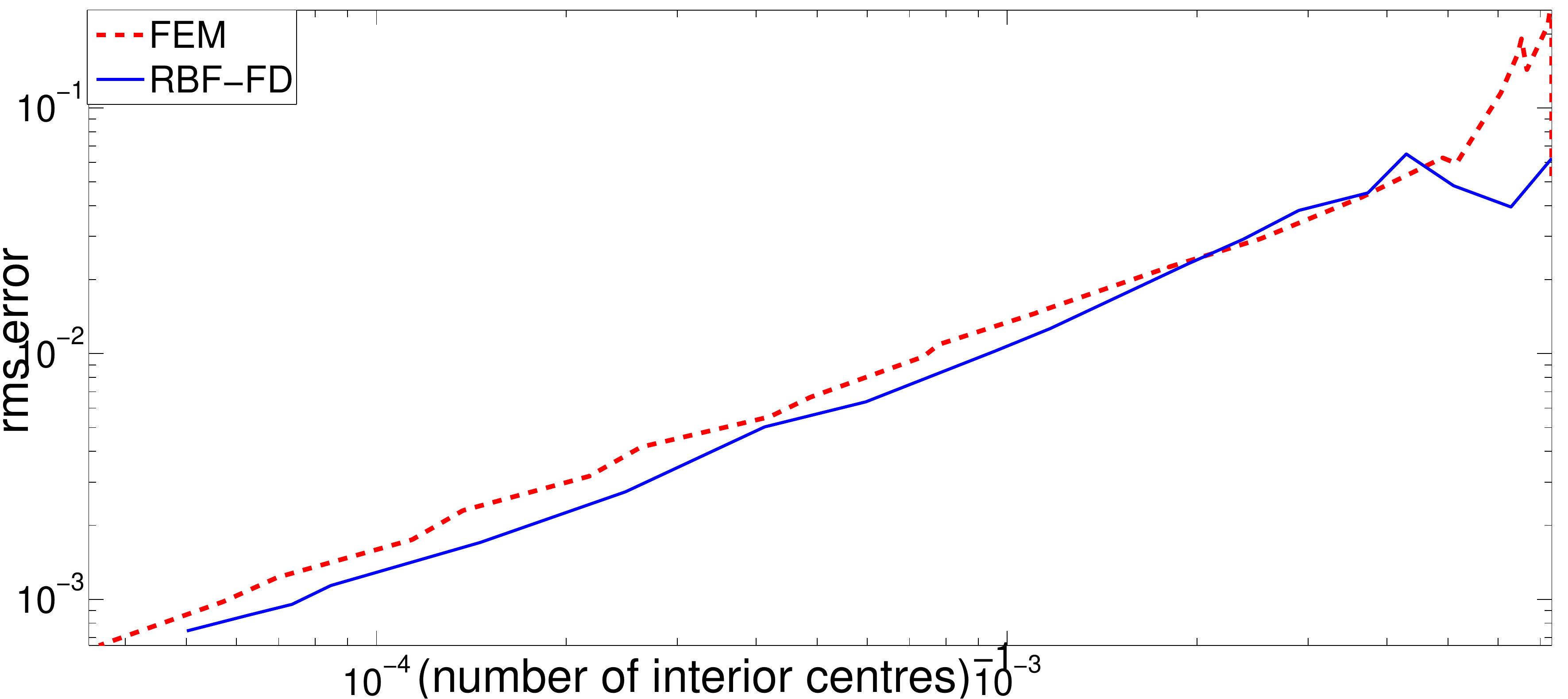}} \quad
 \subfigure[Errors on grid]%
{\includegraphics[width=6cm,height=4cm]{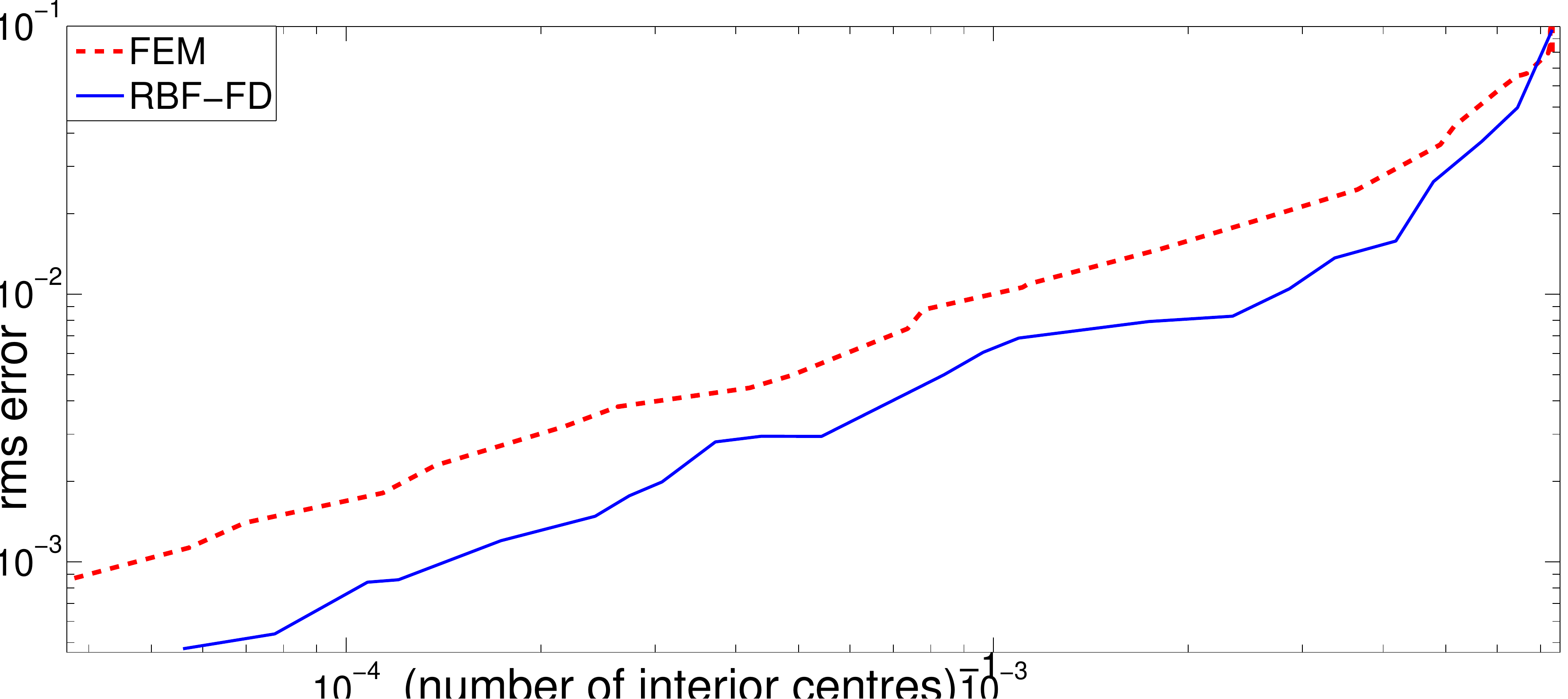}} 
 \subfigure[RBF-FD error]%
{\includegraphics[width=6cm,height=4cm]{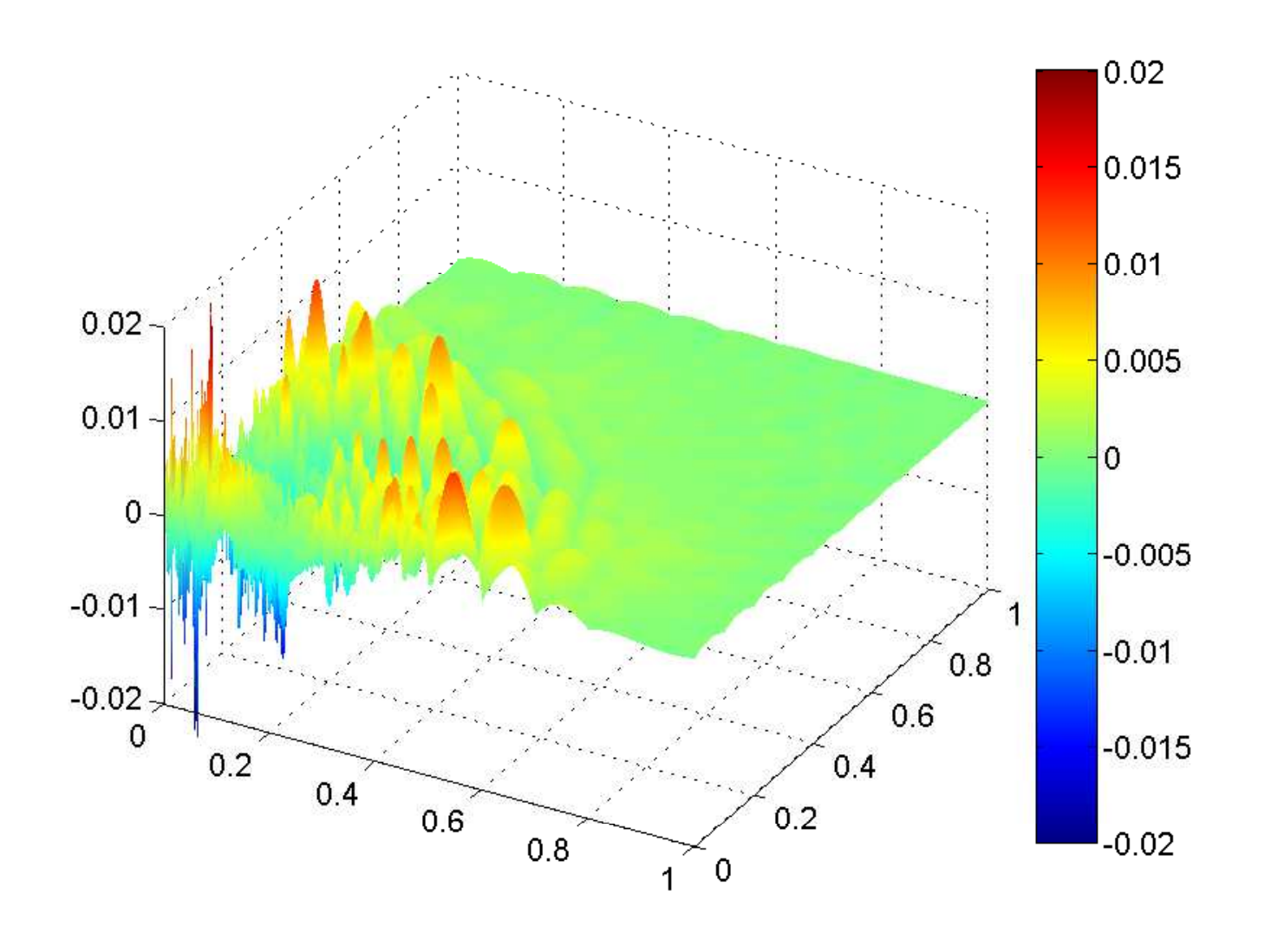}} 
 \subfigure[FEM error]%
{\includegraphics[width=6cm,height=4cm]{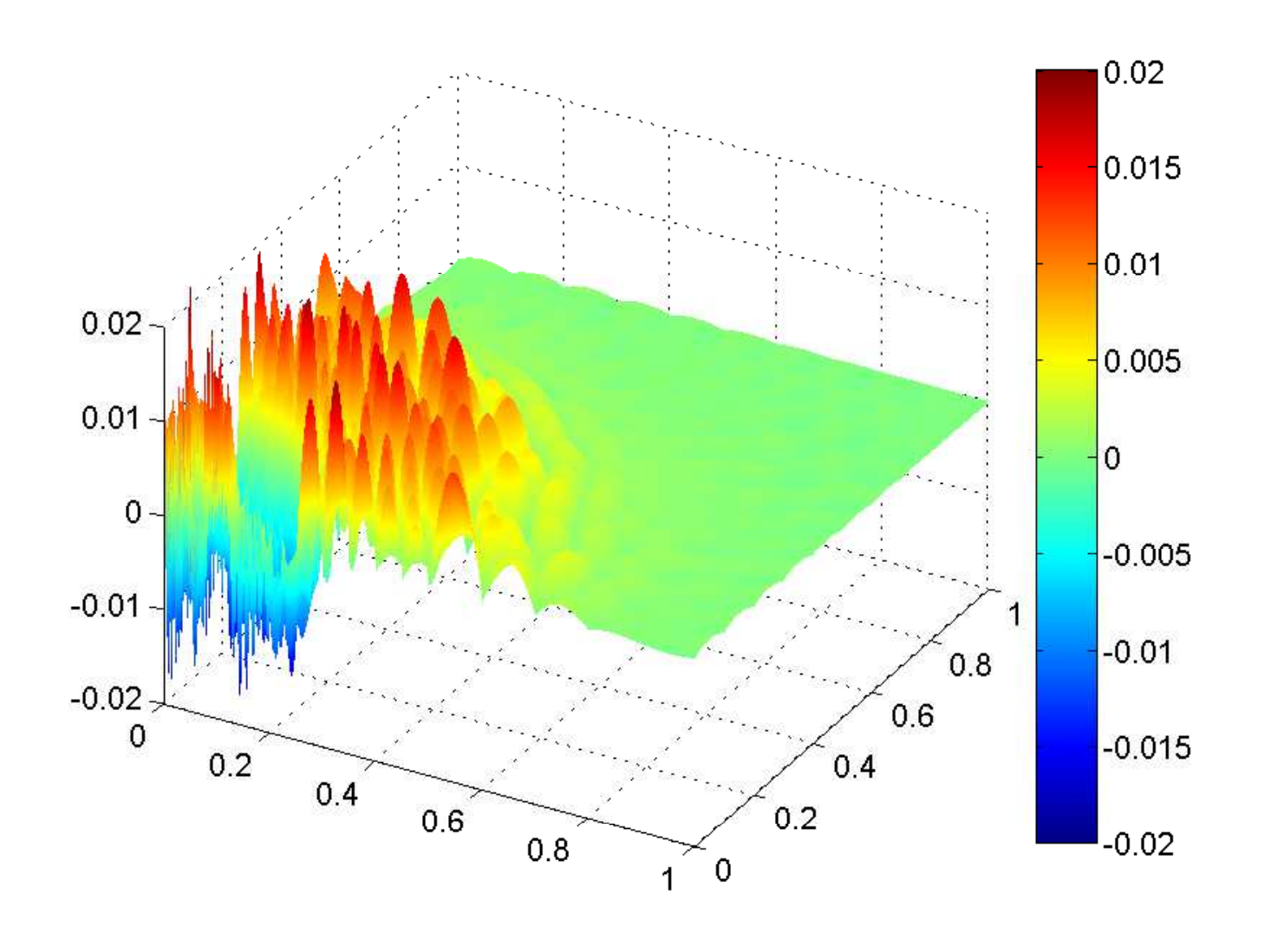}} 
 \subfigure[RBF-FD error: zoom]
{\includegraphics[width=6cm,height=4cm]{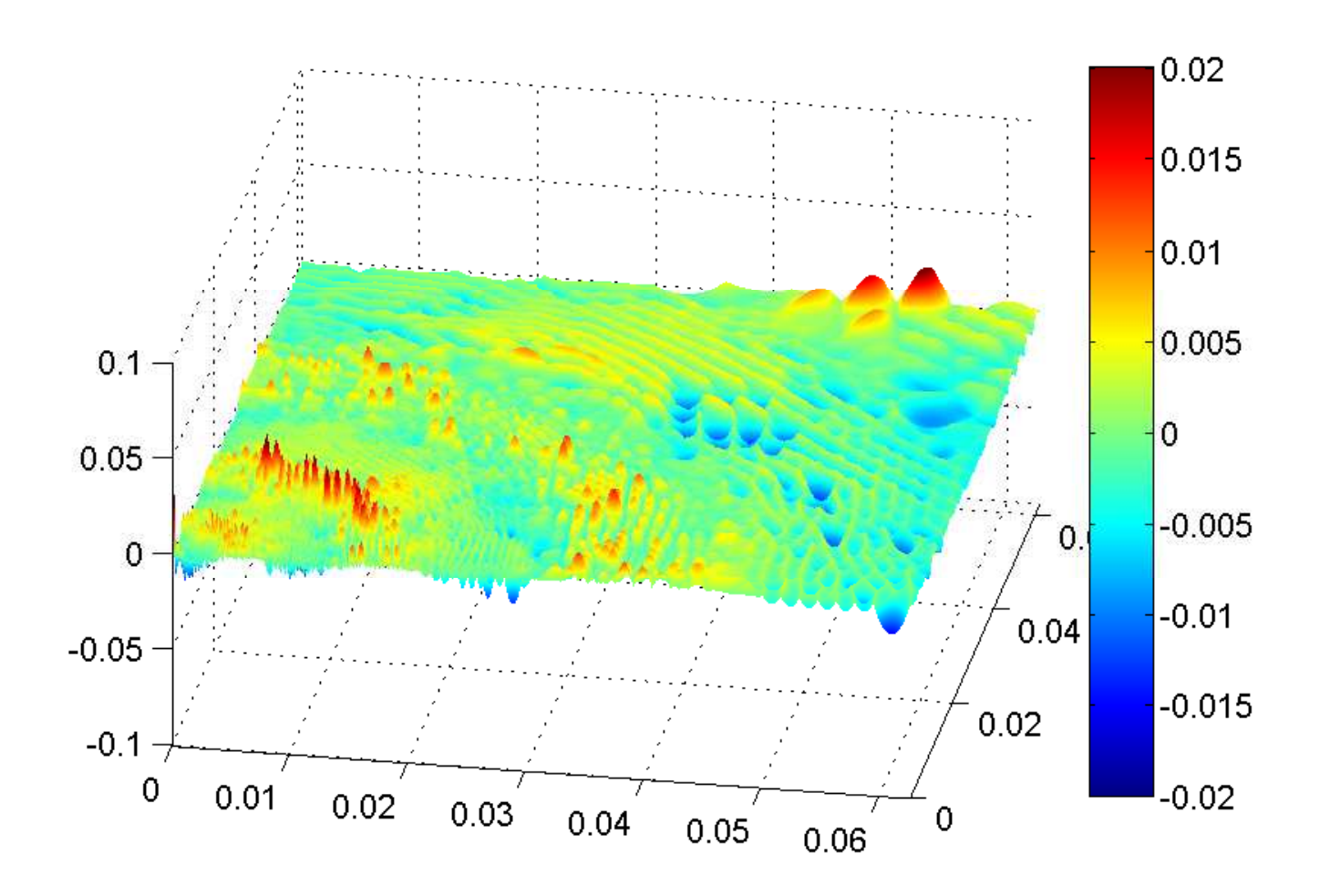}} 
 \subfigure[FEM error: zoom]
{\includegraphics[width=6cm,height=4cm]{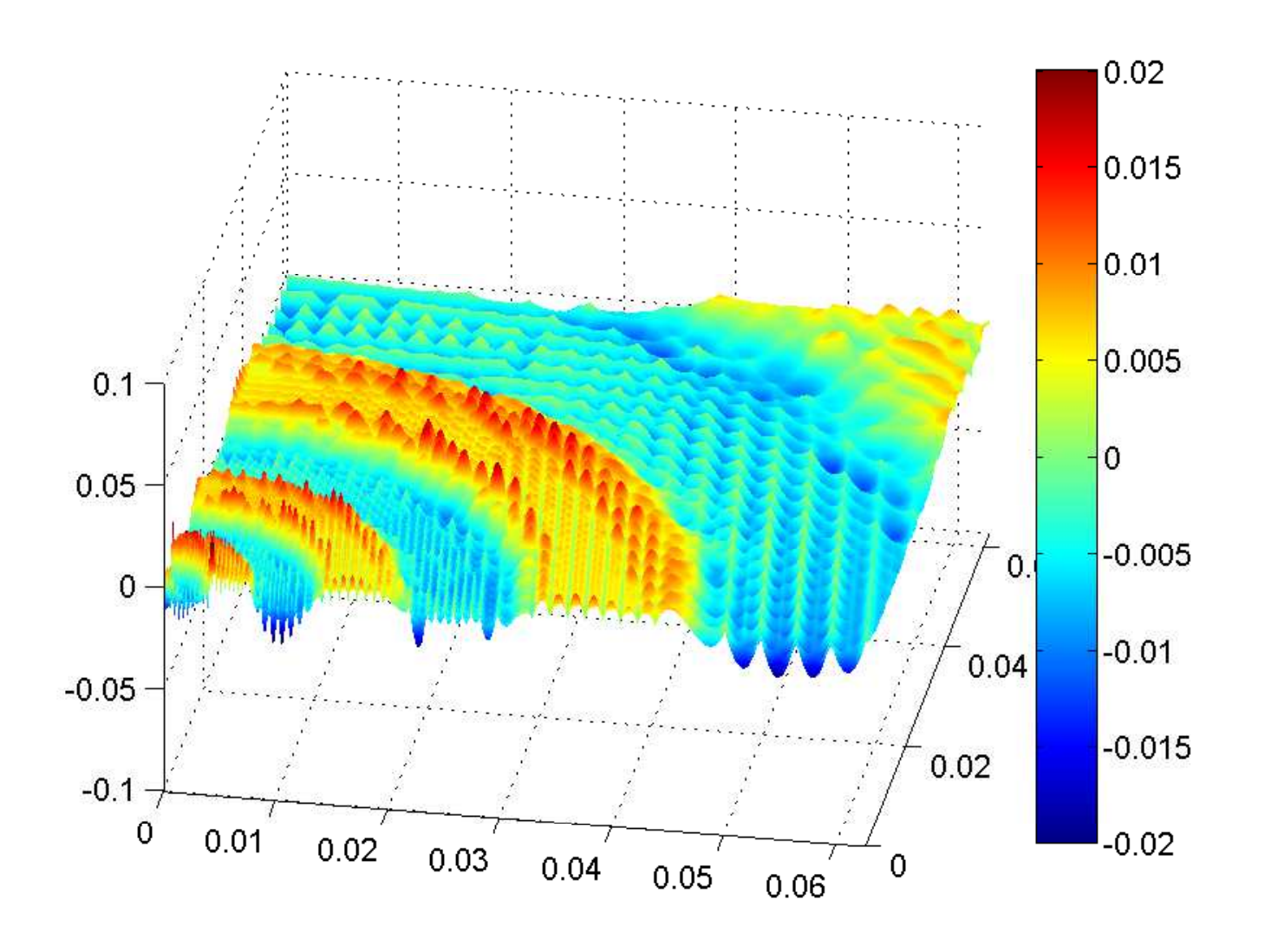}} 
 \subfigure[RBF-FD centers: zoom]%
{\includegraphics[width=6cm,height=4cm]{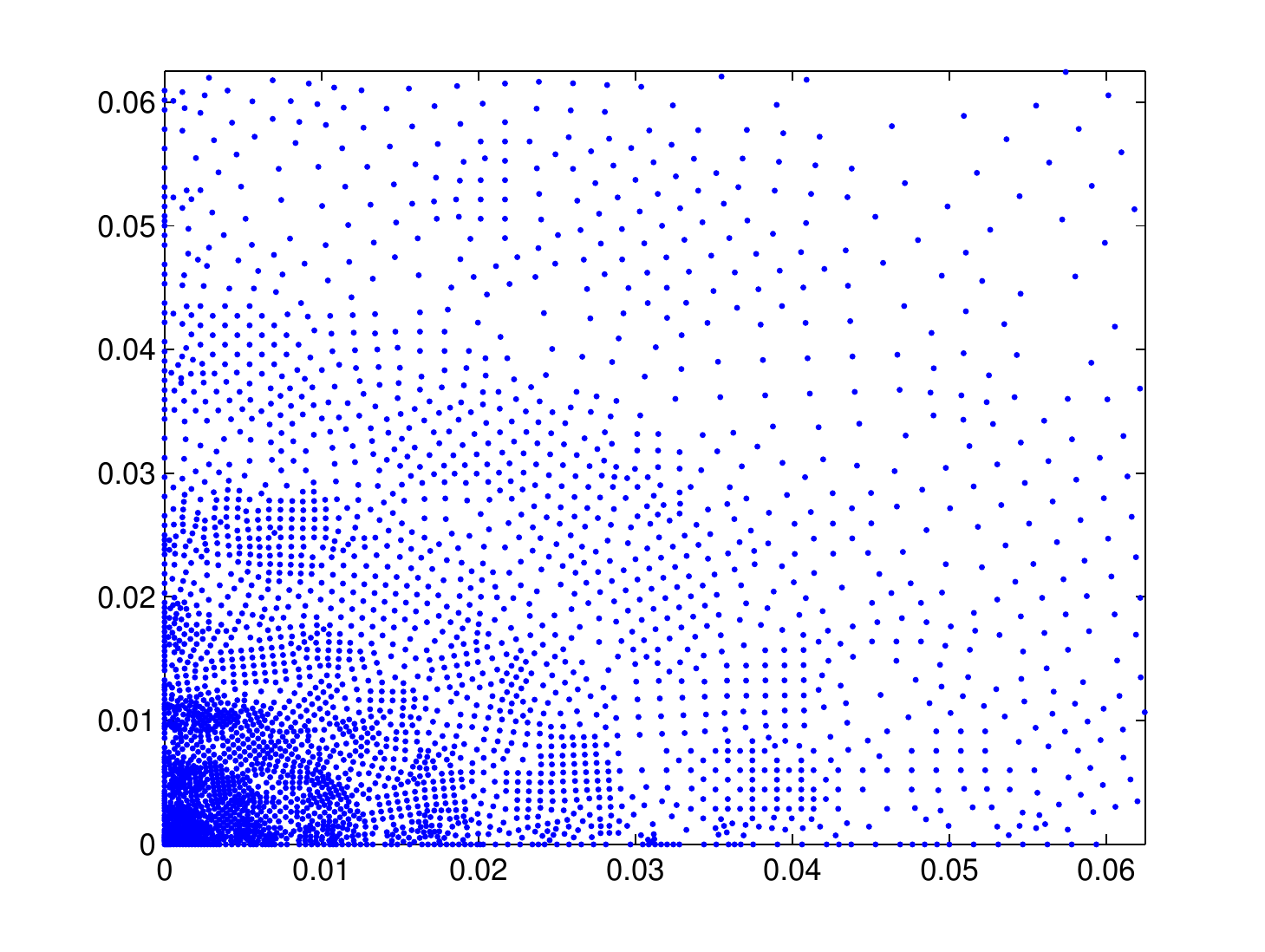}} 
 \subfigure[FEM centers: zoom]%
{\includegraphics[width=6cm,height=4cm]{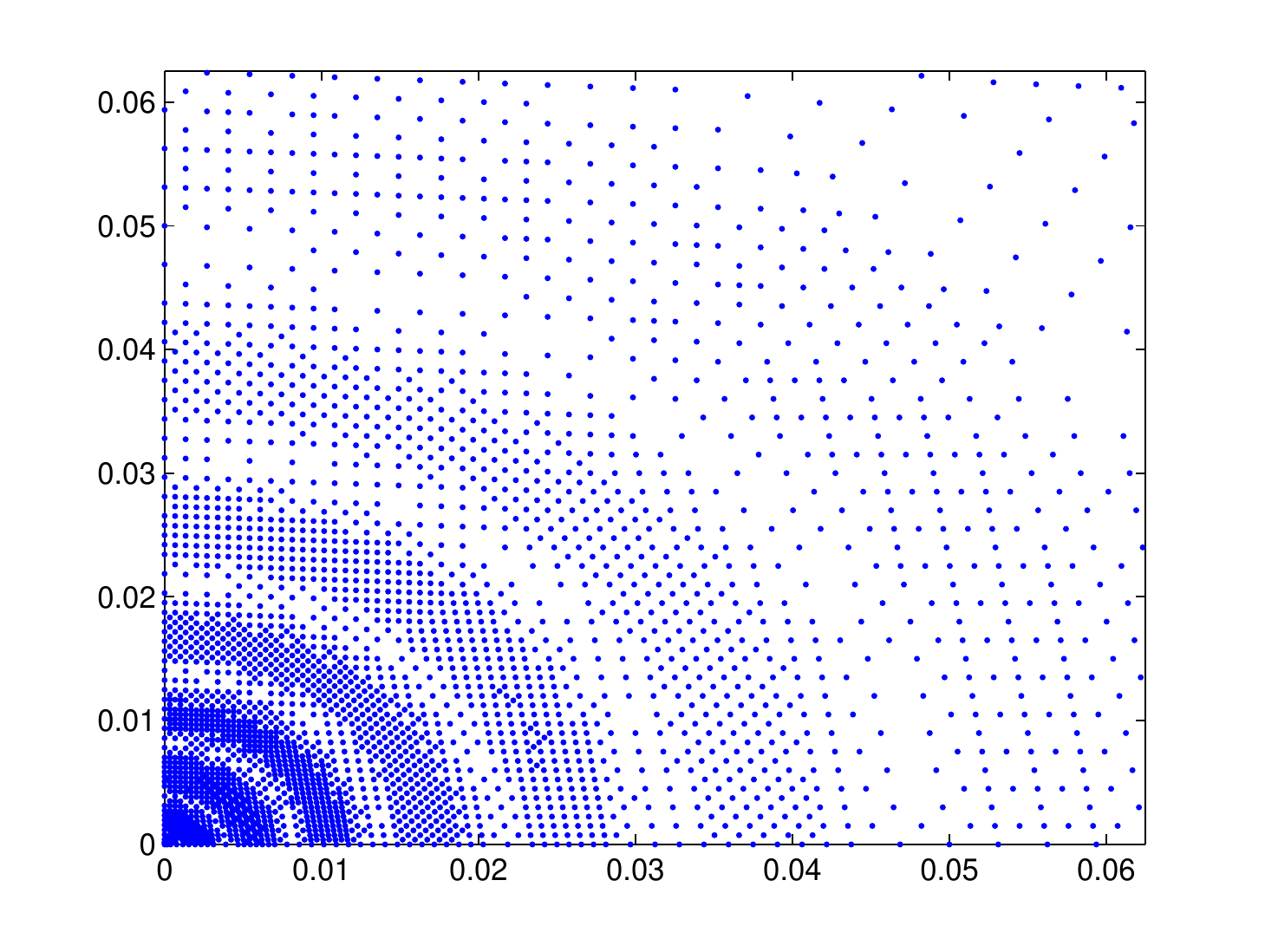}} 
      \end{center}
\caption{Test Problem~\ref{f441} with $\alpha=\frac{1}{10 \pi}$: Errors and centers/vertices.
The plots in  (ceg) are based on the RBF-FD 
solution with 4029 interior centers, whereas those in (dfh) on the FEM 
solution with 3806 interior vertices. 
}
\label{figF441}
\end{figure}

\begin{figure}[htbp!]
 \begin{center}
 \subfigure[Zoom of the exact solution]
{\includegraphics[width=6cm,height=4cm]{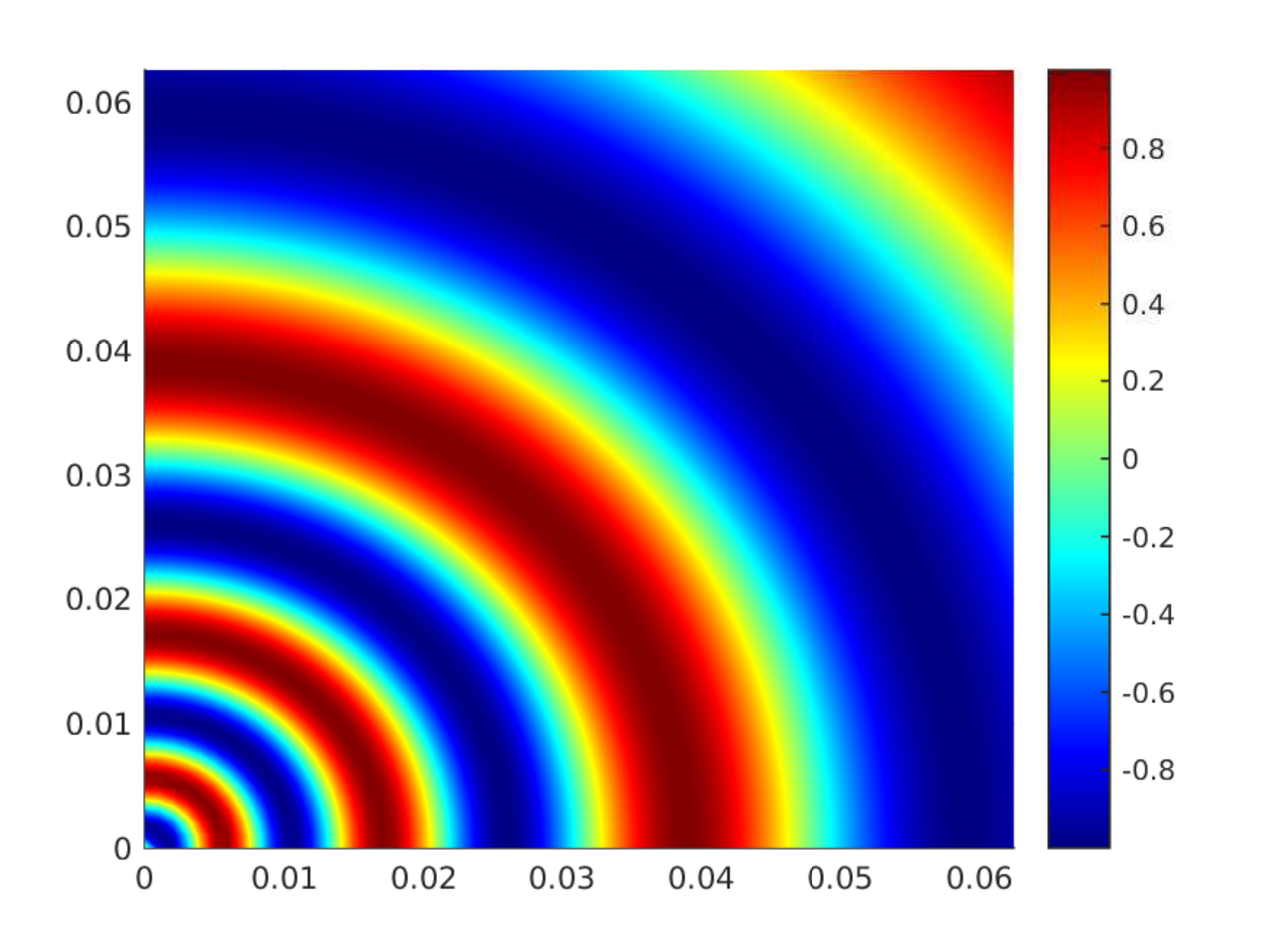}} 
  \subfigure[Centers with error indicator $\eps_0$]%
{\includegraphics[width=6cm,height=4cm]{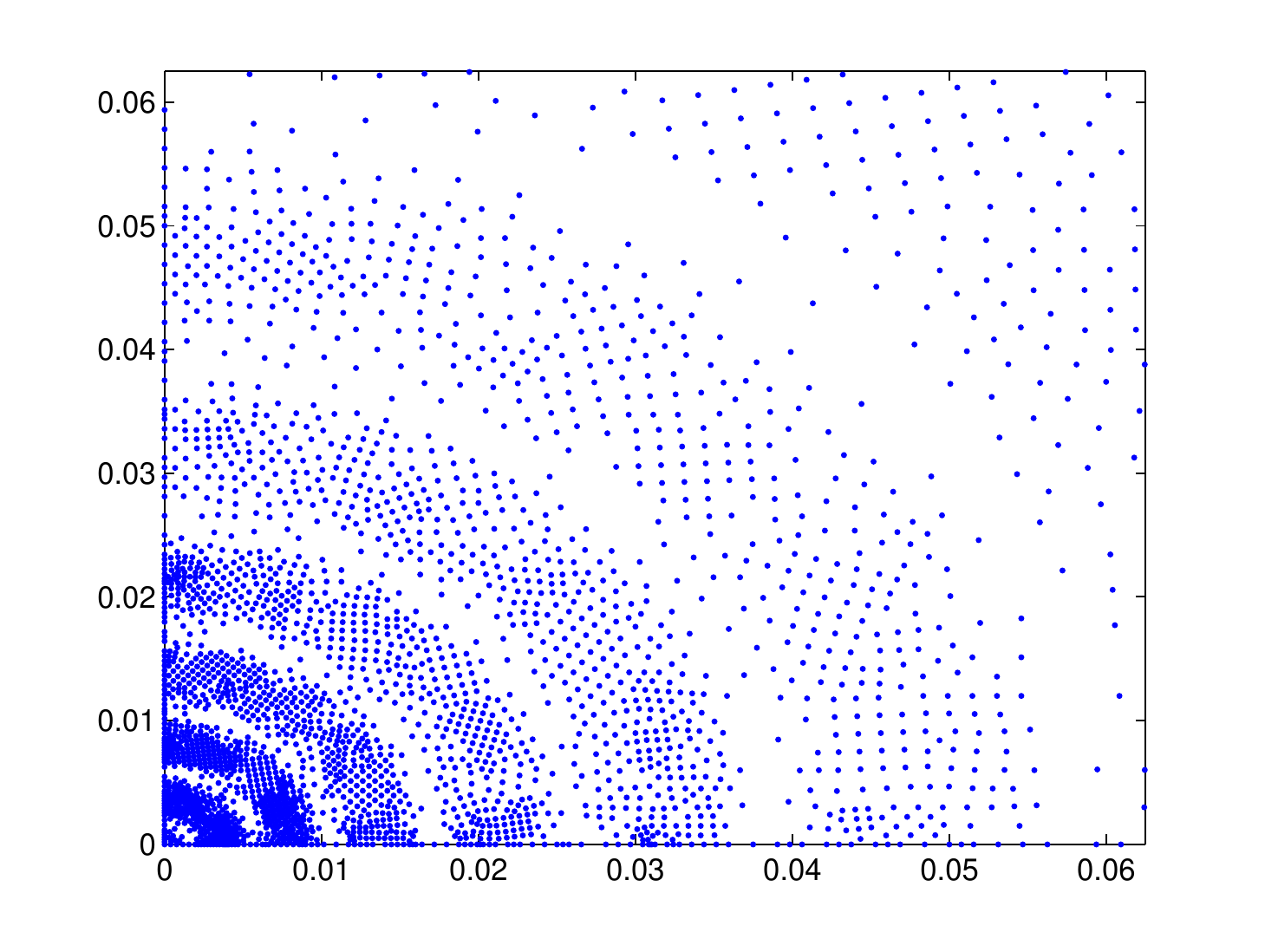}} 
 \subfigure[Error graph]
{\includegraphics[width=6cm,height=4cm]{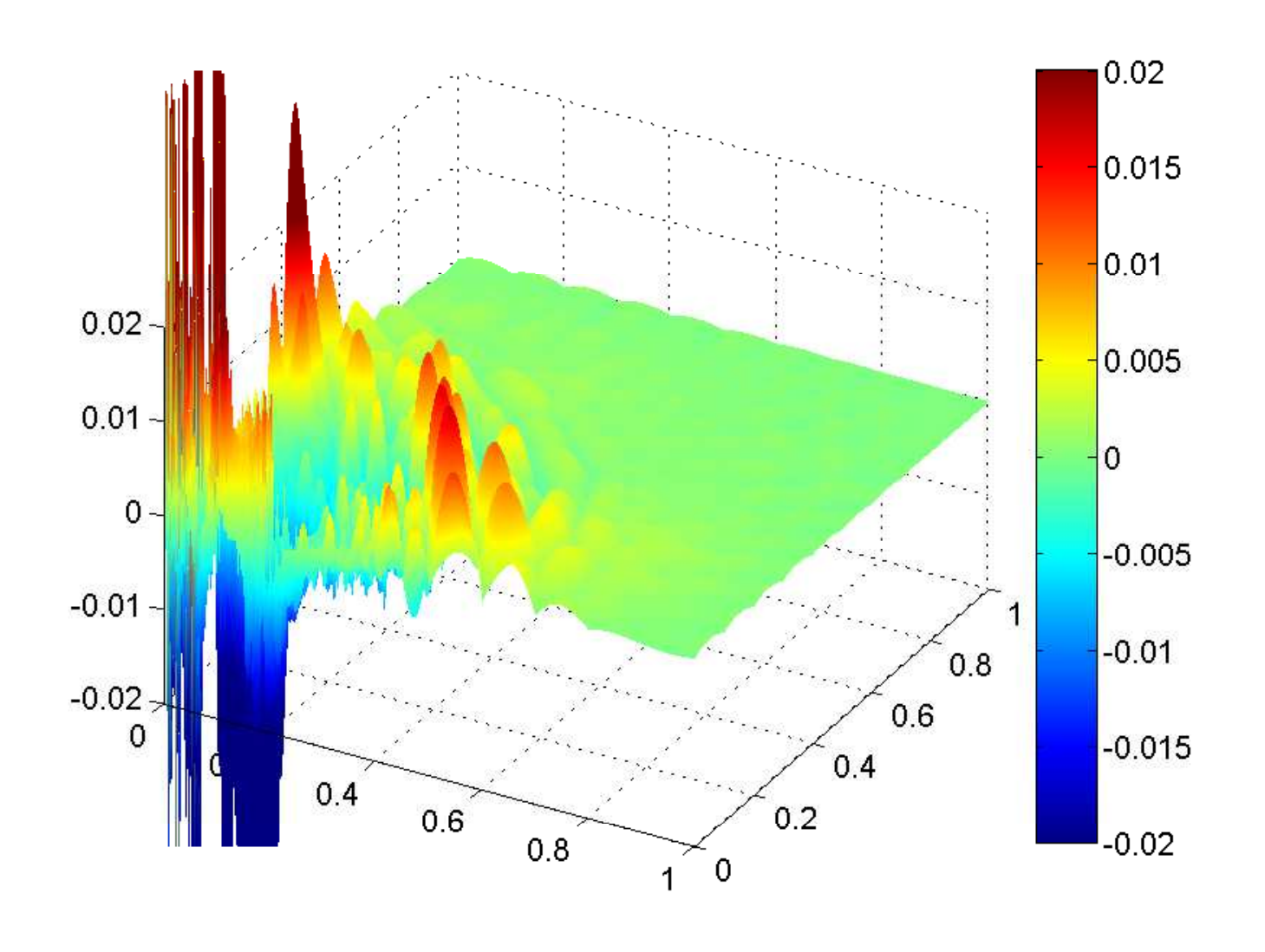}} 
\subfigure[Error graph: zoom]
{\includegraphics[width=6cm,height=4cm]{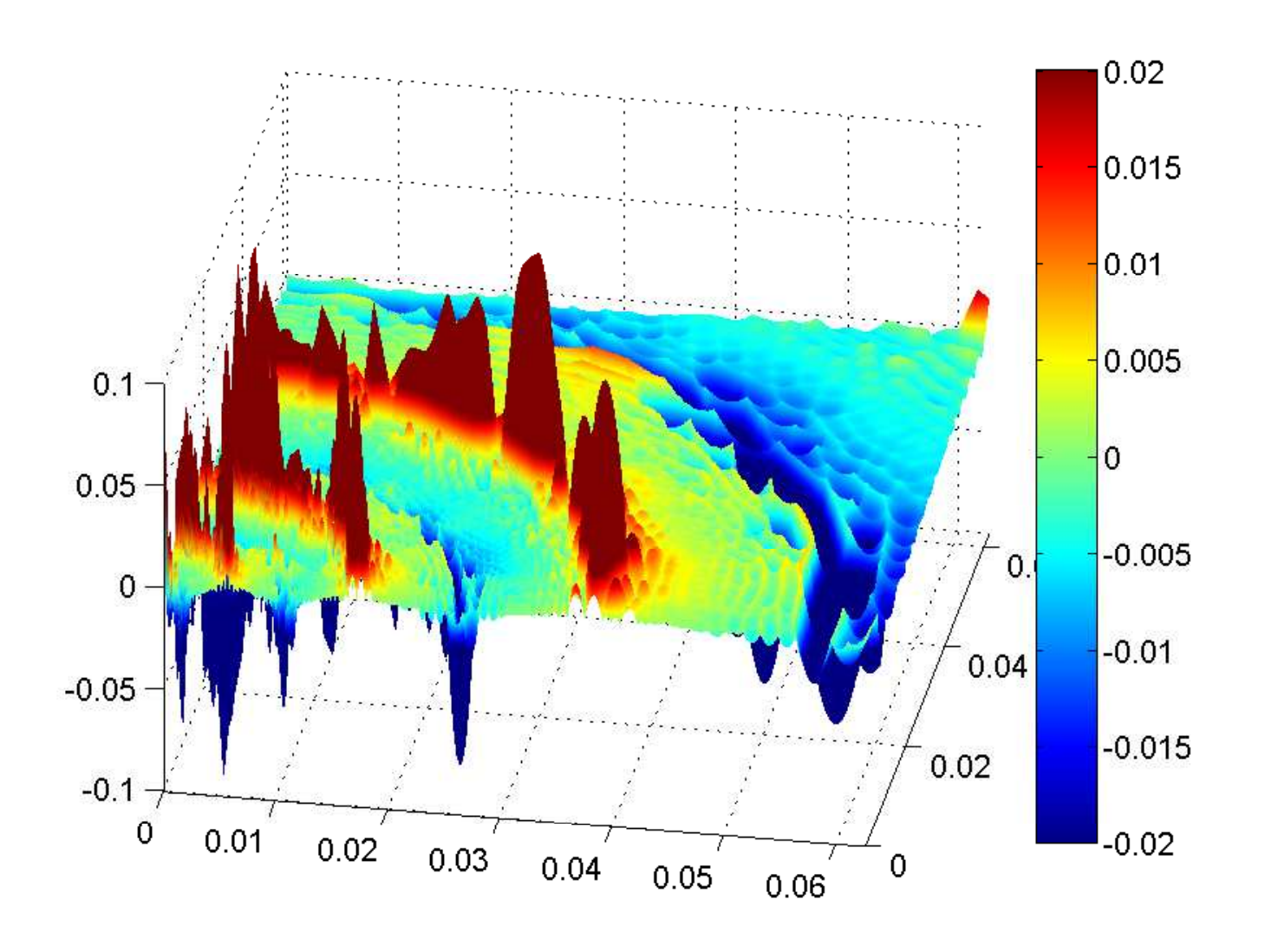}} 
       \end{center}
\caption{Test Problem~\ref{f441} with $\alpha=\frac{1}{10 \pi}$: (a) Exact solution in the 
subregion used in Figure~\ref{figF441}(e-h) and in (bd) of this figure. 
(b-d) Centers and errors for the solution 
with 3679 interior centers obtained by RBF-FD method 
with error indicator $\eps_0(\zeta,\xi)=|\hu_\zeta-\hu_\xi|$.
}
\label{IndicF441}
\end{figure}

Numerical results for Test Problem~\ref{f441} with  $\alpha=\frac{1}{10 \pi}$  are presented 
in Figure~\ref{figF441}. We see that RBF-FD is generally more accurate than FEM in this case. Moreover,
Figures (e) and (f) suggest that the RBF-FD  solution is less susceptible to a bias in the areas of
high oscillation, manifested in a systematic underestimation of the amplitude of the oscillations.
Figure (g) shows that the centers produced by the RBF-FD method reproduce to some extent
the distinctive ring-like pattern seen in the FEM vertices of (h). This pattern is easy to explain by
comparing Figures~\ref{figF441}(gh) with Figure~\ref{IndicF441}(a) which shows a color-mapped image of the exact solution in the
same area. The centers are placed more densely in the highly curved regions at the tops and bottoms of the waves
and neglect the rather flat regions between them. This placement of the centers is known to be advantageous
for piecewise linear approximation and other approximation methods with leading error term related to the
size of the second order derivatives or the curvature of the exact solution \cite{BirmanSolomyak67}. 
This test has been the main
motivation for us to replace the error indication $\eps_0(\zeta,\xi)$ by $\eps_1(\zeta,\xi)$ as explained
in Section~\ref{Refi}. If we apply RBF-FD method of this paper with $\eps_0$ instead of $\eps_1$, then we
get the rings of centers with high density shifted to the flat regions of high gradient as illustrated in
Figure~\ref{IndicF441}(b). As a result, as naturally expected, the overall error is much larger and the accuracy is particularly
poor near the tops and bottoms of waves, see Figure~\ref{IndicF441}(cd).

\begin{figure}[htbp!]
 \begin{center}
 \subfigure[Errors on centers] %
{\includegraphics[width=6cm,height=4cm]{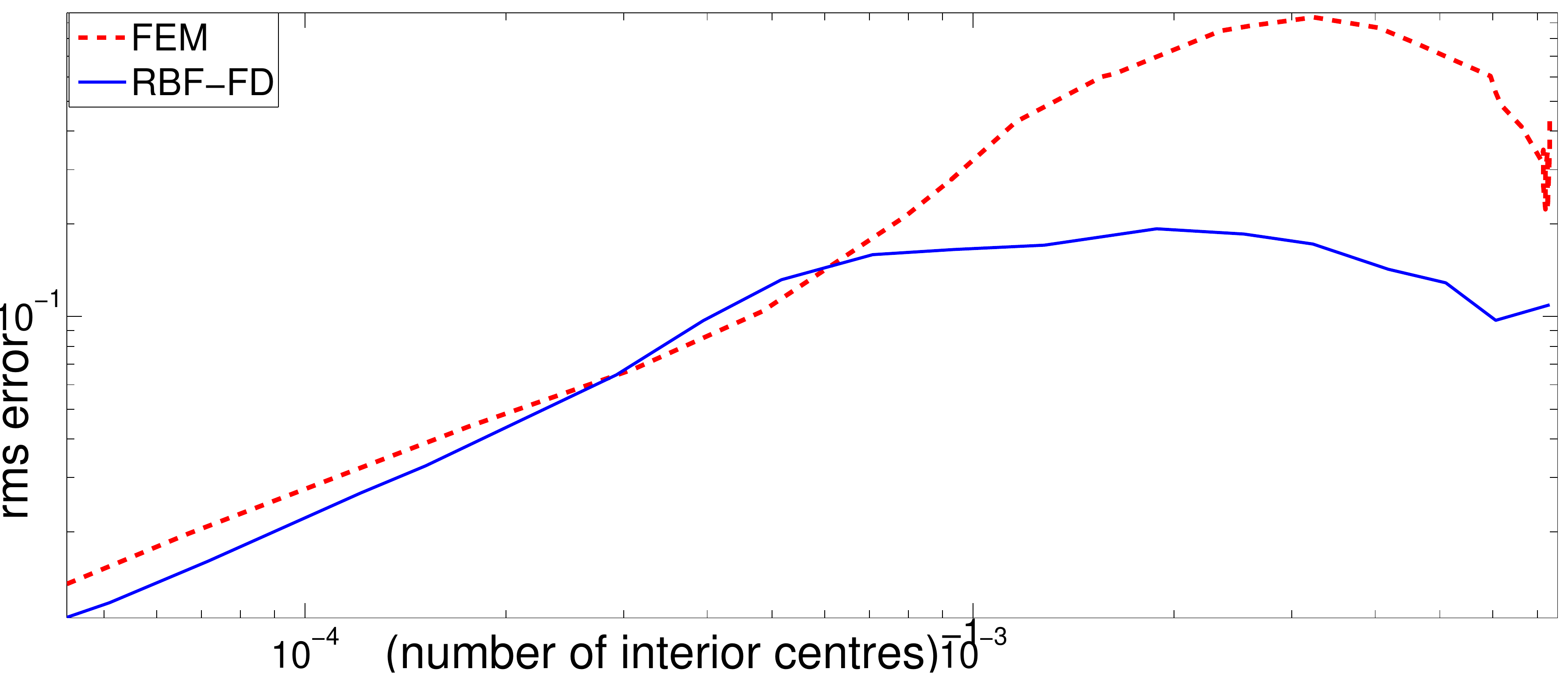}} \quad
 \subfigure[Errors on grid] %
{\includegraphics[width=6cm,height=4cm]{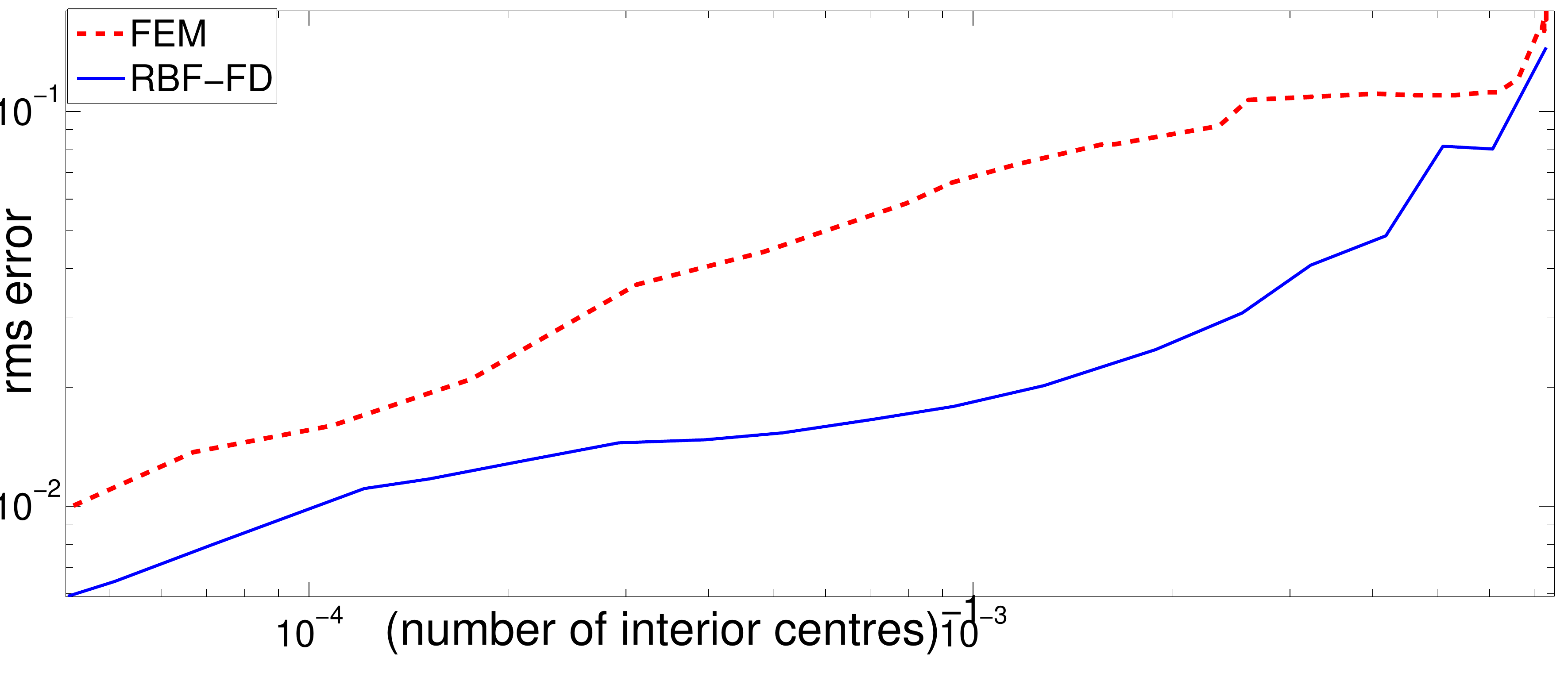}} 
 \subfigure[RBF-FD error]
{\includegraphics[width=6cm,height=4cm]{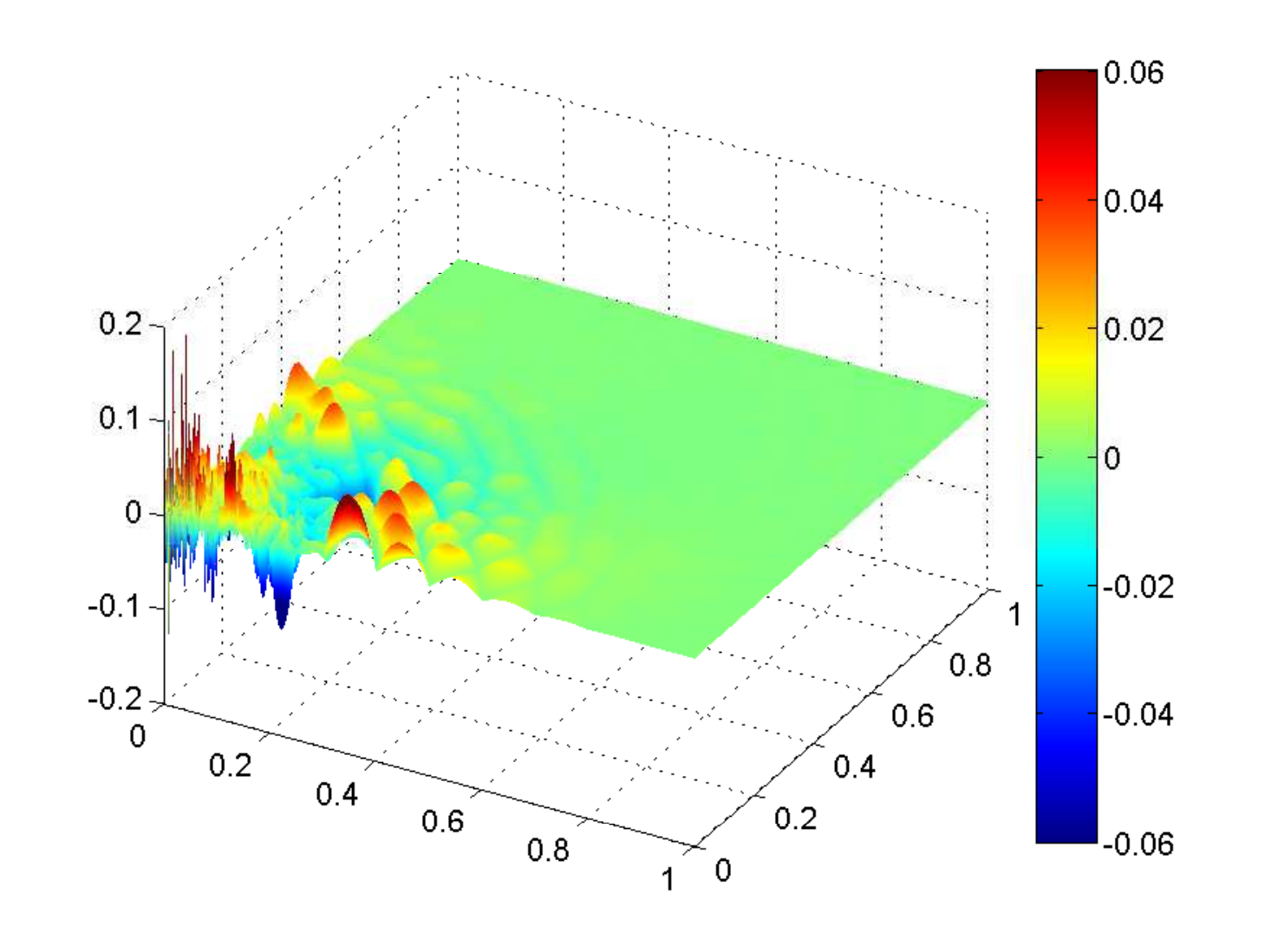}} 
 \subfigure[FEM error ]
{\includegraphics[width=6cm,height=4cm]{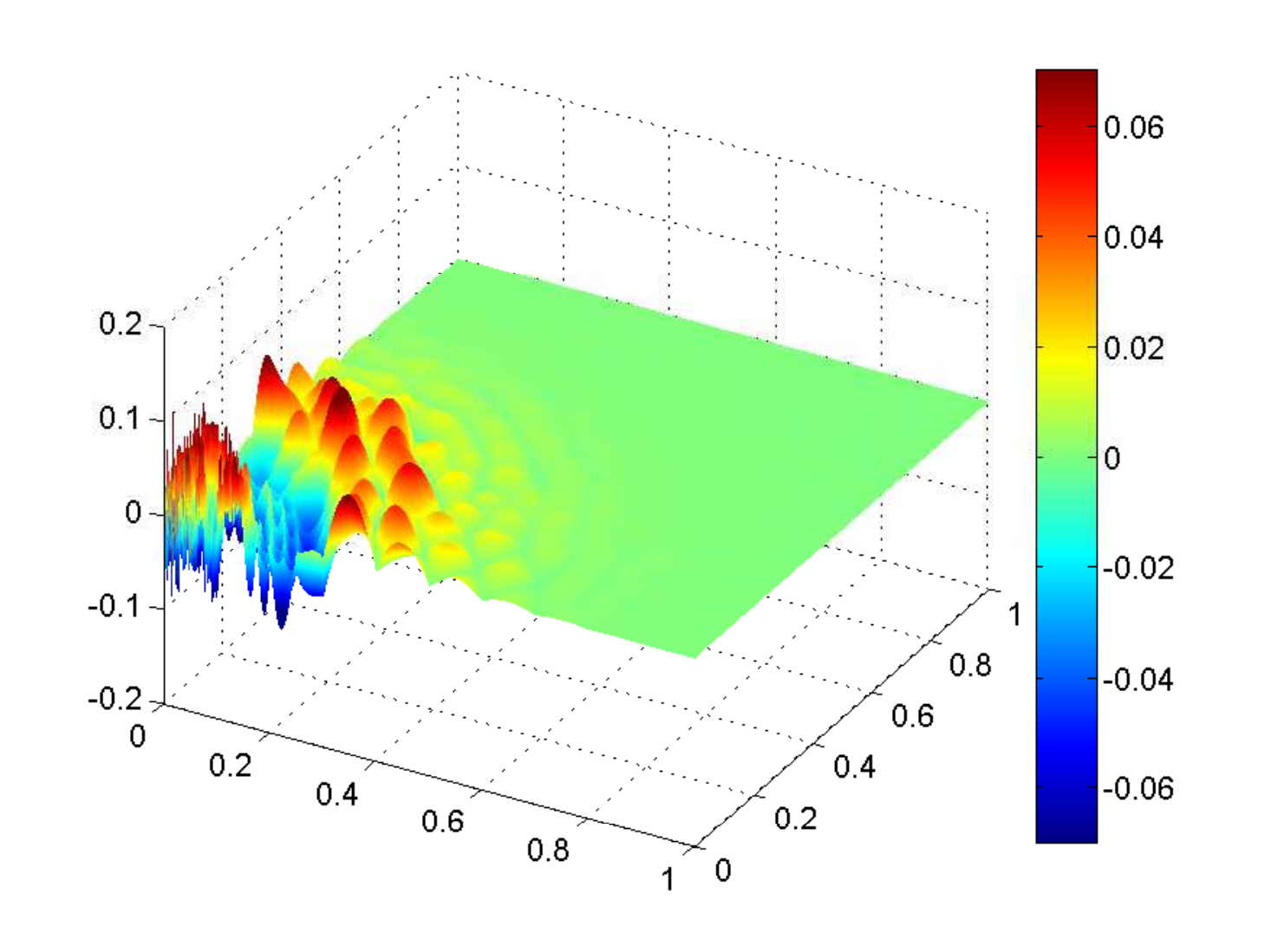}} 
\subfigure[RBF-FD error: zoom]
{\includegraphics[width=6cm,height=4cm]{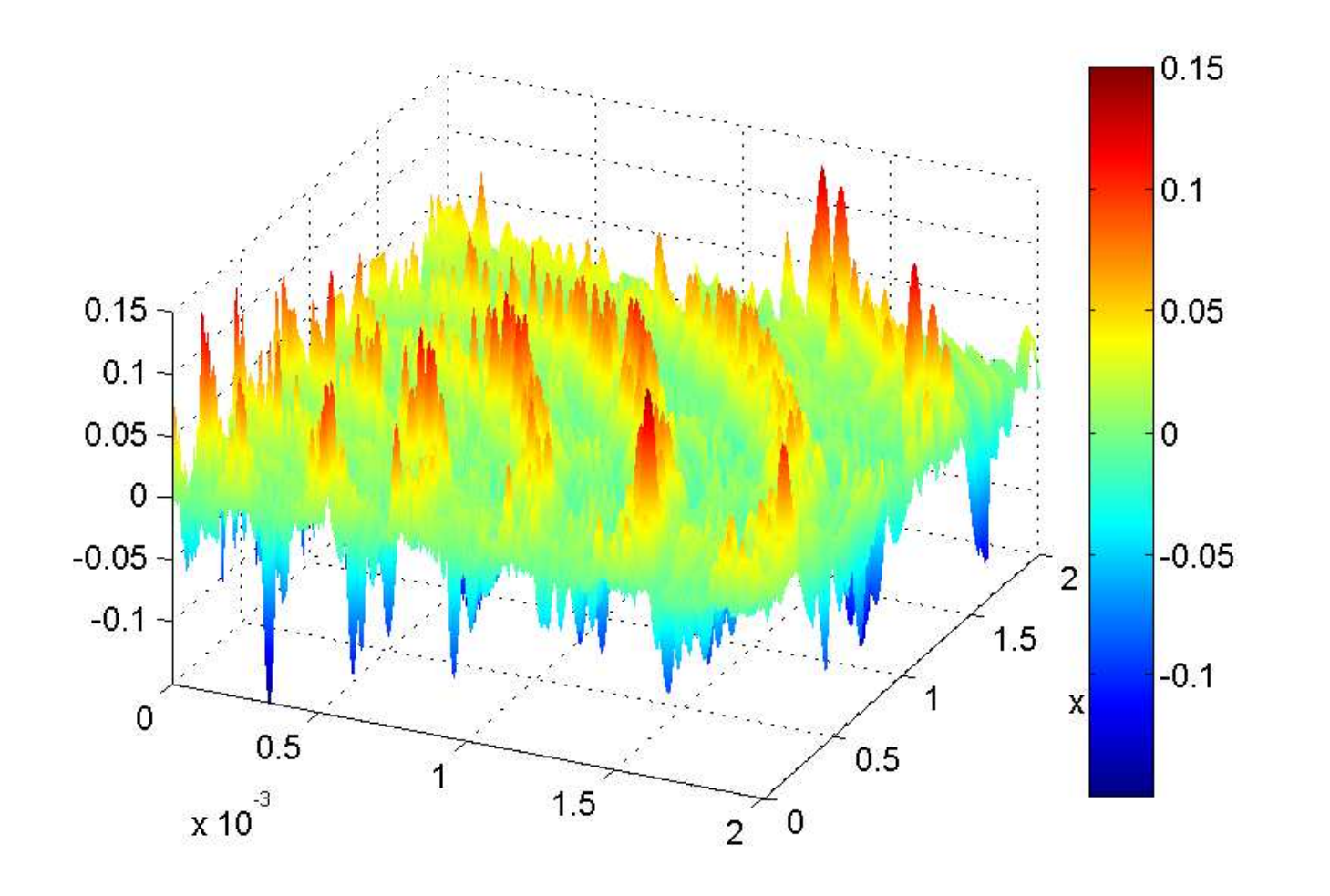}} 
 \subfigure[FEM error: zoom]
{\includegraphics[width=6cm,height=4cm]{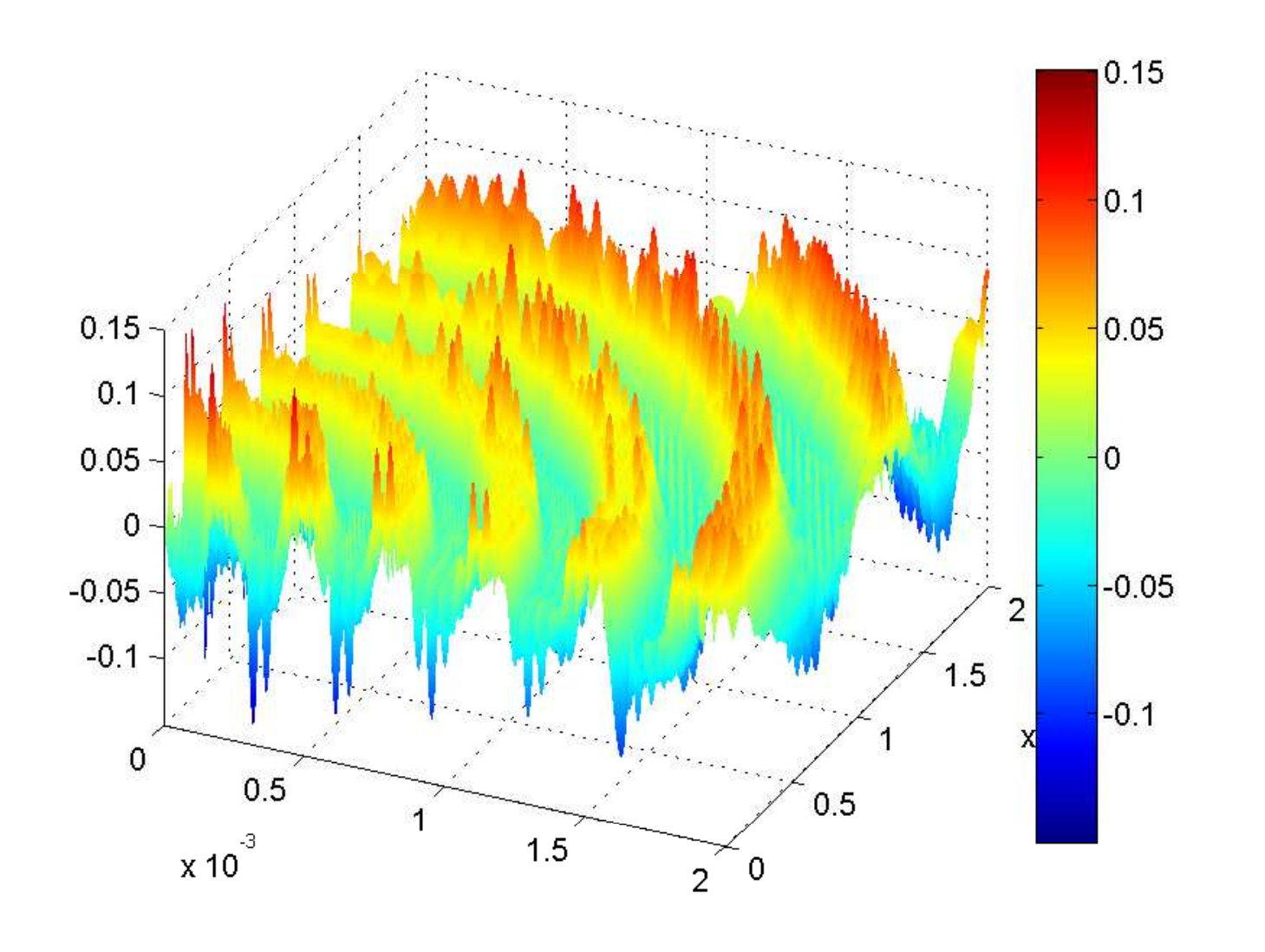}} 
 \subfigure[RBF-FD centers: zoom] %
{\includegraphics[width=6cm,height=4cm]{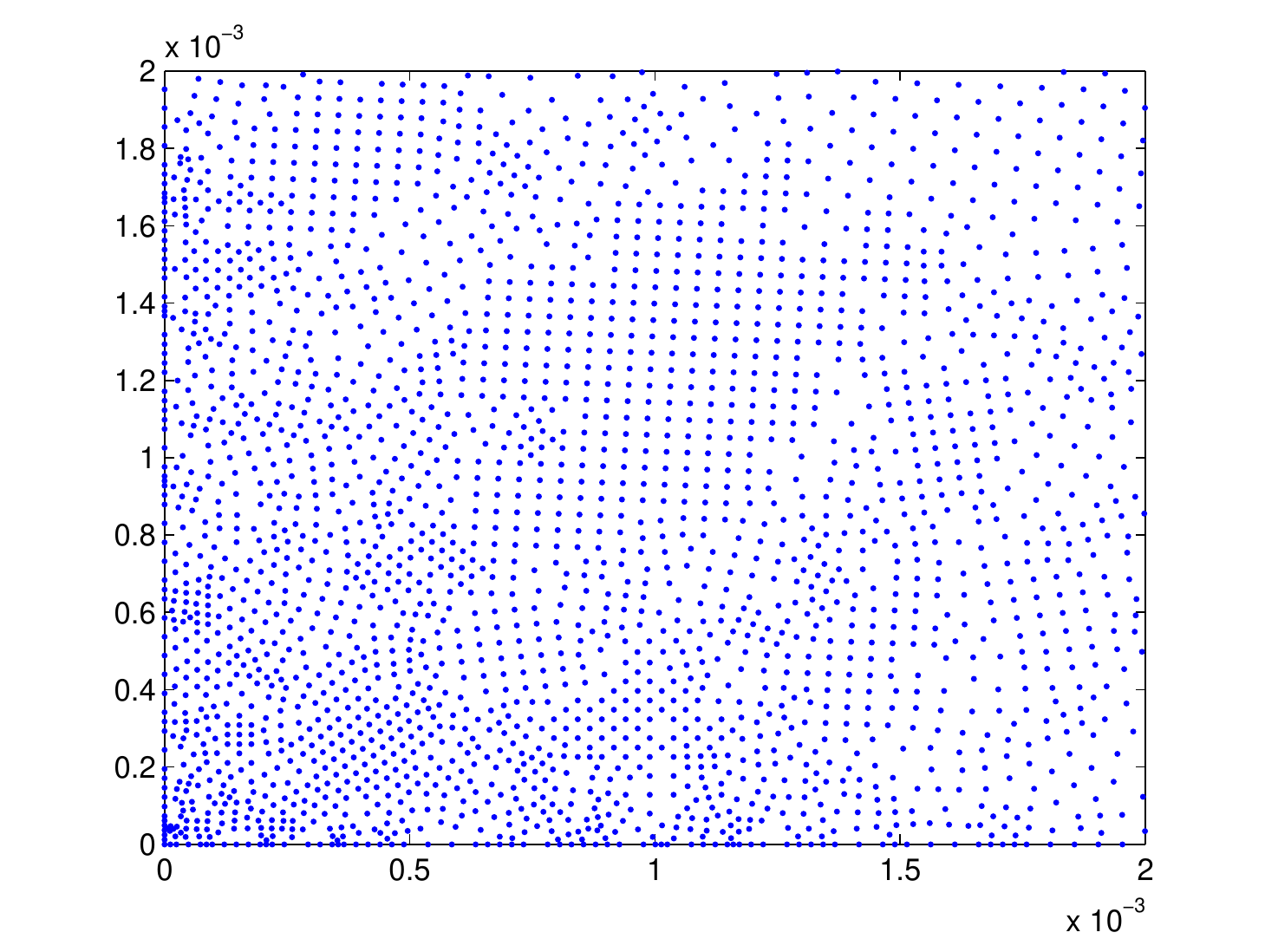}} 
 \subfigure[FEM centers: zoom] %
{\includegraphics[width=6cm,height=4cm]{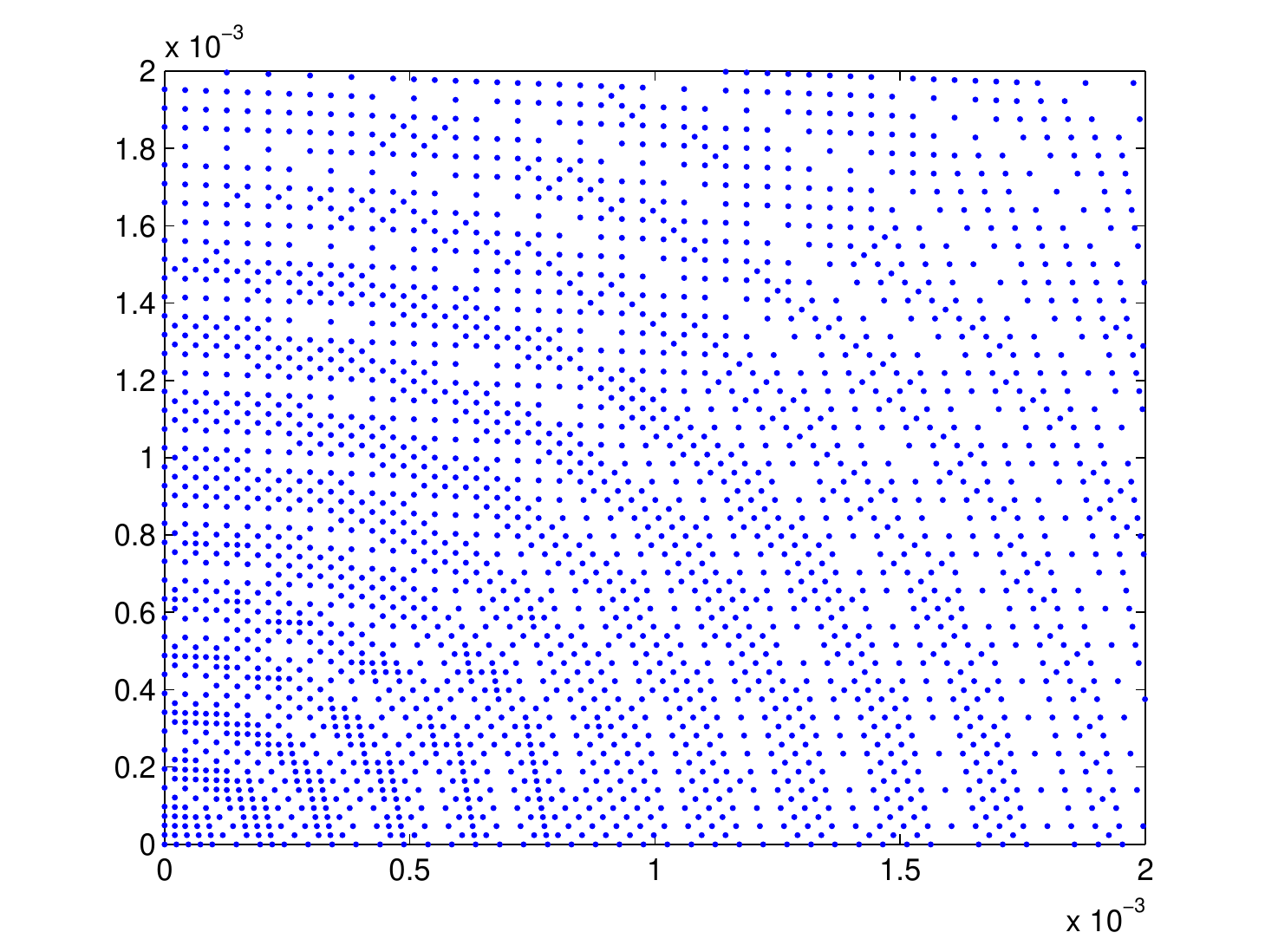}} 
        \end{center}
\caption{Test Problem~\ref{f441} with $\alpha=\frac{1}{50 \pi}$: Errors and centers/vertices.
The plots in  (ceg) are based on the RBF-FD 
solution with 13964 interior centers, whereas those in (dfh) on the FEM 
solution with 14942 interior vertices.
}
\label{mapF442}
\end{figure}

Figure~\ref{mapF442} illustrates the results for Test Problem~\ref{f441} with 
$\alpha=\frac{1}{50 \pi}$. The layout is the same as in Figure~\ref{figF441}, and the results
are similar. There is no ring-like pattern in the locations of the centers in (g) because the 
solution is relatively less resolved than the one in Figure~\ref{figF441}(g). Indeed, in Figures~\ref{mapF442}(gh) we see
about 5 centers across the waves near the origin versus 12 or more in Figures~\ref{figF441}(gh).

\begin{testproblem}\label{f43}\rm 
 \cite[Section 2.4: Peak]{Mitchell2013}
Dirichlet problem (\ref{poi})  for the Poisson equation $\Delta u=f$ in the domain $\Omega=(0,1)^2$, 
where the right hand side $f$ and the boundary conditions are chosen such that the exact solution
is $u(x)=e^{-\alpha \|x-x_0\|^2}$. Following \cite{Mitchell2013}, two 
instances of this problem, with the following values of the parameters $\alpha$ 
(the strength of the peak) and $x_0$ (the location of the peak) will be considered: 
\begin{itemize}
\item[(a)] $\alpha=1000$, $x_0=(0.5,0.5)$,
\item[(b)] $\alpha=100000$, $x_0=(0.51,0.117)$.
\end{itemize}
\end{testproblem}

\begin{figure}[htbp!]
 \begin{center}
\subfigure[\qquad]
{\includegraphics[width=6cm,height=4cm]{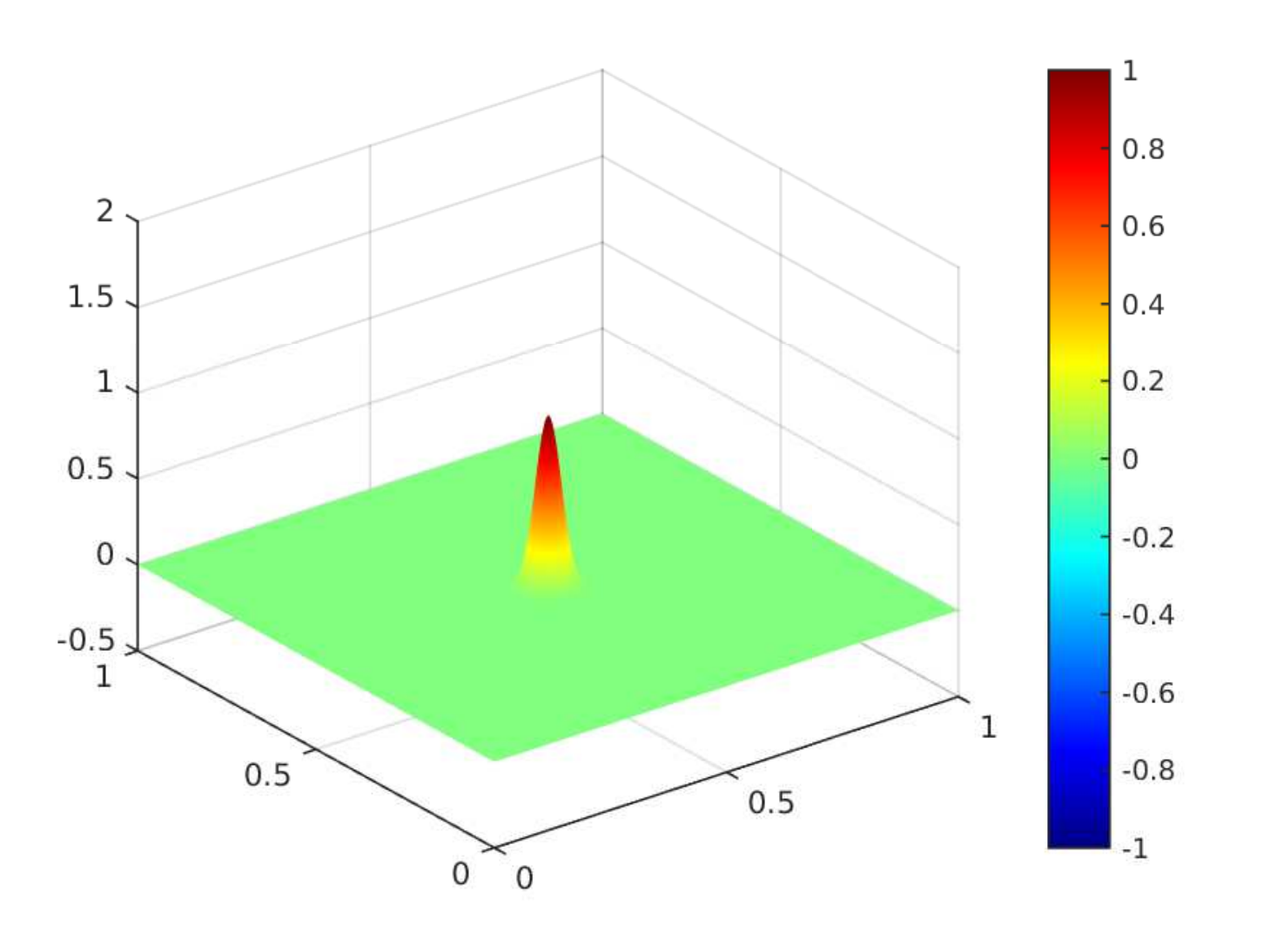}} \quad
 \subfigure[\quad]
{\includegraphics[width=6cm,height=4cm]{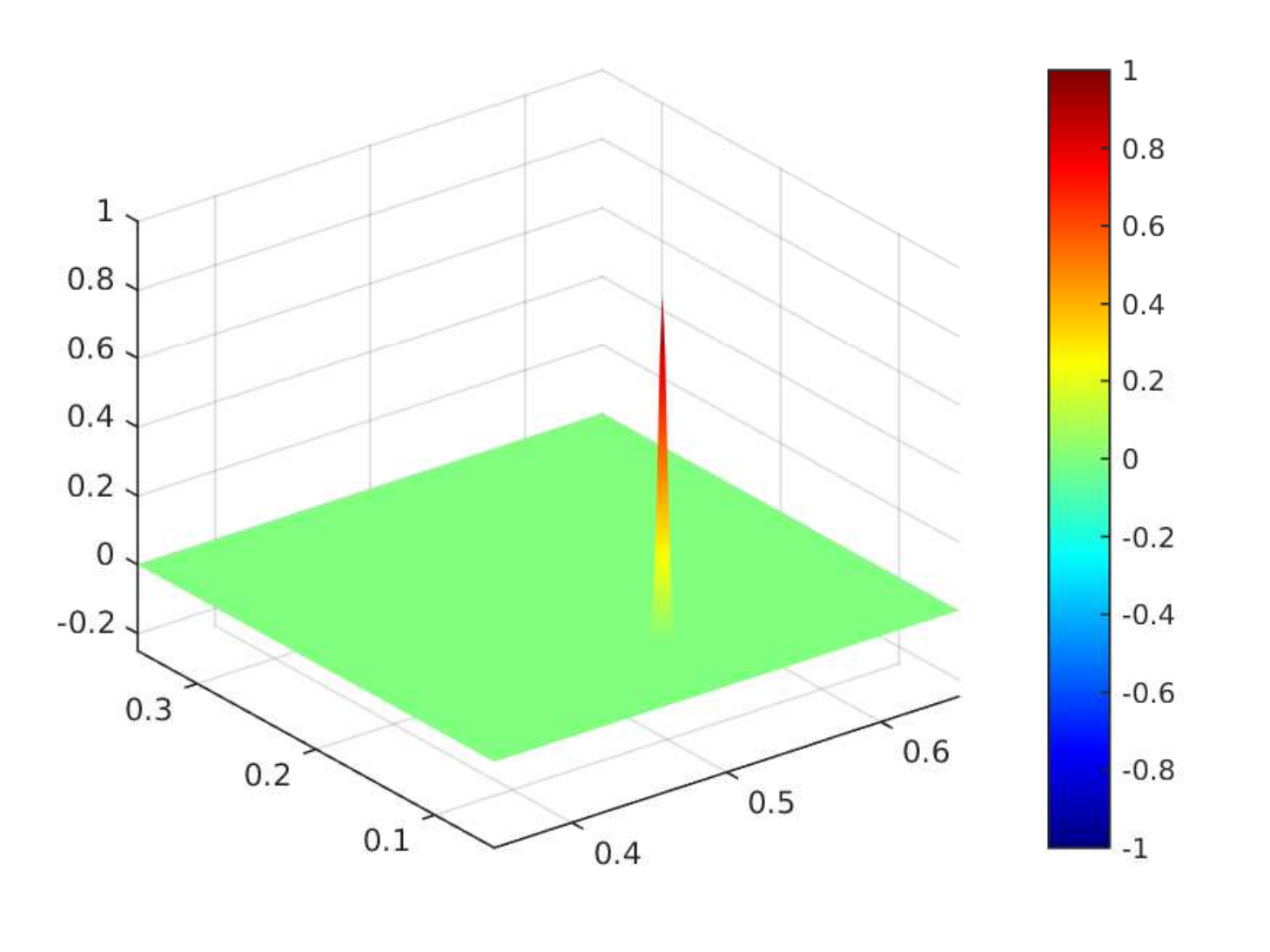}} 
        \end{center}
\caption{Exact solutions for both versions of Test Problem~\ref{f43}.
}
\label{exactF43}
\end{figure}

\begin{figure}[htbp!]
\begin{center} 
 \subfigure[Errors on centers for 5(a)] %
{          \includegraphics[width=6cm,height=4cm]{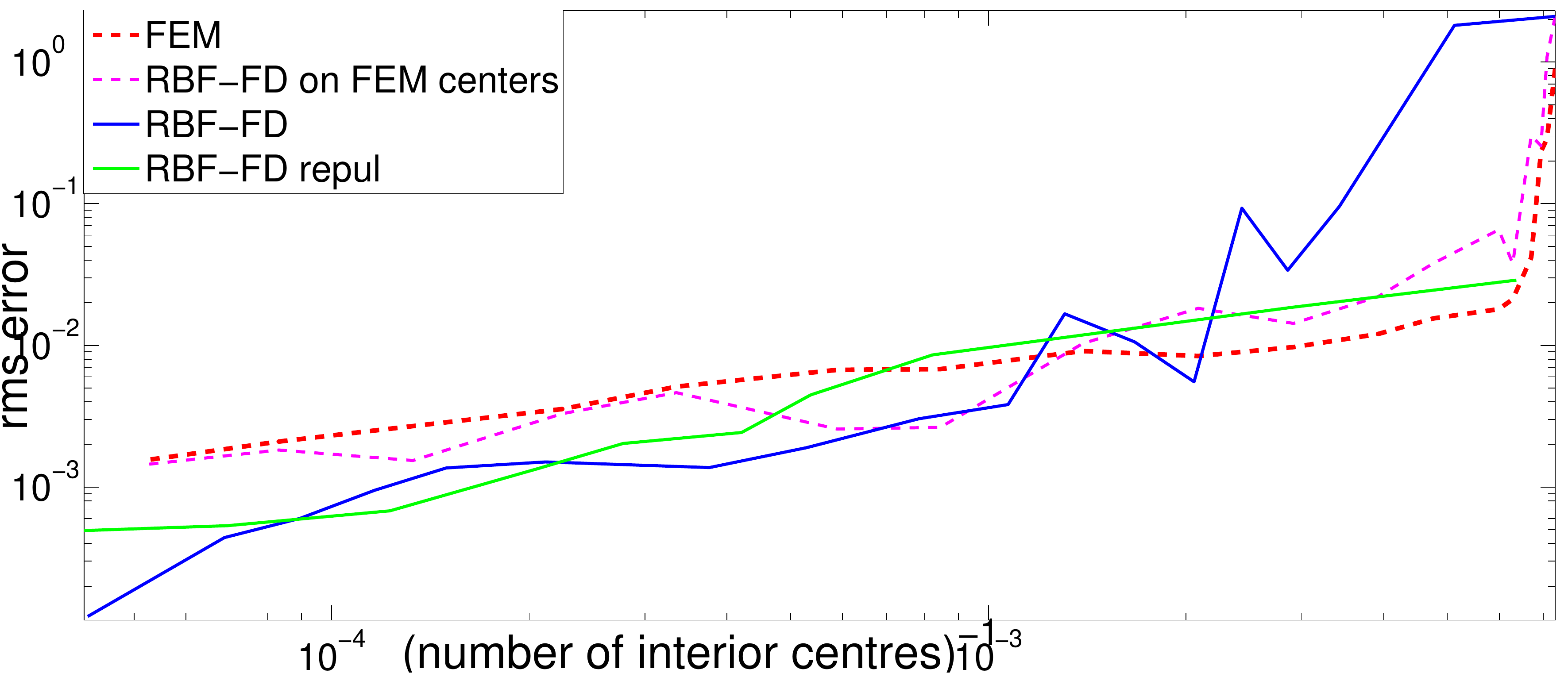}}\qquad
 \subfigure[Errors on grid for 5(a)]
{\includegraphics[width=6cm,height=4cm]{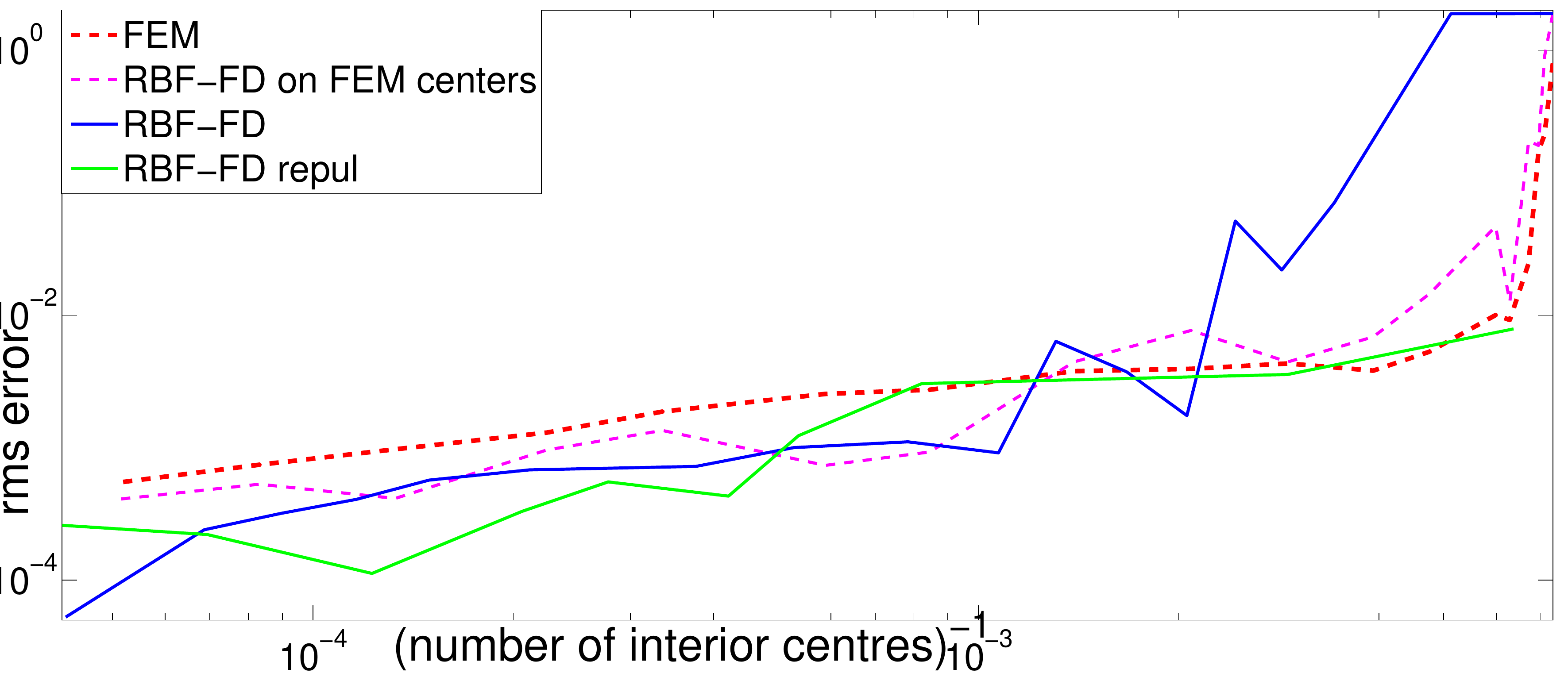}} 
 \subfigure[Errors on centers for 5(b)]
{          \includegraphics[width=6cm,height=4cm]{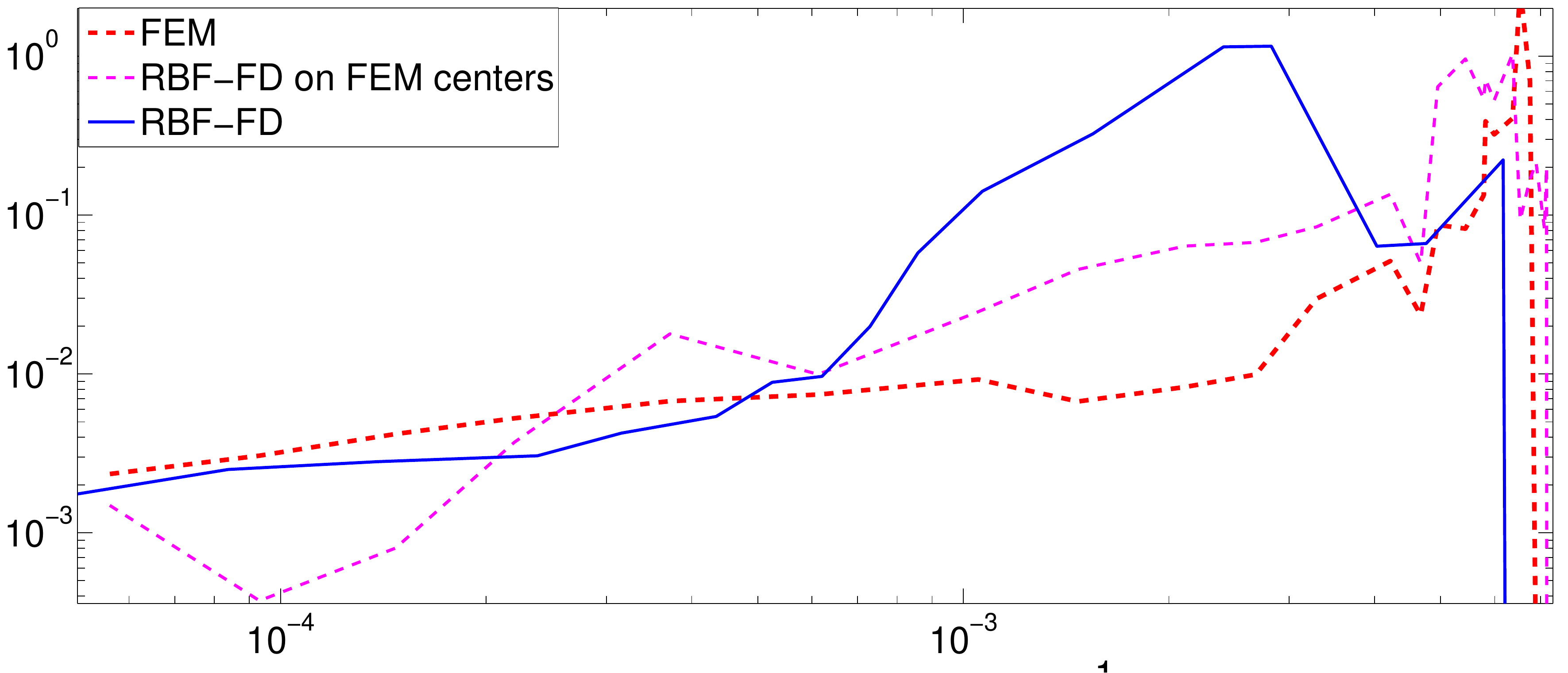}}\qquad
 \subfigure[Errors on grid for 5(b)]
{\includegraphics[width=6cm,height=4cm]{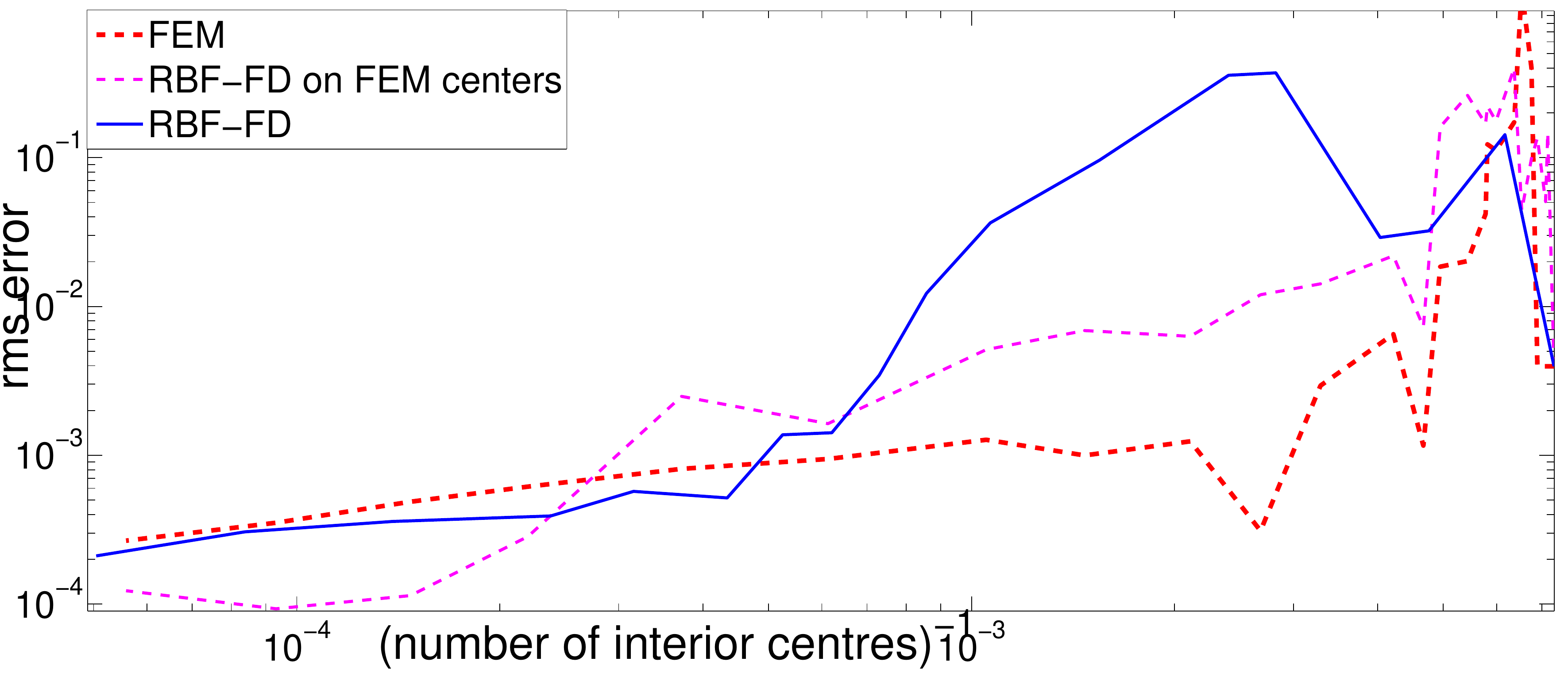}}
\end{center}
\caption{Test Problem~\ref{f43}: Errors with $\alpha = 1000$, $x_0=(0.5,0.5)$ (top) and 
$\alpha = 100000$, $x_0=(0.51,0.117)$ (bottom). {\tt RBF-FD}:  the method of this paper,
{\tt FEM}: finite element method with  piecewise linear shape functions, 
{\tt RBF-FD~on~FEM~centers}: RBF-FD method with stencil support selection
by Algorithm~\ref{alg1} on centers
generated by the adaptive finite element method, {\tt RBF-FD~repul}: RBF-FD method 
with stencil support selection by Algorithm~\ref{alg1} on centers
generated by \emph{a priori} refinement using {\tt DistMesh}.}
\label{errors43New}
\end{figure}

The exact solutions are presented in Figure~\ref{exactF43}, and the
error curves in Figure~\ref{errors43New}. In contrast to the previous test problems
there is a significant difference in the behavior of the FEM and RBF-FD methods.
The errors of the RBF-FD method are smaller on denser sets of 
centers, but higher for a number of initial refinements.
Therefore we also computed  RBF-FD  solutions on the centers/vertices generated by
the adaptive finite element method and included the respective curves in the
same figure. The errors in this case are much closer to those by the FEM. In addition, in 
the case (a) %
we also produced smoothly distributed centers using
{\tt distmesh2d}, which  gives results closer to those by the FEM as well. We again
used different values for the initial edge length $h_0$ and
the edge length function in the form $\sigma(r)=r^\beta$ to obtain sets
with different numbers of centers, and picked the constellations of $h_0,\beta$
shown in Table~\ref{h0Pow5a} which produce results with relatively small errors of the RBF-FD
method.

\begin{table}[htbp!]
\begin{center}\renewcommand{\arraystretch}{1.2}\small
\begin{tabular}{|c|c|c|c|c|c|c|c|c|c|c|}
\hline
\hline
$\#\Xi_{\rm int}$ &157    &   343     & 1216  &    1863 &     2375   &   3600    &  4850   &   8150  &   14411   &  24191 \\ 
\hline
$h_0$&  .00036& .001 &.00055 &.000325&  .00045 &.0006& .00045&  .00035&  .00045& .0005   \\
\hline
$\beta$ & 0.675 & 0.65 & 0.525 &0.575 & 0.525 &0.475&  0.475 &0.475& 0.4& 0.375\\
 \hline
\end{tabular}
\caption{ 
Parameters of {\tt distmesh2d.m} used to produce the \emph{a priori} sets of centers 
for  Test Problem~\ref{f43} with $\alpha = 1000$, $x_0=(0.5,0.5)$: $\#\Xi_{\rm int}$ is the number of interior centers,  $h_0$
the initial edge length, and  $\beta$ the power of the scaled edge 
length function in the form $\sigma(r)=r^{\beta}$. The sets with 
$\#\Xi_{\rm int}=343$ and 1863 are
illustrated in Figure~\ref{F431rep_centers}.}
\label{h0Pow5a}
\end{center}
\end{table}

\begin{figure}[htbp!]
 \begin{center}
 \subfigure[FEM solution (343)]
{\includegraphics[width=6cm,height=4cm]{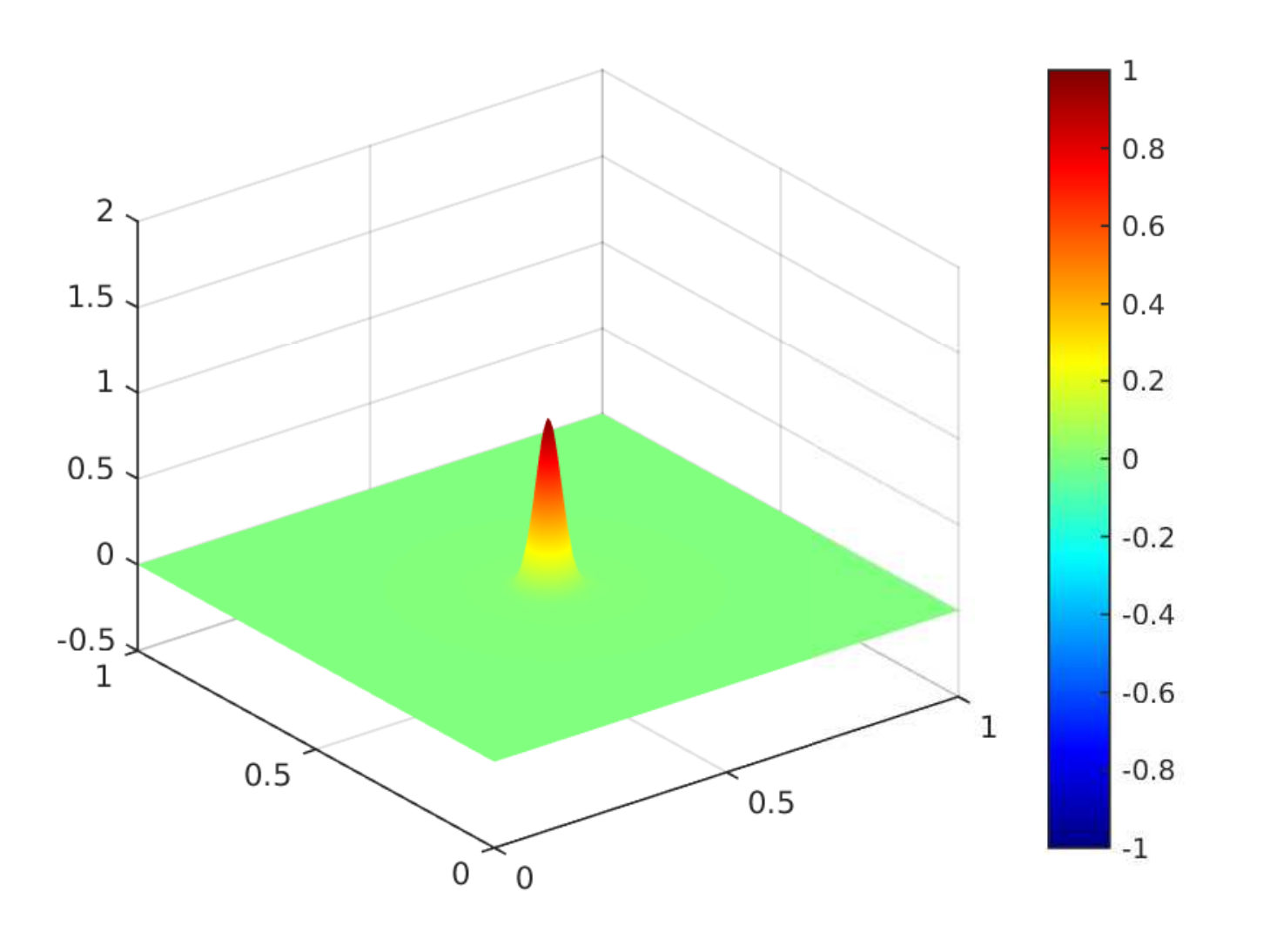}} \qquad 
 \subfigure[RBF-FD solution (350)]
{\includegraphics[width=6cm,height=4cm]{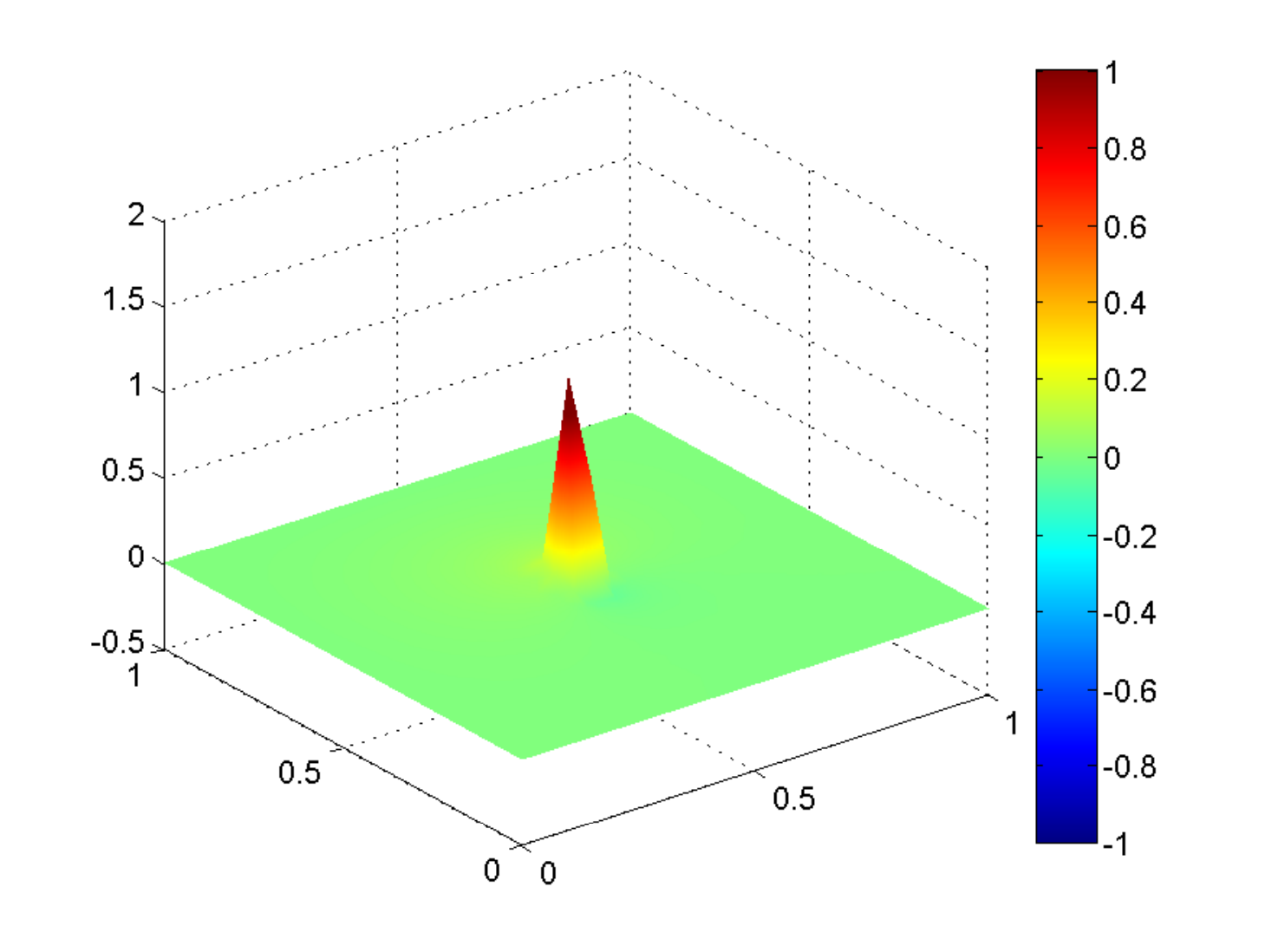}}  \qquad
\subfigure[RBF-FD on FEM centers (343)]
{\includegraphics[width=6cm,height=4cm]{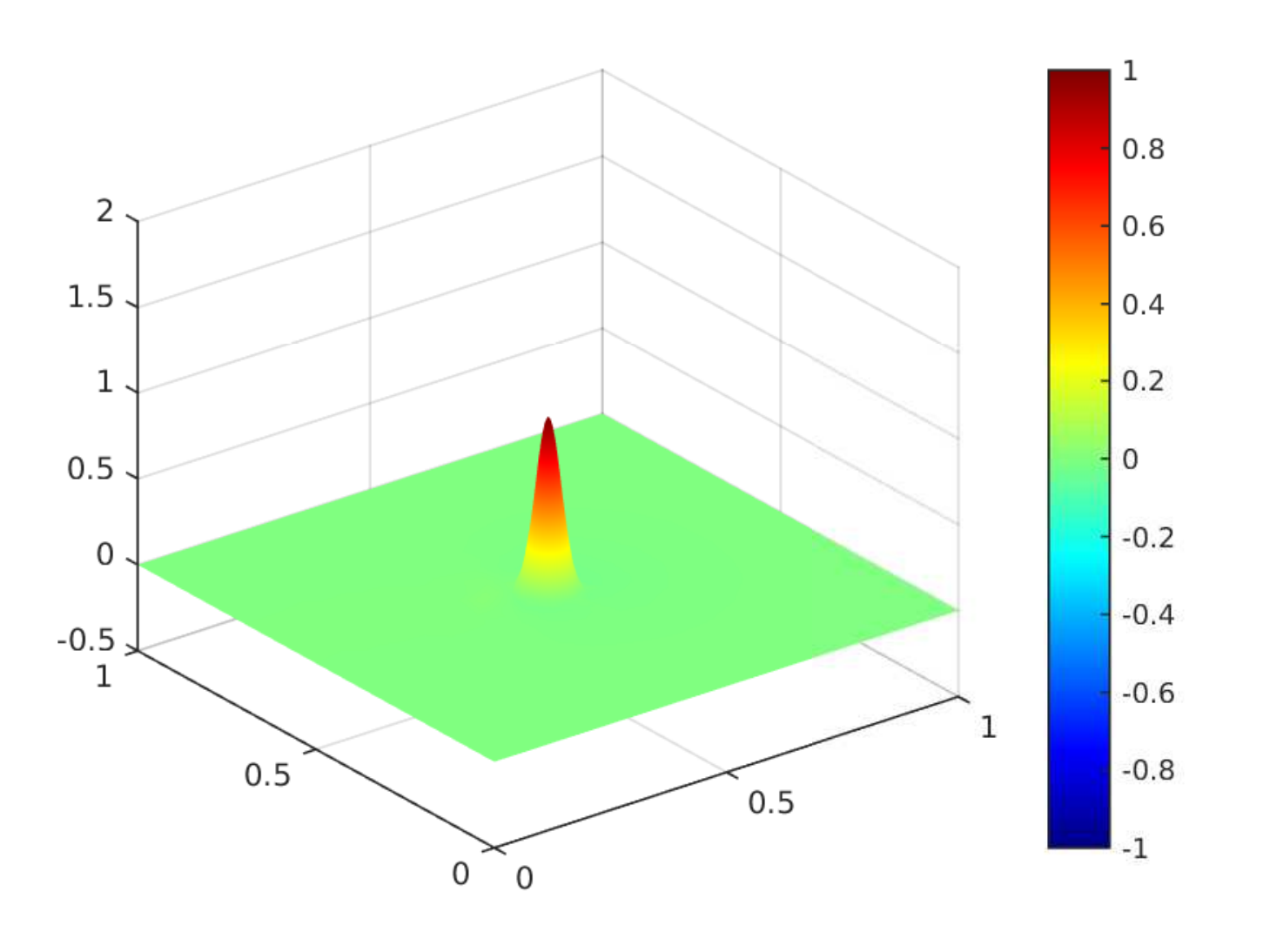}}  \qquad
 \subfigure[RBF-FD on smooth centers (343)]
{\includegraphics[width=6cm,height=4cm]{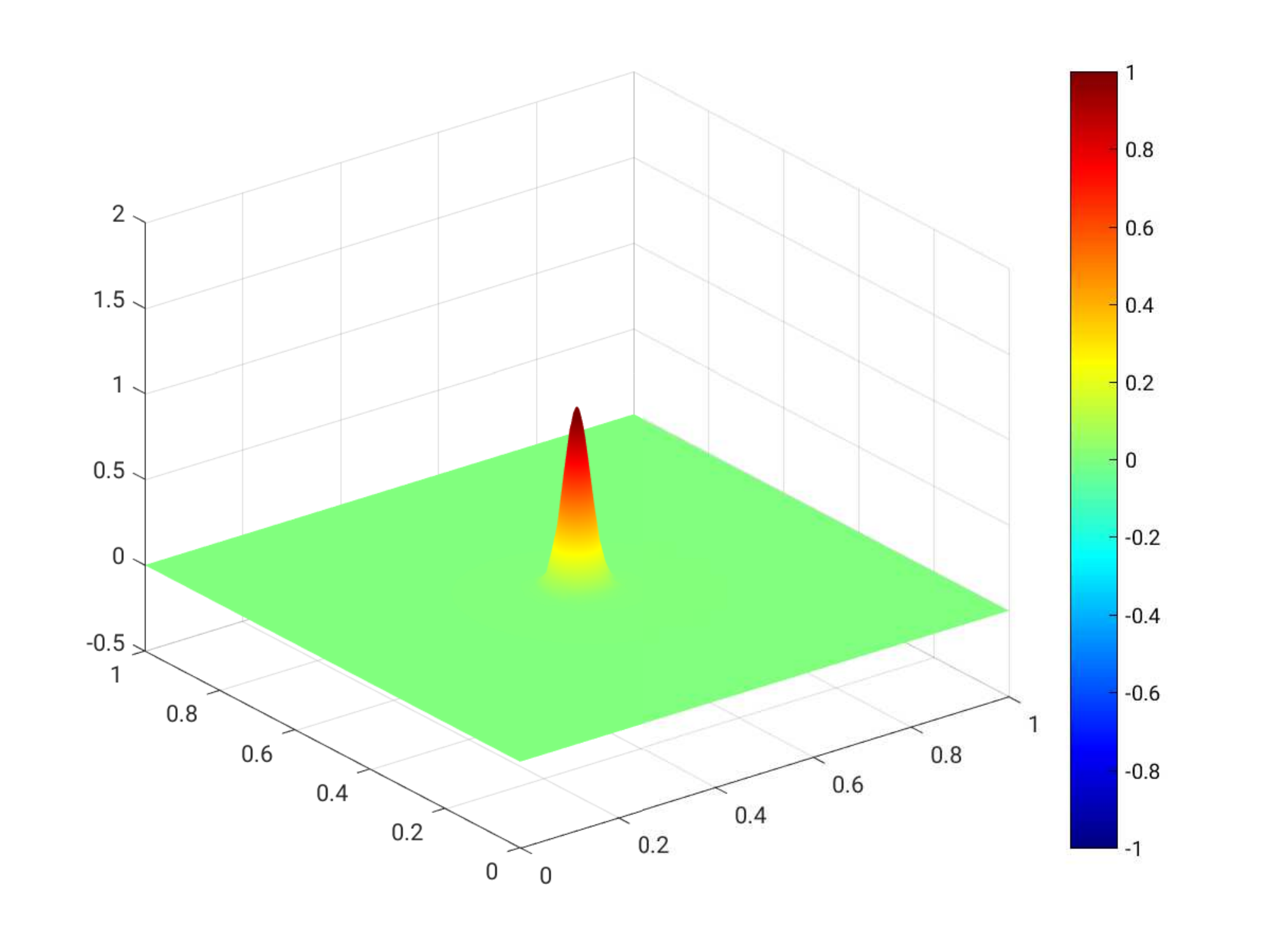}}  \qquad
 \subfigure[FEM error (1699)]
{\includegraphics[width=6cm,height=4cm]{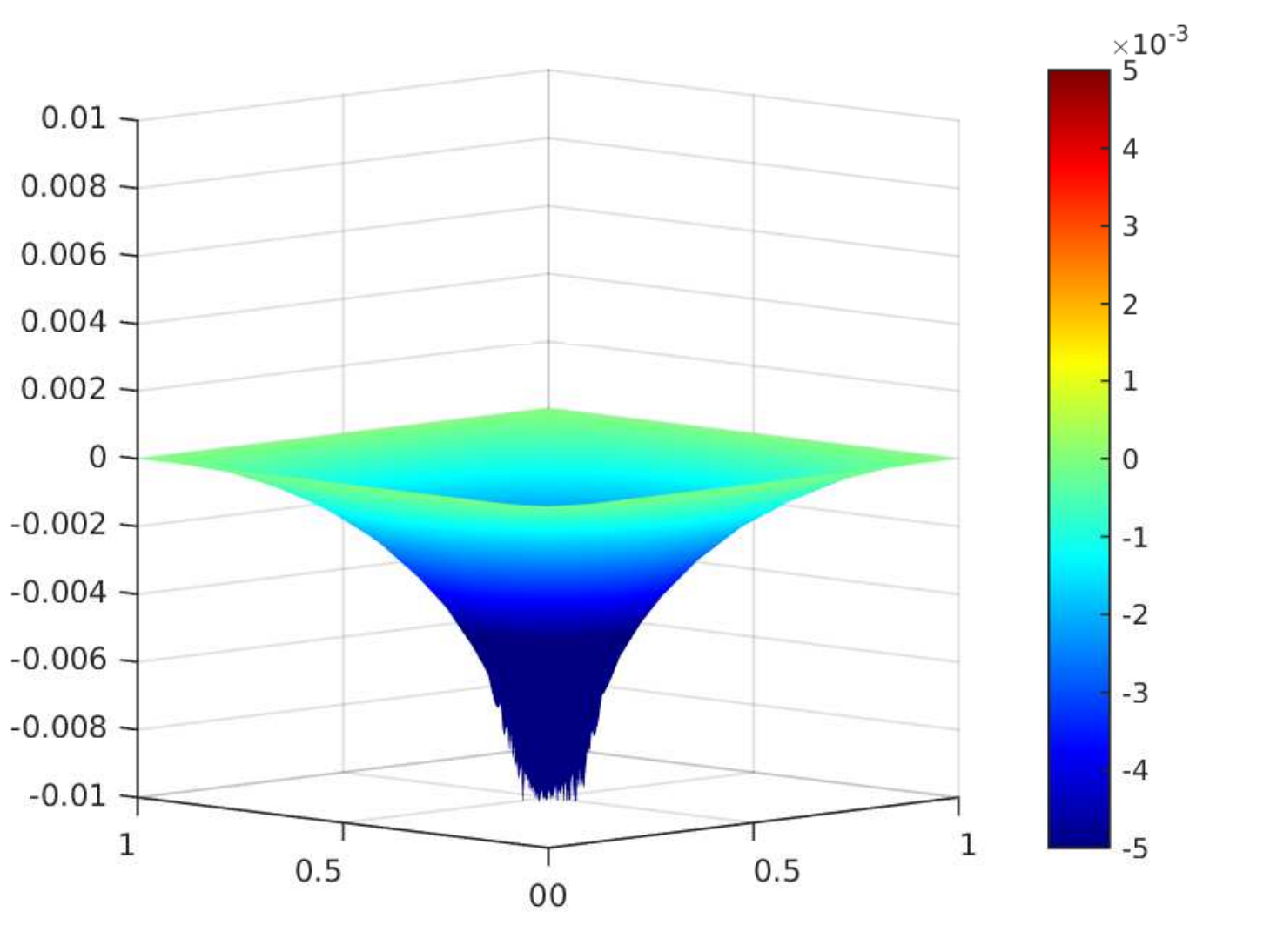}}  \qquad 
  \subfigure[RBF-FD error (1893)] %
{\includegraphics[width=6cm,height=4cm]{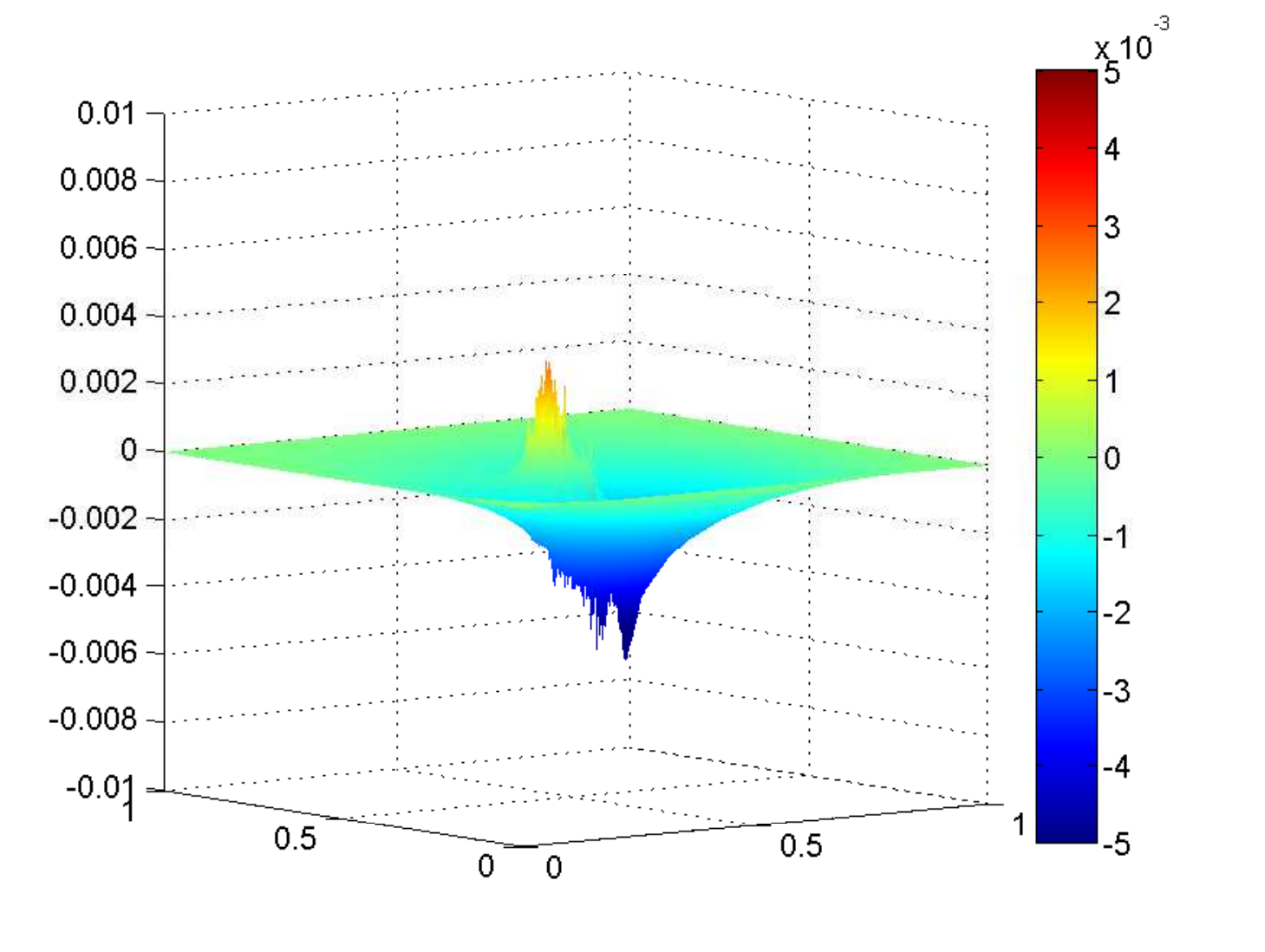}}  \qquad
 \subfigure[RBF-FD on FEM centers  (1699)]
{\includegraphics[width=6cm,height=4cm]{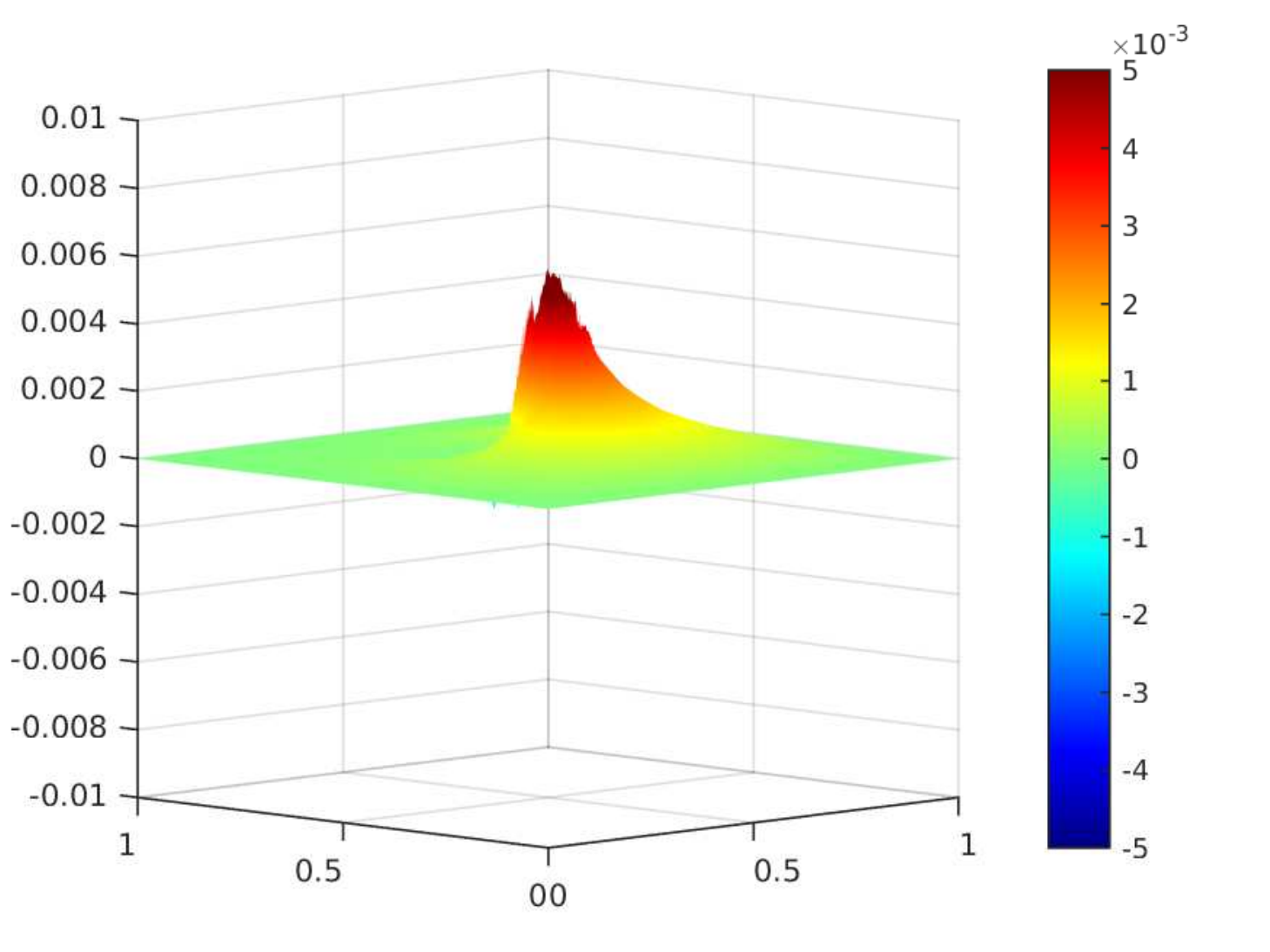}}  \qquad
\subfigure[RBF-FD on smooth centers  (1863)]
{\includegraphics[width=6cm,height=4cm]{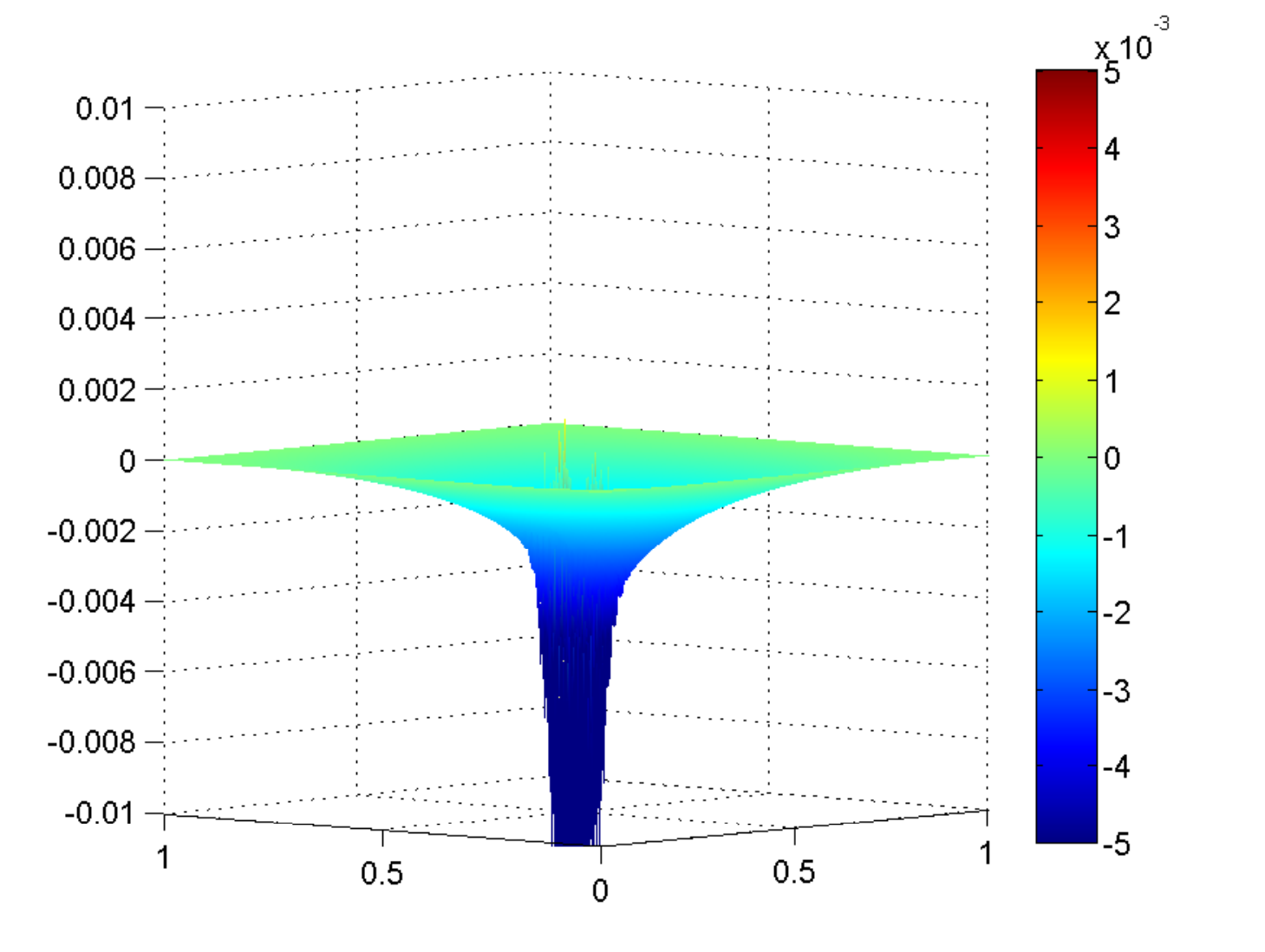}}  \qquad
        \end{center}
\caption{Test Problem~\ref{f43}  with $\alpha = 1000$ and $x_0=(0.5,0.5)$:
(a) FEM solution for a triangulation with 343 interior vertices,
(b) RBF-FD solution for 350 interior centers, 
(c) RBF-FD solution  with stencil support selection
by Algorithm~\ref{alg1} on centers at the 343 interior vertices of FEM triangulation,
(d) RBF-FD solution with stencil support selection
by Algorithm~\ref{alg1} on centers
generated by \emph{a priori} refinement using {\tt DistMesh},
(e-h) error plots for the FEM solution for a triangulation with 
1699 interior vertices and three versions of RBF-FD as above, with the number of
interior vertices given in  brackets. 
}
\label{mapF431new}
\end{figure}

\begin{figure}[htbp!]
 \begin{center} 
 \subfigure[RBF-FD centers {(350)}]
{\includegraphics[width=6cm,height=4cm]{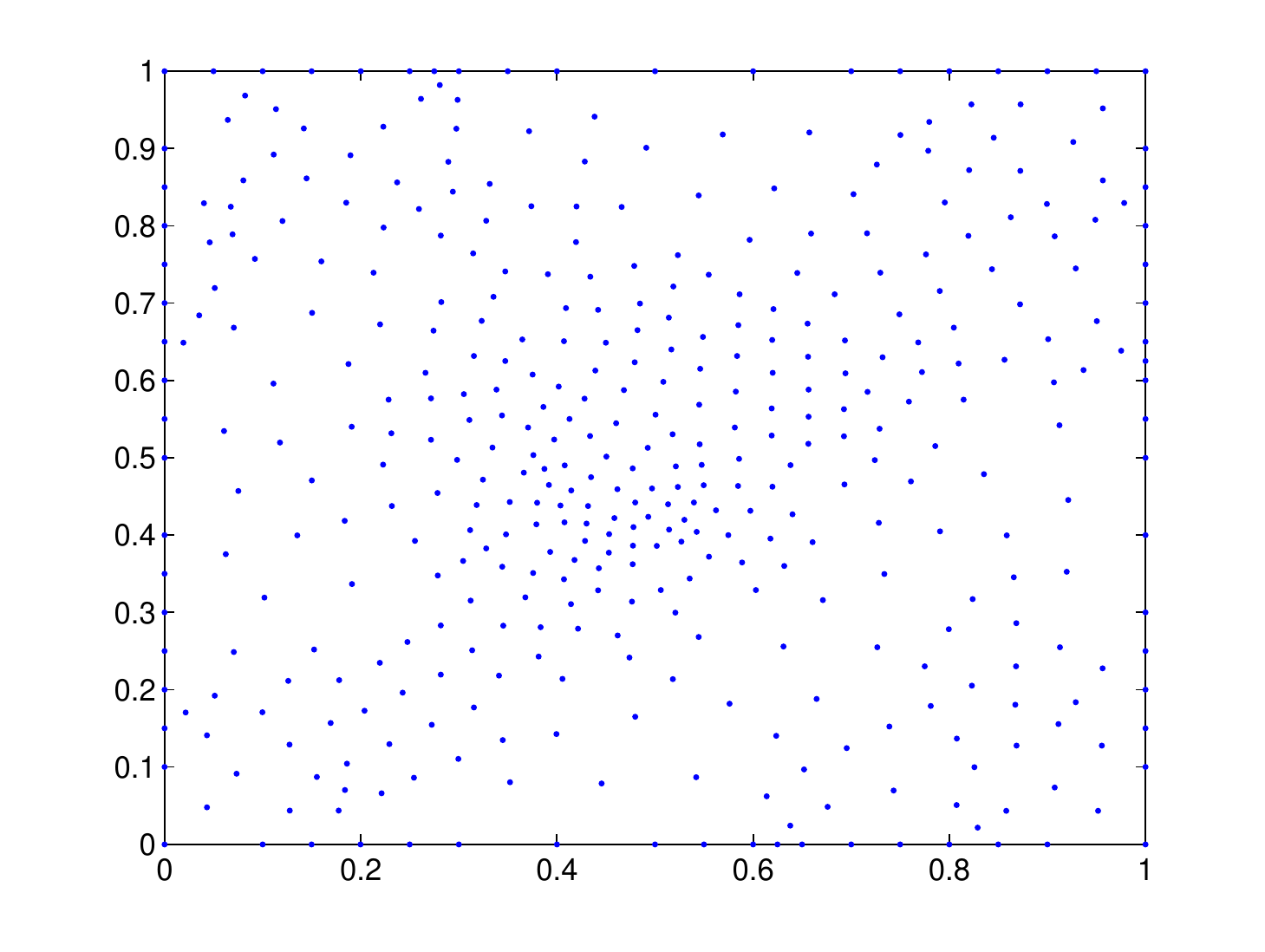}}  \qquad
 \subfigure[FEM  centers (343)]
{\includegraphics[width=6cm,height=4cm]{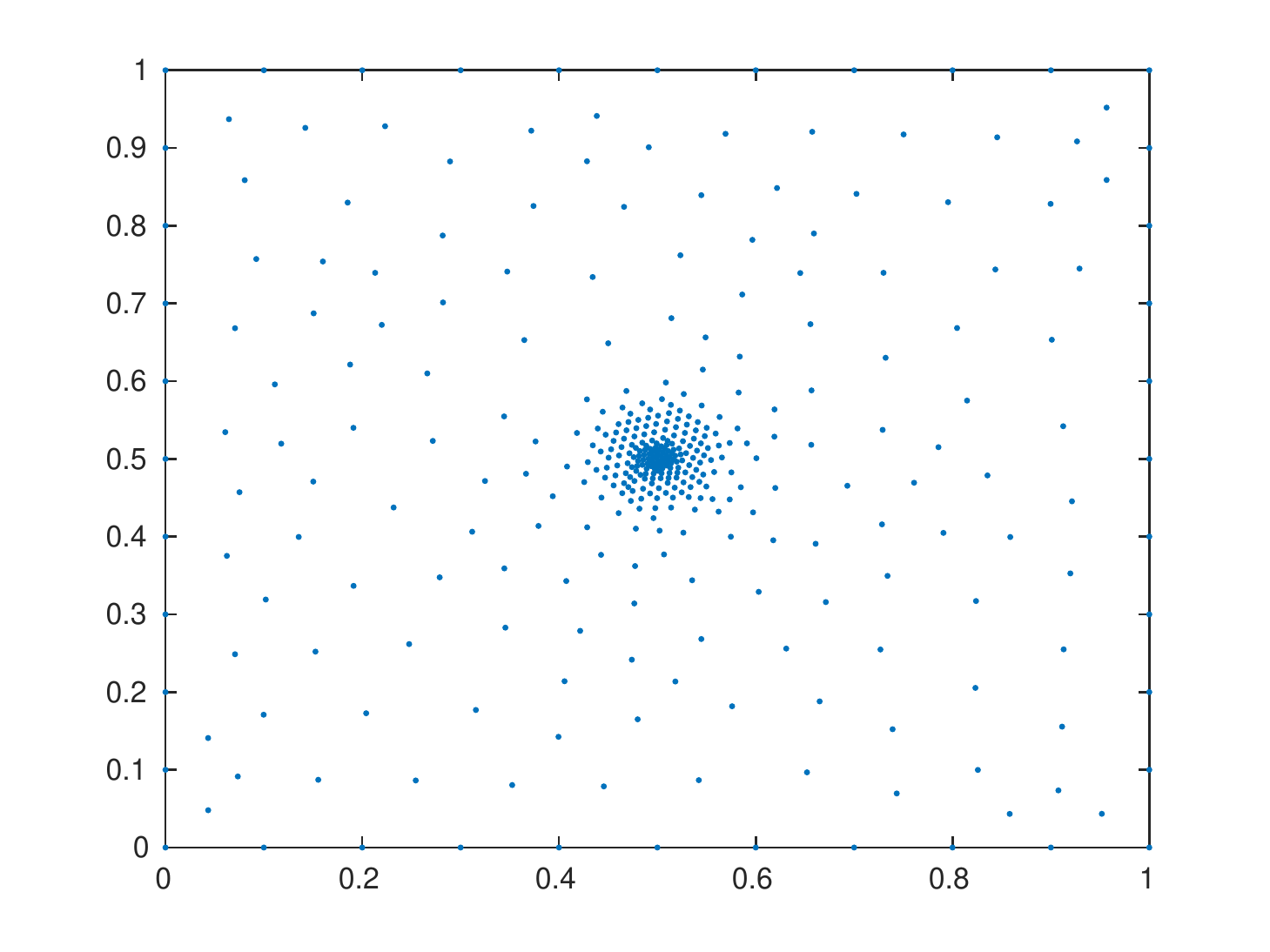}}  \qquad
 \subfigure[RBF-FD  centers  {(765)}]
{\includegraphics[width=6cm,height=4cm]{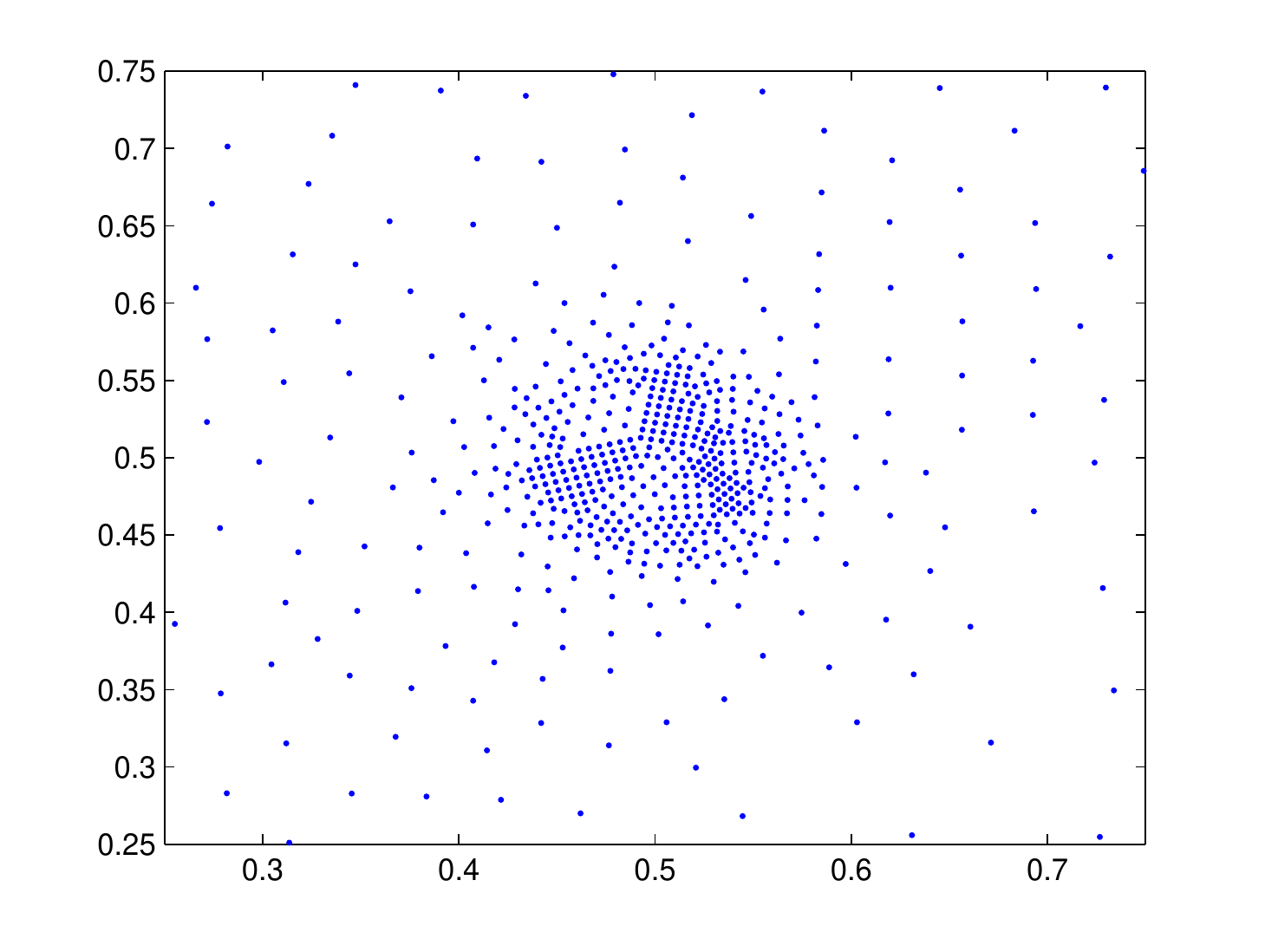}}  \qquad
 \subfigure[FEM centers (718)]
{\includegraphics[width=6cm,height=4cm]{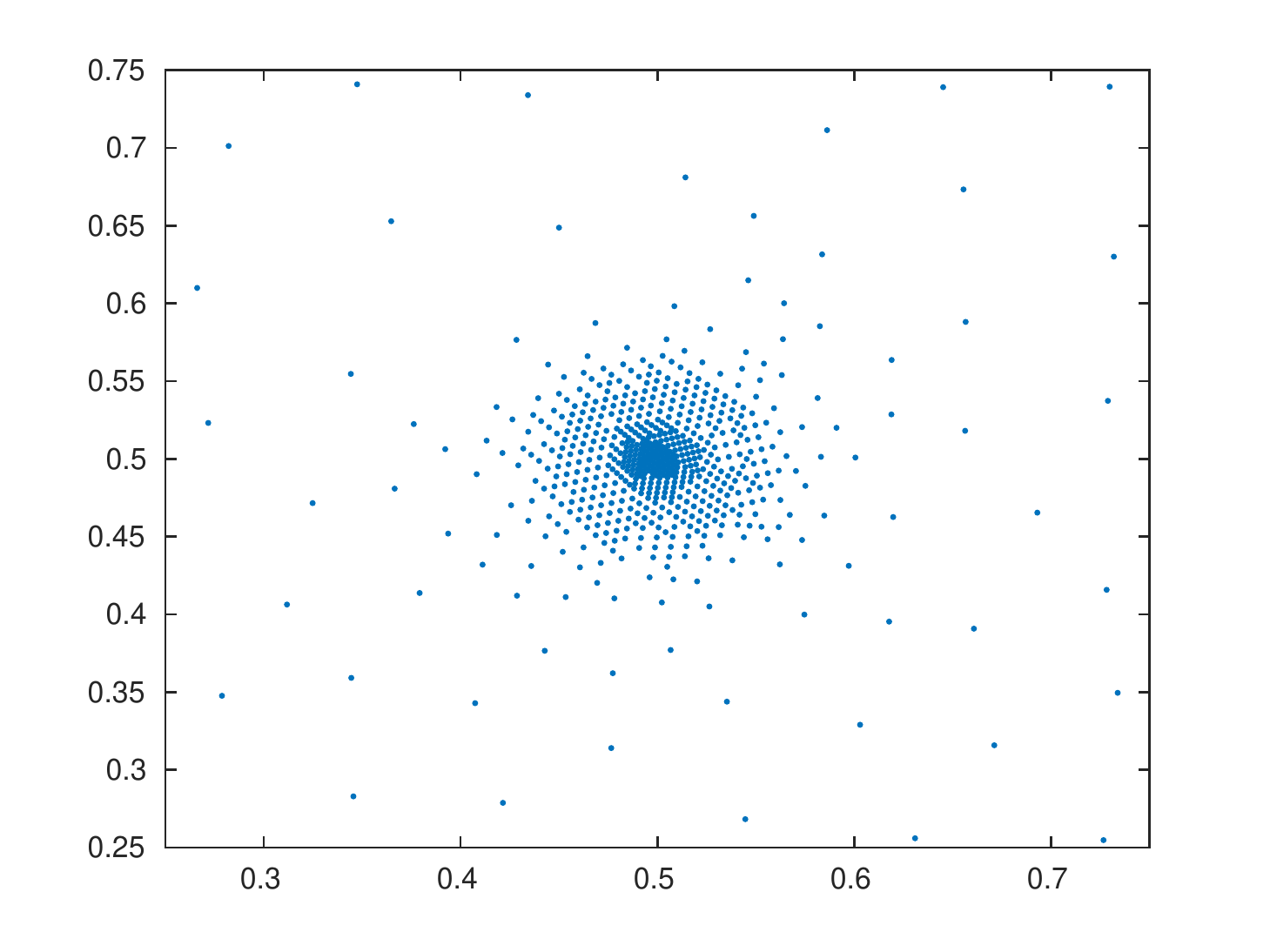}}  \qquad
 \subfigure[RBF-FD centers  {(1893)}]
{\includegraphics[width=6cm,height=4cm]{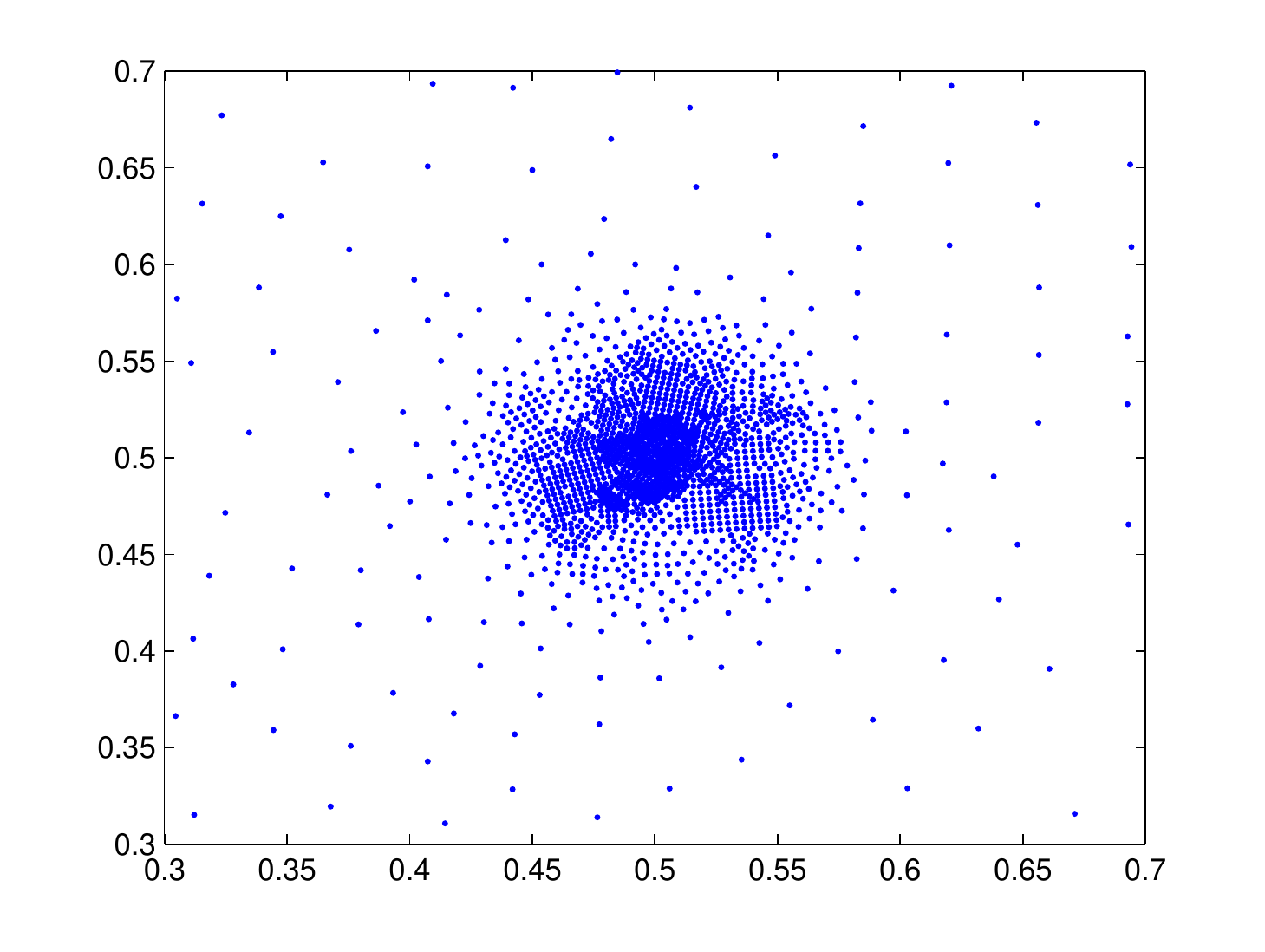}}  \qquad
 \subfigure[FEM centers (1699)]
{\includegraphics[width=6cm,height=4cm]{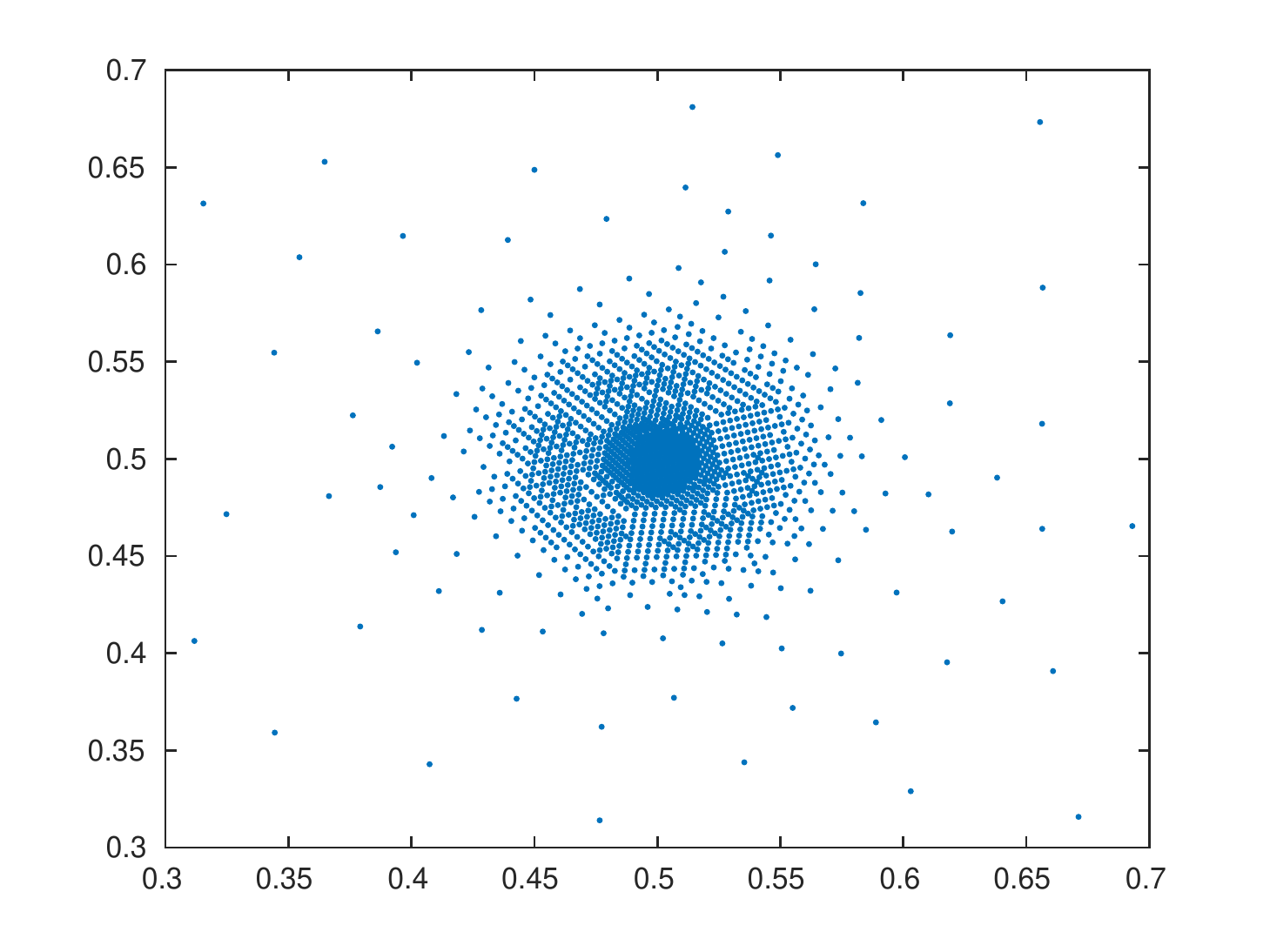}}  \qquad
 \subfigure[RBF centers  {(6689)}]
{\includegraphics[width=6cm,height=4cm]{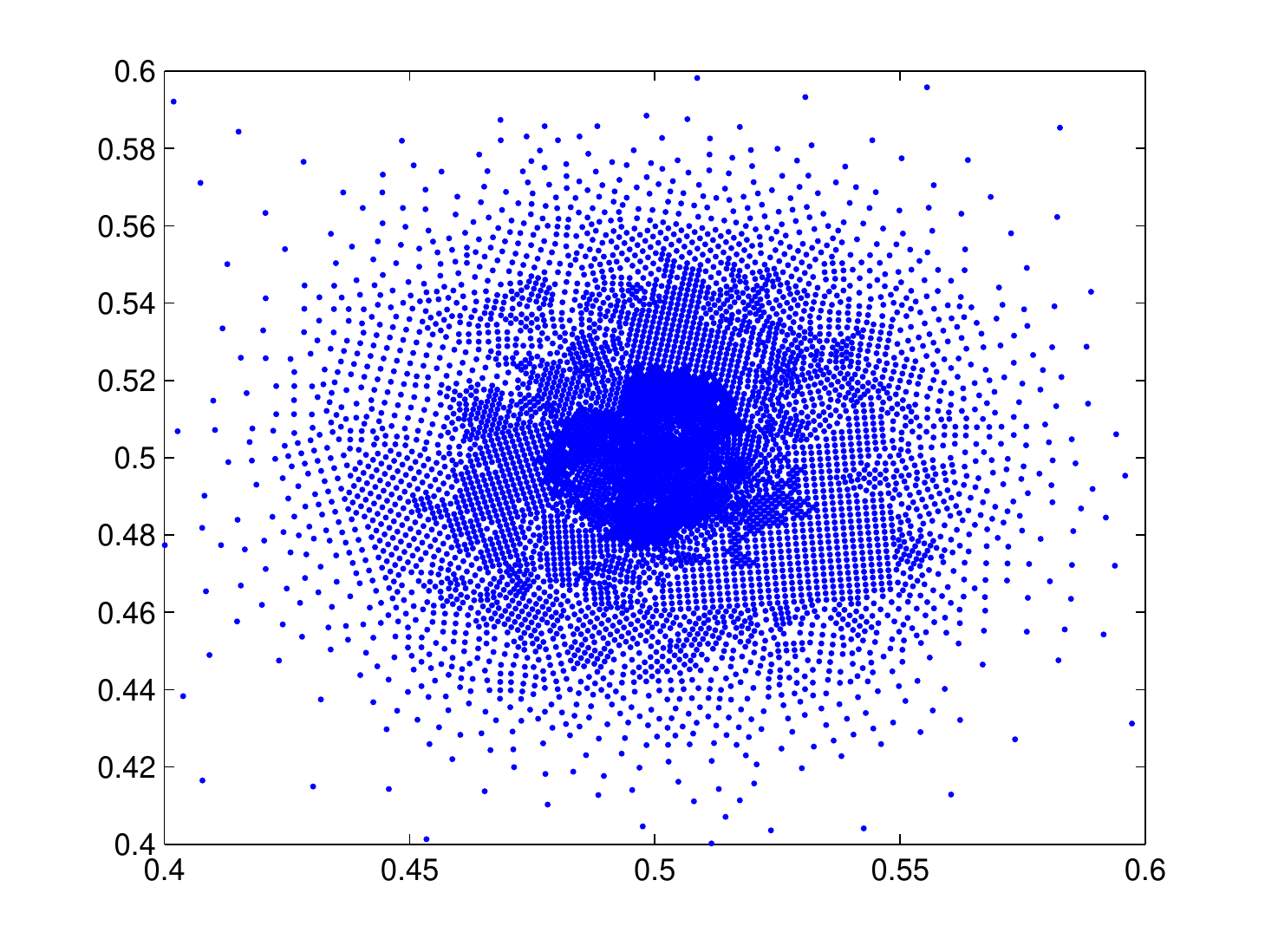}}  \qquad
 \subfigure[FEM centers (7520)]
{\includegraphics[width=6cm,height=4cm]{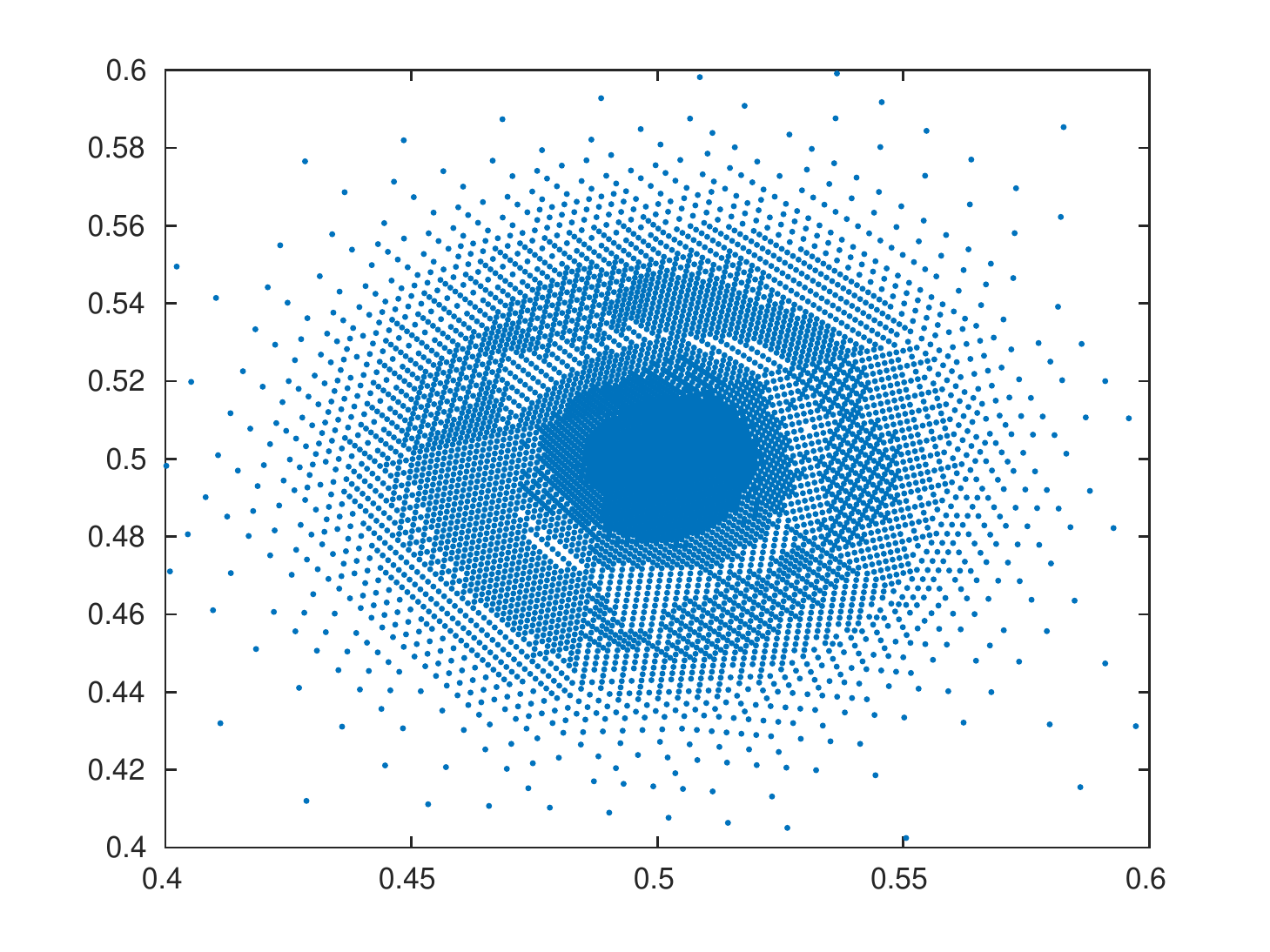}} 
        \end{center}
\caption{Test Problem~\ref{f43}  with $\alpha = 1000$ and $x_0=(0.5,0.5)$: Centers for RBF-FD and FEM 
adaptive refinements.
The number of interior centers/vertices is shown in brackets. The number of refinements of 
an initial set of centers: (a) 3, (b) 8, (c) 7, (d) 10, (e) 10, (f) 14,
(g) {13}, (h) 15.
}
\label{centresF431New}
\end{figure}

\begin{figure}[htbp!]
 \begin{center}
 \subfigure[Repulsion centers {(343)}]
{\includegraphics[width=6cm,height=4cm]{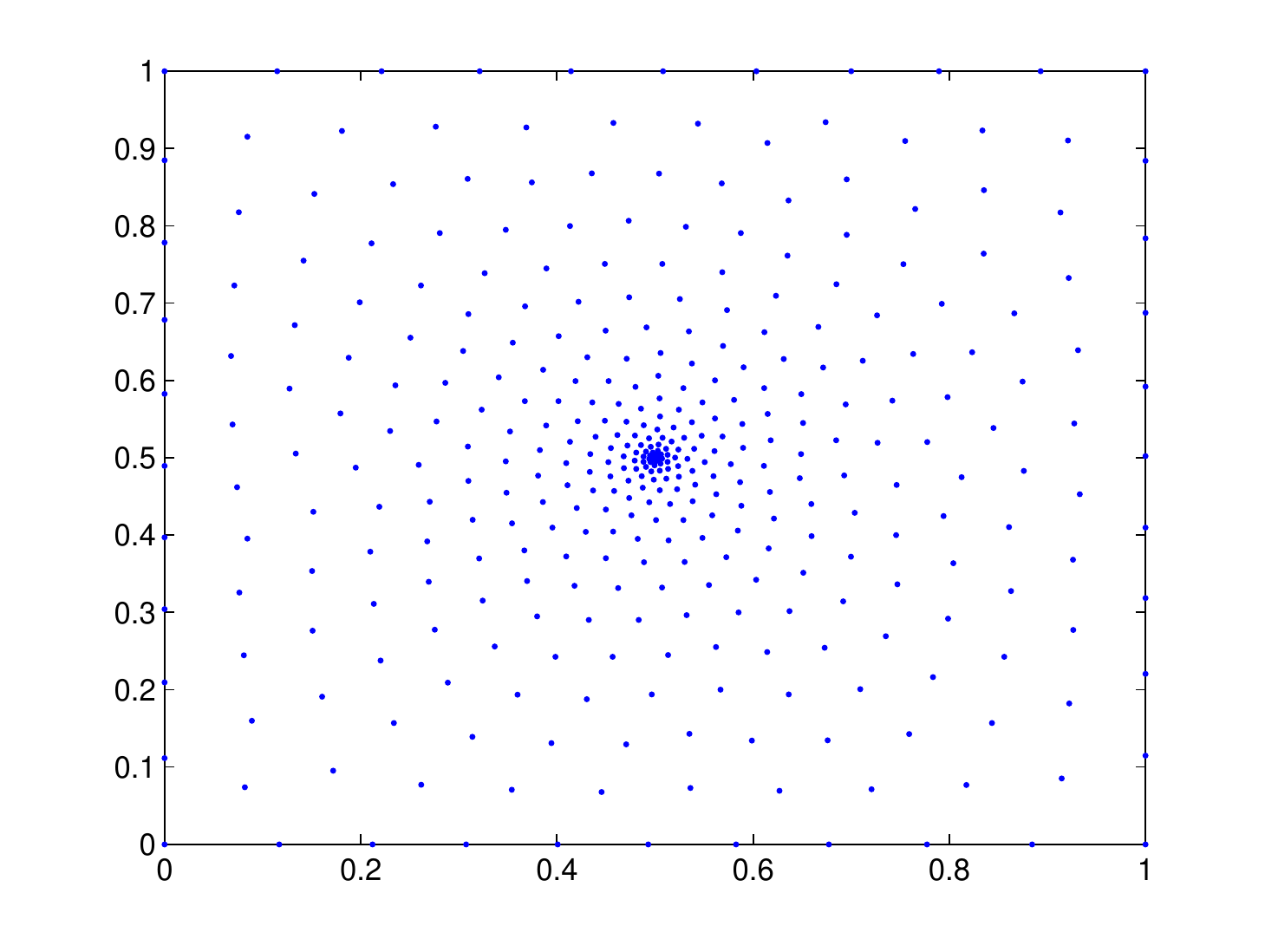}}  \qquad
  \subfigure[Repulsion centers  {(1863)}]
{\includegraphics[width=6cm,height=4cm]{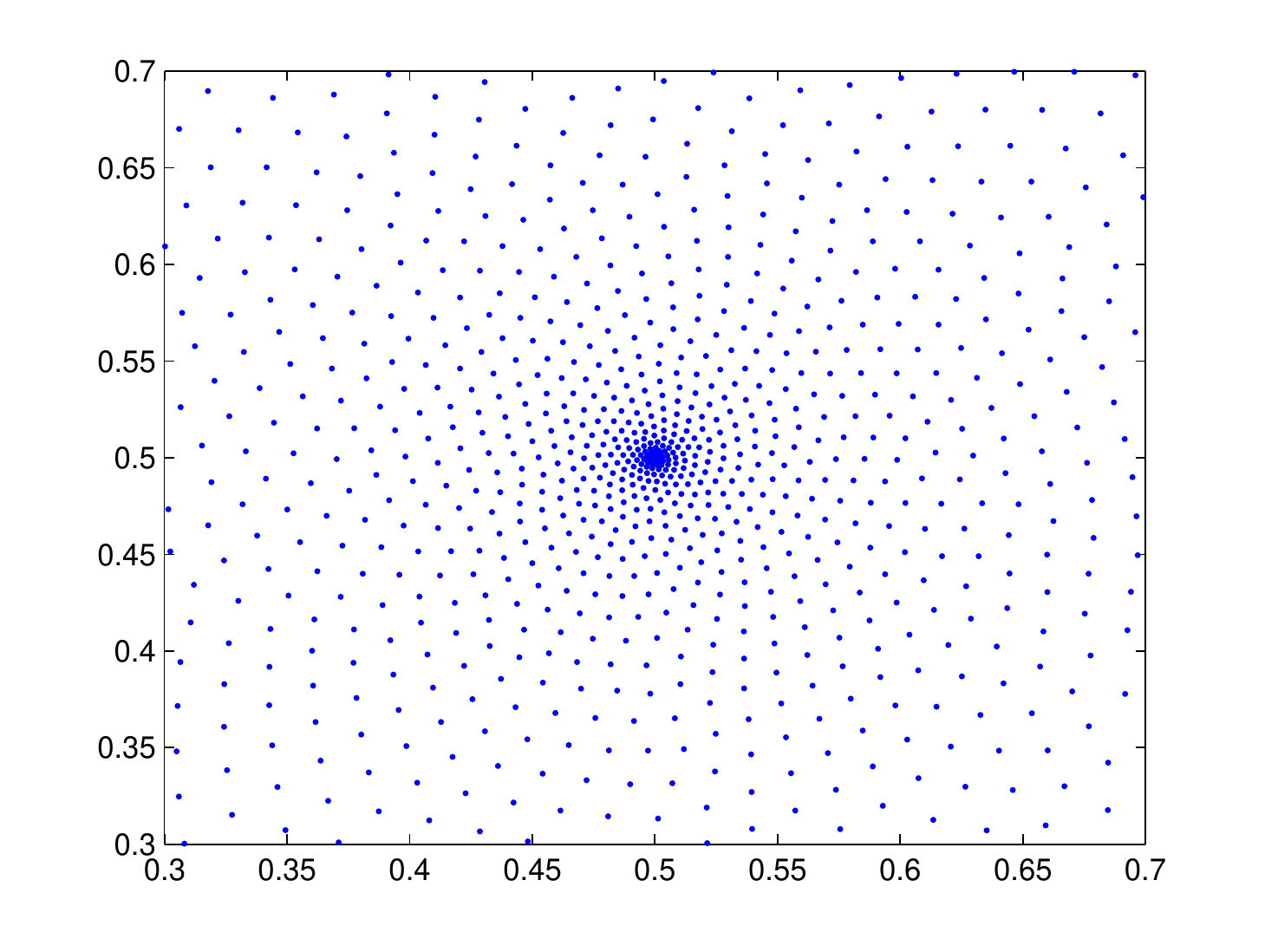}} 
       \end{center}
\caption{Test Problem~\ref{f43}  with $\alpha = 1000$ and $x_0=(0.5,0.5)$: Centers
generated by \emph{a priori} refinement using {\tt DistMesh} with (a) $h_0=0.001$, $\beta=0.65$
and (a) $h_0=0.000325$, $\beta=0.575$.
}
\label{F431rep_centers}
\end{figure}

Some solutions of Test Problem~\ref{f43} with $\alpha = 1000$, 
$x_0=(0.5,0.5)$ by FEM
and the three versions of RBF-FD are illustrated in Figure~\ref{mapF431new}, whereas the centers 
produced by all methods can be seen in Figures~\ref{centresF431New} and \ref{F431rep_centers}.
Figures~\ref{mapF431new}(a-d) compare solutions for a small number of interior
centers/vertices of about 350. We can see in particular that although
the adaptive RBF-FD solution is
visually different from the exact solution in Figure~\ref{exactF43}(a), it does
have the shape of a peak in the correct location. Figures~\ref{mapF431new}(e-h)
compare the errors of the solutions with about 1700--1900 interior centers/vertices.
In this case the RBF-FD-based solutions are more accurate than the FEM solution.
Furthermore, Figure~\ref{centresF431New} compares adaptively generated
centers in the RBF-FD method using Algorithms~\ref{alg1} and \ref{alg2Refi} with
 the  distribution of the 
vertices of the triangulations used by the FEM. In particular, the plots in Figure~\ref{centresF431New}(ab) explain the difference in the
solutions shown in  Figures~\ref{mapF431new}(ab) as they indicate much higher concentration of the FEM vertices near the peak in comparison to
the RBF-FD centers. Note that 350 interior RBF-FD centers are obtained after just 3 refinements, whereas 343 interior vertices
of FEM are generated after 8 refinement steps of the same initial set of centers.
If we compute an RBF-FD solution using
Algorithm~\ref{alg1} on the FEM centers shown in Figure~\ref{centresF431New}(b), then the result is very
close to the FEM solution, see Figure~\ref{mapF431new}(c) and the corresponding points on the curves in 
Figure~\ref{errors43New}(ab). The RBF-FD solution using
Algorithm~\ref{alg1} on the smoothly distributed centers shown in Figure~\ref{F431rep_centers}(a) are also close to the FEM solution,
see  Figure~\ref{mapF431new}(d).
The appropriately zoomed sets of centers/vertices  shown in Figures~\ref{centresF431New}(c-h) demonstrate that RBF-FD centers 
on further refinements are distributed similar to the FEM vertices. Note that the  smoothly distributed centers in 
Figure~\ref{F431rep_centers} do not exhibit the same abrupt change in the density as seen in Figure~\ref{centresF431New}
for the adaptive methods. %
Nevertheless, the error of RBF-FD solutions on these centers is small as well, see  Figures~\ref{errors43New}(ab) and \ref{mapF431new}(h).

Finally, Table~\ref{angl_quo} includes the data on the uniformity measures $v_\mathrm{max},v_\mathrm{aver},c_\mathrm{max},c_\mathrm{aver}$ 
defined at the end of Section~\ref{Selec} for the stencil supports obtained in the experiments for each
test problem solved by the RBF-FD method.  
 We see that the value of $v_\mathrm{aver}$ is around 2.0 for the adaptively generated centers, 
 which means
that for most stencils the angle uniformity quotient $\overline{\alpha}/\underline{\alpha}$ 
is significantly less than the tolerance $v=2.5$ used in
the termination criterion of Step~II.2.ii.b. Nevertheless, there are always stencils with 
$\overline{\alpha}/\underline{\alpha}$ reaching about 4.0 and in some cases up to 7.0, which
means that the termination happens at Step~II.1 in order to prevent a very high non-uniformity of distances.
The average of the distance uniformity quotient \eqref{duq} is also very stable at around 1.3, and 
typical maximum values are about 2.5. This shows altogether that the majority of the stencil supports produced by 
the adaptive RBF-FD method of this paper are remarkably well-balanced geometrically. 
In the experiments with smoothly distributed \emph{a priori} centers generated by {\tt DistMesh} 
(see the rows of the table marked {\tt 1s} and {\tt 6as}) the numbers are significantly smaller,
which is expected.

\begin{table}[htbp!]
\begin{center}\renewcommand{\arraystretch}{1.2}\small
\begin{tabular}{|c||c|c||c|c|}
\hline
 {TP}
 &{$v_\mathrm{max}$}
&{$v_\mathrm{aver}$}
&{$c_\mathrm{max}$}
&{$c_\mathrm{aver}$}\\
\hline
\hline
1& 3.93    & 2.02    & 2.45    & 1.30  \\
\hline
1s & 2.71    & 1.40    & 2.17    & 1.11 \\
\hline
2 & 3.82    & 2.01    & 2.43    & 1.29 \\
\hline
3a & 4.22    & 2.02    & 2.38    & 1.29  \\
\hline
3b  & 3.98    & 2.02    & 2.52    & 1.30\\
\hline
3c  &   5.01    & 2.02    & 2.39    & 1.31 \\
\hline
3d  & 6.97    & 2.02    & 2.52    & 1.31   \\
\hline
4  &7.00    & 2.02    & 2.46    & 1.30    \\
\hline
5a & 5.22    & 2.02    & 2.69    & 1.29  \\
\hline
5b &   4.98    & 2.02    & 2.47    & 1.28     \\
\hline
6a &  4.29    & 2.00    & 2.41    & 1.28  \\
 \hline
6as& 2.85    & 1.39    & 2.24    & 1.11 \\
\hline
6b &   4.39    & 2.00    & 2.46    & 1.28  \\
 \hline

\end{tabular}
\caption{The uniformity measures $v_\mathrm{max},v_\mathrm{aver},c_\mathrm{max},c_\mathrm{aver}$ 
(defined at the end of Section~\ref{Selec}) for the
stencil supports obtained in each test problem. The column marked `TP' contains the test problem number and optional letters.
1s: Test Problem~\ref{f4} with smoothly distributed centers of {\tt DistMesh}.
3a, 3b, 3c, 3d: Test Problem~\ref{f42} with $\omega = \pi+0.01$, $5 \pi/4$, $7 \pi/4$ 
and $2\pi$, respectively. 5a, 5b: Test Problem~\ref{f441} with 
$\alpha = \frac{1}{10 \pi}$ and $\frac{1}{50 \pi}$, respectively. 6a, 6b: versions (a) and (b) of 
Test Problem~\ref{f43}. 6as: version (a) of Test Problem~\ref{f43} with  smoothly distributed centers of {\tt DistMesh}.
}
\label{angl_quo}
\end{center}
\end{table}

\section{Conclusion}
In this paper we suggested  an adaptive meshless RBF-FD method for the 
numerical solution of elliptic problems with point singularities. Both ingredients of the method, 
Algorithm~\ref{alg1} for stencil support selection and Algorithm~\ref{alg2Refi} for adaptive refinement
have been tested on several benchmark problems with different types of point singularities, including
the  problems of this type suggested in  \cite{Mitchell2013} for testing adaptive refinement methods. The numerical results show that the method generates reasonably
distributed adaptive centers similar to the distribution of the centers/vertices from the mesh-based adaptive
finite element method, and the errors are nearly equidistributed. For Test Problems 1--5 the size of the 
error is close to that of the FEM approximation from the same number of centers. For Test
Problem 6 whose solution features an exponential peak in the interior of the domain, the refinement regime of our method significantly differs from that of the
adaptive finite element method used for comparison, which results in a higher accuracy on fine sets of
centers at the expense of a lower accuracy on the initial coarse refinements. When finite element centers
are used instead of those generated by Algorithm~\ref{alg2Refi}, the errors are in a close agreement
with the errors of the FEM method.

The method of this paper significantly improves our earlier adaptive meshless RBF-FD method 
introduced in \cite{DavyOanh11} which suffered from certain deterioration of the quality of the centers
after a number of successive refinements. This undesirable behavior has disappeared, and this has not even
happened at the expense of the efficiency, which has in fact also improved because we were able to 
remove a post-processing step in the predecessor of Algorithm~\ref{alg2Refi} in 
\cite{DavyOanh11} whose goal was to reduce the  deterioration. A single most important new feature of 
Algorithm~\ref{alg2Refi} is an error indicator of Zienkiewicz-Zhu type used instead 
of a simple gradient estimate employed in \cite{DavyOanh11}. This leads to a significant improvement 
of the results for Test Problem~4 with a highly oscillatory solution, as demonstrated in the error plots of
Figures~\ref{figF441} and \ref{IndicF441}. 
The new stencil support selection algorithm (Algorithm~\ref{alg1}) produces geometrically better 
well-balanced stencils thanks to the introduction of a second termination criterion \eqref{term2}.

We expect that 
adaptive algorithms may be further improved by developing more sophisticated error indicators.
Further research is needed in order to achieve a similar competitive performance of meshless RBF-FD
methods for problems with line or curve singularities, boundary layers and wave fronts, as well as for time depending problems with
evolving singularities.

\bibliography{meshless}
\bibliographystyle{plain} %

\end{document}